# ALGEBRAIC STRUCTURES USING SUPER INTERVAL MATRICES

W. B. Vasantha Kandasamy
Florentin Smarandache

2011

# CONTENTS

# PREFACE

In this book authors for the first time introduce the notion of super interval matrices using the special intervals of the form [0, a], a belongs to $Z^+ \cup \{0\}$ or $Z_n$ or $Q^+ \cup \{0\}$ or $R^+ \cup \{0\}$.

The advantage of using super interval matrices is that one can build only one vector space using m × n interval matrices, but in case of super interval matrices we can have several such spaces depending on the partition on the interval matrix.

This book has seven chapters. Chapter one is introductory in nature, just introducing the super interval matrices or interval super matrices. In chapter two essential operations on super interval matrices are defined. Further in this chapter algebraic structures are defined on these super interval matrices using these operation.

Using these super interval matrices semirings and semivector spaces are defined in chapter three. This chapter gives around 90 examples. In chapter four two types of super interval semilinear algebras are introduced.

Super fuzzy interval matrices are introduced in chapter five. This chapter has two sections, in section one super fuzzy interval matrices are introduced using the fuzzy interval [0, 1]. In section two special fuzzy linear algebras using super fuzzy interval matrices are defined and described. Chapter six suggests the probable applications to interval analysis. The final chapter suggests around 110 problems some of which are at research level.

    We thank Dr. K.Kandasamy for proof reading and being extremely supportive.

<div style="text-align:right">
W.B.VASANTHA KANDASAMY  
FLORENTIN SMARANDACHE
</div>

**Chapter One**

# INTRODUCTION

In this chapter we just indicate how we build super interval matrices using intervals of the form [0, a] where a ∈ $R^+ \cup \{0\}$ or $Q^+ \cup \{0\}$ or $Z^+ \cup \{0\}$ or $Z_n$ or intervals of the form [a, b] with a ≤ b, a, b ∈ R or Z or Q.

Consider a m × n matrix

$$A = \begin{bmatrix} a_{11} & a_{12} & \cdots & a_{1n} \\ a_{21} & a_{22} & \cdots & a_{2n} \\ \vdots & \vdots & & \vdots \\ a_{m_1} & a_{m_2} & \cdots & a_{m_n} \end{bmatrix}$$

where $a_{ij}$ = [0, a], a ∈ $R^+ \cup \{0\}$ or $Z^+ \cup \{0\}$ or $Z_n$ or $Q^+ \cup \{0\}$ or [a, b] = $a_{ij}$ with a ≤ b, a, b ∈ R or Q or Z or $Z_n$.

We partition the m × n matrix A and get the super interval matrix. We in this book only use intervals of the form [0, a]; a ∈ $Z^+ \cup \{0\}$ or $Z_n$ or $R^+ \cup \{0\}$ or $Q^+ \cup \{0\}$.

Just we wish to indicate given a 3 × 2 matrix

$$\begin{bmatrix} a_1 & a_2 \\ a_3 & a_4 \\ a_5 & a_6 \end{bmatrix} = A$$

$a_i$'s intervals, we can partition A as

$$A_1 = \left[\begin{array}{c|c} a_1 & a_2 \\ a_3 & a_4 \\ a_5 & a_6 \end{array}\right], A_2 = \left[\begin{array}{cc} a_1 & a_2 \\ \hline a_3 & a_4 \\ a_5 & a_6 \end{array}\right],$$

$$A_3 = \left[\begin{array}{cc} a_1 & a_2 \\ a_3 & a_4 \\ \hline a_5 & a_6 \end{array}\right], A_4 = \left[\begin{array}{c|c} a_1 & a_2 \\ a_3 & a_4 \\ \hline a_5 & a_6 \end{array}\right],$$

$$A_5 = \left[\begin{array}{cc} a_1 & a_2 \\ \hline a_3 & a_4 \\ \hline a_5 & a_6 \end{array}\right], A_6 = \left[\begin{array}{c|c} a_1 & a_2 \\ \hline a_3 & a_4 \\ a_5 & a_6 \end{array}\right] \text{ and } A_7 = \left[\begin{array}{c|c} a_1 & a_2 \\ \hline a_3 & a_4 \\ \hline a_5 & a_6 \end{array}\right].$$

Thus using one 3 × 2 matrix A we can get 7 super interval matrices. Hence such study is not only innovative very useful in applications.

For operations on super matrices refer [17, 47].

**Chapter Two**

# INTERVAL SUPERMATRICES

In this chapter we for the first time introduce the notion of interval super matrices and derive several of their related properties.

We will illustrate this by some examples.

**DEFINITION 2.1:** *Let $X = ([a_1, b_1]\ [a_2, b_2]\ [a_3, b_3]\ |\ [a_4, b_4],\ [a_5, b_5]|\ \ldots\ |\ [a_{n-1}, b_{n-1}]\ [a_n, b_n])$ where $a_i, b_i \in R$ or $Z$ or $Q$ or $Z_m$; $1 \leq i \leq n$; $a_i < b_i$. We define X to be a super interval row matrix.*

Thus we can say in a usual super row matrix $X = (x_1\ x_2\ x_3\ |\ x_4\ x_5\ |\ x_6\ x_7\ |\ \ldots\ |\ x_{n-2}\ x_{n-1}\ x_n)$ if we replace each $x_i$ by an interval $[a_i, b_i]$, $1 \leq i \leq n$, we get the super interval row matrix.

We will illustrate this situation by some examples.

***Example 2.1***: Let $X = ([0, 2]\ |\ [1, 5]\ [2, 7]\ [-9, 13]\ |\ [-8, 2]\ [3, 12]\ [-1\ +1]\ [3, 20])$ be the row interval super matrix.

***Example 2.2***: Let $Y = ([3, 5]\ [0, 12]\ |\ [3, 7],\ [2, 4]\ [7, 9]|\ [-2, 3])$ be the row interval super matrix.

***Example 2.3***: Let A = ([-2, 1] [3, 7] [1, 4] | [5, 9], [0, 4]) be the row interval super matrix.

Now we can perform both the operations of usual addition and multiplication; only when the intervals are of special type we can define multiplication.

Now we proceed onto describe addition of super interval row matrices.

***Example 2.4***: Let X = ([2, 4] [0, 3] | [7, 10] [11, 13] [0, 5]) and Y = ([0, 3] [0, 1] | [2, 8], [10, 11] [2, 4]) be two interval super row matrices. X + Y = ([2, 4] [0, 3] | [7, 10] [11, 13] [0, 5]) + ([0, 3] [0, 1] | [2, 8] [10, 11] [2, 4]) = ([2, 4] + [0, 3] [0, 3] + [0, 1] | | [7, 10] + [2, 8] [11, 13] + [10, 11] [0, 5] + [2, 4]) = ([2, 7] [0, 4] | [9, 18] [21, 24] [2, 9]) is again an interval super row matrix.

Now we just define same type of super row interval matrices.
We say two super interval row matrices A and B are of same type (1) If both A and B have same number of components (2) Both A and B are partitioned in the same way as in case of usual super row matrices.

We will illustrate them by examples. However we wish to state that all 1 × n interval super row matrices need be compatible under addition only when they are of the same type otherwise they cannot be added.

***Example 2.5***: Let X = ([0, 3] [1, 4] | [2, 4] [0, 7] [6, 7] | [0, 8]) and Y = ([2, 4] [2, 5] | [0, 3] [0, 2] [6, 10] | [7, 10]) be two super row interval matrices of same type; we can add them as follows: X+Y = ([0, 3] [1, 4] | [2, 4] [0, 7] [6, 7] [0, 8]) + ([2, 4] [2, 5] | [0, 3] [0, 2] [6, 10] | [7, 10]) = ([0, 3] + [2, 4] [1, 4] + [2, 5] | [2, 4] + [0, 3] [0, 7] + [0, 2] [6, 7] [6, 10] | [0, 8] + [7, 10]) = ([2, 7] [3, 9] | [2, 7] [0, 9] [12, 17] | [7, 18]). We see X + Y is also the row interval super matrix of the same type as that of X and Y.

Thus we see if X and Y are super row interval matrices of the same type then so is their sum X + Y.

Further if X = ([3, 5] | [7, 8] [0, 2]) and Y = ([0, 7] [3, 5] | [6, 8]) are two super row interval matrices. We see both X and Y are 1 × 3 interval row matrices. However X and Y are not of same type, hence we cannot add X with Y. Thus X + Y is undefined. However X and Y as usual row interval matrices can be added.

In view of this we have the following theorem.

**THEOREM 2.1**: *Let T = {collection of all 1 × n row interval super matrices of the same type}. T is a semigroup under addition.*

Proof is straight forward and hence is left as an exercise to the reader.

Now we proceed onto give an example.

***Example 2.6***: Let X = ([0, 3] [1, 9] [3, 6] | [0, 7] [-2, 1] | [3, 5] [-2, 0] [7, 10]) and Y = ([3, 6], [-1, 0] [-3, 2] | [7, 8] [2, 4] | [-3, 2] [2, 7] [3, 5]) be two interval super row matrices X + Y = ([0, 3] [1, 9] [3, 6] | [0, 7] [-2, 1] | [3, 5] [-2, 0] [7, 10]) + ([3, 6] [-1, 0] [-3, 2] | [7, 8] [2, 4] | [-3, 2] [2, 7] [3, 5]) = ([0, 3] + [3, 6] [1, 9] + [-1, 0] [3, 6] + [-3, 2] | [0, 7] + [7, 8] [-2, 1] + [2, 4] | [3, 5] + [-3, 2] [-2, 0] + [2, 7] [7, 10] + [3, 5]) = ([3, 9] [0, 9] [0, 8] | [7, 15] [0, 5] | [0, 7] [0, 7] [10, 15]).

Suppose we take special type of intervals of the form [0, a] in the super interval row matrix we can define multiplication of two super interval row matrices A and B if A and B are of same type and intervals in A and B are of the form [0, a].

Let A = ([0, $a_1$] [0, $a_2$] [0, $a_3$] | [0, $a_4$] [0, $a_5$] | [0, $a_6$] [0, $a_7$] [0, $a_8$] [0, $a_9$] [0, $a_{10}$]) where $a_i \geq 0$ for i = 1, 2, …, 10 and B = ([0, $b_1$] [0, $b_2$] [0, $b_3$] | [0, $b_4$] [0, $b_5$] | [0, $b_6$] [0, $b_7$] [0, $b_8$] [0, $b_9$] [0, $b_{10}$]), $b_i \geq 0$ be two super interval matrices of the same type.
We define A + B = ([0, $a_1+b_1$] [0, $a_2+b_2$] [0, $a_3+b_3$] | [0, $a_4+b_4$] [0, $a_5+b_5$] | [0, $a_6+b_6$] [0, $a_7+b_7$] [0, $a_8+b_8$] [0, $a_9+b_9$] [0, $a_{10}+b_{10}$]); A+B is again a same type of super interval row matrix.
Now A × B = ([0, $a_1$] [0, $a_2$] [0, $a_3$] | [0, $a_4$] [0, $a_5$] | [0, $a_6$] [0, $a_7$] [0, $a_8$] [0, $a_9$] [0, $a_{10}$]) ([0, $b_1$] [0, $b_2$] [0, $b_3$] | [0, $b_4$] [0, $b_5$] | [0, $b_6$]

[0, $b_7$] [0, $b_8$] [0, $b_9$] [0, $b_{10}$]) = ([0, $a_1b_1$] [0, $a_2b_2$] [0, $a_3b_3$] | [0, $a_4b_4$] [0, $a_5b_5$] | [0, $a_6b_6$] [0, $a_7b_7$] [0, $a_8b_8$] [0, $a_9b_9$] [0, $a_{10}b_{10}$]) is again a super interval row matrix of the same type.

Let A be any row interval matrix. If A is vertically partitioned atleast once A becomes the super row interval matrix.

Consider A = ([$b_1$ $a_1$] [$b_2$, $a_2$] … [$b_n$, $a_n$]) be a row interval matrix $b_i \le a_i$; $1 \le i \le n$. Now $\overline{A}$ = ([$b_1$ $a_1$] [$b_2$, $a_2$] | [$b_3$ $a_3$] … [$b_n$, $a_n$]) is the super row interval matrix.

$\overline{A}$ = ([$b_1$ $a_1$] [$b_2$, $a_2$] [$b_3$ $a_3$] | [$b_4$ $a_4$] [$b_5$ $a_5$] [$b_6$ $a_6$] [$b_7$ $a_7$] … [$b_n$, $a_n$]). We can get for a given interval row matrix several super interval row matrices the number of such interval super row matrices depends on n.

When n = 2 we get only one super row interval matrix viz if A = ([$a_1$ $b_1$] [$a_2$ $b_2$]) then $\overline{A}$ = ([$a_1$ $b_1$] | [$a_2$ $b_2$]) is the only super row interval matrix obtained from A.

When n = 3, let A = ([$a_1$ $b_1$] [$a_2$ $b_2$] [$a_3$ $b_3$]) be a 1 × 3 row interval matrix. The number of ways in which A can be partitioned are as follows.

$\overline{A}$ = ([$a_1$ $b_1$] | [$a_2$ $b_2$] | [$a_3$ $b_3$])
$\overline{A}$ = ([$a_1$ $b_1$] | [$a_2$ $b_2$] [$a_3$ $b_3$]) and
$\overline{A}$ = ([$a_1$ $b_1$] [$a_2$ $b_2$] | [$a_3$ $b_3$])

are the three super row interval matrices obtained from A. Suppose n = 4, let X = ([$a_1$ $b_1$] [$a_2$ $b_2$] [$a_3$ $b_3$] [$a_4$ $b_4$]) be a 1 × 4 row interval matrix. The super row interval matrices obtained from A are

$\overline{A}$ = ([$a_1$ $b_1$] | [$a_2$ $b_2$] | [$a_3$ $b_3$] | [$a_4$ $b_4$])
$\overline{A}$ = ([$a_1$ $b_1$] [$a_2$ $b_2$] | [$a_3$ $b_3$] | [$a_4$ $b_4$])
$\overline{A}$ = ([$a_1$ $b_1$] | [$a_2$ $b_2$] [$a_3$ $b_3$] | [$a_4$ $b_4$])
$\overline{A}$ = ([$a_1$ $b_1$] | [$a_2$ $b_2$] | [$a_3$ $b_3$] [$a_4$ $b_4$])
$\overline{A}$ = ([$a_1$ $b_1$] [$a_2$ $b_2$] | [$a_3$ $b_3$] | [$a_4$ $b_4$])
$\overline{A}$ = ([$a_1$ $b_1$] [$a_2$ $b_2$] | [$a_3$ $b_3$] [$a_4$ $b_4$]) and
$\overline{A}$ = ([$a_1$ $b_1$] [$a_2$ $b_2$] [$a_3$ $b_3$] | [$a_4$ $b_4$]).

We get from one 1 × 4 row interval matrix obtain seven super row interval matrices. That A can be partitioned into super row interval matrices in seven ways. The problem of finding the number of super row interval matrices given an interval row matrix of say order $1 \times n$ is left as an exercise to the reader.

Now having seen how to obtain super row interval matrices from a row interval matrix we proceed onto define and give a few properties related with column super interval matrix or super column interval matrix.

**DEFINITION 2.2**: *Let* $A = \begin{bmatrix} [x_1, y_1] \\ \vdots \\ [x_n, y_n] \end{bmatrix}$ *be a column interval matrix* $x_i \leq y_i$; $1 \leq i \leq n$, $x_i$, $y_i \in R$, $Q$ *or* $Z$ *or* $Z_m$; *by horizontally partitioning the column we obtain the super column interval matrix.*

*Or equivalently we can say in a super column matrix,*

$$X = \begin{bmatrix} x_1 \\ x_2 \\ \hline x_3 \\ \vdots \\ \hline \vdots \\ \hline x_n \end{bmatrix}$$

*if we replace each* $x_i$ *by* $[a_i, b_i]$; $a_i \leq b_i$, $1 \leq i \leq n$ *then we obtain a super column interval matrix.*

We will illustrate this situation by some examples.

*Example 2.7*: Let

$$A = \begin{bmatrix} [0,3] \\ \hline [6,7] \\ \hline [1,2] \\ \hline [3,4] \\ \hline [9,18] \\ \hline [1,5] \end{bmatrix}$$

be a super column interval matrix.

***Example 2.8***: Let

$$B = \begin{bmatrix} [1,5] \\ \overline{[3,7]} \\ [0,9] \\ \overline{[-2,1]} \\ [-3,0] \\ [3,9] \end{bmatrix}$$

be a super column interval matrix.

Let

$$A = \begin{bmatrix} [a_1,b_1] \\ [a_2,b_2] \end{bmatrix}$$

be a column interval matrix. We get only one super column interval matrix $\overline{A}$ from A by partitioning A.

$$\overline{A} = \begin{bmatrix} [a_1,b_1] \\ [a_2,b_2] \end{bmatrix}.$$

Suppose

$$A = \begin{bmatrix} [a_1,b_1] \\ [a_2,b_2] \\ [a_3,b_3] \end{bmatrix}$$

be a column interval matrix.

We obtain the following super interval column matrices using A.

$$\overline{A} = \begin{bmatrix} [a_1,b_1] \\ \overline{[a_2,b_2]} \\ [a_3,b_3] \end{bmatrix}, \overline{A} = \begin{bmatrix} [a_1,b_1] \\ [a_2,b_2] \\ \overline{[a_3,b_3]} \end{bmatrix}$$

$$\text{and } \overline{A} = \begin{bmatrix} [a_1,b_1] \\ \overline{[a_2,b_2]} \\ \overline{[a_3,b_3]} \end{bmatrix}.$$

Let

$$A = \begin{bmatrix} [0,1] \\ [2,3] \\ [5,6] \\ [3,4] \\ [-3,1] \\ [0,3] \\ [-3,0] \end{bmatrix}$$

be a column interval matrix. Some of the super column interval matrices obtained from A are as follows.

$$\overline{A} = \begin{bmatrix} [0,1] \\ \overline{[2,3]} \\ \overline{[5,6]} \\ [3,4] \\ \overline{[-3,1]} \\ \overline{[0,3]} \\ \overline{[-3,0]} \end{bmatrix}, \overline{A} = \begin{bmatrix} [0,1] \\ [2,3] \\ \overline{[5,6]} \\ [3,4] \\ \overline{[-3,1]} \\ \overline{[0,3]} \\ \overline{[-3,0]} \end{bmatrix}, \overline{A} = \begin{bmatrix} [0,1] \\ [2,3] \\ [5,6] \\ [3,4] \\ \overline{[-3,1]} \\ \overline{[0,3]} \\ \overline{[-3,0]} \end{bmatrix}, \overline{A} = \begin{bmatrix} [0,1] \\ [2,3] \\ [5,6] \\ [3,4] \\ [-3,1] \\ \overline{[0,3]} \\ \overline{[-3,0]} \end{bmatrix},$$

$$\overline{A} = \begin{bmatrix} [0,1] \\ \overline{[2,3]} \\ \overline{[5,6]} \\ [3,4] \\ \overline{[-3,1]} \\ \overline{[0,3]} \\ [-3,0] \end{bmatrix}, \overline{A} = \begin{bmatrix} [0,1] \\ [2,3] \\ \overline{[5,6]} \\ [3,4] \\ \overline{[-3,1]} \\ \overline{[0,3]} \\ [-3,0] \end{bmatrix}, \overline{A} = \begin{bmatrix} [0,1] \\ [2,3] \\ [5,6] \\ [3,4] \\ \overline{[-3,1]} \\ \overline{[0,3]} \\ [-3,0] \end{bmatrix}, \overline{A} = \begin{bmatrix} [0,1] \\ [2,3] \\ [5,6] \\ [3,4] \\ [-3,1] \\ \overline{[0,3]} \\ [-3,0] \end{bmatrix},$$

$$\bar{A} = \begin{bmatrix} [0,1] \\ [2,3] \\ \overline{[5,6]} \\ [3,4] \\ \overline{[-3,1]} \\ [0,3] \\ \overline{[-3,0]} \end{bmatrix}, \bar{A} = \begin{bmatrix} [0,1] \\ [2,3] \\ \overline{[5,6]} \\ [3,4] \\ [-3,1] \\ \overline{[0,3]} \\ [-3,0] \end{bmatrix}, \bar{A} = \begin{bmatrix} [0,1] \\ [2,3] \\ \overline{[5,6]} \\ [3,4] \\ [-3,1] \\ [0,3] \\ \overline{[-3,0]} \end{bmatrix}, \bar{A} = \begin{bmatrix} [0,1] \\ [2,3] \\ \overline{[5,6]} \\ [3,4] \\ [-3,1] \\ [0,3] \\ [-3,0] \end{bmatrix}$$

and so on.

Thus with a 7 × 1 column interval matrix itself we have so many super interval column matrices.

Thus for a given n × 1 interval column matrix it is an interesting problem to find the number of super column interval matrices.

Now we say two super column interval matrices A and B are of same type;

(1) If they are of same natural order that is both A and B are n × 1 column interval matrices.
(2) They are exactly partitioned in the same way.

All n × 1 column interval super matrices need not be of the same type. For take

$$A = \begin{bmatrix} [a_1,b_1] \\ \overline{[a_2,b_2]} \\ [a_3,b_3] \\ [a_4,b_4] \\ [a_5,b_5] \end{bmatrix} \text{ and } B = \begin{bmatrix} [a_1,b_1] \\ [a_2,b_2] \\ \overline{[a_3,b_3]} \\ [a_4,b_4] \\ [a_5,b_5] \end{bmatrix}$$

A and B are not of the same type though both are 5 × 1 interval column matrices.

Let

$$A = \begin{bmatrix} [x_1, y_1] \\ [x_2, y_2] \\ [x_3, y_3] \\ [x_4, y_4] \\ [x_5, y_5] \end{bmatrix} \text{ and } B = \begin{bmatrix} [a_1, b_1] \\ [a_2, b_2] \\ [a_3, b_3] \\ [a_4, b_4] \\ [a_5, b_5] \end{bmatrix}$$

be two super column interval matrices of the same type A + B is defined and A + B is once again the super column interval matrix of the same type given by

$$A + B = \begin{bmatrix} [x_1 + a_1, y_1 + b_1] \\ [x_2 + a_2, y_2 + b_2] \\ [x_3 + a_3, y_3 + b_3] \\ [x_4 + a_4, y_4 + b_4] \\ [x_5 + a_5, y_5 + b_5] \end{bmatrix}.$$

Now as in case of usual interval column matrices of same order we cannot add super interval column matrices. Addition of super column interval matrices is possible only when they are of the same type. In view of this we have the following theorem.

**THEOREM 2.2**: *Super column interval matrices of same type is a semigroup under addition.*

Proof is straight forward and is left as an exercise for the reader.

Now in super column interval matrices we do not define multiplication as even in usual column matrices multiplication is not compatible.
Now having seen the concept of super row interval matrix and super column interval matrix we proceed on to define and illustrate super interval matrices.
Suppose A be an interval matrix and if we partition A horizontally and or vertically the resultant matrix A is defined as the super interval matrix. We can equivalently define a super interval matrix as a super matrix in which every entry in it is an interval. In case of super interval matrices the submatrices are also interval matrices.

We will illustrate this by some examples.

*Example 2.9*: Let

$$A = \begin{bmatrix} [0,6] & [1,2] & [5,9] \\ [3,9] & [8,9] & [8,15] \\ [4,8] & [10,11] & [-7,3] \\ [-3,1] & [-1,0] & [2,7] \\ [-4,0] & [4,5] & [3,8] \end{bmatrix}$$

be an interval matrix.

Suppose we draw a line between the $3^{rd}$ and $4^{th}$ row the resultant matrix denoted by $\overline{A}$ is a super interval matrix,

$$\overline{A} = \begin{bmatrix} [0,6] & [1,2] & [5,9] \\ [3,9] & [8,9] & [8,15] \\ [4,8] & [10,11] & [-7,3] \\ \hline [-3,1] & [-1,0] & [2,7] \\ [-4,0] & [4,5] & [3,8] \end{bmatrix}.$$

Suppose we vertically partition $\overline{A}$ between the first and the second column we obtain a super interval matrix given by

$$B = \begin{bmatrix} [0,6] & [1,2] & [5,9] \\ [3,9] & [8,9] & [8,15] \\ [4,8] & [10,11] & [-7,3] \\ \hline [-3,1] & [-1,0] & [2,7] \\ [-4,0] & [4,5] & [3,8] \end{bmatrix}.$$

Thus given an interval matrix A we can obtain several super interval matrices obtained from A. Hence when we adopt to super interval matrices we have several choices to work with.
   Suppose

$$A = \begin{bmatrix} [0,8] & [-1,5] \\ [-7,2] & [6,9] \end{bmatrix}$$

be an interval matrix. How many super interval matrices can be constructed using A.

We get

$$A_1 = \begin{bmatrix} [0,8] & | & [-1,5] \\ [-7,2] & | & [6,9] \end{bmatrix}, A_2 = \begin{bmatrix} [0,8] & [-1,5] \\ [-7,2] & [6,9] \end{bmatrix}$$

and

$$A_3 = \begin{bmatrix} [0,8] & | & [-1,5] \\ \hline [-7,2] & | & [6,9] \end{bmatrix}.$$

We have for this given interval matrix A there are 3 super interval matrices constructed using A.

Thus $A_1 = [B_1 \ B_2]$ where

$$B_1 = \begin{bmatrix} [0,8] \\ [-7,2] \end{bmatrix} \text{ and } B_2 = \begin{bmatrix} [-1,5] \\ [6,9] \end{bmatrix}$$

are just interval matrices so $B_1$ and $B_2$ are subinterval matrices of $A_1$. $B_1$ and $B_2$ are interval submatrices of the super interval matrix $A_1$.

Now $A_2 = \begin{bmatrix} C_1 \\ C_2 \end{bmatrix}$ where $C_1$ and $C_2$ are interval matrices given by $C_1$ = ([0, 8] [-1, 5]) and $C_2$ = ([-7, 2] [6, 9]). $C_1$ and $C_2$ are interval submatrices of $A_2$.

Also $A_3 = \begin{bmatrix} D_1 & D_2 \\ D_3 & D_4 \end{bmatrix}$ where $D_1$, $D_2$, $D_3$ and $D_4$ are interval submatrices. In this case we see $D_1$, $D_2$, $D_3$ and $D_4$ are just intervals. Thus we can also say a super interval matrix is an interval matrix whose entries are just interval matrices. The components of a super interval matrix are also called as interval submatrices.

Suppose

$$A = \begin{bmatrix} [3,6] & [1,5] & [0,3] & [6,7] \\ [-1,0] & [-5,1] & [0,8] & [0,3] \\ [2,7] & [3,8] & [0,1] & [7,10] \end{bmatrix}$$

be any interval matrix.

We will show how we can obtain super interval matrices using A.

$$A_1 = \begin{bmatrix} [3,6] & [1,5] & [0,3] & [6,7] \\ [-1,0] & [-5,1] & [0,8] & [0,3] \\ [2,7] & [3,8] & [0,1] & [7,10] \end{bmatrix}$$

$$A_2 = \begin{bmatrix} [3,6] & [1,5] & [0,3] & [6,7] \\ [-1,0] & [-5,1] & [0,8] & [0,3] \\ [2,7] & [3,8] & [0,1] & [7,10] \end{bmatrix}$$

$$A_3 = \begin{bmatrix} [3,6] & [1,5] & [0,3] & [6,7] \\ [-1,0] & [-5,1] & [0,8] & [0,3] \\ [2,7] & [3,8] & [0,1] & [7,10] \end{bmatrix},$$

$$A_4 = \begin{bmatrix} [3,6] & [1,5] & [0,3] & [6,7] \\ [-1,0] & [-5,1] & [0,8] & [0,3] \\ [2,7] & [3,8] & [0,1] & [7,10] \end{bmatrix}$$

$$A_5 = \begin{bmatrix} [3,6] & [1,5] & [0,3] & [6,7] \\ [-1,0] & [-5,1] & [0,8] & [0,3] \\ [2,7] & [3,8] & [0,1] & [7,10] \end{bmatrix}$$

$$A_6 = \begin{bmatrix} [3,6] & [1,5] & [0,3] & [6,7] \\ [-1,0] & [-5,1] & [0,8] & [0,3] \\ [2,7] & [3,8] & [0,1] & [7,10] \end{bmatrix}$$

$$A_7 = \begin{bmatrix} [3,6] & [1,5] & [0,3] & [6,7] \\ [-1,0] & [-5,1] & [0,8] & [0,3] \\ [2,7] & [3,8] & [0,1] & [7,10] \end{bmatrix}$$

$$A_8 = \begin{bmatrix} [3,6] & [1,5] & [0,3] & [6,7] \\ [-1,0] & [-5,1] & [0,8] & [0,3] \\ [2,7] & [3,8] & [0,1] & [7,10] \end{bmatrix}$$

$$A_9 = \begin{bmatrix} [3,6] & [1,5] & [0,3] & [6,7] \\ \hline [-1,0] & [-5,1] & [0,8] & [0,3] \\ \hline [2,7] & [3,8] & [0,1] & [7,10] \end{bmatrix}$$

$$A_{10} = \begin{bmatrix} [3,6] & [1,5] & [0,3] & [6,7] \\ [-1,0] & [-5,1] & [0,8] & [0,3] \\ [2,7] & [3,8] & [0,1] & [7,10] \end{bmatrix}$$

and so on.

The main importance of these super interval matrices is that for a given interval matrix A we can get several super interval matrices. Also simultaneous working of finite element analysis can be made.

Suppose we take yet another interval matrix.

$$X = \begin{bmatrix} [-3,2] & [1,0] & [9,10] \\ [8,9] & [8,10] & [9,11] \\ [8,11] & [-3,2] & [0,2] \\ [3,4] & [-9,0] & [0,1] \\ [0,5] & [-8,1] & [10,11] \\ [-7,0] & [0,5] & [-3,2] \end{bmatrix}.$$

To find some of the super interval matrices obtained from X. However we leave it as an exercise for the reader to find all super interval matrices obtained from X.

$$X_1 = \begin{bmatrix} [-3,2] & [1,0] & [9,10] \\ [8,9] & [8,10] & [9,11] \\ \hline [8,11] & [-3,2] & [0,2] \\ [3,4] & [-9,0] & [0,1] \\ [0,5] & [-8,1] & [10,11] \\ [-7,0] & [0,5] & [-3,2] \end{bmatrix},$$

$$X_2 = \begin{bmatrix} [-3,2] & [1,0] & [9,10] \\ [8,9] & [8,10] & [9,11] \\ \hline [8,11] & [-3,2] & [0,2] \\ [3,4] & [-9,0] & [0,1] \\ [0,5] & [-8,1] & [10,11] \\ [-7,0] & [0,5] & [-3,2] \end{bmatrix},$$

$$X_3 = \begin{bmatrix} [-3,2] & [1,0] & [9,10] \\ [8,9] & [8,10] & [9,11] \\ [8,11] & [-3,2] & [0,2] \\ \hline [3,4] & [-9,0] & [0,1] \\ [0,5] & [-8,1] & [10,11] \\ [-7,0] & [0,5] & [-3,2] \end{bmatrix}$$

and so on.

We leave it as an exercise for the reader to find the number of super interval matrices using the interval matrix X.

Addition of super interval matrices are possible only when they are of the same type.

Consider

$$A = \begin{bmatrix} [0,3] & [0,5] \\ \hline [8,9] & [-1,0] \end{bmatrix} \text{ and } B = \begin{bmatrix} [9,18] & [-7,0] \\ [0,2] & [-8,1] \end{bmatrix},$$

A and B are not of the same type so A + B is not defined.

Consider

$$A = \begin{bmatrix} [0,3] & [-2,0] \\ [9,10] & [-1,0] \\ [3,4] & [3,8] \end{bmatrix} \text{ and } B = \begin{bmatrix} [9,11] & [2,4] \\ [3,4] & [1,4] \\ [0,1] & [-3,0] \end{bmatrix}$$

two super interval matrices of same type.

$$A + B = \begin{bmatrix} [0,3] & [-2,0] \\ [9,10] & [-1,0] \\ [3,4] & [3,8] \end{bmatrix} + \begin{bmatrix} [9,11] & [2,4] \\ [3,4] & [1,4] \\ [0,1] & [-3,0] \end{bmatrix}$$

$$= \begin{bmatrix} [9,14] & [0,4] \\ [12,14] & [0,4] \\ [3,5] & [0,8] \end{bmatrix}.$$

A + B is a super interval matrix of same type. For multiplication of two interval matrices we need to follow certain rules which will be described later.

We will just illustrate super interval diagonal matrix.

*Example 2.10*:

$$P = \begin{bmatrix} [0,1] & [-1,0] & [2,5] & 0 & 0 \\ [7,8] & [9,12] & [1,2] & 0 & 0 \\ \hline 0 & 0 & 0 & [7,9] & [-1,0] \\ 0 & 0 & 0 & [3,4] & [7,8] \\ 0 & 0 & 0 & [0,2] & [-5,0] \end{bmatrix}$$

$$P = \begin{bmatrix} M_1 & 0 \\ \hline 0 & M_2 \end{bmatrix},$$

P is super interval diagonal matrix or P is an interval super diagonal matrix.

**Example 2.11**: Let

$$I = \begin{bmatrix} [1,1] & [0,0] & [0,0] \\ \hline [0,0] & [1,1] & [0,0] \\ [0,0] & [0,0] & [1,1] \end{bmatrix}$$

be a super interval diagonal matrix which is the super interval identity matrix or interval super identity matrix or super identity interval matrix.

Now we will just illustrate how the transpose of a super interval matrix is obtained.

Let

$$V = \begin{bmatrix} [0,3] \\ [4,5] \\ \hline [-7,0] \\ [-2,4] \\ [1,3] \\ \hline [9,10] \\ [-3,1] \\ [4,5] \\ \hline [6,7] \end{bmatrix}$$

be a super column interval matrix. Transpose of V denoted by

$$V^t = [[0, 3]\ [4, 5]\ |\ [-7, 0]\ [-2, 4]\ [1, 3]\ |\ [9, 10]\ [-3, 1]\ [4, 5]\ |\ [6, 7]].$$

Suppose

$$V = \begin{bmatrix} V_1 \\ \hline V_2 \\ \hline V_3 \\ \hline V_4 \end{bmatrix}$$

then

$$V^t = \begin{bmatrix} V_1^t \mid V_2^t \mid V_3^t \mid V_4^t \end{bmatrix}$$

where

$$V_1 = \begin{bmatrix} [0,3] \\ [4,5] \end{bmatrix}, V_2 = \begin{bmatrix} [-7,0] \\ [-2,4] \\ [1,3] \end{bmatrix}, V_3 = \begin{bmatrix} [9,10] \\ [-3,1] \\ [4,5] \end{bmatrix}$$

and $V_4 = ([6, 7])$.

$$V_4^t = [6, 7] = V_4.$$
$$V_3^t = ([9, 10] \, [-3, 1] \, [4, 5]),$$
$$V_2^t = ([-7, 0] \, [-2, 4] \, [1, 3]) \text{ and}$$
$$V_1^t = ([0, 3] \, [4, 5]).$$

***Example 2.12***: Let P = ([0, 3] [7, 8] [-1, 0] | [4,5] [-3, 2] | [-7, 1] | [8, 9] [$\sqrt{2}, \sqrt{17}$ ], [3, 12], [-1, 0]) be a super row interval matrix.
P = (P$_1$ | P$_2$ | P$_3$ | P$_4$) where P$_1$ = ([0, 3] [7, 8] [-1, 0]), P$_2$ = ([4, 5], [-3, 2]), P$_3$ = [-7, 1] and P$_4$ = ([8, 9] [$\sqrt{2}, \sqrt{17}$ ] [3, 12] [-1, 0]).

$$P^t = (P_1 \mid P_2 \mid P_3 \mid P_4)^t = \begin{bmatrix} P_1^t \mid P_2^t \mid P_3^t \mid P_4^t \end{bmatrix}^t$$

where

$$P_1^t = \begin{bmatrix} \dfrac{P_1^t}{P_2^t} \\ \dfrac{P_2^t}{P_3^t} \\ P_4^t \end{bmatrix}$$

where

$$P_1^t = \begin{bmatrix} [0,3] \\ [7,8] \\ [-1,0] \end{bmatrix}, P_2^t = \begin{bmatrix} [4,5] \\ [3,2] \end{bmatrix}, P_3^t = [-7, 1] \text{ and}$$

$$P_4^t = \begin{bmatrix} [8,9] \\ [\sqrt{2},\sqrt{17}] \\ [3,12] \\ [-1,0] \end{bmatrix}.$$

Now we will find the transpose of any general super interval matrix.

*Example 2.13*: Let

$$M = \left[\begin{array}{c|cc} [0,3] & [6,8] & [4,5] \\ [0,4] & [-7,4] & [0,2] \\ [-1,0] & [-3,1] & [-1,0] \\ \hline [2,5] & [2,7] & [7,8] \end{array}\right] = \left[\begin{array}{c|c} M_1 & M_2 \\ \hline M_3 & M_4 \end{array}\right]$$

be a super interval matrix.
Here

$$M_1 = \begin{bmatrix} [0,3] \\ [0,4] \\ [-1,0] \end{bmatrix}, M_2 = \begin{bmatrix} [6,8] & [4,5] \\ [-7,4] & [0,2] \\ [-3,1] & [-1,0] \end{bmatrix},$$

$M_3 = [2, 5]$ and $M_4 = ([2, 7] [7, 8])$.

To find the transpose of M.

$$M^t = \left[\begin{array}{ccc|c} [0,3] & [0,4] & [-1,0] & [2,5] \\ [6,8] & [-7,4] & [-3,1] & [2,7] \\ [4,5] & [0,2] & [-1,0] & [7,8] \end{array}\right] = \left[\begin{array}{c|c} M_1^t & M_2^t \\ \hline M_3^t & M_4^t \end{array}\right];$$

where
$$M_3^t = [2, 5] = M_3, \ M_1^t = ([0, 3] \ [0, 4] \ [-1, 0]),$$

$$M_4^t = \begin{pmatrix} [2,7] \\ [7,8] \end{pmatrix} \text{ and } M_2^t = \begin{pmatrix} [6,8] & [-7,4] & [-3,1] \\ [4,5] & [0,2] & [-1,0] \end{pmatrix}.$$

*Example 2.14*: Let

$$G = \begin{bmatrix} [0,1] & [6,7] & [2,3] \\ [8,9] & [3,8] & [-7,0] \\ \hline [-7,3] & [8,9] & [-1,2] \\ [4,8] & [-5,1] & [3,5] \\ \hline [-1,0] & [-7,8] & [8,12] \end{bmatrix} = \begin{bmatrix} G_1 & G_2 \\ G_3 & G_4 \\ G_5 & G_6 \end{bmatrix}$$

where

$$G_1 = \begin{bmatrix} [0,1] & [6,7] \\ [8,9] & [3,8] \end{bmatrix}, G_2 = \begin{bmatrix} [2,3] \\ [-7,0] \end{bmatrix},$$

$$G_3 = \begin{bmatrix} [-7,3] & [8,9] \\ [4,8] & [-5,1] \end{bmatrix}, G_4 = \begin{bmatrix} [-1,2] \\ [3,5] \end{bmatrix},$$

$$G_5 = ([-1, 0]\ [-7, 8]) \text{ and } G_6 = [8, 12]$$

be a super interval matrix.

$$G^t = \begin{bmatrix} [0,1] & [6,7] & [2,3] \\ [8,9] & [3,8] & [-7,0] \\ \hline [-7,3] & [8,9] & [-1,2] \\ [4,8] & [-5,1] & [3,5] \\ \hline [-1,0] & [-7,8] & [8,12] \end{bmatrix}^t$$

$$= \begin{bmatrix} [0,1] & [8,9] & [-7,3] & [4,8] & [-1,0] \\ [6,7] & [3,8] & [8,9] & [-5,1] & [-7,8] \\ [2,3] & [-7,0] & [-1,2] & [3,5] & [8,12] \end{bmatrix}$$

$$= \begin{bmatrix} G_1^t & G_2^t & G_3^t \\ G_4^t & G_5^t & G_6^t \end{bmatrix}.$$

We give another example of a super interval diagonal matrix.

*Example 2.15*:

$$A = \begin{bmatrix} [0,4] & [6,8] & [-1,0] & 0 & 0 & 0 & 0 & 0 & 0 & 0 & 0 \\ [9,11] & [-8,4] & [7,9] & 0 & 0 & 0 & 0 & 0 & 0 & 0 & 0 \\ [-5,1] & [3,7] & [3,4] & 0 & 0 & 0 & 0 & 0 & 0 & 0 & 0 \\ 0 & 0 & 0 & [3,5] & [7,9] & 0 & 0 & 0 & 0 & 0 & 0 \\ 0 & 0 & 0 & [-1,0] & [8,11] & 0 & 0 & 0 & 0 & 0 & 0 \\ 0 & 0 & 0 & 0 & 0 & [0,9] & [3,4] & [-1,0] & [1,4] & 0 & 0 \\ 0 & 0 & 0 & 0 & 0 & [0,2] & [0,7] & [-3,1] & [0,6] & 0 & 0 \\ 0 & 0 & 0 & 0 & 0 & [-1,0] & [6,8] & [-5,0] & [1,2] & 0 & 0 \\ 0 & 0 & 0 & 0 & 0 & [3,4] & [3,5] & [0,7] & [1,2] & 0 & 0 \\ 0 & 0 & 0 & 0 & 0 & 0 & 0 & 0 & 0 & [1,2] & [3,7] \\ 0 & 0 & 0 & 0 & 0 & 0 & 0 & 0 & 0 & [0,5] & [-1,0] \end{bmatrix}$$

is a super interval diagonal matrix.

*Example 2.16*: Let

$$A = \begin{bmatrix} [0,3] & [7,8] & 0 & 0 & 0 \\ [-1,0] & [3,8] & 0 & 0 & 0 \\ 0 & 0 & [7,9] & 0 & 0 \\ 0 & 0 & 0 & [0,4] & [-1,2] \\ 0 & 0 & 0 & [7,8] & [9,10] \\ 0 & 0 & 0 & [3,4] & [5,7] \end{bmatrix}$$

be the super interval diagonal matrix.

Now

$$A^t = \begin{bmatrix} [0,3] & [7,8] & 0 & 0 & 0 \\ [-1,0] & [3,8] & 0 & 0 & 0 \\ \hline 0 & 0 & [7,9] & 0 & 0 \\ \hline 0 & 0 & 0 & [0,4] & [-1,2] \\ 0 & 0 & 0 & [7,8] & [9,10] \\ 0 & 0 & 0 & [3,4] & [5,7] \end{bmatrix}^t$$

$$= \begin{bmatrix} [0,3] & [-1,0] & 0 & 0 & 0 & 0 \\ [7,8] & [3,8] & 0 & 0 & 0 & 0 \\ \hline 0 & 0 & [7,9] & 0 & 0 & 0 \\ \hline 0 & 0 & 0 & [0,4] & [7,8] & [3,4] \\ 0 & 0 & 0 & [-1,2] & [9,10] & [5,7] \end{bmatrix}.$$

Next we see when ever we multiply two intervals in general the resultant may or may not belong to the usual form of intervals.

For consider [-1, 2] and [-3, 0]. Their product is [-1, 2] [-3, 0] = [3, 0] and [3, 0] is not the usual form of intervals.

We over come this problem in two ways.

We use only intervals of the form [0, a], $0 \leq a$ is the first method.

The second method is define the notion of increasing and decreasing interval.

We say [a, b] if $a \leq b$ to be an increasing interval [a, b] in which a = b is a degenerate interval.

[a, b] such that a > b is defined as the decreasing interval [-4, -9] is a decreasing interval. [5, 0] is a decreasing interval when we multiply them we get [-4, 9] × [5, 0] = [-20, 0] which is an increasing interval. Likewise we get the product of two increasing intervals to be a decreasing interval. Thus when we consider [-5, 0] and [-7, 1], two increasing intervals we see [-5, 0] × [-7, 1] = [35, 0] which is a decreasing interval. Hence if our collection of special intervals include a decreasing interval, increasing interval and a degenerate one we see this collection works well.

This sort of matrices will not be discussed in this book.

Now we consider only super matrices with intervals of the form [0, a] where a ∈ $Z^+ \cup \{0\}$ or $R^+ \cup \{0\}$ or $Q^+ \cup \{0\}$ or $Z_n$ (n < ∞).

We can call them as positive interval super matrices or positive super interval matrices. Before we proceed onto define operations on them we will illustrate this situation by some examples.

*Example 2.17*: Let A = ([0, 5] [0, 9] | [0, 12] [0, 14] [0, 98] | [0, 1012]) be a row positive interval super matrix with entries from $Z^+ \cup \{0\}$.

*Example 2.18*: Let B = ([0, $\sqrt{2}$ ] | [0, $\sqrt{7}$ ] [0, 5/3] [0, 3$\sqrt{3}$ ] [0, 3 $\sqrt{19}$ ] | [0, 100 + 5 $\sqrt{17}$ ] [0, 120$\sqrt{43}$ + $\sqrt{7}$ ]) be a positive row interval super matrix with entries from $R^+ \cup \{0\}$.

*Example 2.19*: Let C = ([0, 7/2] [0, 9/5] [0, 12/13] | [0, 8] [0, 9/19] [0, 4] | [0, 3/7] [0, 81] [0, 40]) be a positive row interval super matrix.

*Example 2.20*: Let S = ([0, 3] | [0, 2] | [0, 7] [0, 12] [0, 15] | [0, 16] [0, 14] | [0, 10] [0, 7]) be a row interval super matrix with entries from $Z_{17}$.

Now having seen positive row interval matrix we now proceed onto give examples of positive column interval super matrix.

*Example 2.21*: Let

$$P = \begin{bmatrix} [0,9] \\ \left[0,\sqrt{8}\right] \\ \hline \left[0,\sqrt{7}/2\right] \\ \left[0,5\sqrt{2}\right] \\ \hline \left[0,17+\sqrt{3}\right] \end{bmatrix}$$

be a positive row interval super matrix with entries from $R^+ \cup \{0\}$.

*Example 2.22*: Let

$$R = \begin{bmatrix} \overline{[0,3]} \\ [0,2] \\ [0,1] \\ \overline{[0,42]} \\ [0,146] \\ [0,27] \\ \overline{[0,1001]} \\ [0,12] \end{bmatrix}$$

be a positive column interval super matrix with entries from $Z^+ \cup \{0\}$.

*Example 2.24*: Let

$$W = \begin{bmatrix} [0,1] \\ [0,2] \\ [0,3] \\ [0,5] \\ \overline{[0,6]} \\ [0,4] \\ [0,5] \\ \overline{[0,1]} \\ [0,6] \\ [0,3] \end{bmatrix}$$

be a positive super column interval matrix with entries from $Z_7$.

*Example 2.23*: Let $\quad T = \begin{bmatrix} [0,3] \\ [0,7/2] \\ [0,8/9] \\ \overline{[0,12/13]} \\ [0,19] \\ \overline{[0,8]} \\ \overline{[0,1]} \end{bmatrix}$

be a positive column interval super matrix from $Q^+ \cup \{0\}$.

Now we will give examples of positive super interval matrices before we proceed onto define operations using them.

**Example 2.25**: Let

$$W = \begin{bmatrix} [0,3] & [0,7] \\ [0,8] & [0,10] \\ \hline [0,1] & [0,40] \\ [0,14] & [0,11] \end{bmatrix}$$

be a positive super interval matrix with entries from $Z_{45}$.

**Example 2.26**: Let

$$P = \begin{bmatrix} [0,\sqrt{8}] & [0,1] & [0,41] & [0,4] & [0,3+\sqrt{17}] & [0,11] & [0,4] \\ [0,\sqrt{11}] & [0,9] & [0,14] & [0,\sqrt{5}+4] & [0,\sqrt{19}] & [0,2] & [0,7] \\ \hline [0,3] & [0,\sqrt{43}] & [0,\sqrt{17}/2] & [0,\sqrt{3}+4\sqrt{7}] & [0,\sqrt{23}] & [0,15] & [0,9] \end{bmatrix}$$

be a positive super interval matrix with entries from $R^+ \cup \{0\}$.

**Example 2.27**: Let

$$A = \begin{bmatrix} [0,3] & [0,8] & [0,9] & [0,1] \\ \hline [0,0] & [0,5] & [0,1021] & [0,120] \\ [0,16] & [0,14] & [0,401] & [0,2] \\ \hline [0,91] & [0,101] & [0,3] & [0,5] \end{bmatrix}$$

be a positive super interval square matrix with entries from $Z^+ \cup \{0\}$.

**Example 2.28**: Let

$$X = \begin{bmatrix} [0,1] & [0,7] & [0,1] & [0,4] & [0,3] & [0,1] \\ \hline [0,2] & [0,6] & [0,6] & [0,5] & [0,4] & [0,2] \\ [0,3] & [0,5] & [0,2] & [0,2] & [0,5] & [0,3] \\ \hline [0,4] & [0,3] & [0,6] & [0,1] & [0,1] & [0,0] \\ [0,5] & [0,4] & [0,3] & [0,5] & [0,5] & [0,6] \\ [0,6] & [0,0] & [0,6] & [0,2] & [0,1] & [0,0] \end{bmatrix}$$

be a positive interval super square matrix with entries from $Z_8$.

*Example 2.29*: Let

$$P = \begin{bmatrix} \begin{bmatrix} 0,\sqrt{2} \end{bmatrix} & [0,0] & \begin{bmatrix} 0,4+\sqrt{2} \end{bmatrix} \\ \begin{bmatrix} 0,\sqrt{17} \end{bmatrix} & [0,1] & [0,6] \\ \begin{bmatrix} 0,4+\sqrt{3} \end{bmatrix} & [0,7] & \begin{bmatrix} 0,\sqrt{3} \end{bmatrix} \end{bmatrix}$$

be a super positive interval square matrix with entries from $R^+ \cup \{0\}$.

Having seen examples of positive super interval matrices we will now illustrate how the product of a interval super row matrix with its transpose is determined.

*Example 2.30*: Let A = ([0, 3] | [0, 5] [0, 1] [0, 12] | [0, 10] [0, 1] [0, 0] [0, 2]) super row positive interval matrix with entries from $Z^+ \cup \{0\}$.

$A^t$ = ([0, 3] | [0, 5] [0, 1] [0, 12] | [0, 10] [0, 1] [0, 0] [0, 2])$^t$

$$= \begin{bmatrix} [0,3] \\ [0,5] \\ [0,1] \\ [0,12] \\ [0,10] \\ [0,1] \\ [0,0] \\ [0,2] \end{bmatrix}$$

is the transpose of A which is a positive column interval super matrix.

   Now we will find the product

$$AA^t = ([0, 3] \mid [0, 5]\ [0, 1]\ [0, 12] \mid [0, 10]\ [0, 1]\ [0, 0]\ [0, 2]) \times \begin{bmatrix} [0,3] \\ \overline{[0,5]} \\ [0,1] \\ [0,12] \\ \overline{[0,10]} \\ [0,1] \\ [0,0] \\ [0,2] \end{bmatrix}$$

$$= ([0, 3] \times [0, 3] + ([0, 5]\ [0, 1]\ [0, 12]) \times \begin{bmatrix} [0,5] \\ [0,1] \\ [0,12] \end{bmatrix} +$$

$$([0, 10]\ [0, 1]\ [0, 0]\ [0, 2]) \times \begin{bmatrix} [0,10] \\ [0,1] \\ [0,0] \\ [0,2] \end{bmatrix}$$

$$= ([0, 9] + ([0, 5] \times [0, 5] + [0, 1] \times [0, 1] + [0, 12] \times [0, 12]) + ([0, 10] \times [0, 10] + [0, 1]\ [0, 1] + [0, 0]\ [0, 0] + [0, 2] \times [0, 2])$$

$$= ([0, 9] + [0, 25] + [0, 1] + [0, 144]) + ([0, 100] + [0, 1] + [0, 0] + [0, 4])$$

$$= [0, 179] + [0, 105]$$

$$= [0, 284].$$

Now always the product of a positive super row interval matrix with its transpose can be found and it will always be an interval which is also positive and is of the form [0, a].

***Example 2.31***: Let P = ([0, 4] [0, 2] | [0, 7] [0, 8] [0, 0] [0, 1] | [0, 5] [0, 3] [0, 9]) be a super row positive interval matrix with entries from $Z_{10}$. To find $AA^t$.

$$A^t = \begin{bmatrix} [0,4] \\ [0,2] \\ \hline [0,7] \\ [0,8] \\ [0,0] \\ [0,1] \\ \hline [0,5] \\ [0,3] \\ [0,9] \end{bmatrix}.$$

Now

$$\begin{aligned}
A \cdot A^t &= ([0,4]\ [0,2] \mid [0,7]\ [0,8]\ [0,0]\ [0,1] \mid [0,5]\ [0,3] \\
&\quad [0,9]) \begin{bmatrix} [0,4] \\ [0,2] \\ \hline [0,7] \\ [0,8] \\ [0,0] \\ [0,1] \\ \hline [0,5] \\ [0,3] \\ [0,9] \end{bmatrix} \\
&= (([0,4]\ [0,2]) \begin{bmatrix} [0,4] \\ [0,2] \end{bmatrix} + ([0,7]\ [0,8]\ [0,0]\ [0,1]) \\
&\quad \begin{bmatrix} [0,7] \\ [0,8] \\ [0,0] \\ [0,1] \end{bmatrix} + ([0,5]\ [0,3]\ [0,9]) \begin{bmatrix} [0,5] \\ [0,3] \\ [0,9] \end{bmatrix} \\
&= ([0,16] + [0,4]) + ([0,49] + [0,64] + [0,0] + [0,1]) + \\
&\quad ([0,25] + [0,9] + [0,81])\ (\bmod\ 10) \\
&= [0,0] + [0,4] + [0,5]\ (\bmod\ 10) \\
&= [0,9].
\end{aligned}$$

**Example 2.32**: To find $A^t$ with A where A = ([0, 3] [0, 2] | [0, 7] | [0, 1] [0, 8] [0, 5]). Now

$$A^t = \begin{bmatrix} [0,3] \\ [0,2] \\ \hline [0,7] \\ \hline [0,1] \\ [0,8] \\ [0,5] \end{bmatrix}$$

is the transpose of A.
To find

$$A^t \times A = \begin{bmatrix} [0,3] \\ [0,2] \\ \hline [0,7] \\ \hline [0,1] \\ [0,8] \\ [0,5] \end{bmatrix} \times ([0, 3]\ [0, 2] \mid [0, 7] \mid [0, 1]\ [0, 8]\ [0, 5]).$$

$$= \begin{bmatrix} \begin{bmatrix} [0,3] \\ [0,2] \end{bmatrix}([0,3][0,2]) & \begin{bmatrix} [0,3] \\ [0,2] \end{bmatrix}[0,7] & \begin{bmatrix} [0,3] \\ [0,2] \end{bmatrix}([0,1][0,8][0,5]) \\ \hline [0,7]([0,3][0,2]) & [0,7][0,7] & [0,7]([0,1][0,8][0,5]) \\ \hline \begin{bmatrix} [0,1] \\ [0,8] \\ [0,5] \end{bmatrix}([0,3][0,2]) & \begin{bmatrix} [0,1] \\ [0,8] \\ [0,5] \end{bmatrix}[0,7] & \begin{bmatrix} [0,1] \\ [0,8] \\ [0,5] \end{bmatrix}([0,1][0,8][0,5]) \end{bmatrix}$$

$$= \begin{bmatrix} [0,9] & [0,6] & [0,21] & [0,3] & [0,24] & [0,15] \\ [0,6] & [0,4] & [0,14] & [0,2] & [0,16] & [0,10] \\ \hline [0,21] & [0,14] & [0,49] & [0,7] & [0,56] & [0,35] \\ \hline [0,3] & [0,2] & [0,7] & [0,1] & [0,8] & [0,5] \\ [0,24] & [0,16] & [0,56] & [0,8] & [0,64] & [0,40] \\ [0,15] & [0,10] & [0,35] & [0,5] & [0,40] & [0,25] \end{bmatrix}.$$

***Example 2.33***: Let

$$S = \begin{bmatrix} [0,2] \\ [0,1] \\ [0,5] \\ \hline [0,7] \\ [0,1] \\ [0,11] \\ [0,10] \end{bmatrix}$$

be a positive row interval super matrix.

$$S^t = ([0, 2]\ [0, 1]\ [0, 5]\ |\ [0, 7]\ [0, 1]\ [0, 11]\ [0, 10])$$

Now

$$SS^t = \begin{bmatrix} [0,2] \\ [0,1] \\ [0,5] \\ \hline [0,7] \\ [0,1] \\ [0,11] \\ [0,10] \end{bmatrix} ([0, 2]\ [0, 1]\ [0, 5]\ |\ [0, 7]\ [0, 1][0, 11]\ [0, 10])$$

$$= \begin{bmatrix} \begin{bmatrix} [0,2] \\ [0,1] \\ [0,5] \end{bmatrix}([0,2][0,1][0,5]) & \begin{bmatrix} [0,2] \\ [0,7] \\ [0,5] \end{bmatrix}([0,7][0,1][0,11][0,10]) \\ \hline \begin{bmatrix} [0,7] \\ [0,1] \\ [0,11] \\ [0,10] \end{bmatrix}([0,2][0,1][0,5]) & \begin{bmatrix} [0,7] \\ [0,1] \\ [0,11] \\ [0,10] \end{bmatrix}([0,7][0,1][0,11][0,10]) \end{bmatrix}$$

$$= \begin{bmatrix} \begin{pmatrix} [0,4] & [0,2] & [0,10] \\ [0,2] & [0,1] & [0,5] \\ [0,10] & [0,5] & [0,25] \end{pmatrix} & \begin{pmatrix} [0,14] & [0,2] & [0,22] & [0,20] \\ [0,7] & [0,1] & [0,11] & [0,10] \\ [0,35] & [0,5] & [0,55] & [0,50] \end{pmatrix} \\ \begin{bmatrix} [0,14] & [0,7] & [0,35] \\ [0,2] & [0,1] & [0,5] \\ [0,22] & [0,11] & [0,55] \\ [0,20] & [0,10] & [0,50] \end{bmatrix} & \begin{pmatrix} [0,49] & [0,7] & [0,77] & [0,70] \\ [0,7] & [0,1] & [0,11] & [0,10] \\ [0,77] & [0,11] & [0,121] & [0,110] \\ [0,70] & [0,10] & [0,110] & [0,100] \end{pmatrix} \end{bmatrix}.$$

We have just shown how product of a row super interval positive matrix with its transpose is defined. Likewise the product of a column positive super interval matrix with its transpose also is illustrated.

Now we will exhibit how the product of a positive column interval super matrix with a positive row interval super matrix is carried out.

*Example 2.34*: Let

$$P = \begin{bmatrix} [0,2] \\ \overline{[0,1]} \\ [0,0] \\ \overline{[0,1]} \\ [0,2] \\ [0,0] \end{bmatrix}$$

be a super positive column interval matrix and Q = ([0, 3] [0, 1] [0, 2] | [0, 1] [0, 2] | [0,3] [0, 4] [0, 0] [0, 0] | [0, 5]) be a positive row interval super matrix. To find the product

$$PQ = \begin{bmatrix} [0,2] \\ \overline{[0,1]} \\ [0,0] \\ \overline{[0,1]} \\ [0,2] \\ [0,0] \end{bmatrix} \; ([0, 3] \; [0, 1] \; [0, 2] \; | \; [0, 1] \; [0, 2] \; | \; [0,3] \; [0, 4] \; [0, 0]$$

$$[0, 0] \; | \; [0, 5])$$

$$= \begin{bmatrix} [0,2]([0,3][0,1][0,2]) & [0,2]([0,1][0,2][0,3][0,4]0,0) & [0,2][0,5] \\ \begin{bmatrix} [0,1] \\ [0,0] \end{bmatrix}([0,3][0,1][0,2]) & \begin{bmatrix} [0,1] \\ [0,0] \end{bmatrix}([0,1][0,2][0,3][0,4]0,0) & \begin{bmatrix} [0,1] \\ [0,0] \end{bmatrix}[0,5] \\ \hline \begin{bmatrix} [0,1] \\ [0,2] \\ [0,0] \end{bmatrix}([0,3][0,1][0,2]) & \begin{bmatrix} [0,1] \\ [0,2] \\ [0,0] \end{bmatrix}([0,1][0,2][0,3][0,4]0,0) & \begin{bmatrix} [0,1] \\ [0,2] \\ [0,0] \end{bmatrix}[0,5] \end{bmatrix}$$

$$= \begin{bmatrix} [0,6] & [0,2] & [0,4] & [0,2] & [0,4] & [0,6] & [0,8] & 0 & 0 & [0,10] \\ [0,3] & [0,1] & [0,2] & [0,1] & [0,2] & [0,3] & [0,4] & 0 & 0 & [0,5] \\ [0,0] & [0,0] & [0,0] & 0 & 0 & 0 & 0 & 0 & 0 & 0 \\ \hline [0,3] & [0,1] & [0,2] & [0,1] & [0,2] & [0,3] & [0,4] & 0 & 0 & [0,5] \\ [0,6] & [0,2] & [0,4] & [0,2] & [0,4] & [0,6] & [0,8] & 0 & 0 & [0,10] \\ 0 & 0 & 0 & 0 & 0 & 0 & 0 & 0 & 0 & [0,0] \end{bmatrix}.$$

Now we will illustrate by some examples the product of a positive row interval super matrix with a positive row interval super matrix.

**Example 2.35**: Let V = ([0, 1] [0, 6] [0, 2] | [0, 3] | [0, 5] [0, 6] [0, 1] [0, 4]) be a positive row interval super matrix.
Let

$$W = \begin{bmatrix} [0,5] \\ [0,2] \\ [0,3] \\ \hline [0,1] \\ \hline [0,2] \\ [0,8] \\ [0,2] \\ [0,5] \end{bmatrix}$$

be a positive super column interval matrix. To find the product

$$V \cdot W = ([0, 1]\ [0, 6]\ [0, 2]\ |\ [0, 3]\ |\ [0, 5]\ [0, 6]\ [0, 1]\ [0, 4]) \begin{bmatrix} [0,5] \\ [0,2] \\ [0,3] \\ \overline{[0,1]} \\ \overline{[0,2]} \\ [0,8] \\ [0,2] \\ [0,5] \end{bmatrix}$$

$$= ([0, 1]\ [0, 6]\ [0, 2]) \begin{bmatrix} [0,5] \\ [0,2] \\ [0,3] \end{bmatrix} + ([0, 3])\ ([0,1]) + ([0, 5]\ [0, 6]$$

$$[0, 1]\ [0, 4]) \begin{bmatrix} [0,2] \\ [0,8] \\ [0,2] \\ [0,5] \end{bmatrix}$$

$$= ([0, 5] + [0, 12] + [0, 6]) + [0, 3] + ([0, 10] + [0, 48] + [0, 2] + [0, 20])$$

$$= [0, 23] + [0, 3] + [0, 80] = [0, 106].$$

***Example 2.36***: Let

$$X = \begin{bmatrix} [0,1] & [0,1] & [0,3] & [0,1] & [0,1] & [0,1] \\ [0,2] & [0,3] & [0,1] & [0,2] & [0,0] & [0,1] \\ [0,3] & [0,4] & [0,2] & [0,0] & [0,1] & [0,0] \\ [0,4] & [0,2] & [0,4] & [0,1] & [0,0] & [0,0] \end{bmatrix}$$

and

$$Y = \begin{bmatrix} [0,2] & [0,1] \\ \overline{[0,3]} & \overline{[0,1]} \\ [0,1] & [0,0] \\ \overline{[0,2]} & \overline{[0,1]} \\ [0,1] & [0,2] \\ [0,1] & [0,1] \end{bmatrix}$$

be two positive super interval matrices.

To find the product

$$X.Y = \begin{bmatrix} [0,1] & [0,1] & [0,3] & [0,1] & [0,1] & [0,1] \\ [0,2] & [0,3] & [0,1] & [0,2] & [0,0] & [0,1] \\ [0,3] & [0,4] & [0,2] & [0,0] & [0,1] & [0,0] \\ [0,4] & [0,2] & [0,4] & [0,1] & [0,0] & [0,0] \end{bmatrix} \begin{bmatrix} [0,2] & [0,1] \\ [0,3] & [0,1] \\ [0,1] & [0,0] \\ [0,2] & [0,1] \\ [0,1] & [0,2] \\ [0,1] & [0,1] \end{bmatrix}$$

$$= \begin{bmatrix} [0,1] \\ [0,2] \\ [0,3] \\ [0,4] \end{bmatrix} ([0,2]\ [0,1]) + \begin{bmatrix} [0,1] & [0,3] \\ [0,3] & [0,1] \\ [0,4] & [0,2] \\ [0,2] & [0,4] \end{bmatrix} \begin{bmatrix} [0,3] & [0,1] \\ [0,1] & [0,0] \end{bmatrix}$$

$$+ \begin{bmatrix} [0,1] & [0,1] & [0,1] \\ [0,2] & [0,0] & [0,1] \\ [0,0] & [0,1] & [0,0] \\ [0,1] & [0,0] & [0,0] \end{bmatrix} \begin{bmatrix} [0,2] & [0,1] \\ [0,1] & [0,2] \\ [0,1] & [0,1] \end{bmatrix}$$

$$= \begin{bmatrix} [0,2] & [0,1] \\ [0,4] & [0,2] \\ [0,6] & [0,3] \\ [0,8] & [0,4] \end{bmatrix} + \begin{bmatrix} [0,6] & [0,1] \\ [0,10] & [0,3] \\ [0,14] & [0,4] \\ [0,10] & [0,2] \end{bmatrix} + \begin{bmatrix} [0,4] & [0,4] \\ [0,5] & [0,3] \\ [0,1] & [0,2] \\ [0,2] & [0,1] \end{bmatrix}$$

$$= \begin{bmatrix} [0,12] & [0,6] \\ [0,19] & [0,8] \\ [0,21] & [0,9] \\ [0,20] & [0,7] \end{bmatrix}.$$

**Example 2.37**: Let
$$X = ([0, 5] \mid [0, 0]\ [0, 1] \mid [0, 1]\ [0, 1]\ [0, 1])$$

and

$$Y = \begin{bmatrix} [0,1] & [0,1] & [0,2] & [0,1] \\ \hline [0,1] & [0,0] & [0,2] & [0,4] \\ [0,0] & [0,1] & [0,0] & [0,3] \\ \hline [0,1] & [0,1] & [0,0] & [0,0] \\ [0,1] & [0,0] & [0,1] & [0,1] \\ [0,0] & [0,1] & [0,0] & [0,1] \end{bmatrix}$$

be two positive super interval matrices.

To find

$$XY = ([0, 5] \mid [0, 0]\ [0, 1] \mid [0, 1]\ [0, 1]\ [0, 1]) \times \begin{bmatrix} [0,1] & [0,1] & [0,2] & [0,1] \\ \hline [0,1] & [0,0] & [0,2] & [0,4] \\ [0,0] & [0,1] & [0,0] & [0,3] \\ \hline [0,1] & [0,1] & [0,0] & [0,0] \\ [0,1] & [0,0] & [0,1] & [0,1] \\ [0,0] & [0,1] & [0,0] & [0,1] \end{bmatrix}$$

$$= [0, 5]\ ([0, 1]\ [0, 1]\ [0, 2]\ [0, 1])$$
$$+ ([0, 0]\ [0, 1]) \begin{bmatrix} [0,1] & [0,0] & [0,2] & [0,4] \\ [0,0] & [0,1] & [0,0] & [0,3] \end{bmatrix}$$
$$+ ([0, 1]\ [0, 1]\ [0, 1]) \begin{bmatrix} [0,1] & [0,1] & [0,0] & [0,0] \\ [0,1] & [0,0] & [0,1] & [0,1] \\ [0,0] & [0,1] & [0,0] & [0,1] \end{bmatrix}$$

$$= ([0, 5]\ [0, 5]\ [0, 10]\ [0, 5]) + ([0, 0]\ [0, 1]\ [0, 0]\ [0, 3])$$
$$+ ([0, 2]\ [0, 2]\ [0, 1]\ [0, 2])$$
$$= ([0, 7]\ [0, 8]\ [0, 11]\ [0, 10]).$$

***Example 2.38***: Let

$$X = \begin{bmatrix} [0,1] & [0,2] & [0,1] & [0,1] & [0,2] & [0,3] \\ [0,3] & [0,1] & [0,2] & [0,3] & [0,1] & [0,1] \end{bmatrix}$$

be a interval positive row matrix and

$$Y = \begin{bmatrix} [0,1] & [0,1] & [0,2] & [0,1] \\ \hline [0,1] & [0,0] & [0,2] & [0,4] \\ [0,0] & [0,1] & [0,0] & [0,3] \\ \hline [0,1] & [0,1] & [0,0] & [0,0] \\ [0,1] & [0,0] & [0,1] & [0,1] \\ [0,0] & [0,1] & [0,0] & [0,1] \end{bmatrix}$$

be a positive interval super column matrix.
To find

$$X.Y =$$

$$\begin{bmatrix} [0,1] & [0,2] & [0,1] & [0,1] & [0,2] & [0,3] \\ [0,3] & [0,1] & [0,2] & [0,3] & [0,1] & [0,1] \end{bmatrix} \begin{bmatrix} [0,1] & [0,1] & [0,2] & [0,1] \\ \hline [0,1] & [0,0] & [0,2] & [0,4] \\ [0,0] & [0,1] & [0,0] & [0,3] \\ \hline [0,1] & [0,1] & [0,0] & [0,0] \\ [0,1] & [0,0] & [0,1] & [0,1] \\ [0,0] & [0,1] & [0,0] & [0,1] \end{bmatrix}$$

$$= \begin{bmatrix} [0,1] \\ [0,3] \end{bmatrix} ([0,1]\ [0,1]\ [0,2]\ [0,1])$$

$$+ \begin{bmatrix} [0,2] & [0,1] \\ [0,1] & [0,2] \end{bmatrix} \begin{bmatrix} [0,1] & [0,0] & [0,2] & [0,4] \\ [0,0] & [0,1] & [0,0] & [0,3] \end{bmatrix}$$

$$+ \begin{bmatrix} [0,1] & [0,2] & [0,3] \\ [0,3] & [0,1] & [0,1] \end{bmatrix} \begin{bmatrix} [0,1] & [0,1] & [0,0] & [0,0] \\ [0,1] & [0,0] & [0,1] & [0,1] \\ [0,0] & [0,1] & [0,0] & [0,1] \end{bmatrix}$$

$$= \begin{bmatrix} [0,1] & [0,1] & [0,2] & [0,1] \\ [0,3] & [0,3] & [0,6] & [0,3] \end{bmatrix}$$

$$+ \begin{bmatrix} [0,2] & [0,1] & [0,4] & [0,11] \\ [0,1] & [0,2] & [0,2] & [0,10] \end{bmatrix} + \begin{bmatrix} [0,3] & [0,4] & [0,2] & [0,5] \\ [0,4] & [0,4] & [0,1] & [0,2] \end{bmatrix}$$

$$= \begin{bmatrix} [0,6] & [0,6] & [0,8] & [0,17] \\ [0,8] & [0,9] & [0,9] & [0,15] \end{bmatrix}.$$

Now we will show how the product of super positive interval column vector with a positive super interval row vector is carried out.

*Example 2.39*: Let

$$X = \begin{bmatrix} [0,1] & [0,2] \\ [0,1] & [0,1] \\ \hline [0,1] & [0,2] \\ [0,1] & [0,2] \\ \hline [0,0] & [0,1] \\ \hline [0,1] & [0,2] \end{bmatrix}$$

be a positive super interval column vector.

$$Y = \begin{bmatrix} [0,1] & [0,1] & [0,1] & | & [0,1] & [0,0] & | & [0,1] \\ [0,2] & [0,1] & [0,2] & | & [0,2] & [0,1] & | & [0,2] \end{bmatrix}$$

be a positive super interval row vector.

$$X.Y = \begin{bmatrix} [0,1] & [0,2] \\ [0,1] & [0,1] \\ \hline [0,1] & [0,2] \\ [0,1] & [0,2] \\ \hline [0,0] & [0,1] \\ \hline [0,1] & [0,2] \end{bmatrix} \times \begin{bmatrix} [0,1] & [0,1] & [0,1] & | & [0,1] & [0,0] & | & [0,1] \\ [0,2] & [0,1] & [0,2] & | & [0,2] & [0,1] & | & [0,2] \end{bmatrix}$$

$$= \begin{bmatrix} \begin{bmatrix} [0,1] & [0,2] \\ [0,1] & [0,1] \\ [0,1] & [0,2] \end{bmatrix} \begin{bmatrix} [0,1] & [0,1] & [0,1] \\ [0,2] & [0,1] & [0,2] \end{bmatrix} & \begin{bmatrix} [0,1] & [0,2] \\ [0,1] & [0,1] \\ [0,1] & [0,2] \end{bmatrix} \begin{pmatrix} [0,1] & [0,0] \\ [0,2] & [0,1] \end{pmatrix} & \begin{bmatrix} [0,1] & [0,2] \\ [0,1] & [0,1] \\ [0,1] & [0,2] \end{bmatrix} \begin{bmatrix} [0,1] \\ [0,2] \end{bmatrix} \\ \hline \begin{pmatrix} [0,1] & [0,2] \\ [0,0] & [0,1] \end{pmatrix} \begin{bmatrix} [0,1] & [0,1] & [0,1] \\ [0,2] & [0,1] & [0,2] \end{bmatrix} & \begin{pmatrix} [0,1] & [0,2] \\ [0,0] & [0,1] \end{pmatrix} \begin{pmatrix} [0,1] & [0,0] \\ [0,2] & [0,1] \end{pmatrix} & \begin{pmatrix} [0,1] & [0,2] \\ [0,0] & [0,1] \end{pmatrix} \begin{bmatrix} [0,1] \\ [0,2] \end{bmatrix} \\ \hline ([0,1][0,2]) \begin{bmatrix} [0,1] & [0,1] & [0,1] \\ [0,2] & [0,1] & [0,2] \end{bmatrix} & ([0,1][0,2]) \begin{pmatrix} [0,1] & [0,0] \\ [0,2] & [0,1] \end{pmatrix} & ([0,1][0,2]) \begin{bmatrix} [0,1] \\ [0,2] \end{bmatrix} \end{bmatrix}$$

$$= \begin{bmatrix} [0,5] & [0,3] & [0,5] & [0,5] & [0,2] & [0,5] \\ [0,3] & [0,2] & [0,3] & [0,3] & [0,1] & [0,3] \\ [0,5] & [0,3] & [0,5] & [0,5] & [0,2] & [0,5] \\ \hline [0,5] & [0,3] & [0,5] & [0,5] & [0,2] & [0,5] \\ [0,2] & [0,1] & [0,2] & [0,2] & [0,1] & [0,2] \\ \hline [0,5] & [0,3] & [0,5] & [0,5] & [0,2] & [0,5] \end{bmatrix}.$$

*Example 2.40*: Let

$$X = \begin{bmatrix} [0,3] & [0,1] & [0,6] \\ [0,2] & [0,0] & [0,0] \\ \hline [0,1] & [0,2] & [0,3] \\ \hline [0,6] & [0,3] & [0,0] \\ [0,4] & [0,2] & [0,1] \\ [0,5] & [0,1] & [0,0] \end{bmatrix}$$

be a positive super interval column vector. Let

$$Y = \begin{bmatrix} [0,3] & [0,5] & [0,2] & [0,0] \\ [0,1] & [0,1] & [0,2] & [0,2] \\ [0,0] & [0,3] & [0,1] & [0,0] \end{bmatrix}$$

be a positive super interval row vector.

To find the product X.Y.

$$X.Y = \begin{bmatrix} [0,3] & [0,1] & [0,6] \\ [0,2] & [0,0] & [0,0] \\ \hline [0,1] & [0,2] & [0,3] \\ [0,6] & [0,3] & [0,0] \\ [0,4] & [0,2] & [0,1] \\ [0,5] & [0,1] & [0,0] \end{bmatrix} \cdot \begin{bmatrix} [0,3] & [0,5] & [0,2] & [0,0] \\ [0,1] & [0,1] & [0,2] & [0,2] \\ [0,0] & [0,3] & [0,1] & [0,0] \end{bmatrix}$$

$$= \begin{bmatrix} \begin{pmatrix} [0,3] & [0,1] & [0,6] \\ [0,2] & [0,0] & [0,0] \end{pmatrix} \begin{bmatrix} [0,3] \\ [0,1] \\ [0,0] \end{bmatrix} & \begin{pmatrix} [0,3] & [0,1] & [0,6] \\ [0,2] & [0,0] & [0,0] \end{pmatrix} \begin{pmatrix} [0,5] & [0,2] & [0,0] \\ [0,1] & [0,2] & [0,2] \\ [0,3] & [0,1] & [0,0] \end{pmatrix} \\ \hline ([0,1][0,2][0,3]) \begin{bmatrix} [0,3] \\ [0,1] \\ [0,0] \end{bmatrix} & ([0,1][0,2][0,3]) \begin{pmatrix} [0,5] & [0,2] & [0,0] \\ [0,1] & [0,2] & [0,2] \\ [0,3] & [0,1] & [0,0] \end{pmatrix} \\ \hline \begin{pmatrix} [0,6] & [0,3] & [0,0] \\ [0,4] & [0,2] & [0,1] \\ [0,5] & [0,1] & [0,0] \end{pmatrix} \begin{bmatrix} [0,3] \\ [0,1] \\ [0,0] \end{bmatrix} & \begin{pmatrix} [0,6] & [0,3] & [0,0] \\ [0,4] & [0,2] & [0,1] \\ [0,5] & [0,1] & [0,0] \end{pmatrix} \begin{pmatrix} [0,5] & [0,2] & [0,0] \\ [0,1] & [0,2] & [0,2] \\ [0,3] & [0,1] & [0,0] \end{pmatrix} \end{bmatrix}$$

$$= \begin{bmatrix} [0,10] & [0,34] & [0,14] & [0,2] \\ [0,6] & [0,10] & [0,4] & [0,0] \\ \hline [0,5] & [0,16] & [0,9] & [0,4] \\ \hline [0,21] & [0,33] & [0,18] & [0,6] \\ [0,14] & [0,25] & [0,13] & [0,4] \\ [0,16] & [0,26] & [0,12] & [0,2] \end{bmatrix}.$$

Now we will find the product of a positive interval super matrix X and $X^t$.

***Example 2.41***: Let

$$X = \begin{bmatrix} [0,1] & [0,2] & [0,1] & [0,3] & [0,2] & [0,1] \\ [0,2] & [0,3] & [0,1] & [0,2] & [0,1] & [0,2] \\ \hline [0,1] & [0,4] & [0,2] & [0,3] & [0,2] & [0,2] \\ [0,4] & [0,1] & [0,3] & [0,2] & [0,1] & [0,1] \\ [0,2] & [0,3] & [0,2] & [0,3] & [0,2] & [0,3] \\ \hline [0,3] & [0,4] & [0,1] & [0,1] & [0,4] & [0,2] \\ [0,2] & [0,1] & [0,2] & [0,2] & [0,1] & [0,3] \end{bmatrix}$$

be a positive super interval matrix. We now find $X^t$ and then calculate $X^t.X$.

$$X^t = \begin{bmatrix} [0,1] & [0,2] & [0,1] & [0,4] & [0,2] & [0,3] & [0,2] \\ \hline [0,2] & [0,3] & [0,4] & [0,1] & [0,3] & [0,4] & [0,1] \\ [0,1] & [0,1] & [0,2] & [0,3] & [0,2] & [0,1] & [0,2] \\ \hline [0,3] & [0,2] & [0,3] & [0,2] & [0,3] & [0,1] & [0,2] \\ [0,2] & [0,1] & [0,2] & [0,1] & [0,2] & [0,4] & [0,1] \\ \hline [0,1] & [0,2] & [0,2] & [0,1] & [0,3] & [0,2] & [0,3] \end{bmatrix}$$

is the transpose of X. Now in $X^tX$ the first row of $X^t$ with that of the first column of X is given by

$$([0,1]\ [0,2]\ |\ [0,1]\ [0,4]\ [0,2]\ [0,3]\ |\ [0,2]) \begin{bmatrix} [0,1] \\ [0,2] \\ \hline [0,1] \\ [0,4] \\ [0,2] \\ [0,3] \\ \hline [0,2] \end{bmatrix}$$

$$= \quad ([0, 1]\ [0, 2]) \begin{pmatrix} [0,1] \\ [0,2] \end{pmatrix} + ([0, 1]\ [0, 4]\ [0, 2]\ [0, 3]) \begin{bmatrix} [0,1] \\ [0,4] \\ [0,2] \\ [0,3] \end{bmatrix}$$

$$\begin{aligned} & + [0, 2]\ [0, 2] \\ = \quad & [0, 5] + [0, 30] + [0, 4] \\ = \quad & [0, 39]. \end{aligned}$$

Product of first row of $X^t$ with $2^{nd}$ column of X.

$$([0, 1]\ [0, 2]\ |\ [0, 1]\ [0, 4]\ [0, 2]\ [0, 3]\ |\ [0, 2]) \begin{bmatrix} [0,2] & [0,1] \\ [0,3] & [0,1] \\ \overline{[0,4]} & \overline{[0,2]} \\ [0,1] & [0,3] \\ [0,3] & [0,2] \\ \overline{[0,4]} & \overline{[0,1]} \\ [0,1] & [0,2] \end{bmatrix}$$

$$= \quad ([0, 1]\ [0, 2]) \begin{pmatrix} [0,2] & [0,1] \\ [0,3] & [0,1] \end{pmatrix}$$

$$+ ([0, 1]\ [0, 4]\ [0, 2]\ [0, 3]) \begin{bmatrix} [0,4] & [0,2] \\ [0,1] & [0,3] \\ [0,3] & [0,2] \\ [0,4] & [0,1] \end{bmatrix}$$

$$+ [0, 2]\ ([0, 1]\ [0, 2])$$

$$\begin{aligned} = \quad & ([0, 8]\ [0, 3]) + ([0, 26]\ [0, 21]) + ([0, 2],\ [0, 4]) \\ = \quad & ([0, 36]\ [0, 28]). \end{aligned}$$

The product of first row with third column

$$([0, 1]\ [0, 2]\ |\ [0, 1]\ [0, 4]\ [0, 2]\ [0, 3]\ |\ [0, 2])\begin{bmatrix}[0,3] & [0,2] & [0,1]\\ [0,2] & [0,1] & [0,2]\\ \hline [0,3] & [0,2] & [0,2]\\ [0,2] & [0,1] & [0,1]\\ [0,3] & [0,2] & [0,3]\\ [0,1] & [0,4] & [0,2]\\ \hline [0,2] & [0,1] & [0,3]\end{bmatrix}$$

$$= ([0, 1]\ [0, 2])\begin{pmatrix}[0,3] & [0,2] & [0,1]\\ [0,2] & [0,1] & [0,2]\end{pmatrix}$$

$$+ ([0, 1]\ [0, 4]\ [0, 2]\ [0, 3])\begin{bmatrix}[0,3] & [0,2] & [0,2]\\ [0,2] & [0,1] & [0,1]\\ [0,3] & [0,2] & [0,3]\\ [0,1] & [0,4] & [0,2]\end{bmatrix}$$

$$+ [0, 2]\ ([0, 2]\ [0, 1]\ [0, 3])$$
$$= ([0, 7]\ [0, 4]\ [0, 5]) + ([0, 20]\ [0, 22]\ [0, 18]) +$$
$$([0, 4],\ [0, 2]\ [0, 6])$$
$$= ([0, 31]\ [0, 28]\ [0, 29]).$$

The product of second row of $X^t$ with first column of X.

$$\begin{pmatrix}[0,2] & [0,3] & | & [0,4] & [0,1] & [0,3] & [0,4] & | & [0,1]\\ [0,1] & [0,1] & | & [0,2] & [0,3] & [0,2] & [0,1] & | & [0,2]\end{pmatrix}\begin{bmatrix}[0,1]\\ [0,2]\\ \hline [0,1]\\ [0,4]\\ [0,2]\\ [0,3]\\ \hline [0,2]\end{bmatrix}$$

$$=\begin{pmatrix}[0,2] & [0,3]\\ [0,1] & [0,1]\end{pmatrix}\begin{bmatrix}[0,1]\\ [0,2]\end{bmatrix}+\begin{pmatrix}[0,4] & [0,1] & [0,3] & [0,4]\\ [0,2] & [0,3] & [0,2] & [0,1]\end{pmatrix}\begin{bmatrix}[0,1]\\ [0,4]\\ [0,2]\\ [0,3]\end{bmatrix}+$$

$$\begin{pmatrix} [0,1] \\ [0,2] \end{pmatrix}[0,2]$$

$$= \begin{bmatrix} [0,8] \\ [0,3] \end{bmatrix} + \begin{pmatrix} [0,26] \\ [0,21] \end{pmatrix} + \begin{bmatrix} [0,2] \\ [0,4] \end{bmatrix}$$

$$= \begin{pmatrix} [0,36] \\ [0,28] \end{pmatrix}.$$

The product of second row of $X^t$ with second column of X.

$$\begin{pmatrix} [0,2] & [0,3] & | & [0,4] & [0,1] & [0,3] & [0,4] & | & [0,1] \\ [0,1] & [0,1] & | & [0,2] & [0,3] & [0,2] & [0,1] & | & [0,2] \end{pmatrix} \begin{bmatrix} [0,2] & [0,1] \\ [0,3] & [0,1] \\ \overline{[0,4]} & [0,2] \\ [0,1] & [0,3] \\ [0,3] & [0,2] \\ \overline{[0,4]} & [0,1] \\ [0,1] & [0,2] \end{bmatrix}$$

$$= \begin{pmatrix} [0,2] & [0,3] \\ [0,1] & [0,1] \end{pmatrix} \begin{pmatrix} [0,2] & [0,1] \\ [0,3] & [0,1] \end{pmatrix} + \begin{pmatrix} [0,4] & [0,1] & [0,3] & [0,4] \\ [0,2] & [0,3] & [0,2] & [0,1] \end{pmatrix} \begin{bmatrix} [0,4] & [0,2] \\ [0,1] & [0,3] \\ [0,3] & [0,2] \\ [0,4] & [0,1] \end{bmatrix} +$$

$$\begin{pmatrix} [0,1] \\ [0,2] \end{pmatrix}([0,1][0,2])$$

$$= \begin{pmatrix} [0,13] & [0,5] \\ [0,5] & [0,2] \end{pmatrix} + \begin{pmatrix} [0,42] & [0,21] \\ [0,21] & [0,18] \end{pmatrix} + \begin{pmatrix} [0,1] & [0,2] \\ [0,2] & [0,4] \end{pmatrix}$$

$$= \begin{pmatrix} [0,56] & [0,28] \\ [0,28] & [0,24] \end{pmatrix}.$$

The product of 2$^{nd}$ row of X$^t$ with 3$^{rd}$ column of X.

$$\begin{pmatrix} [0,2] & [0,3] \mid [0,4] & [0,1] & [0,3] & [0,4] \mid [0,1] \\ [0,1] & [0,1] \mid [0,2] & [0,3] & [0,2] & [0,1] \mid [0,2] \end{pmatrix} \begin{bmatrix} [0,3] & [0,2] & [0,1] \\ [0,2] & [0,1] & [0,2] \\ \hline [0,3] & [0,2] & [0,2] \\ [0,2] & [0,1] & [0,1] \\ [0,3] & [0,2] & [0,3] \\ [0,1] & [0,4] & [0,2] \\ \hline [0,2] & [0,1] & [0,3] \end{bmatrix}$$

$$= \begin{pmatrix} [0,2] & [0,3] \\ [0,1] & [0,1] \end{pmatrix} \begin{pmatrix} [0,3] & [0,2] & [0,1] \\ [0,2] & [0,1] & [0,2] \end{pmatrix} + \begin{pmatrix} [0,4] & [0,1] & [0,3] & [0,4] \\ [0,2] & [0,3] & [0,2] & [0,1] \end{pmatrix}.$$

$$\begin{bmatrix} [0,3] & [0,2] & [0,2] \\ [0,2] & [0,1] & [0,1] \\ [0,3] & [0,2] & [0,3] \\ [0,1] & [0,4] & [0,2] \end{bmatrix} + \begin{pmatrix} [0,1] \\ [0,2] \end{pmatrix} ([0,2][0,1][0,3])$$

$$= \begin{pmatrix} [0,12] & [0,7] & [0,8] \\ [0,5] & [0,3] & [0,3] \end{pmatrix} + \begin{pmatrix} [0,27] & [0,31] & [0,26] \\ [0,19] & [0,15] & [0,15] \end{pmatrix} + \begin{pmatrix} [0,2] & [0,1] & [0,3] \\ [0,4] & [0,2] & [0,6] \end{pmatrix}$$

$$= \begin{pmatrix} [0,41] & [0,39] & [0,37] \\ [0,28] & [0,20] & [0,24] \end{pmatrix}.$$

The product of 3$^{rd}$ row of X with first column of X.

$$\begin{pmatrix} [0,3] & [0,2] \mid [0,3] & [0,2] & [0,3] & [0,1] \mid [0,2] \\ [0,2] & [0,1] \mid [0,2] & [0,1] & [0,2] & [0,4] \mid [0,1] \\ [0,1] & [0,2] \mid [0,2] & [0,1] & [0,3] & [0,2] \mid [0,3] \end{pmatrix} \begin{bmatrix} [0,1] \\ [0,2] \\ \hline [0,1] \\ [0,4] \\ [0,2] \\ [0,3] \\ \hline [0,2] \end{bmatrix}$$

$$= \begin{bmatrix} [0,7] \\ [0,4] \\ [0,5] \end{bmatrix} + \begin{bmatrix} [0,20] \\ [0,22] \\ [0,18] \end{bmatrix} + \begin{bmatrix} [0,4] \\ [0,2] \\ [0,6] \end{bmatrix}.$$

$$= \begin{bmatrix} [0,31] \\ [0,28] \\ [0,29] \end{bmatrix}$$

The product of third row of $X^t$ with second column of X.

$$\begin{pmatrix} [0,3] & [0,2] & | & [0,3] & [0,2] & [0,3] & [0,1] & | & [0,2] \\ [0,2] & [0,1] & | & [0,2] & [0,1] & [0,2] & [0,4] & | & [0,1] \\ [0,1] & [0,2] & | & [0,2] & [0,1] & [0,3] & [0,2] & | & [0,3] \end{pmatrix} \begin{bmatrix} [0,2] & [0,1] \\ [0,3] & [0,1] \\ \hline [0,4] & [0,2] \\ [0,1] & [0,3] \\ [0,3] & [0,2] \\ [0,4] & [0,1] \\ \hline [0,1] & [0,2] \end{bmatrix}$$

$$= \begin{bmatrix} [0,12] & [0,5] \\ [0,7] & [0,3] \\ [0,8] & [0,3] \end{bmatrix} + \begin{bmatrix} [0,27] & [0,19] \\ [0,31] & [0,15] \\ [0,26] & [0,15] \end{bmatrix} + \begin{bmatrix} [0,2] & [0,4] \\ [0,1] & [0,2] \\ [0,3] & [0,6] \end{bmatrix}$$

$$= \begin{bmatrix} [0,41] & [0,28] \\ [0,39] & [0,20] \\ [0,37] & [0,24] \end{bmatrix}.$$

The product of $3^{rd}$ row of $X^t$ with third column of X.

$$\begin{pmatrix} [0,3] & [0,2] & | & [0,3] & [0,2] & [0,3] & [0,1] & | & [0,2] \\ [0,2] & [0,1] & | & [0,2] & [0,1] & [0,2] & [0,4] & | & [0,1] \\ [0,1] & [0,2] & | & [0,2] & [0,1] & [0,3] & [0,2] & | & [0,3] \end{pmatrix} \times$$

$$\begin{bmatrix} [0,3] & [0,2] & [0,1] \\ [0,2] & [0,1] & [0,2] \\ \hline [0,3] & [0,2] & [0,2] \\ [0,2] & [0,1] & [0,1] \\ [0,3] & [0,2] & [0,3] \\ [0,1] & [0,4] & [0,2] \\ \hline [0,2] & [0,1] & [0,3] \end{bmatrix}$$

$$= \begin{bmatrix} [0,13] & [0,8] & [0,7] \\ [0,8] & [0,5] & [0,4] \\ [0,7] & [0,4] & [0,5] \end{bmatrix} + \begin{bmatrix} [0,23] & [0,18] & [0,19] \\ [0,18] & [0,25] & [0,19] \\ [0,19] & [0,19] & [0,18] \end{bmatrix} +$$

$$\begin{bmatrix} [0,4] & [0,2] & [0,6] \\ [0,2] & [0,1] & [0,3] \\ [0,6] & [0,3] & [0,9] \end{bmatrix}$$

$$= \begin{bmatrix} [0,40] & [0,28] & [0,32] \\ [0,28] & [0,31] & [0,26] \\ [0,32] & [0,26] & [0,32] \end{bmatrix}.$$

Thus

$$X^t X = \begin{bmatrix} [0,39] & [0,36] & [0,28] & [0,31] & [0,28] & [0,29] \\ [0,36] & [0,56] & [0,28] & [0,41] & [0,39] & [0,37] \\ [0,28] & [0,28] & [0,24] & [0,28] & [0,20] & [0,24] \\ \hline [0,31] & [0,41] & [0,28] & [0,40] & [0,28] & [0,32] \\ [0,28] & [0,39] & [0,20] & [0,28] & [0,31] & [0,26] \\ [0,29] & [0,37] & [0,24] & [0,32] & [0,26] & [0,32] \end{bmatrix}.$$

Now we give the product of X and Y where X and Y are positive super interval matrices.

*Example 2.42*: Let

$$X = \begin{bmatrix} [0,1] & [0,2] & [0,1] & [0,1] & [0,2] & [0,3] \\ [0,3] & [0,1] & [0,2] & [0,3] & [0,1] & [0,1] \\ \hline [0,1] & [0,1] & [0,3] & [0,1] & [0,1] & [0,1] \\ [0,2] & [0,3] & [0,1] & [0,2] & [0,0] & [0,1] \\ [0,3] & [0,4] & [0,2] & [0,0] & [0,1] & [0,0] \\ [0,4] & [0,2] & [0,4] & [0,1] & [0,0] & [0,0] \\ \hline [0,5] & [0,0] & [0,1] & [0,1] & [0,1] & [0,1] \end{bmatrix}$$

and

$$Y = \begin{bmatrix} [0,1] & [0,1] & [0,2] & [0,1] & [0,2] & [0,1] & [0,3] & [0,1] & [0,0] \\ [0,1] & [0,0] & [0,2] & [0,4] & [0,3] & [0,1] & [0,4] & [0,1] & [0,1] \\ \hline [0,0] & [0,1] & [0,0] & [0,3] & [0,1] & [0,0] & [0,1] & [0,2] & [0,1] \\ \hline [0,1] & [0,1] & [0,0] & [0,0] & [0,2] & [0,1] & [0,1] & [0,2] & [0,1] \\ [0,1] & [0,0] & [0,1] & [0,1] & [0,1] & [0,2] & [0,2] & [0,1] & [0,2] \\ [0,0] & [0,1] & [0,0] & [0,1] & [0,1] & [0,1] & [0,1] & [0,0] & [0,0] \end{bmatrix}$$

be two positive super interval matrices. Now we find the major product of X and Y.

The product of first row of X with first column of Y gives

$$\begin{bmatrix} [0,1] & [0,2] & [0,1] & [0,1] & [0,2] & [0,3] \\ [0,3] & [0,1] & [0,2] & [0,3] & [0,1] & [0,1] \end{bmatrix} \begin{bmatrix} [0,1] & [0,1] & [0,2] & [0,1] \\ [0,1] & [0,0] & [0,2] & [0,4] \\ [0,0] & [0,1] & [0,0] & [0,3] \\ [0,1] & [0,1] & [0,0] & [0,0] \\ [0,1] & [0,0] & [0,1] & [0,1] \\ [0,0] & [0,1] & [0,0] & [0,1] \end{bmatrix}$$

$$= \begin{bmatrix} [0,1] \\ [0,3] \end{bmatrix} ([0,1]\ [0,1]\ [0,2]\ [0,1])$$

$$
\begin{aligned}
&+ \begin{bmatrix} [0,2] & [0,1] \\ [0,1] & [0,2] \end{bmatrix} \begin{bmatrix} [0,1] & [0,0] & [0,2] & [0,4] \\ [0,0] & [0,1] & [0,0] & [0,3] \end{bmatrix} \\
&+ \begin{bmatrix} [0,1] & [0,2] & [0,3] \\ [0,3] & [0,1] & [0,1] \end{bmatrix} \begin{bmatrix} [0,1] & [0,1] & [0,0] & [0,0] \\ [0,1] & [0,0] & [0,1] & [0,1] \\ [0,0] & [0,1] & [0,0] & [0,1] \end{bmatrix}
\end{aligned}
$$

$$
= \begin{bmatrix} [0,1] & [0,1] & [0,2] & [0,1] \\ [0,3] & [0,3] & [0,6] & [0,3] \end{bmatrix} + \begin{bmatrix} [0,2] & [0,1] & [0,4] & [0,11] \\ [0,1] & [0,2] & [0,2] & [0,10] \end{bmatrix} +
$$

$$
\begin{bmatrix} [0,3] & [0,4] & [0,2] & [0,5] \\ [0,4] & [0,4] & [0,1] & [0,2] \end{bmatrix}
$$

$$
= \begin{bmatrix} [0,6] & [0,6] & [0,8] & [0,17] \\ [0,8] & [0,9] & [0,9] & [0,15] \end{bmatrix}.
$$

Now

$$
\begin{bmatrix} [0,1] & [0,2] & [0,1] & [0,1] & [0,2] & [0,3] \\ [0,3] & [0,1] & [0,2] & [0,3] & [0,1] & [0,1] \end{bmatrix} \begin{bmatrix} [0,2] & [0,1] \\ [0,3] & [0,1] \\ [0,1] & [0,0] \\ [0,2] & [0,1] \\ [0,1] & [0,2] \\ [0,1] & [0,1] \end{bmatrix}
$$

$$
= \begin{pmatrix} [0,1] \\ [0,3] \end{pmatrix} \begin{pmatrix} [0,2] & [0,1] \end{pmatrix} + \begin{pmatrix} [0,2] & [0,1] \\ [0,1] & [0,2] \end{pmatrix} \begin{pmatrix} [0,3] & [0,1] \\ [0,1] & [0,0] \end{pmatrix} +
$$

$$
\begin{pmatrix} [0,1] & [0,2] & [0,3] \\ [0,3] & [0,1] & [0,1] \end{pmatrix} \begin{pmatrix} [0,2] & [0,1] \\ [0,1] & [0,2] \\ [0,1] & [0,1] \end{pmatrix}
$$

$$
= \begin{bmatrix} [0,2] & [0,1] \\ [0,6] & [0,3] \end{bmatrix} + \begin{bmatrix} [0,7] & [0,2] \\ [0,5] & [0,1] \end{bmatrix} + \begin{bmatrix} [0,7] & [0,8] \\ [0,8] & [0,6] \end{bmatrix}
$$

$$
= \begin{bmatrix} [0,16] & [0,11] \\ [0,19] & [0,10] \end{bmatrix}.
$$

Consider the product of first row with third column

$$\begin{bmatrix} [0,1] & [0,2] & [0,1] & [0,1] & [0,2] & [0,3] \\ [0,3] & [0,1] & [0,2] & [0,3] & [0,1] & [0,1] \end{bmatrix} \begin{bmatrix} [0,3] & [0,1] & [0,0] \\ [0,4] & [0,1] & [0,1] \\ [0,1] & [0,2] & [0,1] \\ [0,1] & [0,2] & [0,1] \\ [0,2] & [0,1] & [0,2] \\ [0,1] & [0,0] & [0,0] \end{bmatrix}$$

$$= \begin{bmatrix} [0,1] \\ [0,3] \end{bmatrix}([0,3][0,1][0,0]) + \begin{bmatrix} [0,2] & [0,1] \\ [0,1] & [0,2] \end{bmatrix}\begin{bmatrix} [0,4] & [0,1] & [0,1] \\ [0,1] & [0,2] & [0,1] \end{bmatrix}$$

$$+ \begin{pmatrix} [0,1] & [0,2] & [0,3] \\ [0,3] & [0,1] & [0,1] \end{pmatrix}\begin{bmatrix} [0,1] & [0,2] & [0,1] \\ [0,2] & [0,1] & [0,2] \\ [0,1] & [0,0] & [0,0] \end{bmatrix}$$

$$= \begin{pmatrix} [0,3] & [0,1] & [0,0] \\ [0,9] & [0,3] & [0,0] \end{pmatrix} + \begin{pmatrix} [0,9] & [0,4] & [0,3] \\ [0,6] & [0,5] & [0,3] \end{pmatrix} + \begin{pmatrix} [0,8] & [0,4] & [0,5] \\ [0,6] & [0,7] & [0,5] \end{pmatrix}$$

$$= \begin{pmatrix} [0,20] & [0,9] & [0,8] \\ [0,21] & [0,15] & [0,8] \end{pmatrix}.$$

The product of second row of X with the second column of Y

$$\begin{bmatrix} [0,1] & [0,1] & [0,3] & [0,1] & [0,1] & [0,1] \\ [0,2] & [0,3] & [0,1] & [0,2] & [0,0] & [0,1] \\ [0,3] & [0,4] & [0,2] & [0,0] & [0,1] & [0,0] \\ [0,4] & [0,2] & [0,4] & [0,1] & [0,0] & [0,0] \end{bmatrix} \begin{bmatrix} [0,2] & [0,1] \\ [0,3] & [0,1] \\ [0,1] & [0,0] \\ [0,2] & [0,1] \\ [0,1] & [0,2] \\ [0,1] & [0,1] \end{bmatrix}$$

$$= \begin{bmatrix} [0,1] \\ [0,2] \\ [0,3] \\ [0,4] \end{bmatrix} ([0,2]\ [0,1]) + \begin{bmatrix} [0,1] & [0,3] \\ [0,3] & [0,1] \\ [0,4] & [0,2] \\ [0,2] & [0,4] \end{bmatrix} \begin{bmatrix} [0,3] & [0,1] \\ [0,1] & [0,0] \end{bmatrix} +$$

$$\begin{bmatrix} [0,1] & [0,1] & [0,1] \\ [0,2] & [0,0] & [0,1] \\ [0,0] & [0,1] & [0,0] \\ [0,1] & [0,0] & [0,0] \end{bmatrix} \begin{bmatrix} [0,2] & [0,1] \\ [0,1] & [0,2] \\ [0,1] & [0,1] \end{bmatrix}$$

$$= \begin{bmatrix} [0,2] & [0,1] \\ [0,4] & [0,2] \\ [0,6] & [0,3] \\ [0,8] & [0,4] \end{bmatrix} + \begin{bmatrix} [0,6] & [0,1] \\ [0,10] & [0,3] \\ [0,14] & [0,4] \\ [0,10] & [0,2] \end{bmatrix} + \begin{bmatrix} [0,4] & [0,4] \\ [0,5] & [0,3] \\ [0,1] & [0,2] \\ [0,2] & [0,1] \end{bmatrix}$$

$$= \begin{bmatrix} [0,12] & [0,6] \\ [0,19] & [0,8] \\ [0,21] & [0,9] \\ [0,20] & [0,7] \end{bmatrix}.$$

The product of the 2$^{nd}$ row with last column of Y.

$$\begin{bmatrix} [0,1] & [0,1] & [0,3] & [0,1] & [0,1] & [0,1] \\ [0,2] & [0,3] & [0,1] & [0,2] & [0,0] & [0,1] \\ [0,3] & [0,4] & [0,2] & [0,0] & [0,1] & [0,0] \\ [0,4] & [0,2] & [0,4] & [0,1] & [0,0] & [0,0] \end{bmatrix} \begin{bmatrix} [0,3] & [0,1] & [0,0] \\ [0,4] & [0,1] & [0,1] \\ [0,1] & [0,2] & [0,1] \\ [0,1] & [0,2] & [0,1] \\ [0,2] & [0,1] & [0,2] \\ [0,1] & [0,0] & [0,0] \end{bmatrix}$$

$$= \begin{bmatrix}[0,1]\\[0,2]\\[0,3]\\[0,4]\end{bmatrix}([0,3]\ [0,1]\ [0,0]) + \begin{bmatrix}[0,1] & [0,3]\\ [0,3] & [0,1]\\ [0,4] & [0,2]\\ [0,2] & [0,4]\end{bmatrix}\begin{bmatrix}[0,4] & [0,1] & [0,1]\\ [0,1] & [0,2] & [0,1]\end{bmatrix}$$

$$+ \begin{pmatrix}[0,1] & [0,1] & [0,1]\\ [0,2] & [0,0] & [0,1]\\ [0,0] & [0,1] & [0,0]\\ [0,1] & [0,0] & [0,0]\end{pmatrix}\begin{pmatrix}[0,1] & [0,2] & [0,1]\\ [0,2] & [0,1] & [0,2]\\ [0,1] & [0,0] & [0,0]\end{pmatrix}$$

$$= \begin{bmatrix}[0,3] & [0,1] & [0,0]\\ [0,6] & [0,2] & [0,0]\\ [0,9] & [0,3] & [0,0]\\ [0,12] & [0,4] & [0,0]\end{bmatrix} + \begin{bmatrix}[0,7] & [0,7] & [0,4]\\ [0,13] & [0,5] & [0,4]\\ [0,18] & [0,8] & [0,6]\\ [0,12] & [0,10] & [0,6]\end{bmatrix} + \begin{bmatrix}[0,4] & [0,3] & [0,3]\\ [0,3] & [0,4] & [0,2]\\ [0,2] & [0,1] & [0,2]\\ [0,1] & [0,2] & [0,1]\end{bmatrix}$$

$$= \begin{bmatrix}[0,14] & [0,11] & [0,7]\\ [0,22] & [0,11] & [0,6]\\ [0,29] & [0,12] & [0,8]\\ [0,25] & [0,16] & [0,7]\end{bmatrix}.$$

The product of 3$^{\text{rd}}$ row of X with first column of Y

$$([0,5]\ |\ [0,0]\ [0,1]\ |\ [0,1]\ [0,1]\ [0,1])\begin{vmatrix}[0,1] & [0,1] & [0,2] & [0,1]\\ [0,1] & [0,0] & [0,2] & [0,4]\\ [0,0] & [0,1] & [0,0] & [0,3]\\ [0,1] & [0,1] & [0,0] & [0,0]\\ [0,1] & [0,0] & [0,1] & [0,1]\\ [0,0] & [0,1] & [0,0] & [0,1]\end{vmatrix}$$

$$= [0,5]\,([0,1]\ [0,1]\ [0,2]\ [0,1]) +$$
$$([0,0]\ [0,1])\begin{pmatrix}[0,1] & [0,0] & [0,2] & [0,4]\\ [0,0] & [0,1] & [0,0] & [0,3]\end{pmatrix} +$$

$$([0, 1] \; [0, 1] \; [0, 1]) \begin{bmatrix} [0,1] & [0,1] & [0,0] & [0,0] \\ [0,1] & [0,0] & [0,1] & [0,1] \\ [0,0] & [0,1] & [0,0] & [0,1] \end{bmatrix}$$

$$= \quad ([0, 5] \; [0, 5] \; [0, 10] \; [0, 5]) + ([0, 0] \; [0, 1] \; [0, 0] \; [0, 3]) + \\ ([0, 2] \; [0, 2] \; [0, 1] \; [0, 2])$$

$$= \quad ([0, 7] \; [0, 8] \; [0, 11] \; [0, 10]).$$

The product of $3^{\text{rd}}$ row of X with the second column of Y.

$$([0, 5] \mid [0, 0] \; [0, 1] \mid [0, 1] \; [0, 1] \; [0, 1]) \begin{bmatrix} [0,2] & [0,1] \\ \overline{[0,3]} & \overline{[0,1]} \\ [0,1] & [0,0] \\ \overline{[0,2]} & \overline{[0,1]} \\ [0,1] & [0,2] \\ [0,1] & [0,1] \end{bmatrix}$$

$$= \quad [0, 5] \, ([0, 2] \; [0, 1])$$

$$+ \; ([0, 0] \; [0, 1]) \begin{pmatrix} [0,3] & [0,1] \\ [0,1] & [0,0] \end{pmatrix}$$

$$+ \; ([0, 1] \; [0, 1] \; [0, 1]) \begin{bmatrix} [0,2] & [0,1] \\ [0,1] & [0,2] \\ [0,1] & [0,1] \end{bmatrix}$$

$$= \quad ([0, 10] \; [0, 5]) + ([0, 1] \; [0, 0]) + ([0, 4] \; [0, 4])$$

$$= \quad ([0, 15] \; [0, 9]).$$

Likewise second row of X with first column of Y is found.

Thus we get XY =

$$\begin{bmatrix} [0,6] & [0,6] & [0,8] & [0,17] & [0,16] & [0,11] & [0,20] & [0,9] & [0,8] \\ [0,8] & [0,8] & [0,9] & [0,15] & [0,19] & [0,10] & [0,21] & [0,15] & [0,8] \\ [0,4] & [0,6] & [0,5] & [0,16] & [0,12] & [0,6] & [0,14] & [0,11] & [0,7] \\ [0,7] & [0,6] & [0,10] & [0,18] & [0,19] & [0,8] & [0,22] & [0,11] & [0,6] \\ [0,8] & [0,5] & [0,15] & [0,26] & [0,21] & [0,9] & [0,29] & [0,12] & [0,8] \\ [0,7] & [0,9] & [0,12] & [0,24] & [0,20] & [0,7] & [0,25] & [0,16] & [0,7] \\ [0,7] & [0,8] & [0,11] & [0,10] & [0,15] & [0,9] & [0,20] & [0,10] & [0,4] \end{bmatrix}$$

Now having seen how the products are carried out in positive interval super matrices now we proceed onto give the properties enjoyed by this collection.

Recall a positive super interval matrix is of same type
1. They should be of same order.
2. They should be partitioned in the same way, that is they are identically partitioned.

When we say a m × n matrix it can be a positive super row vector (matrix) if it is partitioned only vertically (m=1) or any positive super interval matrix. Thus ([0, 2] | [0, 7] [0, 12] [0, 16] | [0, 7] | [0, 5] [0, 4]) is a 1 × 7 positive super row interval matrix where as

$$\begin{pmatrix} [0,2] & [0,1] & [0,11] & [0,0] & [0,9] & [0,1] \\ [0,9] & [0,2] & [0,12] & [0,5] & [0,7] & [0,2] \\ [0,4] & [0,10] & [0,3] & [0,7] & [0,5] & [0,11] \end{pmatrix}$$

is a positive super interval row vector.
Likewise

$$\begin{bmatrix} [0,3] \\ [0,2] \\ [0,14] \\ \hline [0,5] \\ [0,17] \\ [0,3] \\ \hline [0,0] \\ [0,9] \\ [0,12] \end{bmatrix}$$

is a positive super interval column matrix where as

$$\begin{bmatrix} [0,7] & [0,9] & [0,0] & [0,2] \\ [0,1] & [0,4] & [0,8] & [0,0] \\ \hline [0,9] & [0,11] & [0,1] & [0,5] \\ [0,10] & [0,1] & [0,13] & [0,4] \\ [0,12] & [0,14] & [0,1] & [0,0] \\ \hline [0,1] & [0,2] & [0,3] & [0,4] \end{bmatrix}$$

positive super interval column vector for in this 6 × 4 positive interval matrix only horizontal partition is carried out and there is not even a single vertical partition.

Now

$$\begin{bmatrix} [0,2] & [0,7] & [0,8] & [0,9] & [0,4] & [0,1] \\ \hline [0,0] & [0,2] & [0,3] & [0,4] & [0,5] & [0,6] \\ [0,1] & [0,3] & [0,10] & [0,1] & [0,0] & [0,11] \\ [0,9] & [0,0] & [0,0] & [0,2] & [0,3] & [0,0] \\ \hline [0,11] & [0,4] & [0,4] & [0,4] & [0,0] & [0,1] \\ [0,16] & [0,0] & [0,2] & [0,0] & [0,2] & [0,3] \end{bmatrix}$$

is a positive super interval square matrix.

We see

$$\begin{bmatrix} [0,1] & [0,2] & [0,3] & [0,4] & [0,5] & [0,6] \\ [0,7] & [0,8] & [0,9] & [0,0] & [0,1] & [0,2] \\ [0,3] & [0,4] & [0,5] & [0,6] & [0,7] & [0,8] \\ \hline [0,9] & [0,11] & [0,1] & [0,3] & [0,5] & [0,6] \\ [0,10] & [0,0] & [0,2] & [0,4] & [0,7] & [0,0] \end{bmatrix}$$

is a positive super interval rectangular matrix or a positive super 5 × 6 interval matrix.

**THEOREM 2.3**: *Let A = {collection of all m × n positive interval super matrices of same type} with entries from $R^+ \cup \{0\}$ or $Z^+ \cup \{0\}$ or $Q^+ \cup \{0\}$. A is a semigroup under addition.*

The proof is direct and hence is left as an exercise to the reader.

**THEOREM 2.4**: *Let M = {collection of all m × n positive super interval matrices of same type} with entries from $Z_n$ (n < ∞). M is a group under addition.*

*Example 2.43*: Let

$$A = \left\{ \begin{bmatrix} [0,a_1] & [0,a_2] & [0,a_3] \\ [0,a_4] & [0,a_5] & [0,a_6] \\ \hline [0,a_7] & [0,a_8] & [0,a_9] \\ [0,a_{10}] & [0,a_{11}] & [0,a_{12}] \end{bmatrix} \right.$$

where $a_i \in Z_{43}$; $1 \le i \le 12$} be the collection of all positive interval super matrices. A is a group under addition. Of course A is not even compatible with respect to multiplication.

*Example 2.44*: Let

$$A = \left\{ \begin{bmatrix} [0,a_1] & [0,a_2] & [0,a_3] & [0,a_4] \\ [0,a_5] & [0,a_6] & [0,a_7] & [0,a_8] \\ [0,a_9] & [0,a_{10}] & [0,a_{11}] & [0,a_{12}] \end{bmatrix} \right.$$

$a_i \in Z^+ \cup \{0\}$, $1 \le i \le 12$} be the collection of all positive interval vectors. A is a semigroup under addition. Of course multiplication cannot be defined on A.

*Example 2.45*: Let

$$P = \left\{ \begin{bmatrix} [0,a_1] & [0,a_2] & [0,a_3] \\ [0,a_4] & [0,a_5] & [0,a_6] \\ \hline [0,a_7] & [0,a_8] & [0,a_9] \\ [0,a_{10}] & [0,a_{11}] & [0,a_{12}] \\ [0,a_{13}] & [0,a_{14}] & [0,a_{15}] \\ [0,a_{16}] & [0,a_{17}] & [0,a_{18}] \end{bmatrix} \right.$$

$a_i \in Q^+ \cup \{0\}$, $1 \leq i \leq 8\}$ be the set of all positive super interval column vectors. P is a semigroup under addition clearly P is not a group as inverse do not exist.

*Example 2.46*: Let

$$S = \left\{ \left[ \begin{array}{cc|c} [0,a_1] & [0,a_2] & [0,a_3] \\ \hline [0,a_4] & [0,a_5] & [0,a_6] \\ [0,a_7] & [0,a_8] & [0,a_9] \end{array} \right] \right.$$

$a_i \in Z_{480}$, $1 \leq i \leq 9\}$ be a positive super interval matrix. S is a group under addition and S is not in general compatible with respect to multiplication. To define multiplication we need to consider only those positive interval super matrices such that $X = X^t$ but how much such study or definition is practical is to be analysed.

*Example 2.47*: Let

$$G = \left\{ \left[ \begin{array}{cc} [0,a_1] & [0,a_2] \\ [0,a_3] & [0,a_4] \\ [0,a_5] & [0,a_6] \\ \hline [0,a_7] & [0,a_8] \\ [0,a_9] & [0,a_{10}] \\ [0,a_{11}] & [0,a_{12}] \end{array} \right] \right| a_i \in Z_8, 1 \leq i \leq 12 \right\}$$

be a positive super interval column vector G is a group under addition. G is of finite order and G is abelian.
Also

$$H = \left\{ \left[ \begin{array}{cc} [0,a] & [0,a] \\ [0,a] & [0,a] \\ [0,a] & [0,a] \\ \hline [0,a] & [0,a] \\ [0,a] & [0,a] \\ [0,a] & [0,a] \end{array} \right] \right| a \in Z_8 \right\} \subseteq G$$

is a super interval subgroup of G of order 8.

Infact G has several such super interval subgroups.

$$T = \left\{ \left[ \begin{array}{cc} [0,a] & 0 \\ [0,a] & 0 \\ [0,a] & 0 \\ \hline 0 & [0,a] \\ \hline [0,a] & 0 \\ 0 & [0,a] \end{array} \right] \;\middle|\; a \in Z_8 \right\}$$

is again a interval super subgroup of G of order 8. However T and H are in no way related. By varying the portions of the super interval matrix

$$P = \begin{bmatrix} [0,a_1] & [0,a_2] \\ [0,a_3] & [0,a_4] \\ [0,a_5] & [0,a_6] \\ [0,a_7] & [0,a_8] \\ [0,a_9] & [0,a_{10}] \\ [0,a_{11}] & [0,a_{12}] \end{bmatrix}$$

and getting the collection of same type of super interval matrices we get several super interval matrix groups of same order built using $Z_8$. This is one of the convenient method of getting many groups of same order but distinctly different in their structure.

Further in case of a super interval matrix group we can form the transpose. Thus for every super interval matrix group we have a unique super matrix interval group associated with it, which is the transpose of it. Only under addition the transpose of a super row interval matrix will be a group.

We will call the transpose of a super row interval matrix group as transpose group super row interval matrix or a super row interval matrix group transpose. From the context one can easily find the underlying structure.

We will only illustrate this situation by some examples.

***Example 2.48***: Let

$$G = \left\{ \begin{bmatrix} [0,a_1] \\ \overline{[0,a_2]} \\ [0,a_3] \\ \overline{[0,a_4]} \\ [0,a_5] \\ [0,a_6] \end{bmatrix} \middle| \; a_i \in Z_{42}, 1 \le i \le 6 \right\}$$

be super row interval matrices of the same type, G under addition is a group. Now $G^T = ([0, a_1] \mid [0, a_2] \; [0, a_3] \; [0, a_4] \mid [0, a_5] \; [0, a_6]) \mid a_i \in Z_{42}$; $1 \le i \le 6 \}$ be the transpose group of G.

$G^T$ is also group known as the transpose super interval super matrix group.

***Example 2.49***: Let

$$H = \left\{ \begin{bmatrix} [0,a_1] & [0,a_2] \\ [0,a_3] & [0,a_4] \\ [0,a_5] & [0,a_6] \\ \overline{[0,a_7]} & \overline{[0,a_8]} \\ [0,a_9] & [0,a_{10}] \\ [0,a_{11}] & [0,a_{12}] \\ [0,a_{13}] & [0,a_{14}] \\ [0,a_{15}] & [0,a_{16}] \end{bmatrix} \middle| \; a_i \in Z_9; 1 \le i \le 16 \right\}$$

be the super interval row matrix group under addition of the specified type of super interval row matrices described in H.
Now

$$H^T = \begin{pmatrix} [0,a_1] & [0,a_3] & [0,a_5] & [0,a_7] & [0,a_9] & [0,a_{11}] & [0,a_{13}] & [0,a_{15}] \\ [0,a_2] & [0,a_4] & [0,a_6] & [0,a_8] & [0,a_{10}] & [0,a_{12}] & [0,a_{14}] & [0,a_{16}] \end{pmatrix}$$

where $a_i \in Z_9$ ; $1 \le i \le 16\}$ is the transpose group of H called the transpose super interval row vector matrix group.

Now we give some more examples.

*Example 2.50*: Let

$$P = \left\{ \left[ \begin{array}{cc|c} [0,a_1] & [0,a_2] & [0,a_3] \\ [0,a_4] & [0,a_5] & [0,a_6] \\ \hline [0,a_7] & [0,a_8] & [0,a_9] \\ [0,a_{10}] & [0,a_{11}] & [0,a_{12}] \end{array} \right] \,\middle|\, a_i \in Z_{142},\, 1 \leq i \leq 12 \right\}$$

be a collection of super interval matrix group of the same type. P is a group.
Now

$$P^T = \left\{ \left[ \begin{array}{cc|cc} [0,a_1] & [0,a_4] & [0,a_7] & [0,a_{10}] \\ [0,a_2] & [0,a_5] & [0,a_8] & [0,a_{11}] \\ \hline [0,a_3] & [0,a_6] & [0,a_9] & [0,a_{12}] \end{array} \right] \,\middle|\, a_i \in Z_{142},\, 1 \leq i \leq 12 \right\}$$

is the set of all transpose of the super interval matrices of the same type in P, Thus $P^T$ is a transpose group of super interval matrices of the group P. $P^T$ is also a group under addition of super interval matrices.

*Example 2.51*: Let

$$G = \left\{ \left[ \begin{array}{c|ccc|cc} [0,a_1] & [0,a_2] & [0,a_3] & [0,a_4] & [0,a_5] \\ [0,a_6] & [0,a_7] & [0,a_8] & [0,a_9] & [0,a_{10}] \\ \hline [0,a_{11}] & [0,a_{12}] & [0,a_{13}] & [0,a_{14}] & [0,a_{15}] \\ [0,a_{16}] & [0,a_{17}] & [0,a_{18}] & [0,a_{19}] & [0,a_{20}] \\ \hline [0,a_{21}] & [0,a_{22}] & [0,a_{23}] & [0,a_{24}] & [0,a_{25}] \end{array} \right] \,\middle|\, a_i \in Z_{20},\, 1 \leq i \leq 25 \right\}$$

be the collection of all $5 \times 5$ super interval matrices of this type described, G is a group called the super interval $5 \times 5$ matrix group under addition. Clearly G is not a group under multiplication.
Consider

$$G^T = \left\{ \begin{bmatrix} [0,a_1] & [0,a_6] & [0,a_{11}] & [0,a_{16}] & [0,a_{21}] \\ [0,a_2] & [0,a_7] & [0,a_{12}] & [0,a_{17}] & [0,a_{22}] \\ [0,a_3] & [0,a_8] & [0,a_{13}] & [0,a_{18}] & [0,a_{23}] \\ [0,a_4] & [0,a_9] & [0,a_{14}] & [0,a_{19}] & [0,a_{24}] \\ \hline [0,a_5] & [0,a_{10}] & [0,a_{15}] & [0,a_{20}] & [0,a_{25}] \end{bmatrix} \right.$$

where $a_i \in Z_{20}$, $1 \le i \le 25$} is the transpose group of the super interval square matrix under G addition.

Now having seen examples of super interval matrix groups under addition built using $Z_n$, $n < \infty$. We can proceed onto construct super interval matrix semigroups using $Z^+ \cup \{0\}$ or $R^+ \cup \{0\}$ or $Q^+ \cup \{0\}$ or $C^+ \cup \{0\}$.

We see if super interval matrices are built using $Z^+ \cup \{0\}$ or $R^+ \cup \{0\}$ or $Q^+ \cup \{0\}$ or $C^+ \cup \{0\}$ under addition they do not form a group.

We will illustrate this situation by some examples.

***Example 2.52***: Let G = {([0, $a_1$] [0, $a_2$] | [0, $a_3$] [0, $a_4$] [0, $a_5$] [0, $a_6$] | [0, $a_7$] [0, $a_8$] [0, $a_9$] / $a_i \in Z^+ \cup \{0\}$ $1 \le i \le 9$} be a collection of super interval row matrices. Clearly every super interval row matrix in G is of the same type.

We see G is a semigroup under interval row matrix addition. Further G is not a group. G is of infinite order. G is a commutative monoid.

G has no zero divisors G will be known as the super row interval matrix semigroup.

***Example 2.53***: Let H = {([0, $a_1$] [0, $a_2$]| [0, $a_3$] [0, $a_4$] [0, $a_5$] | $a_i \in Q^+ \cup \{0\}$ $1 \le i \le 5$} be a collection of super interval row matrices of the same type. H is only a super row interval matrix semigroup and is not a group. H is also of infinite order.

***Example 2.54***: Let P = {([0, $a_1$] | [0, $a_2$] [0, $a_3$] | [0, $a_4$] [0, $a_5$] | $a_i \in Q^+ \cup \{0\}$; $1 \le i \le 5$} be a collection of super interval row matrices of the same type. P under addition is not a group. P is only a super interval row matrix semigroup of infinite order.

We see the super interval row matrices in P and H are of same order but they are of different types. Thus we have two distinct semigroups under addition got using H and P.

Infact we can get several semigroups of same order (say $(1 \times n)$) but distinct super row interval matrix semigroups.

**Example 2.55**: Let

$$G = \left\{ \begin{pmatrix} [0,a_1] & [0,a_3] & | & [0,a_5] & | & [0,a_7] & [0,a_8] & [0,a_9] & | & [0,a_{10}] & [0,a_{11}] \\ [0,a_2] & [0,a_4] & | & [0,a_6] & | & [0,a_{12}] & [0,a_{13}] & [0,a_{14}] & | & [0,a_{15}] & [0,a_{16}] \end{pmatrix} \right.$$

$a_i \in R^+ \cup \{0\}; 1 \le i \le 16\}$ be a collection of super row interval vectors of the same type. Clearly G is a semigroup under addition known as the super interval row vector semigroup. G is of infinite order and G is commutative. G is not a group as G has no nontrivial zero divisors still no element in G has an additive inverse.

**Example 2.56**: Let

$$T = \left\{ \begin{pmatrix} [0,a_1] & [0,a_4] & | & [0,a_7] & | & [0,a_{10}] & [0,a_{13}] & [0,a_{16}] \\ [0,a_2] & [0,a_5] & | & [0,a_8] & | & [0,a_{11}] & [0,a_{14}] & [0,a_{17}] \\ [0,a_3] & [0,a_6] & | & [0,a_9] & | & [0,a_{12}] & [0,a_{15}] & [0,a_{18}] \end{pmatrix} \right.$$

$a_i \in R^+ \cup \{0\}; 1 \le i \le 18\}$ be the collection of all super row interval of same type. Clearly V under addition of super interval row vectors is a semigroup under addition.

**Example 2.57**: Let

$$T = \left\{ \begin{pmatrix} [0,a_1] & [0,a_4] & | & [0,a_7] & [0,a_{10}] & [0,a_{13}] \\ [0,a_2] & [0,a_5] & | & [0,a_8] & [0,a_{11}] & [0,a_{14}] \\ [0,a_3] & [0,a_6] & | & [0,a_9] & [0,a_{12}] & [0,a_{15}] \end{pmatrix} \middle| a_i \in Q^+ \cup \{0\}; 1 \le i \le 15 \right\}$$

be the collection super row interval vectors of same type. T is a semigroup under super row interval vector addition.

For take
$$X = \begin{pmatrix} [0,2] & [0,3] & [0,0] & [0,5] & [0,7] \\ [0,0] & [0,1] & [0,1/2] & [0,3] & [0,2] \\ [0,1] & [0,4] & [0,1] & [0,0] & [0,10] \end{pmatrix}$$

and

$$Y = \begin{pmatrix} [0,1] & [0,0] & [0,3] & [0,2] & [0,6] \\ [0,2] & [0,1] & [0,1/2] & [0,10] & [0,5] \\ [0,10] & [0,6] & [0,0] & [0,2] & [0,4] \end{pmatrix}$$

be elements in T.

To find $X + Y$. $X + Y$ under component wise addition of row super interval vectors is given by

$$X + Y = \begin{pmatrix} [0,3] & [0,3] & [0,3] & [0,7] & [0,13] \\ [0,2] & [0,2] & [0,1] & [0,13] & [0,7] \\ [0,11] & [0,10] & [0,1] & [0,2] & [0,14] \end{pmatrix}.$$

We see $X + Y$ is in T. Thus T is a semigroup known as the super row interval vector semigroup.

Now we proceed onto describe super interval column matrix semigroup and super interval column vector semigroup by examples.

*Example 2.58*: Let

$$G = \left\{ \begin{bmatrix} [0,a_1] & [0,a_2] & [0,a_3] \\ [0,a_4] & [0,a_5] & [0,a_6] \\ [0,a_7] & [0,a_8] & [0,a_9] \\ [0,a_{10}] & [0,a_{11}] & [0,a_{12}] \\ [0,a_{13}] & [0,a_{14}] & [0,a_{15}] \\ [0,a_{16}] & [0,a_{17}] & [0,a_{18}] \\ [0,a_{19}] & [0,a_{20}] & [0,a_{21}] \end{bmatrix} \right.$$

where $a_i \in Z^+ \cup \{0\}$; $1 \leq i \leq 21\}$ be the collection of super interval column vectors of the same type. G is a semigroup under addition called the super interval column vector semigroup. Clearly G is of infinite order.

**Example 2.59**: Let

$$V = \left\{ \begin{bmatrix} [0,a_1] & [0,a_2] \\ [0,a_3] & [0,a_4] \\ [0,a_5] & [0,a_6] \\ \hline [0,a_7] & [0,a_8] \\ [0,a_9] & [0,a_{10}] \\ \hline [0,a_{11}] & [0,a_{12}] \\ [0,a_{13}] & [0,a_{14}] \\ [0,a_{15}] & [0,a_{16}] \end{bmatrix} \middle| a_i \in Q^+ \cup \{0\}; 1 \le i \le 16 \right\}$$

be the collection of super interval column vectors of the same time. Clearly V is a semigroup under super interval column vector addition.

**Example 2.60**: Let

$$P = \left\{ \begin{bmatrix} [0,a_1] & [0,a_2] & [0,a_3] \\ [0,a_4] & [0,a_5] & [0,a_6] \\ \hline [0,a_7] & [0,a_8] & [0,a_9] \\ [0,a_{10}] & [0,a_{11}] & [0,a_{12}] \end{bmatrix} \right.$$

where $a_i \in R^+ \cup \{0\}$, $1 \le i \le 12$} be a super interval column vector of the same type. P is a super interval column vector semigroup of infinite order under addition.

**Example 2.61**: Let

$$V = \left\{ \begin{bmatrix} [0,a_1] & [0,a_2] & [0,a_3] & [0,a_4] \\ [0,a_5] & [0,a_6] & [0,a_7] & [0,a_8] \\ \hline [0,a_9] & [0,a_{10}] & [0,a_{11}] & [0,a_{12}] \end{bmatrix} \right.$$

where $a_i \in C^+ \cup \{0\}$, $1 \le i \le 12$} be the collection of super interval column vectors of same type. Clearly V is a super interval column vector semigroup under addition of infinite order.

***Example 2.62:*** Let

$$S = \left\{ \begin{bmatrix} [0,a_1] & [0,a_2] & [0,a_3] & [0,a_4] \\ \hline [0,a_5] & [0,a_6] & [0,a_7] & [0,a_8] \\ [0,a_9] & [0,a_{10}] & [0,a_{11}] & [0,a_{12}] \end{bmatrix} \right.$$

where $a_i \in C^+ \cup \{0\}$, $1 \leq i \leq 12\}$ be the collection of super interval column vectors of same type. We see S is a semigroup under addition.

Further S and V are distinct semigroups and they contain different types of super interval column vectors, though they are built using the same set $C^+ \cup \{0\}$.

Now we will proceed onto give examples of super interval matrix of same type.

***Example 2.63:*** Let

$$V = \left\{ \begin{bmatrix} [0,a_1] & [0,a_2] & [0,a_3] \\ [0,a_4] & [0,a_5] & [0,a_6] \\ [0,a_7] & [0,a_8] & [0,a_9] \\ \hline [0,a_{10}] & [0,a_{11}] & [0,a_{12}] \end{bmatrix} \right.$$

where $a_i \in Q^+ \cup \{0\}$, $1 \leq i \leq 12\}$ be a collection of super interval matrices of the same type. V under addition is a semigroup. V is a super interval matrix semigroup under addition.

***Example 2.64:*** Let

$$V = \left\{ \begin{bmatrix} [0,a_1] & [0,a_2] & [0,a_3] & [0,a_4] \\ [0,a_5] & [0,a_6] & [0,a_7] & [0,a_8] \\ [0,a_9] & [0,a_{10}] & [0,a_{11}] & [0,a_{12}] \\ [0,a_{13}] & [0,a_{14}] & [0,a_{15}] & [0,a_{16}] \end{bmatrix} \right.$$

where $a_i \in R^+ \cup \{0\}$, $1 \leq i \leq 16\}$ be the collection of all super interval matrices of the specified type. Under addition V is a semigroup of infinite order.

**Example 2.65**: Let

$$W = \left\{ \begin{bmatrix} [0,a_1] & [0,a_2] & [0,a_3] & [0,a_4] \\ [0,a_5] & [0,a_6] & [0,a_7] & [0,a_8] \\ [0,a_9] & [0,a_{10}] & [0,a_{11}] & [0,a_{12}] \\ [0,a_{13}] & [0,a_{14}] & [0,a_{15}] & [0,a_{16}] \end{bmatrix} \right.$$

where $a_i \in R^+ \cup \{0\}$, $1 \le i \le 16\}$ be the collection of all super interval matrices of the specified type. Clearly W is a semigroup of infinite order under addition. However both V and W are different.

Thus we see one can get many semigroups of super interval matrices using the same n × m interval matrix by partitioning them differently. This is one of the advantages of using super interval matrices instead of interval matrices. Now we can for every semigroup of super interval matrices define the concept of transpose semigroup of super interval matrices or super interval matrices of transpose semigroup under addition.

We will illustrate this situation by some examples.

**Example 2.66**: Let

$$V = \left\{ \begin{bmatrix} [0,a_1] \\ [0,a_2] \\ \overline{[0,a_3]} \\ [0,a_4] \\ [0,a_5] \\ [0,a_6] \\ \overline{[0,a_7]} \end{bmatrix} \middle| a_i \in Z^+ \cup \{0\}, 1 \le i \le 7 \right\}$$

be a collection of super interval column matrices of same type. Under addition V is a semigroup called the super interval column matrix semigroup. Clearly V is of infinite order.

Now $V^t = \{([0, a_1] [0, a_2] \mid [0, a_3] [0, a_4] [0, a_5] [0, a_6] \mid [0, a_7])$ $\mid a_i \in Z^+ \cup \{0\}; 1 \le i \le 7\}$ be the transpose of V that is set all super row interval matrices of the same type $V^t$ is also a semigroup under addition. $V^t$ is called the super interval matrix transpose semigroup of V.

*Example 2.67*: Let $W = \{([0, a_1] [0, a_2] \mid [0, a_3] [0, a_4] [0, a_5]) \mid a_i \in Z^+ \cup \{0\}; 1 \le i \le 5\}$ be the collection of super row interval matrices of the same type. W under addition is a semigroup. Now

$$W^t = \left\{ \begin{bmatrix} [0,a_1] \\ [0,a_2] \\ \hline [0,a_3] \\ [0,a_4] \\ [0,a_5] \end{bmatrix} \middle| a_i \in Z^+ \cup \{0\}, 1 \le i \le 5 \right\}$$

is the transpose of W is the collection of all super interval column matrices of same type and W is also a semigroup under addition know as the transpose semigroup of W.

*Example 2.68*: Let

$$G = \left\{ \begin{bmatrix} [0,a_1] & [0,a_2] & [0,a_3] & [0,a_4] & [0,a_5] & [0,a_6] \\ [0,a_7] & [0,a_8] & . & . & . & [0,a_{12}] \\ \hline [0,a_{13}] & [0,a_{14}] & . & . & . & [0,a_{18}] \\ [0,a_{19}] & [0,a_{20}] & . & . & . & [0,a_{24}] \\ [0,a_{25}] & [0,a_{26}] & . & . & . & [0,a_{30}] \\ [0,a_{31}] & [0,a_{32}] & . & . & . & [0,a_{36}] \\ \hline [0,a_{37}] & [0,a_{38}] & . & . & . & [0,a_{42}] \\ [0,a_{43}] & [0,a_{44}] & . & . & . & [0,a_{48}] \end{bmatrix} \right.$$

$a_i \in Q^+ \cup \{0\}; 1 \le i \le 48\}$ be the collection of super interval matrices of the same type. Clearly G is a semigroup under super interval matrix addition of infinite order. Now consider the transpose semigroup $G^T$ of G.

$$G^T = \left\{ \begin{bmatrix} [0,a_1] & [0,a_7] & [0,a_{13}] & \dots & [0,a_{37}] & [0,a_{43}] \\ [0,a_2] & [0,a_8] & [0,a_{14}] & \dots & [0,a_{38}] & [0,a_{44}] \\ \vdots & \vdots & \vdots & \vdots & \vdots & \vdots \\ [0,a_6] & [0,a_{12}] & [0,a_{18}] & \dots & [0,a_{42}] & [0,a_{48}] \end{bmatrix} \right.$$

where $a_i \in Q^+ \cup \{0\}; 1 \le i \le 48\}$ is the transpose semigroup of G.

Thus given a semigroup we have a transpose semigroup to exist provided the super interval matrix semigroup under consideration is under the operation addition otherwise in general the semigroup may not have transpose semigroup. We now proceed onto show how in general a multiplicative semigroup does not in general have transpose.

**Example 2.69**: Let G = {([0, $a_1$] [0, $a_2$] | [0, $a_3$] [0, $a_4$] [0, $a_5$] [0, $a_6$] | [0, $a_7$]) | $a_i \in Z^+ \cup \{0\}; 1 \le i \le 7$} be a collection of super interval semigroup of the stipulated type. Now G under multiplication is a semigroup. The transpose of G given by

$$G^t = \left\{ \begin{bmatrix} [0,a_1] \\ [0,a_2] \\ \overline{[0,a_3]} \\ [0,a_4] \\ [0,a_5] \\ \overline{[0,a_6]} \\ [0,a_7] \end{bmatrix} \middle| a_i \in Z^+ \cup \{0\}, 1 \le i \le 7 \right\}$$

is not closed under multiplication thus $G^t$ is not a semigroup under multiplication.

If on $G^t$ we can define addition we call $G^t$ a dual transpose semigroup of G under the operation +. We cannot say for every semigroup we may have dual transpose semigroup. Only in some cases we have dual transpose in some case we have dual transpose as well as transpose.

*Example 2.70*: Let

$$V = \left\{ \begin{bmatrix} [0,a_1] \\ [0,a_2] \\ \hline [0,a_3] \\ [0,a_4] \\ \hline [0,a_5] \\ [0,a_6] \\ [0,a_7] \end{bmatrix} \middle| a_i \in Q^+ \cup \{0\}, 1 \le i \le 7 \right\}$$

be a collection of super interval column matrices of the same type. V under addition is a semigroup.

$V^t$ = Let G = {([0, $a_1$] [0, $a_2$] | [0, $a_3$] [0, $a_4$] | [0, $a_5$] [0, $a_6$] [0, $a_7$]) | $a_i \in Q^+ \cup \{0\}$; $1 \le i \le 7$} can be both a dual transpose semigroup super row interval matrix under multiplication as well as transpose semigroup super interval row matrix under addition.

*Example 2.71*: Let

$$P = \left\{ \begin{bmatrix} [0,a_1] & [0,a_2] & [0,a_3] & [0,a_4] \\ [0,a_5] & [0,a_6] & [0,a_7] & [0,a_8] \\ \hline [0,a_9] & [0,a_{10}] & [0,a_{11}] & [0,a_{12}] \\ [0,a_{13}] & [0,a_{14}] & [0,a_{15}] & [0,a_{16}] \\ [0,a_{17}] & [0,a_{18}] & [0,a_{19}] & [0,a_{20}] \\ \hline [0,a_{21}] & [0,a_{22}] & [0,a_{23}] & [0,a_{24}] \end{bmatrix} \right.$$

where $a_i \in R^+ \cup \{0\}$; $1 \leq i \leq 24\}$ be a collection of super interval matrices of the same type. Now P under super interval matrix addition is a semigroup.

$$P^t = \left\{ \begin{bmatrix} [0,a_1] & [0,a_5] & [0,a_9] & [0,a_{13}] & [0,a_{17}] & [0,a_{21}] \\ [0,a_2] & [0,a_6] & [0,a_{10}] & [0,a_{14}] & [0,a_{18}] & [0,a_{22}] \\ [0,a_3] & [0,a_7] & [0,a_{11}] & [0,a_{15}] & [0,a_{19}] & [0,a_{23}] \\ [0,a_4] & [0,a_8] & [0,a_{12}] & [0,a_{16}] & [0,a_{20}] & [0,a_{24}] \end{bmatrix} \right.$$

where $a_i \in R^+ \cup \{0\}$; $1 \leq i \leq 24\}$ be the transpose of P. $P^t$ is only a transpose semigroup under addition and $P^t$ or P are not semigroup under multiplication as product is not compatible on P. This P has no dual transpose semigroup.

Now we can define interval groups under multiplication only for certain special type of super interval matrices we will illustrate this situation by some examples.

***Example 2.72***: Let G = {([0, $a_1$], [0, $a_2$] | … | [0, $a_9$]) | $a_i \in Q^+$; $1 \leq i \leq 9$} be the collection of all super row interval row matrices of the same type. Under component wise special multiplication G is a group.

We just illustrate the product which is new to super row matrices.
Let X = ([0, 3] [0, 1] | [0, 7] [0, 2] [0, 1] [0, 5] [0, 1/8] [0, 1/5] | [0, 3]) and Y = ([0, 1/3] [0, 2] | [0, 1/7] [0, 1/2] [0, 1] [0, 5] [0, 8] [0, 5] | [0, 3]) be two super row matrices in G. X.Y = ([0, 1] [0, 2] | [0, 1] [0, 1] [0, 1] [0, 25] [0, 1] [0, 1] | [0, 9]) is in G. It is easily verified G under this new multiplication operation is a group known as the group of super row interval matrices.

***Example 2.73***: Let P = {((([0, $a_1$] [0, $a_2$] [0, $a_3$] | [0, $a_4$] | [0, $a_5$] [0, $a_6$]) | $a_i \in R^+$; $1 \leq i \leq 6$} be the collection of all super row interval matrices of the same type. P under this special multiplication is a

group. Clearly ([0, 1] [0, 1] [0, 1] | [0, 1] | [0, 1] [0, 1]) in P acts as the multiplication identity.

However if we replace $R^+$ by $Z^+$ or $R^+ \cup \{0\}$ or $Q^+ \cup \{0\}$ the structure fails to be a group. Now we proceed onto define substructures in them.

Suppose G is a super interval row matrix collection of same type and G is a group of super interval matrices of the same type then we define a proper subset $H \subseteq G$ to be a subgroup if H under the operations of G is a group.

We will illustrate this situation by some examples.

***Example 2.74***: Let G = {([0, $a_1$] [0, $a_2$] | [0, $a_3$] [0, $a_4$] [0, $a_5$]) | $a_i \in R^+$; $1 \leq i \leq 5$} be a super interval row matrices of same type. G under the special multiplication is a group.

Consider H = {([0, $a_1$] [0, $a_2$] | [0, $a_3$] [0, $a_4$] [0, $a_5$]) | $a_i \in Q^+$; $1 \leq i \leq 5$} $\subseteq$ G be a proper subset of G. H under the operations of G is a group. Thus G is a Smarandache special definite group.

***Example 2.75***: Let X = {([0, $a_1$] [0, $a_2$] | [0, $a_3$] [0, $a_4$]) | $a_i \in Z_7 \setminus \{0\}$, $1 \leq i \leq 4$} be the collection of all super interval matrices of same type. X is a group.

Clearly X has a subgroup given by Y = ([0, $a_1$] [0, $a_2$] | [0, $a_3$] [0, $a_4$]) | $a_i \in \{1, 6\} \subseteq Z_7 \setminus \{0\}$; $1 \leq i \leq 4$} $\subseteq$ X is a subgroup. However X has no proper semigroup. Thus X is not a S-special definite group.

However all the examples of groups which are group of super interval matrices of same type happen to be commutative.

***Example 2.76***: Let P = {([0, $a_1$] [0, $a_2$] [0, $a_3$] | [0, $a_4$] [0, $a_5$] | [0, $a_6$] | [0, $a_7$] [0, $a_8$]) | $a_i \in Z_{13} \setminus \{0\}$; $1 \leq i \leq 8$} be a collection of super row interval matrices of same type. P is not a Smarandache special definite group.

In view of this we have a theorem the proof of which is left as an exercise to the reader.

**THEOREM 2.5**: *Let $G = \{([0, a_1] [0, a_2] | [0, a_3] \ldots [0, a_r] | [0, a_{r+1}] \ldots [0, a_{n-2}] | [0, a_{n-1}] [0, a_n])$ where $a_i \in Z_p \setminus \{0\}$; p a prime $1 \leq i \leq n\}$ be a collection of super row interval matrices of same type. G is not a Smarandache special definite group but only a group of super row interval matrices under component wise multiplication.*

***Example 2.77***: Let $G = \{([0, a_1] [0, a_2] | [0, a_3] [0, a_4] [0, a_5] | [0, a_6]) | a_i \in Z^+ \cup \{0\}; 1 \leq i \leq 6\}$ be the collection of super row interval vectors of same type. G is only a semigroup of super row interval matrices under new multiplication. Infact this semigroup is not even a Smarandache semigroup. This semigroup has however many subsemigroups.

Further G is also a semigroup under addition and even under addition G is not a Smarandache semigroup.

***Example 2.78***: Let $G = \{([0, a_1] [0, a_2] [0, a_3] [0, a_4] | [0, a_5] [0, a_6] | [0, a_7] [0, a_8]) | a_i \in Q^+ \cup \{0\}$ where $1 \leq i \leq 8\}$ be a collection of super row interval matrices of same type. G is a semigroup multiplication.

Further G is a Smarandache semigroup under multiplication as $H = \{([0, a_1] [0, a_2] [0, a_3] [0, a_4] | [0, a_5] [0, a_6] | [0, a_7] [0, a_8]) | a_i \in Q^+ ; 1 \leq i \leq 8\} \subseteq G$ is a group under multiplication. Thus G is a Smarandache semigroup under multiplication.

It is still interesting to note that G is not a Smarandache semigroup under addition. Thus we see a lot of difference between these super interval row matrices under addition and multiplication.

***Example 2.79***: Let $G = \{([0, a_1] | [0, a_2] [0, a_3] [0, a_4] [0, a_5] | [0, a_6] [0, a_7]) | a_i \in Z_{12}; 1 \leq i \leq 7\}$ be the set of all super interval row

vectors of the same type. G is a semigroup under multiplication. G is a Smarandache semigroup.

We can have the notion of Smarandache subsemigroup. Here some examples of them are given.

**Example 2.80**: Let G = {([0, $a_1$] [0, $a_2$] | [0, $a_3$] | [0, $a_4$] [0, $a_5$] [0, $a_6$] [0, $a_7$] [0, $a_8$]) | $a_i \in Z_{13}$; $1 \le i \le 8$} be a collection of super interval row vectors of the same type. G is a semigroup super row interval matrix of finite order under new multiplication.

Consider H = {([0, a] [0, a] | [0, a] | [0, a] [0, a] [0, a] | [0, a] [0, a]) | a $\in Z_{13} \setminus \{0\}$} $\subseteq$ G, H is a group. If we consider P = {([0, a] [0, a] | [0, a] | [0, a] [0, a] [0, a] | [0, a] [0, a]) | a $\in Z_{13}$ } $\subseteq$ G is a subsemigroup of G. Since H $\subseteq$ P $\subseteq$ G. P is a S-subsemigroup of G. Infact G is also a S-semigroup.

**Example 2.81**: Let G = {([0, $a_1$] [0, $a_2$] | [0, $a_3$] [0, $a_4$] [0, $a_5$] [0, $a_6$] | [0, $a_7$]) | $a_i \in Z_{14}$; $1 \le i \le 7$} be the collection of all super row vectors of same type.

G is a group under addition. Infact G is not a S-definite group.

**Example 2.82**: Let

$$H = \left\{ \begin{bmatrix} [0,a_1] & [0,a_2] \\ [0,a_3] & [0,a_4] \\ \hline [0,a_5] & [0,a_6] \\ [0,a_7] & [0,a_8] \\ [0,a_9] & [0,a_{10}] \\ \hline [0,a_{11}] & [0,a_{12}] \end{bmatrix} \middle| a_i \in Z_{25}, 1 \le i \le 12 \right\}$$

be a collection of super interval column matrices of same type. H is a group of super interval matrices H is not a S-definite group.

Transpose of H is also a group which is not a S-definite group under addition.

Interested reader can study the substructures like subgroups and normal subgroups of group of super interval column (row) vectors under addition. Since the operation addition is commutative we see all subgroups are normal subgroups. However we can define another innovative type of multiplication as in case of the special multiplication of super row interval matrices.

This multiplication we call as extended multiplication. Here multiplication of super interval square matrices of same type alone can be defined.

This is carried out in the following way.

Let A and B two super interval square matrices of same type. We define the extended product on AB as follows.

We as in case of usual square matrices determine the product AB and partition it as the same type for A and B · AB will be the square matrix of same order.

We will proceed onto describe them with examples.
Suppose
$$A = \begin{bmatrix} [0,8] & [0,3] \\ [0,4] & [0,1] \end{bmatrix} \text{ and } B = \begin{bmatrix} [0,1] & [0,2] \\ [0,3] & [0,4] \end{bmatrix}$$

be two 2 × 2 interval super matrices of same type.
Then
$$AB = \begin{bmatrix} [0,8] & [0,3] \\ [0,4] & [0,1] \end{bmatrix} \begin{bmatrix} [0,1] & [0,2] \\ [0,3] & [0,4] \end{bmatrix}$$

$$= \begin{bmatrix} [0,8][0,1]+[0,3][0,3] & [0,8][0,2]+[0,3][0,4] \\ [0,4][0,1]+[0,1][0,3] & [0,4][0,2]+[0,1][0,4] \end{bmatrix}$$

$$= \begin{bmatrix} [0,17] & [0,28] \\ [0,7] & [0,12] \end{bmatrix}.$$

This product is defined as the extendent product. In this product we do not bother about the partition but partition it later.

This is infact partly artificial but such product gives compatibility and can be seen as a nice algebraic structure. This type of product can be used in special situations when interactions among the elements takes place without respecting the partition but later fall under the same partition after producted. Clearly under this extended product also in general AB ≠ BA.

Suppose

$$A = \begin{bmatrix} [0,10] & [0,1] \\ \hline [0,4] & [0,2] \end{bmatrix} \text{ and } B = \begin{bmatrix} [0,2] & [0,7] \\ \hline [0,5] & [0,3] \end{bmatrix}$$

be two super square interval matrices of same type.

Then we fine the usual product AB and then partition AB as the same type as A or B

$$AB = \begin{bmatrix} [0,10] & [0,1] \\ \hline [0,4] & [0,2] \end{bmatrix} \begin{bmatrix} [0,2] & [0,7] \\ \hline [0,5] & [0,3] \end{bmatrix}$$

$$= \begin{bmatrix} [0,10][0,2]+[0,1][0,5] & [0,10][0,7]+[0,1][0,3] \\ \hline [0,4][0,2]+[0,2][0,5] & [0,4][0,7]+[0,2][0,3] \end{bmatrix}$$

$$= \begin{bmatrix} [0,25] & [0,73] \\ \hline [0,18] & [0,34] \end{bmatrix}.$$

Suppose

$$A = \begin{bmatrix} [0,9] & [0,1] & [0,2] & [0,0] \\ \hline [0,1] & [0,3] & [0,4] & [0,1] \\ [0,4] & [0,2] & [0,5] & [0,4] \\ [0,1] & [0,1] & [0,0] & [0,2] \end{bmatrix} \text{ and }$$

$$B = \begin{bmatrix} [0,0] & [0,1] & [0,2] & [0,3] \\ \hline [0,1] & [0,0] & [0,1] & [0,0] \\ [0,2] & [0,1] & [0,0] & [0,5] \\ [0,3] & [0,0] & [0,3] & [0,0] \end{bmatrix}$$

be two 9 × 4 super interval matrices of same type. We find AB the usual product extended to these super matrices and then partition AB of the same type.

Thus

$$AB = \begin{bmatrix} [0,5] & [0,11] & [0,19] & [0,37] \\ [0,14] & [0,5] & [0,8] & [0,23] \\ [0,24] & [0,9] & [0,22] & [0,37] \\ [0,7] & [0,1] & [0,9] & [0,3] \end{bmatrix}$$

is again a super interval 4 × 4 matrix of the same type.

These extended product on super interval matrices of same type will have a group or semigroup structure depending on the set over which they are defined.

This extended type of multiplication may be useful when interaction takes place overlooking the partition. When they are to be treated only as submatrix components this product is not suitable. Depending on the situation this extended product can be used. We can extend the notion of super linear algebra in case of super interval matrices. For in the first pace we use super interval matrices in which the intervals are of the special form. These collection of special intervals in general cannot be a field. Only in case the intervals are from $Z_p$, p a prime they form a field. So the super interval linear algebra will have meaning only when defined over $Z_p$ or intervals from $Z_p$.

They can only form a semifield when intervals are taken from $Q^+ \cup \{0\}$ or $Z^+ \cup \{0\}$ or $R^+ \cup \{0\}$. So we will discuss about super interval semivector spaces and super interval semilinear algebra and super interval linear algebras when $Z_p$ is used. This will form the following chapter.

**Chapter Three**

# SEMIRINGS AND SEMIVECTOR SPACES USING SUPER INTERVAL MATRICES

Here we will just define the notion of semivector spaces using super interval matrices. For more about general interval semirings or semifields please refer [46, 50].

**DEFINITION 3.1**: *Let $V = \{([0, a_1] \mid [0, a_2] [0, a_3] [0, a_4] \mid ... \mid [0, a_{n-2}] [0, a_{n-1}] [0, a_n]) \mid a_i \in Z^+ \cup \{0\}\}$ be a collection of super interval row matrices of the same type. Clearly V is a semigroup with identity under addition and V is also a semigroup under componentwise new multiplication. Thus $(V, +, \times)$ is a semiring known as the semiring of super row interval matrices.*

We will give examples of them as well as the working of their algebraic structure.

**Example 3.1**: Let $S = \{([0, a_1] [0, a_2] \mid [0, a_3] [0, a_4] [0, a_5]) \mid a_i \in Z^+ \cup \{0\}; 1 \leq i \leq 5\}$. S is a collection of super interval matrices of same type. Clearly (S, +) is the commutative monoid under addition.

(S, .) is a semigroup under multiplication. We see + distributes over '.' as if x, y, z ∈ S. We see (x + y).z = xz + yz.

We will just show how this works. x = ([0, 3] [0, 2] | [0, 1] [0, 5] [0, 1]) and y = ([0, 8] [0, 1] | [0, 3] [0, 1] [0, 4]) be in S. Consider x + y = ([0, 11] [0, 3] | [0, 4] [0, 6] [0, 5]) is in S. Also x + y = y + x. Further 0 = (0 0 | 0 0 0) acts as the additive identity in S and S is a commutative semigroup. Now if for same x, y in S. We see x.y = ([0, 24] [0, 2] | [0, 3] [0, 5] [0, 4]) is in S. Now x.y = y.x Thus (S, .) is also commutative semigroup.

Suppose z = ([0, 7] [0, 4] | [0, 2] [0, 1] [0, 4]) is in S. We see

$$
\begin{aligned}
(x + y)z &= ([0, 11]\ [0, 3]\ |\ [0, 4]\ [0, 6]\ [0, 5]) \times \\
&\quad ([0, 7]\ [0, 4]\ |\ [0, 2]\ [0, 1]\ [0, 4]) \\
&= ([0, 77]\ [0, 12]\ |\ [0, 8]\ [0, 6]\ [0, 20]). \qquad \text{I}
\end{aligned}
$$

Consider

$$
\begin{aligned}
xz + yz &= ([0, 21]\ [0, 8]\ |\ [0, 2]\ [0, 5]\ [0, 4]) + \\
&\quad ([0, 56]\ [0, 4]\ |\ [0, 6]\ [0, 1]\ [0, 16]) \\
&= ([0, 77]\ [0, 12]\ |\ [0, 8]\ [0, 6]\ [0, 20]) \qquad \text{II}
\end{aligned}
$$

Clearly I and II are equal hence distributive. Thus (S, +, .) is the semiring of super row interval matrices.

**Example 3.2**: Let S = {([0, $a_1$] [0, $a_2$] [0, $a_3$] [0, $a_4$] | [0, $a_5$] … [0, $a_{11}$]) | $a_i \in Q^+ \cup \{0\}$; $1 \le i \le 11$} be the collection of super row interval matrices of same type. S is a semiring of super row interval matrices.

It is pertinent to mention here that using the super interval row matrices we get several semirings of 1 × 11 super row interval matrices; on the contrary if we are using just interval row matrices of order 1 × 11 we will get only one semiring. We get several semirings using interval 1 × 11 row super matrices.

**Example 3.3**: Let P = {([0, $a_1$] [0, $a_2$] | [0, $a_3$] | [0, $a_4$] | [0, $a_5$] [0, $a_6$]) | $a_i \in Z_{18}$; $1 \leq i \leq 6$} be the collection of super interval row matrices of same type. P is a semiring we get several semiring of different types using the $1 \times 6$ interval row matrix. All of them are of finite order.

**Example 3.4**: Let W = {([0, $a_1$] [0, $a_2$] [0, $a_3$] [0, $a_4$] | [0, $a_5$] [0, $a_6$] [0, $a_7$] | [0, $a_8$] [0, $a_9$] [0, $a_{10}$]) | $a_i \in Z_{43}$; $1 \leq i \leq 10$} be the collection of super interval row matrices of same type. W is a semiring of finite order.

Now we have seen several semirings of super interval row matrices; now we proceed onto define the subsemiring and ideal of a semiring of super interval row matrices.

**DEFINITION 3.2**: *Let S = {([0, $a_1$] | [0, $a_2$] ... [0, $a_r$] | [0, $a_{r+1}$] ... [0, $a_t$] | [0, $a_{t+1}$] ... [0, $a_n$]) | $a_i \in R^+ \cup \{0\}$ or $Q^+ \cup \{0\}$ or $Z^+ \cup \{0\}$ or $Z_n$; $n < \infty$} be a super interval row matrix of same type. S is a semiring of super interval row matrices. Suppose $P \subseteq S$; if P itself is a semiring under the operations of S then we call P to be a subsemiring of super interval row matrices of S of same type.*

We will illustrate this by examples.

**Example 3.5**: Let W = {([0, $a_1$] [0, $a_2$] [0, $a_3$] | [0, $a_4$] | [0, $a_5$] | [0, $a_6$]) | $a_i \in Z^+ \cup \{0\}$; $1 \leq i \leq 6$} be the semiring of super row interval matrices of same type. Consider H = {([0, $a_1$] [0, $a_2$] [0, $a_3$] | [0, $a_4$] | [0, $a_5$] | [0, $a_6$]) | $a_i \in 3Z^+ \cup \{0\}$; $1 \leq i \leq 6$} $\subseteq$ W; H is a subsemiring of super row interval matrices of W.

**Example 3.6**: Let V = {([0, $a_1$] | [0, $a_2$] [0, $a_3$] | [0, $a_4$] [0, $a_5$] [0, $a_6$] | [0, $a_7$] [0, $a_8$] [0, $a_9$]) | $a_i \in R^+ \cup \{0\}$; $1 \leq i \leq 9$} be a semiring of super interval row matrices of same type. Consider H = {([0, $a_1$] | [0, $a_2$] [0, $a_3$] | [0, $a_4$] [0, $a_5$] [0, $a_6$] | [0, $a_7$] [0, $a_8$] [0, $a_9$]) | $a_i \in Q^+$

∪ {0}; 1 ≤ i ≤ 9} ⊆ V; H is a subsemiring of super interval row matrices of same type.

**Example 3.7**: Let S = {([0, $a_1$] | [0, $a_2$] | [0, $a_3$] [0, $a_4$] [0, $a_5$]) | $a_i$ ∈ $Z_{24}$; 1 ≤ i ≤ 5} be the semiring of super interval row matrices of same type. Consider V = {([0, $a_1$] | [0, $a_2$] | [0, $a_3$] [0, $a_4$] [0, $a_5$]) | $a_i$ ∈ {0, 2, 4, 6, 8, 10, …, 20, 22} ⊆ $Z_{24}$; 1 ≤ i ≤ 5} ⊆ S; V is a subsemiring of super interval row matrices of finite order.

**Example 3.8**: Let S = {([0, $a_1$] [0, $a_2$] [0, $a_3$] [0, $a_4$] [0, $a_5$] | [0, $a_6$]) | $a_i$ ∈ $Z_{17}$; 1 ≤ i ≤ 6} be the semiring of super interval row matrices of finite order. We see S has no subsemirings.

We have a class of semirings which have no subsemirings.

**THEOREM 3.1**: *Let S = {([0, $a_1$] [0, $a_2$] | [0, $a_3$] [0, $a_4$] [0, $a_5$] | [0, $a_6$] [0, $a_7$] … | [0, $a_t$] [0, $a_{t+1}$] | … | [0, $a_{n-2}$] [0, $a_{n-1}$] [0, $a_n$]) | $a_i$ ∈ $Z_p$; p a prime, $a_1$ = $a_2$ = … = $a_n$ = a; 1 ≤ i ≤ n} be a semiring of super interval matrices of finite order. S has no subsemirings.*

The proof is left as an exercise to the reader.
   We will define now the notion of ideals in the semiring of super interval row matrices.

**DEFINITION 3.3**: *Let S = {([0, $a_1$] [0, $a_2$] | [0, $a_3$] [0, $a_4$] [0, $a_5$] [0, $a_6$] | … | [0, $a_{n-4}$] [0, $a_{n-3}$] [0, $a_{n-2}$] [0, $a_{n-1}$] [0, $a_n$]) | $a_i$ ∈ $Z_n$ or $Z^+$ ∪ {0}, $Q^+$ ∪ {0}, $R^+$ ∪ {0}} be the semiring of super interval row matrix of same type. Let I ⊆ S; I is an ideal of S if I is a subsemiring of S and for every x ∈ S and every i ∈ I; xi and ix are in I.*

We will illustrate this by some examples.

**Example 3.9**: Let S = {([0, $a_1$] [0, $a_2$] | [0, $a_3$] [0, $a_4$] [0, $a_5$] | [0, $a_6$] [0, $a_7$]) | $a_i$ ∈ $Z^+$ ∪ {0}; 1 ≤ i ≤ 7} be a semiring of super interval

row matrices of same type. Consider I = {([0, $a_1$] [0, $a_2$] | [0, $a_3$] [0, $a_4$] [0, $a_5$] | [0, $a_6$] [0, $a_7$]) | $a_i$ ∈ $3Z^+$ ∪ {0}; 1 ≤ i ≤ 7} ⊆ S, I is an ideal of S.

Infact S has infinite number of ideals.

**Example 3.10**: Let V = {([0, a] | [0, a] | [0, a] | [0, a] [0, a] [0, a] [0, a]) where a ∈ $Z_7$} be a semiring of super interval row matrices. V has no ideals and V has no subsemirings.

**Example 3.11**: Let W = {([0, $a_1$] [0, $a_2$] | [0, $a_3$] | [0, $a_4$] [0, $a_5$] | [0, $a_6$] [0, $a_7$]) | $a_i$ ∈ $Z_7$; 1 ≤ i ≤ 7} be a semiring of super interval row matrices. P = {[0, $a_1$] [0, $a_2$] | 0 | 0 0 0 0) | $a_1$, $a_2$ ∈ $Z_7$} ⊆ W, P is an ideal of W.

Infact W has several ideals. Consider T = {([0, a] [0, a] | [0, a] | [0, a] [0, a] | [0, a] [0, a]) | a ∈ $Z_7$} ⊆ W ; T is only a subsemiring. We see T is not an ideal of W.

**Example 3.12**: Let S = {([0, $a_1$] [0, $a_2$] | [0, $a_3$] [0, $a_4$] [0, $a_5$]) | $a_i$ ∈ $R^+$ ∪ {0}; 1 ≤ i ≤ 5} be a semiring of super interval row matrices. R = {([0, $a_1$] 0 | 0 0 [0, $a_2$]) | $a_1$ $a_2$ ∈ $R^+$ ∪ {0}} ⊆ S, T is an ideal of S, ideal of the semiring of super interval row matrices.

**Example 3.13**: Let P = {([0, $a_1$] [0, $a_2$] [0, $a_3$] | [0, $a_4$]) | $a_i$ ∈ $R^+$ ∪ {0}; 1 ≤ i ≤ 4} be a semiring of super interval row matrices. W = {([0, $a_1$] | [0, $a_2$] [0, $a_3$] | [0, $a_4$]) | $a_i$ ∈ $Q^+$ ∪ {0}; 1 ≤ i ≤ 4} ⊆ P, W is only a subsemiring of P and is not an ideal of P.

Thus we have the following theorem the proof of which is direct.

**THEOREM 3.2**: *Let S = {([0, $a_1$] [0, $a_2$] [0, $a_3$] | [0, $a_4$] [0, $a_5$] | ... | [0, $a_{n-1}$] [0, $a_n$]) | $a_i$ ∈ $R^+$ ∪ {0}, $Q^+$ ∪ {0}, $Z_n$, $Z^+$ ∪ {0}} be the semiring of super interval row matrices. Every ideal of S is a subsemiring of S but every subsemiring in general is not an ideal of S.*

We can prove the later part by examples.

***Example 3.14***: Let W = {([0, $a_1$] [0, $a_2$] | [0, $a_3$] [0, $a_4$] [0, $a_5$] [0, $a_6$] | [0, $a_7$]) | $a_i$ ∈ $Z^+$ ∪ {0}; 1 ≤ i ≤ 7} be the semiring of super interval row matrices. T = {([0, a] [0, a] | [0, a] [0, a] [0, a] [0, a] | [0, a]) | a ∈ $Z^+$ ∪ {0}} ⊆ W is only a subsemiring but is not an ideal of W. Thus every subsemiring in general is not an ideal.

Now we can define the notion of super interval row matrix zero divisors, idempotents, units and nilpotents in these, semirings, we leave this task to the reader.

We will illustrate them by some examples.

***Example 3.15***: Let V = {([0, $a_1$] [0, $a_2$] [0, $a_3$] | [0, $a_4$] [0, $a_5$]) | $a_i$ ∈ $Z_{12}$; 1 ≤ i ≤ 5} be a semiring super interval row matrix same type.
  Let x = ([0, 8] [0, 4] [0, 2] | [0, 6] [0, 9]) only y = ([0, 3] [0,3] [0, 6] | [0, 6] [0, 4]) in V. We see x.y = ( 0 0 0 | 0 0).
  Consider x = ([0, 5] [0, 7] [0, 11] | [0, 1] [0, 7]) we see y = ([0, 5] [0, 7] [0, 11] | [0, 1] [0, 7]) in S in such that x.y = ([0, 1] [0, 1] [0, 1] | [0, 1] [0, 1]). Thus x is the inverse of y.
  Let x = ([0, 4] [0, 9] [0, 4] | [0, 1] [0, 9]) in S is such that $x^2$ = x. Thus x has super row matrix interval idempotent in S.
  Thus we see S has zero divisors units and idempotents.

***Example 3.16***: Let W = {([0, $a_1$] [0, $a_2$] | [0, $a_3$] [0, $a_4$] [0, $a_5$] | [0, $a_6$] [0, $a_7$]) | $a_i$ ∈ $Z^+$ ∪ {0}; 1 ≤ i ≤ 7} be the semiring of super interval matrices of same type. W has zero divisors but no idempotents or units.

***Example 3.17***: Let V = {([0, $a_1$] [0, $a_2$] | [0, $a_3$] [0, $a_4$] | [0, $a_5$]) | $a_i$ ∈ $R^+$ ∪ {0}; 1 ≤ i ≤ 5} be the semiring of super interval row matrices. V has ideals, subsemirings and subsemirings which are not ideals. R has only zero divisors.

Now we will proceed onto define strict semiring. We say a super interval row matrix semiring is a strict semiring S if $a + b = (0)$ for a, b $\in$ S is possible if and only if $a = b = 0$.

We will give examples of strict semirings built using super interval row matrices of same type.

***Example 3.18***: Let $V = \{([0, a_1] [0, a_2] [0, a_3] [0, a_4] | [0, a_5] | [0, a_6] [0, a_7] [0, a_8]) | a_i \in Z^+ \cup \{0\}; 1 \le i \le 8\}$ be a semiring of super interval row matrices. V is a strict semiring.

***Example 3.19***: Let $S = \{([0, a_1] | [0, a_2] [0, a_3] | [0, a_4] | [0, a_5] [0, a_6] | [0, a_7]) | a_i \in R^+ \cup \{0\}; 1 \le i \le 7\}$ be a super interval row matrix semiring. S is also a strict semiring. However S has zero divisors.

***Example 3.20***: Let $T = \{([0, a_1] [0, a_2] | [0, a_3] [0, a_4] [0, a_5] | [0, a_6]) | a_i \in Z_{240}; 1 \le i \le 6\}$ be a semiring of super interval row matrices of same type. T is not a strict semiring for $a + b = (0)$ exists even if $a \ne (0)$ and $b \ne (0)$; a, b $\in$ T.

***Example 3.21***: Let $W = \{([0, a_1] [0, a_2] | [0, a_3] [0, a_4] [0, a_5] | [0, a_6] | [0, a_7]) | a_i \in Z_{23}; 1 \le i \le 7\}$ be the semiring of super interval row matrices of same type. Clearly W is not a strict semiring for take $x = ([0, 20] [0, 10] | [0, 15] [0, 21] [0, 22] [0, 3] | [0, 15])$ and $y = ([0, 3] [0, 13] | [0, 8] [0, 2] [0, 1] [0, 20] | [0, 8])$ in W. $x+y = (0\ 0\ |\ 0\ 0\ 0\ |\ 0)$. Thus W is not a strict semiring.

**THEOREM 3.3**: *Let $V = \{([0, a_1] [0, a_2] | [0, a_3] | [0, a_4] [0, a_5] [0, a_6] | ... | [0, a_{n-4}] ... [0, a_n]) | a_i \in Z_n, n < \infty; 1 \le i \le n\}$ be the semiring of super row interval matrices. V is not a strict semiring.*

The proof is straight forward and hence is left as an exercise to the reader.

**THEOREM 3.4**: *Let $V = \{([0, a_1] [0, a_2] | [0, a_3] | [0, a_4] [0, a_5] | [0, a_6] ... [0, a_r] | ... [0, a_{r+1}] ... [0, a_n]) | a_i \in Z^+ \cup \{0\}$ or $R^+ \cup$*

*{0} or $Q^+ \cup \{0\}\}$ be a semiring of super interval row matrices of same type. V is a strict semiring.*

Now we proceed onto define the notion of semified. We say a semiring of super interval row matrices is a semifield

(a) If S is a strict semiring
(b) S has no zero divisors and
(c) S is a commutative semiring.

We give examples of such semifields.

***Example 3.22***: Let S = {([0, $a_1$] [0, $a_2$] [0, $a_3$] | [0, $a_4$] … [0, $a_r$] | [0, $a_{r+1}$] … [0, $a_n$]) | $a_i \in Z^+$; $1 \le i \le n$} $\cup$ {(0 0 0 | 0 … 0 | 0 … 0)} be the semiring of super row interval matrices. S is a semifield.

***Example 3.23***: Let P = {([0, $a_1$] [0, $a_2$] | [0, $a_3$] [0, $a_4$] [0, $a_5$] | [0, $a_6$]) | $a_i \in Q^+$; $1 \le i \le 6$} $\cup$ (0 0 | 0 0 0 | 0) be the semifield.

***Example 3.24***: Let W = {([0, $a_1$] | [0, $a_2$] [0, $a_3$] | [0, $a_4$] [0, $a_5$]) | $a_i \in Q^+$; $1 \le i \le 5$} $\cup$ {(0 0 | 0 0 | 0 0)} be the semifield.

***Example 3.25***: Let P = {([0, $a_1$] | [0, $a_2$] [0, $a_3$] [0, $a_4$] | [0, $a_5$] [0, $a_6$]) | $a_i \in Z^+ \cup \{0\}$ $1 \le i \le 6$} is a semiring which is a strict semiring but P is not a semifield as P has zero divisors.

***Example 3.26***: Let V = {([0, $a_1$] [0, $a_2$] | [0, $a_3$] [0, $a_4$] [0, $a_5$] [0, $a_6$] | [0, $a_7$]) | $a_i \in Z_n$, $1 \le i \le 7$} is a semiring which is not a strict semiring and also has zero divisors. Hence is not a semifield.

We see we have semifields. We will built semivector spaces using two types of semifields.
    (i)    Interval semifields and
    (ii)   Usual semifields.

**DEFINITION 3.4**: *Let V be an abelian super interval semigroup under addition with zero of one fixed type. S be a semifield of characteristic zero (ie. $Z^+ \cup \{0\}$ or $Q^+ \cup \{0\}$ or $R^+ \cup \{0\}$). If V is a semivector space over S then we define V to be a super interval row matrix semi vector space over S of type I.*

For more about semivector spaces please refer [46, 50].

We will illustrate this situation by some examples.

*Example 3.27*: Let $V = \{([0, a_1] | [0, a_2] [0, a_3] | [0, a_4] [0, a_5] [0, a_6]) | a_i \in R^+ \cup \{0\}, 1 \leq i \leq 6\}$ be a semigroup of interval row matrices of same type under addition. V is a commutative monoid as $(0 | 0\,0 | 0\,0\,0) \in V$ acts as the additive identity.

Consider $S = Z^+ \cup \{0\}$ the semifield. V is a semivector space of super row interval matrices over S of type I.

Consider $x = ([0, 3] | [0, \sqrt{2}], [0, \sqrt{5}] | [0, 4] [0, 0] [0, 7\sqrt{3}+5])$ in V. Take $s = 12$ in S. We see $sx = 12\,([0, 3] | [0, \sqrt{2}] [0, \sqrt{5}] | [0, 4] [0, 0] [0, 7\sqrt{3}+5]) = (0, 36] | [0, 12\sqrt{2}] [0, 12\sqrt{5}] | [0, 48] [0, 0] [0, 84\sqrt{3} + 60]) \in V$.

*Example 3.28*: Let $P = \{([0, a_1] | [0, a_2] | [0, a_3] | [0, a_4] [0, a_5] [0, a_6] [0, a_7]) | a_i \in Z^+ \cup \{0\}, 1 \leq i \leq 7\}$ be an abelian semigroup of super interval row matrices with additive identity $(0 | 0 | 0 | 0\,0\,0\,0)$.

Consider $S = Z^+ \cup \{0\}$ the semifield. P is a semivector space of super interval row matrices over S of type I.

*Example 3.29*: Let $V = \{([0, a_1] [0, a_2] [0, a_3] | [0, a_4] | [0, a_5] [0, a_6] [0, a_7] [0, a_8]) | a_i \in Q^+ \cup \{0\}, 1 \leq i \leq 8\}$ be an abelian semigroup of super interval row matrices with same type under addition.

V is a semivector space of super interval row matrices over $Q^+ \cup \{0\}$ of type I. Further V is also a semivector space of super interval row matrices over $Z^+ \cup \{0\}$ of type I. However V is not a semivector space of super interval row matrices over $R^+ \cup \{0\}$ of type I for if $x = \sqrt{3} \in R^+ \cup \{0\}$ and $v = ([0, 1] [0, 3] [0, 5/2] \mid [0, 5/9] \mid [0, 2] [0, 9/2] [0, 8] [0, 11/2])$ in V then $x \, v \notin V$ (easily verified).

Hence V is not a semivector space of super row interval matrices over $R^+ \cup \{0\}$ of type I.

Thus the concept of semivector spaces of super row interval matrices depend on the semifield over which they are defined. Further every additive abelian semigroup of super interval row vectors need not be a semivector space over every semifield.

***Example 3.30***: Let $V = \{([0, a_1] [0, a_2] [0, a_3] [0, a_4] \mid [0, a_5] [0, a_6] \mid [0, a_7]) \mid a_i \in 3Z^+ \cup \{0\}, 1 \le i \le 7\}$ be an additive abelian semigroup with identity; $(0\ 0\ 0\ 0 \mid 0\ 0 \mid 0)$. V is a semivector space of super interval matrices over the semifield $S = Z^+ \cup \{0\}$ of type I. Clearly V is a semivector space of super interval matrices of type I over $R^+ \cup \{0\}$ or $Q^+ \cup \{0\}$. (It is left for the reader to verify).

***Example 3.31***: Let $P = \{([0, a_1] \mid [0, a_2] \mid [0, a_3] \mid [0, a_4])\}$ be an additive semigroup of super interval row matrices where $a_i \in 5Z^+ \cup \{0\}$. P is a semivector space of super interval row matrices of type I over $Z^+ \cup \{0\}$ and not over $Q^+ \cup \{0\}$ or $R^+ \cup \{0\}$.

Now we will proceed onto define the substructure.

**DEFINITION 3.5**: *Let V be a semigroup of super row interval matrices of same type under addition with identity. Suppose V is a semivector space of super row interval matrices over S type I. Let $G \subseteq V$ (a proper subset of V). If G itself is a semivector space of super row interval matrices over S of type I then we define G to be a subsemivector space (semivector subspace) of super row interval matrices of type I over S of V.*

We will illustrate this situation by some examples.

**Example 3.32**: Let $V = \{([0, a_1] [0, a_2] \mid [0, a_3] [0, a_4] \mid [0, a_5] \mid [0, a_6]) \mid a_i \in Z^+ \cup \{0\}, 1 \leq i \leq 6\}$ be a semigroup of super interval row matrices of same type. V is a semivector space of super interval row matrices of type I over the semifield $S = Z^+ \cup \{0\}$. Take $W = \{([0, a_1] [0, a_2] \mid 0\ 0 \mid 0 \mid [0, a_3]) \mid a_1, a_2, a_3 \in Z^+ \cup \{0\}\} \subseteq V$; W is a semivector subspace of super interval row matrices of type I over the semifield S. Infact V has infinitely many sub-semivector spaces of super row interval vectors of type I over S.

**Example 3.33**: Let $V = \{([0, a_1] [0, a_2] \mid [0, a_3] \ldots [0, a_9] \mid [0, a_{10}] [0, a_{11}] \mid [0, a_{12}]) \mid a_i \in Q^+ \cup \{0\}, 1 \leq i \leq 12\}$ be a semivector space of super interval row matrices of type I over the semifield $S = Z^+ \cup \{0\}$. Consider $W = \{[0, a_1] [0, a_2] \mid [0, a_3] \ldots [0, a_9] \mid [0, a_{10}] [0, a_{11}] \mid [0, a_{12}]) \mid a_i \in Z^+ \cup \{0\}, 1 \leq i \leq 12\} \subseteq V$ is a semivector subspace of super interval row matrices of type I over the semifield S. Infact V has infinite number of semivector subspaces of super interval row matrices of type I over S.

Now we will illustrate other semivector spaces of super interval matrices over a semifields of type I.

**Example 3.34**: Let

$$V = \left\{ \begin{bmatrix} [0,a_1] \\ [0,a_2] \\ \hline [0,a_3] \\ [0,a_4] \\ \hline [0,a_5] \\ [0,a_6] \end{bmatrix} \middle| a_i \in Q^+ \cup \{0\}, 1 \leq i \leq 6 \right\}$$

be a semigroup of super interval column matrices with

$$\begin{bmatrix} 0 \\ \hline 0 \\ 0 \\ \hline 0 \\ 0 \\ 0 \end{bmatrix}$$

as its additive identity. $S = Z^+ \cup \{0\}$ is a semifield. V is a semivector space of super interval column matrices of type I over S. Suppose

$$x = \begin{bmatrix} [0,5] \\ \hline [0,3/2] \\ [0,7] \\ \hline [0,10] \\ [0,2] \\ [0,5] \end{bmatrix} \text{ and } s = 10 \in S.$$

We see 
$$xs = sx = \begin{bmatrix} [0,50] \\ \hline [0,15] \\ [0,70] \\ \hline [0,100] \\ [0,20] \\ [0,50] \end{bmatrix} \in V.$$

*Example 3.35*: Let

$$V = \left\{ \begin{bmatrix} [0,a_1] \\ \hline [0,a_2] \\ [0,a_3] \\ \hline [0,a_4] \\ [0,a_5] \\ [0,a_6] \\ \hline [0,a_7] \\ [0,a_8] \end{bmatrix} \middle| a_i \in Z^+ \cup \{0\}, 1 \le i \le 8 \right\}$$

be a semivector space of super column vectors of type I over the semifield $Z^+ \cup \{0\}$. Infact V cannot be a semivector space over other semifields.

*Example 3.36*: Let

$$V = \left\{ \begin{bmatrix} [0,a_1] \\ [0,a_2] \\ [0,a_3] \\ \vdots \\ [0,a_{10}] \\ [0,a_{12}] \\ \vdots \\ [0,a_{16}] \end{bmatrix} \middle| a_i \in Q^+ \cup \{0\}, 1 \le i \le 16 \right\}$$

be a semivector space of super column matrices of same type over the semifield $S = Q^+ \cup \{0\}$ of type I.

*Example 3.37*: Let

$$V = \left\{ \begin{bmatrix} [0,a_1] \\ [0,a_2] \\ [0,a_3] \\ [0,a_4] \end{bmatrix} \middle| a_i \in R^+ \cup \{0\}, 1 \le i \le 4 \right\}$$

be a semivector space of super column matrices of same type over $Z^+ \cup \{0\}$ of type I.

Infact V in example 3.36 can be a semivector space over $Z^+ \cup \{0\}$ but not over $R^+ \cup \{0\}$. Likewise V in example 3.37 can be a semivector space over $Q^+ \cup \{0\}$ also.

We will give some examples of semivector subspaces of super column matrices of type I.

*Example 3.38*: Let

$$V = \left\{ \begin{bmatrix} \overline{[0,a_1]} \\ [0,a_2] \\ [0,a_3] \\ \overline{[0,a_4]} \\ [0,a_5] \\ [0,a_6] \\ [0,a_7] \\ [0,a_8] \end{bmatrix} \middle| a_i \in Z^+ \cup \{0\}, 1 \le i \le 8 \right\}$$

be a semivector space of super interval column matrices over the semifield $S = Z^+ \cup \{0\}$ of type I.

Consider

$$W = \left\{ \begin{bmatrix} \overline{[0,a_1]} \\ 0 \\ [0,a_2] \\ \overline{0} \\ 0 \\ 0 \\ 0 \\ [0,a_3] \end{bmatrix} \middle| a_1, a_2, a_3 \in Z^+ \cup \{0\} \right\} \subseteq V;$$

W is a semi vector subspace of V of type I over the semifield S.

*Example 3.39*: Let

$$V = \left\{ \begin{bmatrix} \overline{[0,a_1]} \\ [0,a_2] \\ [0,a_3] \\ [0,a_4] \\ \overline{[0,a_5]} \\ [0,a_6] \\ [0,a_7] \\ [0,a_8] \end{bmatrix} \middle| a_i \in Q^+ \cup \{0\}, 1 \le i \le 8 \right\}$$

be a semivector space of super column matrices of type I over the semifield $S = Z^+ \cup \{0\}$.

Consider

$$W = \left\{ \begin{bmatrix} [0,a_1] \\ [0,a_2] \\ [0,a_3] \\ \hline [0,a_4] \\ \hline [0,a_5] \\ \hline [0,a_6] \\ [0,a_7] \\ [0,a_8] \end{bmatrix} \middle| a_i \in Z^+ \cup \{0\}, 1 \le i \le 8 \right\} \subseteq V;$$

W is a semivector subspace of V of type I over the semifield S. Infact V has infinitely many semivector subspaces.

*Example 3.40*: Let

$$V = \left\{ \begin{bmatrix} [0,a_1] \\ [0,a_2] \\ \hline [0,a_3] \\ [0,a_4] \\ \vdots \\ [0,a_{10}] \end{bmatrix} \middle| a_i \in R^+ \cup \{0\}, 1 \le i \le 10 \right\}$$

be a semivector space of super column vectors of type I over the semifield $S = Z^+ \cup \{0\}$.

Consider W which is a semivector subspace of super column matrices of type I of V over the semifield S where

$$W = \left\{ \begin{bmatrix} 0 \\ \overline{[0,a_2]} \\ \overline{[0,a_1]} \\ [0,a_4] \\ 0 \\ [0,a_3] \\ 0 \\ [0,a_5] \\ [0,a_6] \\ [0,a_7] \end{bmatrix} \middle| a_i \in Q^+ \cup \{0\}, 1 \le i \le 7 \right\} \subseteq V;$$

*Example 3.41*: Let

$$V = \left\{ \begin{bmatrix} [0,a_1] \\ \vdots \\ [0,a_5] \\ \overline{[0,a_6]} \\ [0,a_7] \\ \overline{[0,a_8]} \\ \vdots \\ [0,a_{12}] \\ \overline{[0,a_{13}]} \\ [0,a_{14}] \end{bmatrix} \middle| a_i \in R^+ \cup \{0\}, 1 \le i \le 14 \right\}$$

be a semivector space of super column vectors of type I over the semifield $S = Z^+ \cup \{0\}$.

Consider

$$P = \left\{ \begin{bmatrix} \begin{bmatrix} [0,a_1] \\ \vdots \\ [0,a_5] \\ \hline 0 \\ 0 \\ \hline 0 \\ 0 \\ 0 \\ 0 \\ \hline [0,a_6] \\ 0 \end{bmatrix} \end{bmatrix} \middle| a_i \in Z^+ \cup \{0\}, 1 \le i \le 6 \right\} \subseteq V$$

is a semivector subspace of V of type I over the semifield $S = Z^+ \cup \{0\}$.

Infact V has infinitely many subsemivector space of type I.

Now we will give examples of type I semivector spaces.

*Example 3.42*: Let

$$V = \left\{ \begin{bmatrix} \begin{bmatrix} [0,a_1] & [0,a_2] \\ \hline [0,a_3] & [0,a_4] \\ [0,a_5] & [0,a_6] \\ \hline [0,a_7] & [0,a_8] \\ [0,a_9] & [0,a_{10}] \\ [0,a_{11}] & [0,a_{12}] \end{bmatrix} \end{bmatrix} \middle| a_i \in Q^+ \cup \{0\}, 1 \le i \le 12 \right\}$$

be a semigroup of super column vectors of same type under addition. V is a semivector space of super column vectors of type I over the semifield $S = Z^+ \cup \{0\}$.

*Example 3.43*: Let

$$G = \left\{ \begin{bmatrix} [0,a_1] & [0,a_2] & | & [0,a_7] & [0,a_8] & [0,a_9] & | & [0,a_{10}] \\ [0,a_3] & [0,a_4] & | & [0,a_{11}] & [0,a_{12}] & [0,a_{13}] & | & [0,a_{14}] \\ [0,a_5] & [0,a_6] & | & [0,a_{15}] & [0,a_{16}] & [0,a_{17}] & | & [0,a_{18}] \end{bmatrix} \middle| a_i \in Z^+ \cup \{0\}, 1 \le i \le 18 \right\}$$

be the collection of super row interval vector of same type. G is a semigroup under addition G is a semivector space of super row interval vectors over the semifield $S = Z^+ \cup \{0\}$ of type I.

*Example 3.44*: Let

$$V = \left\{ \begin{bmatrix} [0,a_1] & | & [0,a_4] & [0,a_7] & | & [0,a_{10}] & [0,a_{13}] & [0,a_{16}] \\ [0,a_2] & | & [0,a_5] & [0,a_8] & | & [0,a_{11}] & [0,a_{14}] & [0,a_{17}] \\ [0,a_3] & | & [0,a_6] & [0,a_9] & | & [0,a_{12}] & [0,a_{15}] & [0,a_{18}] \end{bmatrix} \middle| a_i \in Z^+ \cup \{0\}, 1 \le i \le 18 \right\}$$

be a collection of super row interval vector of type different from one mentioned in example 3.43. V is a semivector space of super interval row vectors of type I over the semifield $S = Z^+ \cup \{0\}$.

Thus using the same 3 × 6 interval matrix we can get several semivector super interval row vectors of type I where as we can in case of 3 × 6 interval matrix construct only one semivector interval row vector space.

*Example 3.45*: Let

$$V = \left\{ \begin{bmatrix} [0,a_1] & | & [0,a_2] & [0,a_{11}] & | & [0,a_{16}] \\ [0,a_3] & | & [0,a_7] & [0,a_{12}] & | & [0,a_{17}] \\ [0,a_4] & | & [0,a_8] & [0,a_{13}] & | & [0,a_{18}] \\ [0,a_5] & | & [0,a_9] & [0,a_{14}] & | & [0,a_{19}] \\ [0,a_6] & | & [0,a_{10}] & [0,a_{15}] & | & [0,a_{20}] \end{bmatrix} \middle| a_i \in Q^+ \cup \{0\}, 1 \le i \le 20 \right\}$$

be the collection of super interval row vectors of the same type. V is a semivector space of super interval row vectors of type I over S = $Z^+ \cup \{0\}$.

We can have several semivector space of super interval row vectors of type I over S.

Now we proceed onto give examples super semivector space of interval matrices of type I and super semivector space of interval square matrices of type I.

*Example 3.46*: Let

$$V = \left\{ \begin{bmatrix} [0,a_1] & [0,a_2] & [0,a_3] \\ [0,a_4] & [0,a_5] & [0,a_6] \\ [0,a_7] & [0,a_8] & [0,a_9] \\ [0,a_{10}] & [0,a_{11}] & [0,a_{12}] \\ [0,a_{13}] & [0,a_{14}] & [0,a_{15}] \end{bmatrix} \middle| a_i \in Z^+ \cup \{0\}, 1 \le i \le 15 \right\}$$

be the collection of 5 × 3 super interval matrices of same type. V is a semigroup under addition. Clearly V is a super semivector space of interval matrices of type I over the semifield S = $Z^+ \cup \{0\}$ (or semivector space of super interval matrices of type I over the semifield S = $Z^+ \cup \{0\}$).

*Example 3.47*: Let

$$W = \left\{ \begin{bmatrix} [0,a_1] & [0,a_5] & [0,a_9] & [0,a_{10}] & [0,a_{12}] & [0,a_{13}] & [0,a_{14}] \\ [0,a_2] & [0,a_6] & [0,a_{15}] & [0,a_{11}] & [0,a_{16}] & [0,a_{17}] & [0,a_{18}] \\ [0,a_3] & [0,a_7] & [0,a_{19}] & [0,a_{20}] & [0,a_{21}] & [0,a_{22}] & [0,a_{23}] \\ [0,a_4] & [0,a_8] & [0,a_{24}] & [0,a_{25}] & [0,a_{26}] & [0,a_{27}] & [0,a_{28}] \end{bmatrix} \middle| a_i \in Q^+ \cup \{0\}, 1 \le i \le 28 \right\}$$

be the collection of all 4 × 7 super interval matrices of a specified type. Clearly W is a semigroup under super interval matrix addition. Now W is a semivector space of super interval 4 × 7 matrices of type I over $Z^+ \cup \{0\} = S$.

We see infact we can have several semivector spaces of super interval 4 × 7 matrices of type I over $Z^+ \cup \{0\}$ by changing the types of super interval matrices which is impossible in case of 4 × 7 interval matrices. This is possible only in case of super interval matrices.

*Example 3.48*: Let

$$V = \left\{ \begin{bmatrix} [0,a_1] & [0,a_2] & [0,a_3] & [0,a_4] \\ [0,a_5] & [0,a_6] & [0,a_7] & [0,a_8] \\ [0,a_9] & [0,a_{10}] & [0,a_{11}] & [0,a_{12}] \\ \hline [0,a_{13}] & [0,a_{14}] & [0,a_{15}] & [0,a_{16}] \\ [0,a_{17}] & [0,a_{18}] & [0,a_{19}] & [0,a_{20}] \\ \hline [0,a_{21}] & [0,a_{22}] & [0,a_{23}] & [0,a_{24}] \\ \hline [0,a_{25}] & [0,a_{26}] & [0,a_{27}] & [0,a_{28}] \end{bmatrix} \middle| a_i \in R^+ \cup \{0\}, 1 \le i \le 28 \right\}$$

be a collection of super interval 7 × 4 matrices of a specified type. Consider V under super interval matrix addition. V is a commutative monoid. Now V is a semivector space of super interval 7 × 4 matrices of type I over the semifield $S = Z^+ \cup \{0\}$.

Infact several semivector spaces can be constructed by changing the partition of the 7 × 4 matrices.

*Example 3.49*: Let

$$V = \left\{ \begin{bmatrix} [0,a_1] & [0,a_2] & [0,a_3] \\ [0,a_4] & [0,a_5] & [0,a_6] \\ [0,a_7] & [0,a_8] & [0,a_9] \\ \hline [0,a_{10}] & [0,a_{11}] & [0,a_{12}] \\ \hline [0,a_{13}] & [0,a_{14}] & [0,a_{15}] \end{bmatrix} \right.$$

where $a_i \in Q^+ \cup \{0\}$ ; $1 \leq i \leq 15\}$ be the collection of $5 \times 3$ super interval matrices of one specified type of partition.

V is semigroup under addition. V is a semivector space of super interval $5 \times 3$ matrices over the semifield $S = Q^+ \cup \{0\}$ of over I.

*Example 3.50*: Let

$$P = \left\{ \begin{bmatrix} [0,a_1] & [0,a_6] & [0,a_{11}] & [0,a_{12}] & [0,a_{13}] \\ [0,a_2] & [0,a_7] & [0,a_{14}] & [0,a_{15}] & [0,a_{16}] \\ \hline [0,a_3] & [0,a_8] & [0,a_{17}] & [0,a_{18}] & [0,a_{19}] \\ \hline [0,a_4] & [0,a_9] & [0,a_{20}] & [0,a_{21}] & [0,a_{22}] \\ [0,a_5] & [0,a_{10}] & [0,a_{23}] & [0,a_{24}] & [0,a_{25}] \end{bmatrix} \middle| a_i \in Z^+ \cup \{0\}, 1 \leq i \leq 25 \right\}$$

be the collection of super interval square matrices of one specifield type. P is a semigroup under addition. Infact P is a semivector space of super interval square matrices of type I over the semifield $S = Z^+ \cup \{0\}$.

Infact by varying the partition on the $5 \times 5$ interval square matrices we get different classes of semivector spaces of type I over S.

*Example 3.51*: Let

$$V = \left\{ \begin{bmatrix} [0,a_1] & [0,a_4] & [0,a_7] \\ \hline [0,a_2] & [0,a_5] & [0,a_8] \\ \hline [0,a_3] & [0,a_6] & [0,a_9] \end{bmatrix} \middle| a_i \in Z^+ \cup \{0\}, 1 \leq i \leq 9 \right\}$$

be a semivector space of super interval square matrices of a specified type over $S = Z^+ \cup \{0\}$ the semifield of type I.

We will study how many such super semivector spaces can be got using $3 \times 3$ interval matrices.

Taking

$$W_1 = \left\{ \left[ \begin{array}{c|c|c} [0,a_1] & [0,a_4] & [0,a_7] \\ \hline [0,a_2] & [0,a_5] & [0,a_8] \\ \hline [0,a_3] & [0,a_6] & [0,a_9] \end{array} \right] \,\middle|\, a_i \in Z^+ \cup \{0\}, 1 \le i \le 9 \right\}$$

is again semivector space of type I over S.

$$W_2 = \left\{ \left[ \begin{array}{c|cc} [0,a_1] & [0,a_4] & [0,a_7] \\ \hline [0,a_2] & [0,a_5] & [0,a_8] \\ \hline [0,a_3] & [0,a_6] & [0,a_9] \end{array} \right] \,\middle|\, a_i \in Z^+ \cup \{0\}, 1 \le i \le 9 \right\}$$

is again a semivector space of type I over S.

$$W_3 = \left\{ \left[ \begin{array}{c|cc} [0,a_1] & [0,a_4] & [0,a_7] \\ {[0,a_2]} & [0,a_5] & [0,a_8] \\ {[0,a_3]} & [0,a_6] & [0,a_9] \end{array} \right] \,\middle|\, a_i \in Z^+ \cup \{0\}, 1 \le i \le 9 \right\}$$

is again semivector space of type I over S.

$$W_4 = \left\{ \left[ \begin{array}{cc|c} [0,a_1] & [0,a_4] & [0,a_7] \\ {[0,a_2]} & [0,a_5] & [0,a_8] \\ {[0,a_3]} & [0,a_6] & [0,a_9] \end{array} \right] \,\middle|\, a_i \in Z^+ \cup \{0\}, 1 \le i \le 9 \right\}$$

is again a semivector space of type I over S.

$$W_5 = \left\{ \left[ \begin{array}{ccc} [0,a_1] & [0,a_4] & [0,a_7] \\ [0,a_2] & [0,a_5] & [0,a_8] \\ [0,a_3] & [0,a_6] & [0,a_9] \end{array} \right] \middle| a_i \in Z^+ \cup \{0\}, 1 \le i \le 9 \right\}$$

is again a semivector space of type I over S.

$$W_6 = \left\{ \left[ \begin{array}{ccc} [0,a_1] & [0,a_4] & [0,a_7] \\ [0,a_2] & [0,a_5] & [0,a_8] \\ \hline [0,a_3] & [0,a_6] & [0,a_9] \end{array} \right] \middle| a_i \in Z^+ \cup \{0\}, 1 \le i \le 9 \right\}$$

is again a semivector space of type I over S.

$$W_7 = \left\{ \left[ \begin{array}{c|c|c} [0,a_1] & [0,a_4] & [0,a_7] \\ [0,a_2] & [0,a_5] & [0,a_8] \\ [0,a_3] & [0,a_6] & [0,a_9] \end{array} \right] \middle| a_i \in Z^+ \cup \{0\}, 1 \le i \le 9 \right\}$$

is again a semivector space of type I over S.

$$W_8 = \left\{ \left[ \begin{array}{ccc} [0,a_1] & [0,a_4] & [0,a_7] \\ \hline [0,a_2] & [0,a_5] & [0,a_8] \\ [0,a_3] & [0,a_6] & [0,a_9] \end{array} \right] \middle| a_i \in Z^+ \cup \{0\}, 1 \le i \le 9 \right\}$$

is again a semivector space of type I over S.

$$W_9 = \left\{ \left[ \begin{array}{cc|c} [0,a_1] & [0,a_4] & [0,a_7] \\ [0,a_2] & [0,a_5] & [0,a_8] \\ [0,a_3] & [0,a_6] & [0,a_9] \end{array} \right] \middle| a_i \in Z^+ \cup \{0\}, 1 \le i \le 9 \right\}$$

is again a semivector space of type I over S.

$$W_{10} = \left\{ \left[ \begin{array}{cc|c} [0,a_1] & [0,a_4] & [0,a_7] \\ \hline [0,a_2] & [0,a_5] & [0,a_8] \\ {[0,a_3]} & [0,a_6] & [0,a_9] \end{array} \right] \middle| a_i \in Z^+ \cup \{0\}, 1 \le i \le 9 \right\}$$

is again a semivector space of type I over S.

$$W_{11} = \left\{ \left[ \begin{array}{c|cc} [0,a_1] & [0,a_4] & [0,a_7] \\ {[0,a_2]} & [0,a_5] & [0,a_8] \\ \hline [0,a_3] & [0,a_6] & [0,a_9] \end{array} \right] \middle| a_i \in Z^+ \cup \{0\}, 1 \le i \le 9 \right\}$$

is again a semivector space of type I over S.

$$W_{12} = \left\{ \left[ \begin{array}{cc|c} [0,a_1] & [0,a_4] & [0,a_7] \\ \hline [0,a_2] & [0,a_5] & [0,a_8] \\ {[0,a_3]} & [0,a_6] & [0,a_9] \end{array} \right] \middle| a_i \in Z^+ \cup \{0\}, 1 \le i \le 9 \right\}$$

is again a semivector space of type I over S.

$$W_{13} = \left\{ \left[ \begin{array}{c|c|c} [0,a_1] & [0,a_4] & [0,a_7] \\ {[0,a_2]} & [0,a_5] & [0,a_8] \\ {[0,a_3]} & [0,a_6] & [0,a_9] \end{array} \right] \middle| a_i \in Z^+ \cup \{0\}, 1 \le i \le 9 \right\}$$

is again a semivector space of type I over S.

$$W_{14} = \left\{ \left[ \begin{array}{c|c|c} [0,a_1] & [0,a_4] & [0,a_7] \\ {[0,a_2]} & [0,a_5] & [0,a_8] \\ {[0,a_3]} & [0,a_6] & [0,a_9] \end{array} \right] \middle| a_i \in Z^+ \cup \{0\}, 1 \le i \le 9 \right\}$$

is a semivector space of type I over S.

Thus we have 15 semivector space of type I using super interval 3 × 3 matrices, which is impossible in case of 3 × 3 interval matrices.

*Example 3.52*: Let

$$V = \left\{ \left[ \begin{array}{cc|cc} [0,a_1] & [0,a_2] & [0,a_3] & [0,a_4] \\ [0,a_5] & \ldots & \ldots & [0,a_8] \\ \hline [0,a_7] & \ldots & \ldots & [0,a_{12}] \\ [0,a_{13}] & \ldots & \ldots & [0,a_{16}] \end{array} \right] \middle| a_i \in Z^+ \cup \{0\}, 1 \le i \le 16 \right\}$$

be a semivector space of super interval square matrices of type I over the semifield $S = Z^+ \cup \{0\}$.

Now having seen examples of semivector spaces using super interval matrices we now proceed onto give examples of subsemivector spaces of type I.

*Example 3.53*: Let

$$V = \left\{ \left[ \begin{array}{cc|cc} [0,a_1] & [0,a_2] & [0,a_3] & [0,a_4] \\ [0,a_5] & [0,a_6] & [0,a_7] & [0,a_8] \\ \hline [0,a_9] & [0,a_{10}] & [0,a_{11}] & [0,a_{12}] \\ [0,a_{13}] & [0,a_{14}] & [0,a_{15}] & [0,a_{16}] \end{array} \right] \middle| a_i \in Z^+ \cup \{0\}, 1 \le i \le 16 \right\}$$

be a semivector space of super interval square matrices of type I over the semifield $S = Z^+ \cup \{0\}$.
Take

$$W = \left\{ \begin{bmatrix} 0 & [0,a_2] & [0,a_3] & [0,a_4] \\ 0 & [0,a_6] & [0,a_7] & [0,a_8] \\ \hline [0,a_1] & 0 & 0 & 0 \\ [0,a_2] & 0 & 0 & 0 \end{bmatrix} \middle| a_i \in Z^+ \cup \{0\}, 1 \le i \le 8 \right\} \subseteq V$$

is a semivector subspace of type I over S of V.

Consider

$$W_1 = \left\{ \begin{bmatrix} [0,a_1] & \ldots & \ldots & [0,a_4] \\ [0,a_5] & \ldots & \ldots & \ldots \\ \hline [0,a_9] & \ldots & \ldots & \ldots \\ [0,a_{12}] & \ldots & \ldots & [0,a_{16}] \end{bmatrix} \middle| a_i \in 5Z^+ \cup \{0\}, 1 \le i \le 16 \right\} \subseteq V$$

is a semivector subspace of V of type I over S.

*Example 3.54*: Let

$$V = \left\{ \begin{bmatrix} [0,a_1] & \ldots & [0,a_6] & [0,a_7] & [0,a_8] & \ldots & [0,a_{15}] \\ [0,a_{16}] & \ldots & [0,a_{21}] & [0,a_{22}] & [0,a_{23}] & \ldots & [0,a_{30}] \\ [0,a_{31}] & \ldots & [0,a_{36}] & [0,a_{37}] & [0,a_{38}] & \ldots & [0,a_{45}] \end{bmatrix} \middle| a_i \in Q^+ \cup \{0\}, 1 \le i \le 45 \right\}$$

be a semivector space of type I over $S = Z^+ \cup \{0\}$.

Consider

$$P = \left\{ \begin{bmatrix} [0,a_1] & \ldots & [0,a_6] & [0,a_7] & [0,a_8] & \ldots & [0,a_{15}] \\ [0,a_{16}] & \ldots & [0,a_{21}] & [0,a_{22}] & [0,a_{23}] & \ldots & [0,a_{30}] \\ [0,a_{31}] & \ldots & [0,a_{36}] & [0,a_{37}] & [0,a_{38}] & \ldots & [0,a_{45}] \end{bmatrix} \middle| a_i \in 3Z^+ \cup \{0\}, 1 \le i \le 45 \right\}$$

$\subseteq V$, P is a semivector subspace of V of type I over the semifield $S = Z^+ \cup \{0\}$.

*Example 3.55*: Let

$$V = \left\{ \begin{bmatrix} [0,a_1] & [0,a_2] & [0,a_3] \\ \underline{[0,a_4] & [0,a_5] & [0,a_6]} \\ [0,a_7] & [0,a_8] & [0,a_9] \\ [0,a_{10}] & [0,a_{11}] & [0,a_{12}] \\ \underline{[0,a_{13}] & [0,a_{14}] & [0,a_{15}]} \\ [0,a_{16}] & [0,a_{17}] & [0,a_{18}] \\ [0,a_{19}] & [0,a_{20}] & [0,a_{21}] \\ [0,a_{22}] & [0,a_{23}] & [0,a_{24}] \\ [0,a_{25}] & [0,a_{26}] & [0,a_{27}] \\ \underline{[0,a_{28}] & [0,a_{29}] & [0,a_{30}]} \end{bmatrix} \middle| a_i \in Z^+ \cup \{0\}, 1 \le i \le 30 \right\}$$

be a semi vector space of super interval column vectors of type I over the semifield $S = Z^+ \cup \{0\}$.

Let T = {all 10 × 3 super interval column vectors of type described in V with entries from $5Z^+ \cup \{0\} \subseteq V$. T is a semivector subspace of V of type I over the semifield $S = Z^+ \cup \{0\}$.

Now we will define subsemifield semivector subspace W of a semivector space V over the semifield S.

**DEFINITION 3.6**: *Let V be a semivector space of super interval matrices of type I over the semifield S. If $W \subseteq V$, be such that W is a semivector space of super interval matrices of type I over a subsemifield F of the semifield S, then we define W to be the subsemifield semivector subspace of type I of V over the subsemifield F of the semifield S.*

We will illustrate this by some simple examples.

*Example 3.56*: Let

$$V = \left\{ \begin{bmatrix} [0,a_1] & [0,a_3] \\ [0,a_2] & [0,a_5] \\ [0,a_4] & [0,a_6] \end{bmatrix} \middle| a_i \in Q^+ \cup \{0\}, 1 \leq i \leq 6 \right\}$$

be a semivector space of super interval matrices of type I over the semifield $S = Q^+ \cup \{0\}$.

Consider

$$W = \left\{ \begin{bmatrix} [0,a_1] & [0,a_3] \\ [0,a_2] & [0,a_5] \\ [0,a_4] & [0,a_6] \end{bmatrix} \middle| a_i \in Z^+ \cup \{0\}, 1 \leq i \leq 6 \right\} \subseteq V;$$

W is a subsemifield semivector subspace of V of type I over the subsemifield $Z^+ \cup \{0\} = F$ of the semifield $S = Q^+ \cup \{0\}$.

Clearly W is a semivector subspace of V over the semifield $S = Q^+ \cup \{0\}$.

If a semifield S has no subsemifield and if V is a semivector space of super interval matrices defined over that S then we say V is a pseudo simple semivector space of super interval matrices of type I over S.

We will illustrate this situation by some examples.

*Example 3.57*: Let

$$V = \left\{ \begin{bmatrix} [0,a_1] & [0,a_2] & [0,a_3] \\ [0,a_4] & [0,a_5] & [0,a_6] \\ [0,a_7] & [0,a_8] & [0,a_9] \\ [0,a_{10}] & [0,a_{11}] & [0,a_{12}] \\ [0,a_{13}] & [0,a_{14}] & [0,a_{15}] \end{bmatrix} \right.$$

where $a_i \in Z^+ \cup \{0\}$, $1 \le i \le 15\}$ be a semivector space of super interval column vectors of type I over the semifield $S = Z^+ \cup \{0\}$.

Clearly S has no subsemifields hence V is a pseudo simple vector space of super interval column vectors of type I over S. We see V has non trivial semivector subspaces of super interval column vector of type I over S.

*Example 3.58*: Let

$$V = \left\{ \begin{bmatrix} [0,a_1] & [0,a_2] & |[0,a_3] & [0,a_4] & [0,a_5]| & [0,a_6] \\ [0,a_7] & [0,a_8] & ... & ... & ... & [0,a_{12}] \\ [0,a_{13}] & [0,a_{14}] & ... & ... & ... & [0,a_{18}] \\ [0,a_{19}] & [0,a_{20}] & ... & ... & ... & [0,a_{24}] \end{bmatrix} \middle| a_i \in Q^+ \cup \{0\}, 1 \le i \le 24 \right\}$$

be a semivector space of super interval matrices of type I over the semifield $S = Z^+ \cup \{0\}$.

Since S has no subsemifields V is a pseudo semivector space of super interval matrices of type I over S.

Now we can characterize pseudo simple semivector spaces of super interval matrices.

**THEOREM 3.5**: *Let V be a semivector space of super interval matrices of type I over the semifield $S = Z^+ \cup \{0\}$. V is a pseudo simple super interval matrices of type I over S.*

The proof is direct and hence is left as an exercise for the reader to prove.

Now we proceed onto define other types of semivector spaces which we choose to call as semivector space of super interval matrices of type II.

We just recall $I(Z^+ \cup \{0\}) = \{[0, a] \mid a \in Z^+ \cup \{0\}.\}$ is a semifield.

Similarly $I(Q^+ \cup \{0\}) = \{[0, a] \mid a \in Q^+ \cup \{0\}\}$ is a semifield and $I(R^+ \cup \{0\}) = \{[0, a] \mid a R^+ \cup \{0\}\}$ is a semifield. Clearly $I(Z^+ \cup \{0\})$ has no subsemifields but $I(Q^+ \cup \{0\}$ and $I(R^+ \cup \{0\})$ has proper subsemifields.

**DEFINITION 3.7**: *Let V be a collection of super interval matrices of same type be a semigroup under addition (that is with same partition).*

*$S = I(Z^+ \cup \{0\})$ or $I(R^+ \cup \{0\})$ or $I(Q^+ \cup \{0\})$ be the semifield of intervals. If V is a semivector space over S then we define V to be a semivector space of super interval matrices of super interval matrices of type II over the same field S.*

We will illustrate this situation by some examples.

***Example 3.59***: Let $V = \{([0, a_1] \mid [0, a_2] [0, a_3] [0, a_4] \mid [0, a_5]) \mid a_i \in Z^+ \cup \{0\}, 1 \leq i \leq 5\}$ be a semivector space of super interval matrices of type II over the interval semifield $S = I(Z^+ \cup \{0\})$.

***Example 3.60***: Let $V = \{([0, a_1] [0, a_2] \mid [0, a_3] \mid [0, a_4] [0, a_5] [0, a_6]) \mid a_i \in Q^+ \cup \{0\}, 1 \leq i \leq 6\}$ be the collection of all super row interval matrices of same type. V is a semigroup under addition. V is a semivector space of super interval row matrices of type II over the semiring $S = I(Z^+ \cup \{0\})$.

For any x = ([0, 3] [0, 1/2] | [0, 7/5] | [0, 1] [0, 4] [0, 5]) and s = [0, 2] we have sx = ([0, 6] [0, 1] | [0, 14/5] | [0, 2] [0, 8] [0, 10]) ∈ V. Thus V is a semivector space of super interval row matrices of type II over S.

**Example 3.61**: Let V = {([0, $a_1$] | [0, $a_2$] [0, $a_3$] | [0, $a_4$] [0, $a_5$] [0, $a_6$] [0, $a_7$] | [0, $a_8$]) | $a_i$ ∈ $Z^+$ ∪ {0}, 1 ≤ i ≤ 8} be the collection of all super interval row matrices of same type. V under addition is a semigroup with identity. V is a semivector space of super row interval matrices of type II over the semifield S = I ($Z^+$ ∪ {0}).

**Example 3.62**: Let V = {(([0, $a_1$] [0, $a_2$] [0, $a_3$] [0, $a_4$] | [0, $a_5$]) | $a_i$ ∈ $Q^+$ ∪ {0}, 1 ≤ i ≤ 5} be a collection of super interval matrices of same type. V is a semigroup under addition. V is a semivector space of super row interval matrices of type II over the semifield I ($Q^+$ ∪ {0}) = S.

We will give examples of other super interval matrices.

**Example 3.63**: Let

$$V = \left\{ \begin{bmatrix} [0,a_1] \\ \overline{[0,a_2]} \\ [0,a_3] \\ \overline{[0,a_4]} \\ [0,a_5] \\ [0,a_6] \\ \overline{[0,a_7]} \\ [0,a_8] \\ [0,a_9] \end{bmatrix} \middle| a_i \in Q^+ \cup \{0\}, 1 \leq i \leq 9 \right\}$$

be a collection of super interval column matrices of same type. V is a semigroup under addition. Choose S = I ($Z^+$ ∪ {0}) semifield.

Clearly V is a semivector space of super column interval matrices of type II over S the interval semifield.

For if

$$x = \begin{bmatrix} \overline{[0,3]} \\ \overline{[0,2]} \\ \overline{[0,1]} \\ \overline{[0,0]} \\ [0,1] \\ \overline{[0,2]} \\ \overline{[0,0]} \\ [0,3] \\ [0,9] \end{bmatrix} \in V \text{ and } s = [0, 4] \text{ then } sx = \begin{bmatrix} \overline{[0,12]} \\ \overline{[0,8]} \\ \overline{[0,4]} \\ \overline{[0,0]} \\ [0,4] \\ \overline{[0,8]} \\ \overline{[0,0]} \\ [0,12] \\ [0,36] \end{bmatrix}.$$

*Example 3.64*: Let

$$P = \left\{ \begin{bmatrix} \overline{[0,a_1]} \\ \overline{[0,a_2]} \\ \overline{[0,a_3]} \\ \overline{[0,a_4]} \\ [0,a_5] \\ \overline{[0,a_6]} \\ \overline{[0,a_7]} \\ \overline{[0,a_8]} \\ [0,a_9] \\ [0,a_{10}] \end{bmatrix} \middle| a_i \in R^+ \cup \{0\}, 1 \le i \le 10 \right\}$$

be the collection of super interval column matrices of same type. P is a semigroup under addition.

Thus P is a semivector space of super interval column matrices over the interval semifield $S = I(R^+ \cup \{0\})$.

*Example 3.65*: Let

$$V = \left\{ \begin{bmatrix} [0,a_1] \\ [0,a_2] \\ \hline [0,a_3] \\ [0,a_4] \\ \hline [0,a_5] \\ [0,a_6] \end{bmatrix} \middle| a_i \in Z^+ \cup \{0\}, 1 \le i \le 6 \right\}$$

be a semivector space of super column interval matrices of type II over the interval semifield $S = I(Z^+ \cup \{0\})$. It is interesting to observe that V cannot be a semivector space over any other interval semifields.

*Example 3.66*: Let

$$V = \left\{ \begin{bmatrix} [0,a_1] & [0,a_2] & [0,a_3] & [0,a_4] & [0,a_5] & [0,a_6] \\ [0,a_7] & [0,a_8] & [0,a_9] & [0,a_{10}] & [0,a_{11}] & [0,a_{12}] \\ [0,a_{13}] & [0,a_{14}] & [0,a_{15}] & [0,a_{16}] & [0,a_{17}] & [0,a_{18}] \\ [0,a_{19}] & [0,a_{20}] & [0,a_{21}] & [0,a_{22}] & [0,a_{23}] & [0,a_{24}] \end{bmatrix} \middle| a_i \in Z^+ \cup \{0\}, 1 \le i \le 24 \right\}$$

be a semivector space of super row vectors over the interval semifield $S = I(Z^+ \cup \{0\})$ of type II.

*Example 3.67*: Let

$$W = \left\{ \begin{bmatrix} [0,a_1] & [0,a_2] & [0,a_3] \\ \hline [0,a_4] & [0,a_5] & [0,a_6] \\ [0,a_7] & [0,a_8] & [0,a_9] \\ [0,a_{10}] & [0,a_{11}] & [0,a_{12}] \\ \hline [0,a_{13}] & [0,a_{14}] & [0,a_{15}] \\ [0,a_{16}] & [0,a_{17}] & [0,a_{18}] \\ [0,a_{19}] & [0,a_{20}] & [0,a_{21}] \end{bmatrix} \middle| a_i \in Z^+ \cup \{0\}, 1 \le i \le 21 \right\}$$

be a semivector space of super column vectors of type II over the interval semifield $S = I(Z^+ \cup \{0\})$.

*Example 3.68*: Let

$$V = \left\{ \begin{bmatrix} [0,a_1] & [0,a_2] & [0,a_3] \\ [0,a_4] & [0,a_5] & [0,a_6] \\ [0,a_7] & [0,a_8] & [0,a_9] \\ [0,a_{10}] & [0,a_{11}] & [0,a_{12}] \\ [0,a_{13}] & [0,a_{14}] & [0,a_{15}] \\ [0,a_{16}] & [0,a_{17}] & [0,a_{18}] \\ [0,a_{19}] & [0,a_{20}] & [0,a_{21}] \\ [0,a_{22}] & [0,a_{23}] & [0,a_{24}] \end{bmatrix} \middle| a_i \in Q^+ \cup \{0\}, 1 \le i \le 24 \right\}$$

be a semivector space of super column vectors of type II over the interval semifield $S = I(Z^+ \cup \{0\})$.

*Example 3.69*: Let

$$V = \left\{ \begin{bmatrix} [0,a_1] & [0,a_6] & [0,a_7] & [0,a_{16}] & [0,a_{17}] & [0,a_{18}] \\ [0,a_2] & [0,a_8] & [0,a_9] & [0,a_{19}] & [0,a_{20}] & [0,a_{21}] \\ [0,a_3] & [0,a_{10}] & [0,a_{11}] & [0,a_{22}] & [0,a_{23}] & [0,a_{24}] \\ [0,a_4] & [0,a_{12}] & [0,a_{13}] & [0,a_{25}] & [0,a_{26}] & [0,a_{27}] \\ [0,a_5] & [0,a_{14}] & [0,a_{15}] & [0,a_{28}] & [0,a_{29}] & [0,a_{30}] \end{bmatrix} \middle| a_i \in R^+ \cup \{0\}, 1 \le i \le 30 \right\}$$

be the semigroup of super row interval vectors under addition. V is a semivector space of super row interval vectors of type II over $S = I(Z^+ \cup \{0\})$.

*Example 3.70*: Let

$$V = \left\{ \begin{bmatrix} [0,a_1] & [0,a_2] & [0,a_3] \\ [0,a_4] & [0,a_5] & [0,a_6] \\ [0,a_7] & [0,a_8] & [0,a_9] \\ [0,a_{10}] & [0,a_{11}] & [0,a_{12}] \\ [0,a_{13}] & [0,a_{14}] & [0,a_{15}] \\ [0,a_{16}] & [0,a_{17}] & [0,a_{18}] \end{bmatrix} \middle| a_i \in Q^+ \cup \{0\}, 1 \le i \le 18 \right\}$$

be a semivector space of super interval matrices over the interval semifield $S = I(Z^+ \cup \{0\})$ of type II.

*Example 3.71*: Let

$$W = \left\{ \begin{bmatrix} [0,a_1] & [0,a_2] & [0,a_3] \\ \hline [0,a_4] & [0,a_5] & [0,a_6] \\ [0,a_7] & [0,a_8] & [0,a_9] \end{bmatrix} \middle| a_i \in Z^+ \cup \{0\}, 1 \leq i \leq 9 \right\}$$

be a semivector space of super interval square matrices over the interval semifield $S = I(Z^+ \cup \{0\})$ of type II.

*Example 3.72*: Let

$$V = \left\{ \begin{bmatrix} [0,a_1] & [0,a_2] \\ \hline [0,a_3] & [0,a_4] \\ [0,a_5] & [0,a_6] \end{bmatrix} \middle| a_i \in Z^+ \cup \{0\}, 1 \leq i \leq 6 \right\}$$

be a semivector space of super interval matrices of type II over the semifield $S = I(Z^+ \cup \{0\})$.

We see we have by this definition get more and more semivector spaces. We now give the number of semivector spaces V given in example 3.72

$$V_1 = \left\{ \begin{bmatrix} [0,a_1] & [0,a_2] \\ [0,a_3] & [0,a_4] \\ [0,a_5] & [0,a_6] \end{bmatrix} \middle| a_i \in Z^+ \cup \{0\}, 1 \leq i \leq 6 \right\}$$

is a semivector space of type II over S.

$$V_2 = \left\{ \begin{bmatrix} [0,a_1] & [0,a_2] \\ \hline [0,a_3] & [0,a_4] \\ [0,a_5] & [0,a_6] \end{bmatrix} \middle| a_i \in Z^+ \cup \{0\}, 1 \leq i \leq 6 \right\}$$

is a semivector space of type II over S.

$$V_3 = \left\{ \begin{bmatrix} [0,a_1] & [0,a_2] \\ \hline [0,a_3] & [0,a_4] \\ [0,a_5] & [0,a_6] \end{bmatrix} \middle| a_i \in Z^+ \cup \{0\}, 1 \le i \le 6 \right\}$$

is a semivector space of type II.

$$V_4 = \left\{ \begin{bmatrix} [0,a_1] & [0,a_2] \\ \hline [0,a_3] & [0,a_4] \\ \hline [0,a_5] & [0,a_6] \end{bmatrix} \middle| a_i \in Z^+ \cup \{0\}, 1 \le i \le 6 \right\}$$

is a semivector space of type II.

$$V_5 = \left\{ \begin{bmatrix} [0,a_1] & | & [0,a_2] \\ [0,a_3] & | & [0,a_4] \\ [0,a_5] & | & [0,a_6] \end{bmatrix} \middle| a_i \in Z^+ \cup \{0\}, 1 \le i \le 6 \right\}$$

is a semivector space of type II V is given in example 3.72.

$$V_6 = \left\{ \begin{bmatrix} [0,a_1] & | & [0,a_2] \\ [0,a_3] & | & [0,a_4] \\ \hline [0,a_5] & | & [0,a_6] \end{bmatrix} \middle| a_i \in Z^+ \cup \{0\}, 1 \le i \le 6 \right\}$$

is a semivector space of type II over $S = I(Z^+ \cup \{0\})$.

$$V_7 = \left\{ \begin{bmatrix} [0,a_1] & | & [0,a_2] \\ \hline [0,a_3] & | & [0,a_4] \\ [0,a_5] & | & [0,a_6] \end{bmatrix} \middle| a_i \in Z^+ \cup \{0\}, 1 \le i \le 6 \right\}$$

is a semivector space of type II over S.

We see using the interval $3 \times 2$ matrix we can construct seven such super semivector spaces. This is one of the main advantages of using super interval matrices in place of usual interval matrices.

*Example 3.73*: Let

$$V = \left\{ \begin{bmatrix} [0,a_1] & [0,a_2] & [0,a_3] & [0,a_4] \\ [0,a_5] & [0,a_6] & [0,a_7] & [0,a_8] \\ [0,a_9] & [0,a_{10}] & [0,a_{11}] & [0,a_{12}] \\ [0,a_{13}] & [0,a_{14}] & [0,a_{15}] & [0,a_{16}] \\ [0,a_{17}] & [0,a_{18}] & [0,a_{19}] & [0,a_{20}] \\ [0,a_{21}] & [0,a_{22}] & [0,a_{23}] & [0,a_{24}] \\ [0,a_{25}] & [0,a_{26}] & [0,a_{27}] & [0,a_{28}] \\ [0,a_{29}] & [0,a_{30}] & [0,a_{31}] & [0,a_{32}] \end{bmatrix} \middle| a_i \in Q^+ \cup \{0\}, 1 \le i \le 32 \right\}$$

is a semivector space of super interval square matrices of type II over the interval semifield $S = I (Z^+ \cup \{0\})$.

Consider

$$H = \left\{ \begin{bmatrix} [0,a_1] & [0,a_2] & 0 & [0,a_{10}] \\ 0 & [0,a_3] & 0 & [0,a_{11}] \\ 0 & [0,a_4] & 0 & [0,a_{12}] \\ 0 & 0 & [0,a_5] & 0 \\ 0 & 0 & [0,a_6] & 0 \\ 0 & [0,a_7] & 0 & 0 \\ [0,a_{13}] & 0 & [0,a_8] & 0 \\ 0 & [0,a_9] & 0 & 0 \end{bmatrix} \middle| a_i \in Q^+ \cup \{0\}, 1 \le i \le 13 \right\} \subseteq V;$$

H is a semivector subspace of V over the interval semifield $S = I (Q^+ \cup \{0\})$ of type II.

*Example 3.74*: Let

$$V = \left\{ \begin{bmatrix} [0,a_1] & [0,a_2] & [0,a_3] & [0,a_4] \\ [0,a_5] & [0,a_6] & [0,a_7] & [0,a_8] \\ [0,a_9] & [0,a_{10}] & [0,a_{11}] & [0,a_{12}] \\ \hline [0,a_{13}] & [0,a_{14}] & [0,a_{15}] & [0,a_{16}] \\ [0,a_{17}] & [0,a_{18}] & [0,a_{19}] & [0,a_{20}] \\ \hline [0,a_{21}] & \ldots & \ldots & [0,a_{36}] \\ [0,a_{22}] & \ldots & \ldots & [0,a_{37}] \\ [0,a_{23}] & \ldots & \ldots & [0,a_{38}] \\ [0,a_{24}] & \ldots & \ldots & [0,a_{39}] \\ [0,a_{25}] & \ldots & \ldots & [0,a_{40}] \end{bmatrix} \middle| a_i \in Q^+ \cup \{0\}, 1 \le i \le 40 \right\}$$

be a semivector space of super column vectors of type II over the semifield $S = I\,(Z^+ \cup \{0\})$.

Consider

$$H = \left\{ \begin{bmatrix} [0,a_1] & \ldots & \ldots & [0,a_4] \\ [0,a_5] & \ldots & \ldots & [0,a_8] \\ [0,a_9] & \ldots & \ldots & [0,a_{12}] \\ \hline [0,a_{13}] & \ldots & \ldots & [0,a_{16}] \\ [0,a_{17}] & \ldots & \ldots & [0,a_{20}] \\ \hline [0,a_{21}] & \ldots & \ldots & [0,a_{36}] \\ \vdots & \ldots & \ldots & \vdots \\ \vdots & \ldots & \ldots & \vdots \\ \vdots & \ldots & \ldots & \vdots \\ [0,a_{25}] & \ldots & \ldots & [0,a_{40}] \end{bmatrix} \middle| a_i \in Z^+ \cup \{0\}, 1 \le i \le 40 \right\} \subseteq V;$$

H is a semivector subspace of super column vectors of type II over the interval semifield $S = I\,(Z^+ \cup \{0\})$.

*Example 3.75*: Let

$$V = \left\{ \begin{bmatrix} [0,a_1] & [0,a_2] & [0,a_3] & [0,a_4] & \ldots & [0,a_{12}] \\ [0,a_{13}] & [0,a_{14}] & [0,a_{15}] & [0,a_{16}] & \ldots & [0,a_{24}] \\ [0,a_{25}] & [0,a_{26}] & [0,a_{27}] & [0,a_{28}] & \ldots & [0,a_{36}] \\ [0,a_{37}] & [0,a_{38}] & [0,a_{39}] & [0,a_{40}] & \ldots & [0,a_{48}] \\ [0,a_{49}] & [0,a_{50}] & [0,a_{51}] & [0,a_{52}] & \ldots & [0,a_{60}] \end{bmatrix} \Big| a_i \in Z^+ \cup \{0\}, 1 \le i \le 60 \right\}$$

be a semivector space of super interval row vectors of type II over the interval semifield $S = I(Z^+ \cup \{0\})$.

Consider

$$W = \left\{ \begin{bmatrix} [0,a_1] & [0,a_2] & [0,a_3] & [0,a_4] & \ldots & [0,a_{12}] \\ \vdots & \vdots & \vdots & \vdots & \ldots & \vdots \\ \vdots & \vdots & \vdots & \vdots & \ldots & \vdots \\ \vdots & \vdots & \vdots & \vdots & \ldots & \vdots \\ [0,a_{49}] & [0,a_{50}] & [0,a_{51}] & [0,a_{52}] & \ldots & [0,a_{60}] \end{bmatrix} \Big| a_i \in 3Z^+ \cup \{0\}, 1 \le i \le 60 \right\}$$

$\subseteq V$. W is a semivector subspace of V over the semifield S of type II.

*Example 3.76*: Let

$$V = \left\{ \begin{bmatrix} [0,a_1] \\ \vdots \\ [0,a_5] \\ \hline [0,a_6] \\ \vdots \\ [0,a_{16}] \\ \hline [0,a_{17}] \\ [0,a_{18}] \\ [0,a_{19}] \\ [0,a_{20}] \end{bmatrix} \Big| a_i \in Z^+ \cup \{0\}, 1 \le i \le 20 \right\}$$

be a semivector space of super column interval matrices of type II over $S = I(Z^+ \cup \{0\})$.

Consider

$$H = \left\{ \begin{bmatrix} [0,a_1] \\ \vdots \\ [0,a_5] \\ \hline [0,a_6] \\ \vdots \\ [0,a_{16}] \\ \hline [0,a_{17}] \\ \vdots \\ \vdots \\ [0,a_{20}] \end{bmatrix} \;\middle|\; a_i \in 5Z^+ \cup \{0\}, 1 \le i \le 20 \right\} \subseteq V;$$

H is a semivector subspace of the super column interval matrices of type II over $S = I(Z^+ \cup \{0\})$.

***Example 3.77***: Let $P = \{([0, a_1] \ldots [0, a_5] \mid [0, a_6] [0, a_7] \mid [0, a_8] [0, a_9] [0, a_{10}] [0, a_{11}] \mid [0, a_{12}]) \mid a_i \in Q^+ \cup \{0\}; 1 \le i \le 12\}$ be a semivector space of super interval row matrix over the interval semifield $S = I(Z^+ \cup \{0\})$ of type II.

Consider $K = \{([0, a_1] \ldots [0, a_5] \mid [0, a_6] [0, a_7] \mid [0, a_8] [0, a_9] [0, a_{10}] [0, a_{11}] \mid [0, a_{12}]) \mid a_i \in 3Z^+ \cup \{0\}; 1 \le i \le 12\} \subseteq P$; K is a semivector subspace of super interval row matrix over the semifield S of type II of P.

***Example 3.78***: Let

$$V = \left\{ \begin{bmatrix} [0,a_1] & [0,a_2] & [0,a_3] \\ [0,a_4] & [0,a_5] & [0,a_6] \\ \hline [0,a_7] & [0,a_8] & [0,a_9] \\ [0,a_{10}] & [0,a_{11}] & [0,a_{12}] \\ [0,a_{13}] & [0,a_{14}] & [0,a_{15}] \\ [0,a_{16}] & [0,a_{17}] & [0,a_{18}] \end{bmatrix} \;\middle|\; a_i \in Q^+ \cup \{0\}, 1 \le i \le 18 \right\}$$

be a semivector space of super interval matrices of type II over the interval semifield $S = I(Z^+ \cup \{0\})$.

Consider

$$H = \left\{ \begin{bmatrix} [0,a_1] & [0,a_2] & [0,a_3] \\ [0,a_4] & [0,a_5] & [0,a_6] \\ \hline [0,a_7] & [0,a_8] & [0,a_9] \\ [0,a_{10}] & [0,a_{11}] & [0,a_{12}] \\ [0,a_{13}] & [0,a_{14}] & [0,a_{15}] \\ [0,a_{16}] & [0,a_{17}] & [0,a_{18}] \end{bmatrix} \middle| a_i \in Z^+ \cup \{0\}, 1 \leq i \leq 18 \right\} \subseteq V;$$

is a semivector subspace of super interval matrices of type II over $S = I(Z^+ \cup \{0\})$.

Now we will proceed onto define basis and some related results of these two types of semivector super interval matrix spaces.

Let S be a semifield and V a semivector space of super interval matrices. If for any $\beta = \sum_i s_i v_i$ ($v_i$ is a super interval matrix in V) and $s_i \in S$. We say $\beta$ in V is a linear combination of $v_i$'s where $v_i$'s are the super interval matrices. We also say $\beta$ is linearly dependent on $v_i$'s if $\beta$ can be expressed as a linear combination of $v_i$'s. We say the relation is a non trivial relation if atleast one of the coefficients $s_i$'s is non zero.

The set of super interval matrices $(v_1 \ldots v_k)$ satisfies a non trivial relation if $v_j$'s is a linear combination of $(v_1, \ldots, v_{j-1}, v_{j+1}, \ldots, v_k)$

We will illustrate this situation by an example.

Suppose $V = \{([0, a_1] \mid [0, a_2] [0, a_3] \mid [0, a_4] [0, a_5] [0, a_6]) \mid a_i \in Z^+ \cup \{0\}\ 1 \leq i \leq 6\}$ be a semivector space of super interval row matrices of type I over the semifield $S = Z^+ \cup \{0\}$.

Consider the collection $\{v_1, v_2, \ldots, v_6\}$ where

$v_1 = ([0, 1] \mid 0\ 0 \mid 0\ 0\ 0)$
$v_2 = (0 \mid [0, 1]\ [0, 1] \mid 0\ 0\ 0)$
$v_3 = (0 \mid 0\ 0 \mid [0, 1]\ 0\ 0)$
$v_4 = (0 \mid 0\ 0 \mid 0\ [0, 1]\ 0)$

$$v_5 = (0 \mid 0\ 0 \mid 0\ 0\ [0, 1]) \text{ and}$$
$$v_6 = ([0, 7] \mid [0, 2]\ [0, 2] \mid [0, 5]\ [0, 9]\ [0, 3]) \text{ in V.}$$

Now $v_6 = 7v_1 + 2v_2 + 5(v_3) + 9(v_4) + 3v_5$.

However no other element can represent $v_i$ interms of the other elements in the set $\{v_1, v_2, \ldots, v_6\} \subseteq V$.

Consider

$$V = \left\{ \begin{bmatrix} [0,a_1] & [0,a_2] & [0,a_3] \\ [0,a_4] & [0,a_6] & [0,a_7] \\ [0,a_5] & [0,a_8] & [0,a_9] \end{bmatrix} \middle| a_i \in Z^+ \cup \{0\}, 1 \le i \le 9 \right\}$$

the super semivector space of super interval matrices over the semifield $S = Z^+ \cup \{0\}$.

Suppose

$$P = \left\{ P_1 = \begin{bmatrix} [0,1] & 0 & 0 \\ 0 & 0 & 0 \\ 0 & 0 & 0 \end{bmatrix}, P_2 = \begin{bmatrix} 0 & 0 & 0 \\ [0,1] & 0 & 0 \\ 0 & 0 & 0 \end{bmatrix}, \right.$$

$$P_3 = \begin{bmatrix} 0 & 0 & 0 \\ 0 & 0 & 0 \\ [0,1] & 0 & 0 \end{bmatrix}, P_4 = \begin{bmatrix} 0 & [0,1] & 0 \\ 0 & 0 & 0 \\ 0 & 0 & 0 \end{bmatrix},$$

$$P_5 = \begin{bmatrix} 0 & 0 & [0,1] \\ 0 & 0 & 0 \\ 0 & 0 & 0 \end{bmatrix}, P_6 = \begin{bmatrix} 0 & 0 & 0 \\ 0 & [0,1] & 0 \\ 0 & 0 & 0 \end{bmatrix},$$

$$P_7 = \begin{bmatrix} 0 & 0 & 0 \\ 0 & 0 & [0,1] \\ 0 & 0 & 0 \end{bmatrix}, P_8 = \begin{bmatrix} 0 & 0 & 0 \\ 0 & 0 & 0 \\ 0 & [0,1] & 0 \end{bmatrix} \text{ and}$$

$$P_9 = \begin{bmatrix} [0,5] & [0,1] & [0,4] \\ \hline [0,7] & [0,6] & [0,7] \\ \hline [0,2] & [0,9] & 0 \end{bmatrix}$$

where $P_i \in V$; $1 \leq i \leq 9$.

It is easily verified

$$P_9 = 5P_1 + 7P_2 + 3P_3 + 1P_4 + 4P_5 + 6P_6 + 7P_7 + 9P_8$$

$$= \begin{bmatrix} [0,5] & [0,1] & [0,4] \\ \hline [0,7] & [0,6] & [0,7] \\ \hline [0,2] & [0,9] & 0 \end{bmatrix}.$$

Let

$$V = \left\{ \begin{bmatrix} [0,a_1] & [0,a_2] \\ \hline [0,a_3] & [0,a_4] \\ \hline [0,a_5] & [0,a_6] \\ \hline [0,a_7] & [0,a_8] \end{bmatrix} \middle| a_i \in Z^+ \cup \{0\}, 1 \leq i \leq 8 \right\}$$

be a semivector space of super interval column vector of type I over the semifield $S = Z^+ \cup \{0\}$.

Consider the set

$$P = \left\{ P_1 = \begin{bmatrix} [0,2] & 0 \\ \hline 0 & 0 \\ \hline 0 & 0 \\ \hline 0 & 0 \end{bmatrix}, P_2 = \begin{bmatrix} 0 & [0,12] \\ \hline 0 & 0 \\ \hline 0 & 0 \\ \hline 0 & 0 \end{bmatrix}, \right.$$

$$\left. P_3 = \begin{bmatrix} 0 & 0 \\ \hline 0 & 0 \\ \hline 0 & 0 \\ \hline [0,1] & [0,1] \end{bmatrix}, P_4 = \begin{bmatrix} [0,4] & [0,60] \\ \hline 0 & 0 \\ \hline 0 & 0 \\ \hline [0,7] & [0,7] \end{bmatrix} \right\} \subseteq V.$$

Clearly $P_4$ is a linear combination of $P_1$, $P_2$ and $P_3$ as
$$P_4 = 2P_1 + 5P_2 + 7P_3.$$

Now having seen examples of a linear combination of super vectors in the semivector space we define a subset P of a super semivector space V to be a linearly dependent set if there exists a relation among them other wise the set P is said to be a linearly independent set.

We will now give some examples of linearly independent set.

Consider the semivector space V = {([0, $a_1$] | [0, $a_2$] [0, $a_3$] [0, $a_4$] | [0, $a_5$] [0, $a_6$]) | $a_i \in Z^+ \cup \{0\}$ $1 \leq i \leq 6$} be a semivector space of super interval row matrices of type I over the semifield S = $Z^+ \cup \{0\}$

Consider the set P = ([0, 2] | 0 0 0 | 0 0), ([0, 5] | 0 0 [0, 2] | 0 0), ([0, 9] | 0 0 0 | [0, 1] 0])} $\subseteq$ V. P is a linearly independent set of super row interval vectors (matrices).

Interested reader can provide more examples.

The following theorem is straight forward and the interested reader can prove.

**THEOREM 3.6**: *Let V be a semivector space of super interval matrices over the semifield S. If $\alpha \in V$ is a linearly dependent set on $\{\beta_i\} \subseteq V$ and if each $\beta_i \in V$ is linearly dependent on $\{\gamma_i\}$ then $\alpha$ is linearly dependent on $\{\gamma_i\} \subseteq V$.*

Now we proceed onto define the notion of set spanned by A $\subseteq$ V.

**DEFINITION 3.8**: *Let V be a semivector space of super interval matrices of type I over the semifield S. For any subset A of V, the set of all linear combination of vectors in A is called the set spanned by A and we denote it by $\langle A \rangle$, and A $\subseteq \langle A \rangle$.*

*Thus if A $\subset$ B then $\langle A \rangle \subseteq \langle B \rangle$.*

We will illustrate this by some examples.

***Example 3.79***: Let V = {([0, $a_1$] [0, $a_2$] | [0, $a_3$] [0, $a_4$] [0, $a_5$] | [0, $a_6$]) | $a_i \in Z^+ \cup \{0\}$; $1 \leq i \leq 6$} be a semivector space of super interval row matrices over the semifield S = $Z^+ \cup \{0\}$.

Consider A = {([0, 1] [0, 2] | 0 0 0 | [0, 1])} $\subseteq$ V. Now <A> = {([0, $a_1$] [0, 2$a_1$] | 0 0 0 | [0, $a_1$]) | a $\in Z^+ \cup \{0\}$} $\subseteq$ V. Clearly A $\subseteq$ <A>. So the set spanned by A is infact a semivector subspace of V over S.

**Example 3.80**: Let

$$V = \left\{ \begin{bmatrix} [0,a_1] \\ [0,a_2] \\ \hline [0,a_3] \\ [0,a_4] \\ \hline [0,a_5] \\ [0,a_6] \\ [0,a_7] \end{bmatrix} \middle| a_i \in Z^+ \cup \{0\}, 1 \leq i \leq 7 \right\}$$

be a super semivector space of super interval column matrices over the semifield $S = Z^+ \cup \{0\}$.

Consider

$$A = \left\{ \begin{bmatrix} [0,1] \\ [0,1] \\ \hline 0 \\ 0 \\ \hline 0 \\ 0 \\ 0 \end{bmatrix}, \begin{bmatrix} 0 \\ 0 \\ \hline 0 \\ [0,1] \\ \hline 0 \\ 0 \\ 0 \end{bmatrix}, \begin{bmatrix} 0 \\ 0 \\ \hline 0 \\ 0 \\ \hline [0,1] \\ [0,1] \\ 0 \end{bmatrix} \right\} \subseteq V.$$

Now

$$\langle A \rangle = \left\{ \begin{bmatrix} [0,a_1] \\ [0,a_1] \\ \hline 0 \\ [0,a_2] \\ \hline [0,a_3] \\ [0,a_3] \\ 0 \end{bmatrix} \right\}$$

where $a_1, a_2, a_3 \in Z^+ \cup \{0\}\}$.

It is easily verified (1) $A \subseteq \langle A \rangle$ and $\langle A \rangle$ is again a super semivector subspace of V over the semifield $S = Z^+ \cup \{0\}$.

*Example 3.81*: Let

$$V = \left\{ \begin{bmatrix} [0,a_1] & [0,a_5] & [0,a_9] & [0,a_{13}] & [0,a_{17}] & [0,a_{21}] \\ [0,a_2] & [0,a_6] & [0,a_{10}] & [0,a_{14}] & [0,a_{18}] & [0,a_{22}] \\ [0,a_3] & [0,a_7] & [0,a_{11}] & [0,a_{15}] & [0,a_{19}] & [0,a_{23}] \\ [0,a_4] & [0,a_8] & [0,a_{12}] & [0,a_{16}] & [0,a_{20}] & [0,a_{24}] \end{bmatrix} \middle| a_i \in Z^+ \cup \{0\}, 1 \le i \le 24 \right\}$$

be a semivector space of super interval row vectors of type I over the semifield $S = Z^+ \cup \{0\}$.

Consider the set

$$A = \left\{ \begin{pmatrix} [0,1] & 0 & 0 & 0 & 0 & [0,1] \\ [0,1] & 0 & 0 & 0 & 0 & 0 \\ 0 & [0,1] & 0 & 0 & 0 & [0,1] \\ 0 & [0,1] & 0 & 0 & 0 & 0 \end{pmatrix}, \begin{pmatrix} 0 & 0 & [0,1] & 0 & 0 & 0 \\ 0 & 0 & 0 & 0 & 0 & 0 \\ 0 & 0 & 0 & 0 & 0 & 0 \\ 0 & 0 & 0 & 0 & [0,1] & 0 \end{pmatrix}, \right.$$

$$\left. \begin{pmatrix} 0 & 0 & 0 & 0 & 0 & 0 \\ 0 & 0 & 0 & [0,1] & 0 & 0 \\ 0 & 0 & [0,1] & 0 & 0 & 0 \\ 0 & 0 & 0 & 0 & 0 & 0 \end{pmatrix} \right\} \subseteq V.$$

Now

$$\langle A \rangle = \left\{ \begin{pmatrix} [0,a_1] & 0 & [0,a_2] & 0 & 0 & [0,a_1] \\ [0,a_1] & 0 & 0 & [0,a_3] & 0 & 0 \\ 0 & [0,a_1] & [0,a_3] & 0 & 0 & [0,a_1] \\ 0 & [0,a_1] & 0 & 0 & [0,a_2] & 0 \end{pmatrix} \middle| a_1, a_2, a_3 \in Z^+ \cup \{0\} \right\}$$

$\subseteq V$ ($A \subseteq \langle A \rangle$).

It is easily verifield $\langle A \rangle$ is a semivector subspace of V of super row interval vectors over the semifield $S = Z^+ \cup \{0\}$.

***Example 3.82:*** Let

$$V = \left\{ \begin{pmatrix} [0,a_1] & [0,a_2] & [0,a_3] \\ \hline [0,a_4] & [0,a_5] & [0,a_6] \\ \hline [0,a_7] & [0,a_8] & [0,a_9] \\ \hline [0,a_{10}] & [0,a_{11}] & [0,a_{12}] \\ [0,a_{13}] & [0,a_{14}] & [0,a_{15}] \\ \hline [0,a_{16}] & [0,a_{17}] & [0,a_{18}] \\ \hline [0,a_{19}] & [0,a_{20}] & [0,a_{21}] \end{pmatrix} \right.$$

where $a_i \in Q^+ \cup \{0\}$, $1 \le i \le 21\}$ be a semivector space of super interval column vectors of type I over the semifield $S = Q^+ \cup \{0\}$.

Consider the set

$$A = \left\{ \begin{pmatrix} 0 & [0,1] & 0 \\ \hline [0,1] & 0 & 0 \\ \hline 0 & 0 & 0 \\ \hline 0 & 0 & [0,1] \\ 0 & 0 & 0 \\ \hline 0 & 0 & 0 \\ \hline 0 & [0,1] & [0,1] \end{pmatrix}, \begin{pmatrix} [0,1] & 0 & 0 \\ \hline 0 & 0 & 0 \\ \hline 0 & 0 & [0,1] \\ \hline 0 & 0 & 0 \\ 0 & [0,1] & 0 \\ \hline [0,1] & 0 & 0 \\ \hline 0 & 0 & 0 \end{pmatrix}, \begin{pmatrix} 0 & 0 & 0 \\ \hline 0 & 0 & 0 \\ \hline 0 & [0,1] & 0 \\ \hline 0 & 0 & 0 \\ 0 & 0 & 0 \\ \hline 0 & 0 & [0,1] \\ \hline 0 & 0 & 0 \end{pmatrix} \right\} \subseteq V$$

be a set in V.
Now

$$\langle A \rangle = \left\{ \begin{pmatrix} [0,a_1] & [0,a_2] & 0 \\ \hline [0,a_2] & 0 & 0 \\ \hline 0 & [0,a_3] & [0,a_1] \\ \hline 0 & 0 & [0,a_2] \\ 0 & [0,a_1] & 0 \\ \hline [0,a_1] & 0 & [0,a_3] \\ \hline [0,a_1] & [0,a_2] & [0,a_3] \end{pmatrix} \middle| a_1, a_2, a_3 \in Q^+ \cup \{0\} \right\} \subseteq V$$

is a semivector subspace of V over the semifield $S = Q^+ \cup \{0\}$.

Now having seen examples of the spanning subset of a semivector space we now proceed onto give a theorem.

**THEOREM 3.7**: *Let V be a semivector space of super interval matrices over the semifield S of type I. Suppose B and C be any two subsets of V such that $B \subseteq C$ then $\langle B \rangle \subseteq \langle C \rangle$.*

The proof is direct and hence is left as an exercise for the reader to prove [38, 46, 50].

**THEOREM 3.8**: *Let V be a semivector space of super interval matrices over the semifield S of type I. Let $A = \{a_1, ..., a_t\}$ be a subset of V. If $a_i \in A$ is dependent on other vectors in A then $\langle A \rangle = \langle A \setminus a_i \rangle$*

This proof is also simple and hence is left as an exercise to the reader [38, 46, 50].

Now we proceed onto define the notion of basis of a super semi vector space over a semifield of type I.

**DEFINITION 3.9**: *A linearly independent set A of semivector space V of super interval matrices over the semifield S is called a basis of V if we that set A can span the semivector space V.*

We will illustrate this situation by some examples.

***Example 3.83***: Let $V = \{([0, a_1] [0, a_2] | [0, a_3] [0, a_4] [0, a_5] | [0, a_6]) | a_i \cup Z^+ \cup \{0\}; 1 \leq i \leq 6\}$ be a semivector space of super interval row matrices over the semifield S of type I. B = {([0, 1] 0 | 0 0 0 | 0), ([0, [0, 1] | 0 0 0 | 0), ([0 0 | [0, 1] 0 0 | 0), (0, 0 | 0, [0, 1], 0 | 0), (0, 0 | 0 0 [0, 1] | 0), (0 0 | 0 0 0 | [0, 1]) $\subseteq$ V is a super interval basis of V over S.

***Example 3.84***: Let $V = \{([0, a_1] | [0, a_2] [0, a_3] | [0, a_4] [0, a_5] [0, a_6] [0, a_7]) | a_i \cup Q^+ \cup \{0\}; 1 \leq i \leq 7\}$ be a semivector space of super interval row matrices over the semifield $S = Z^+ \cup \{0\}$. Consider B = {([0, 1] | 0 0 | 0 0 0 0), ([0 | [0, 1] 0 | 0 0 0 0), ([0 | 0 [0, 1] | 0 0 0 0), (0 | 0 0 | [0, 1] 0 0 0), (0 | 0 0 | 0 [0, 1] 0 0), (0 | 0 0 | 0 0 [0, 1] 0) (0 | 0 0 | 0 0 0 [0, 1])} $\subseteq$ V is a linearly independent subset of V but is not a basis of V over $Z^+ \cup \{0\}$.

However if we take instead of the semifield $S = Z^+ \cup \{0\}$ the semifield $Q^+ \cup \{0\}$ will be a basis of V over $Q^+ \cup \{0\}$; the semifield.

**Example 3.85**: Let $V = \{([0, a_1] [0, a_2] [0, a_3] \mid [0, a_4] [0, a_5] \mid [0, a_6]) \mid a_i \cup Z^+ \cup \{0\}; 1 \leq i \leq 6\}$ be a semivector space of super interval row matrices of type I over the semifield $S = Z^+ \cup \{0\}$.

Consider $B = \{([0, 1] [0, 1] 0 \mid 0 \, 0 \mid 0), ([0 \, 0 \, [0, 1] \mid 0 \, 0 \mid 0), ([0 \, 0 \, 0 \mid [0, 1] [0, 1] \mid 0), (0 \, 0 \, 0 \mid 0 \, 0 \mid [0, 1])\} \subseteq V$. B is a linearly independent subset of V but B is not a basis of V over the semifield $Z^+ \cup \{0\}$.

Thus we see from these examples a linearly independent set is not a basis. Also basis of a semivector space of super interval matrices depends on the semifield over which the semivector space is defined.

Further the following example will show that we can have more than the cardinality of the basis to form a linearly independent set.

**Example 3.86**: Let $V = \{([0, a_1] \mid [0, a_2] [0, a_3] \mid [0, a_4]) \mid a_i \cup Z^+ \cup \{0\}; 1 \leq i \leq 4\}$ be a semivector space over the semifield $S = Z^+ \cup \{0\}$ of type I.

Consider $B = \{([0, 1] \mid [0, 1] [0, 3] \mid 0), (0 \mid [0, 4] [0, 7] \mid 0), (0 \mid 0 \, 0 \mid [0, 5]), ([0, 1] \mid 0 \, 0 \mid [0, 2]), (0 \mid [0, 9] \, 0 \mid [0, 13])\} \subseteq V$. B is a linearly independent subset of V. But $T = \{([0, 1] \mid 0 \, 0 \mid 0), (0 \mid [0, 1] \, 0 \mid 0) \, (0 \mid 0 \, 0 \mid [0, 1]), (0 \mid 0 \, [0, 1] \mid 0)\} \subseteq V$ to be a basis of V over the semifield $S = Z^+ \cup \{0\}$.

Clearly $|T| = 4$ but $|B| = 5$ still B is also a linearly independent set which is not a basis of V and $|B| > |T|$ which is different from usual vector spaces.

**Example 3.87**: Let

$$V = \left\{ \begin{bmatrix} [0,a_1] & [0,a_5] \\ [0,a_2] & [0,a_6] \\ [0,a_3] & [0,a_7] \\ [0,a_4] & [0,a_8] \end{bmatrix} \right.$$

where $a_i \in Z^+ \cup \{0\}$, $1 \leq i \leq 8\}$ be a semivector space of super interval matrices over the semifield $S = Z^+ \cup \{0\}$.

Consider

$$B = \left\{ \begin{bmatrix} [0,1] & 0 \\ 0 & 0 \\ 0 & 0 \\ 0 & 0 \end{bmatrix}, \begin{bmatrix} 0 & 0 \\ 0 & 0 \\ [0,1] & 0 \\ 0 & 0 \end{bmatrix}, \begin{bmatrix} 0 & [0,1] \\ 0 & 0 \\ 0 & 0 \\ 0 & 0 \end{bmatrix}, \begin{bmatrix} 0 & 0 \\ 0 & [0,1] \\ 0 & 0 \\ 0 & 0 \end{bmatrix}, \right.$$

$$\begin{bmatrix} 0 & 0 \\ 0 & 0 \\ [0,1] & 0 \\ 0 & 0 \end{bmatrix}, \begin{bmatrix} 0 & 0 \\ 0 & 0 \\ 0 & [0,1] \\ 0 & 0 \end{bmatrix}, \begin{bmatrix} 0 & 0 \\ 0 & 0 \\ 0 & 0 \\ [0,1] & 0 \end{bmatrix}, \begin{bmatrix} 0 & 0 \\ 0 & 0 \\ 0 & 0 \\ 0 & [0,1] \end{bmatrix} \right\}$$

⊆V, B is a basis of V and |B| = 8. Thus the dimension of V is 8.
Now let

$$T = \left\{ \begin{bmatrix} 0 & [0,8] \\ 0 & 0 \\ 0 & 0 \\ [0,5] & 0 \end{bmatrix}, \begin{bmatrix} [0,9] & [0,31] \\ 0 & 0 \\ 0 & 0 \\ 0 & [0,2] \end{bmatrix}, \right.$$

$$\begin{bmatrix} [0,11] & 0 \\ [0,13] & 0 \\ 0 & 0 \\ 0 & 0 \end{bmatrix}, \begin{bmatrix} 0 & 0 \\ 0 & [0,7] \\ 0 & 0 \\ 0 & [0,5] \end{bmatrix}, \begin{bmatrix} 0 & 0 \\ 0 & 0 \\ 0 & [0,11] \\ [0,7] & 0 \end{bmatrix},$$

$$\begin{bmatrix} 0 & 0 \\ 0 & 0 \\ [0,19] & [0,23] \\ 0 & 0 \end{bmatrix}, \begin{bmatrix} 0 & 0 \\ 0 & 0 \\ 0 & 0 \\ [0,11] & [0,43] \end{bmatrix},$$

$$\left. \begin{bmatrix} 0 & 0 \\ [0,43] & 0 \\ [0,47] & 0 \\ [0,19] & 0 \end{bmatrix}, \begin{bmatrix} 0 & [0,101] \\ 0 & [0,19] \\ 0 & [0,17] \\ 0 & [0,13] \end{bmatrix}, \begin{bmatrix} [0,19] & 0 \\ 0 & [0,23] \\ [0,29] & 0 \\ 0 & [0,47] \end{bmatrix} \right\} \subseteq V,$$

is a linearly independent subset of V.

However T cannot form a basis of V over $S = Z^+ \cup \{0\}$. Further $|T| = 10$. Thus we see one can have more number of linearly independent elements than the basis of V.

***Example 3.88***: Let $V = \{([0, a_1] [0, a_2] [0, a_3] | [0, a_4] [0, a_5] | [0, a_6]) | a_i \cup R^+ \cup \{0\}; 1 \le i \le 6\}$ be a semivector space of super interval row matrices over the semifield $S = Z^+ \cup \{0\}$ of type I. Consider V is of infinite dimension over S. Further V has infinitely many linearly independent elements.

Infact V has $T = \{([0, a_1] [0, a_1] [0, a_1] | [0, a_1] [0, a_1] | [0, a_1]) | a_i \in Z^+ \cup \{0\}\} \subseteq V$ is an infinite set which is a linearly dependent set and is a semivector subspace of V. Clearly T is generated by $\{[0, 1] [0, 1] [0, 1] | [0, 1] [0, 1] | [0, 1]\}$; is a basis of T and dimension of T is one as a semivector subspace of V.

However if V is considered as a semivector space over the semifield $R^+ \cup \{0\} = F$, then V has dimension 6 and the base elements are given by $B = \{([0, 1] 0 0 | 0 0 | 0), (0, [0, 1] 0 | 0 0 | 0), (0 0 [0, 1] | 0 0 | 0), (0 0 0 | [0, 1] 0 | 0), (0 0 0 | 0 [0, 1] | 0), (0 0 0 | 0 0 | [0, 1])\} \subseteq V$ is a basis of V over the semifield $R^+ \cup \{0\}$ and dimension of V is six.

Now we proceed on to define the notion of linear transformation of two super semi vector spaces or semi vector spaces of super interval matrices.

**DEFINITION 3.10**: *Let V and W be two super semivector spaces or semivector spaces of super interval matrices defined over the same semifield S. A map / function $T : V \to W$ is said to be linear transformation of super semi vector spaces if $T(\alpha v + u) = \alpha T(v) + T(u)$ for all $u, v \in V$ and $\alpha \in S$.*

We will illustrate this situation by some simple examples.

**Example 3.89**: Let $V = \{([0, a_1] \mid [0, a_2] [0, a_3] [0, a_4] \mid [0, a_5] [0, a_6]) \mid a_i \in Z^+ \cup \{0\}; 1 \leq i \leq 6\}$ be a semivector space of super interval row matrices over the semifield $S = Z^+ \cup \{0\}$.

Let

$$W = \left\{ \begin{pmatrix} [0,a_1] & [0,a_2] \\ \overline{[0,a_3]} & \overline{[0,a_4]} \\ [0,a_5] & [0,a_6] \\ [0,a_7] & [0,a_8] \\ \overline{[0,a_9]} & \overline{[0,a_{10}]} \\ [0,a_{11}] & [0,a_{12}] \end{pmatrix} \middle| a_i \in Z^+ \cup \{0\}; 1 \leq i \leq 12 \right\}$$

be a semivector space of super interval column vector over the semifield $S = Z^+ \cup \{0\}$. Define $T : V \to W$ by $T([0, a_1] \mid [0, a_2] [0, a_3] [0, a_4] \mid [0, a_5] [0, a_6]) =$

$$\begin{pmatrix} [0,a_1] & [0,a_1] \\ \overline{[0,a_2]} & \overline{[0,a_2]} \\ [0,a_3] & [0,a_3] \\ [0,a_4] & [0,a_4] \\ \overline{[0,a_5]} & \overline{[0,a_5]} \\ [0,a_6] & [0,a_6] \end{pmatrix}$$

is a linear transformation of V to W over S.

**Example 3.90**: Let

$$V = \left\{ \begin{pmatrix} [0,a_1] & [0,a_2] \\ [0,a_3] & [0,a_4] \\ [0,a_5] & [0,a_6] \\ [0,a_7] & [0,a_8] \end{pmatrix} \middle| a_i \in Z^+ \cup \{0\}; 1 \leq i \leq 8 \right\}$$

be a semivector space of super interval matrices over the semifield $Z^+ \cup \{0\} = S$.

Let

$$W = \left\{ \left( \begin{array}{c|ccc} [0,a_1] & [0,a_3] & [0,a_4] & [0,a_5] \\ \hline [0,a_2] & [0,a_6] & [0,a_7] & [0,a_8] \end{array} \right) \middle| a_i \in Z^+ \cup \{0\}; 1 \le i \le 8 \right\}$$

be a semivector space of super interval matrices over the semifield $S = Z^+ \cup \{0\}$.

Define $T : V \to W$ as follows:

$$T \left( \begin{array}{c|c} [0,a_1] & [0,a_2] \\ \hline [0,a_3] & [0,a_4] \\ \hline [0,a_5] & [0,a_6] \\ \hline [0,a_7] & [0,a_8] \end{array} \right)$$

$$= \left( \begin{array}{c|ccc} [0,a_1] & [0,a_3] & [0,a_4] & [0,a_5] \\ \hline [0,a_2] & [0,a_6] & [0,a_7] & [0,a_8] \end{array} \right)$$

for every super interval matrix in V. It is easily verified T is a linear transformation of V to W.

If in the definition we replace W by V itself then the linear transformation is defined to be a linear operator on V.

We will illustrate this situation by some examples.

***Example 3.91***: Let $V = \{([0, a_1] \mid [0, a_2] \; [0, a_3] \; [0, a_4] \mid [0, a_5] \; [0, a_6] \mid [0, a_7]) \mid a_i \cup Z^+ \cup \{0\}; 1 \le i \le 7\}$ be a semivector space over the semifield $S = Z^+ \cup \{0\}$.

Define $T : V \to V$ by $T([0, a_1] \mid [0, a_2] \; [0, a_3] \; [0, a_4] \mid [0, a_5] \; [0, a_6] \mid [0, a_7]) = ([0, 3a_1] \mid [0, a_4] \; [0, a_3] \; [0, a_3] \mid [0, a_5] \; [0, a_6] \mid [0, 5a_7])$. It is easily verified that T is a linear operator on V.

Interested readers can construct examples of linear operators of semivector spaces of super interval matrices of type I over the semifield S.

All results described, illustrated and studied for semivector spaces of super interval matrices over semifields of type I can be easily extended to the case of super semivector spaces of super interval matrices over interval semifields of type II. Further the study of basis of linear independent set dependent set, linear operator and transformation can be defined with appropriate modifications for type II semivector spaces.

**Chapter Four**

# SUPER INTERVAL SEMILINEAR ALGEBRAS

In this chapter we for the first time introduce the notion of semilinear algebra of super interval matrices over semifields of type I (super semilinear algebra of type I) and semilinear algebra of super interval matrices over interval semifields of type II (super semilinear algebra of type II) and study their properties and illustrate them with examples.

**DEFINITION 4.1:** *Let V be a semivector space of super interval matrices defined over the semifield S of type I. If on V we for every pair of elements x, y $\in$ V; x . y is in V where '.' is the product defined on V, then we call V a semilinear algebra of super interval matrices over the semifield S of type I.*

We will illustrate this situation by some examples.

***Example 4.1***: Let V = {([0, $a_1$] [0, $a_2$] | [0, $a_3$] [0, $a_4$] | [0, $a_5$] [0, $a_6$] | [0, $a_7$]) | $a_i$ $\in$ $Z^+$ $\cup$ {0}; 1 $\leq$ i $\leq$ 7} be a semivector space of super interval matrices defined over the semifield S = $Z^+$ $\cup$ {0} of type I.
   Consider x = ([0, 5] [0, 3] | [0, 9] [0, 1] [0, 2] [0, 8] | [0, 6]) and y = ([0, 1] [0, 2] | [0, 3] [0, 5] [0, 3] [0, 1] | [0, 5]) in V. We define the product '.' on V as x.y =([0, 5] [0, 6] | [0, 27] [0, 5] [0, 6] [0, 8] | [0, 30])

∈ V. Thus V is a super semilinear algebra of super interval row matrices over the semifield S of type II.

***Example 4.2***: Let V = {([0, $a_1$] [0, $a_2$] [0, $a_3$] [0, $a_4$] [0, $a_5$] | [0, $a_6$] | [0, $a_7$] [0, $a_8$]) | $a_i$ ∪ $Z^+$ ∪ {0}; 1 ≤ i ≤ 8} be a semivector space over the semifield S = $Z^+$ ∪ {0} of type I. V is a semilinear algebra.

Now we wish to show as in case of every linear algebra is a vector space but every vector space in general is not a linear algebra, likewise every semilinear algebra is a semivector space but a semivector space in general is not a semilinear algebra.

**THEOREM 4.1**: *Let V be a semilinear algebra of super interval matrices over the semifield S of type I then V is a semivector space of super interval matrices over the semifield S of type I.*

But a super vector space of super interval matrices of type I in general is not a semi linear algebra of type I over S.
One way of the proof is obvious by the very definition of semilinear algebra.
To prove the other part we give few examples.

***Example 4.3***: Let

$$V = \left\{ \begin{bmatrix} [0,a_1] & [0,a_2] & [0,a_3] \\ \hline [0,a_4] & [0,a_5] & [0,a_6] \\ [0,a_7] & [0,a_8] & [0,a_9] \\ \hline [0,a_{10}] & [0,a_{11}] & [0,a_{12}] \\ [0,a_{13}] & [0,a_{14}] & [0,a_{15}] \\ [0,a_{16}] & [0,a_{17}] & [0,a_{18}] \end{bmatrix} \right.$$

where $a_i$ ∪ $R^+$ ∪ {0}; 1 ≤ i ≤ 18} be a semivector space of super interval column vectors interval column vectors of type I over the semifield S = $Q^+$ ∪ {0}. Clearly no product can be defined on V so V can never be a super semi linear algebra of super interval column vectors of type I over the semifield S = $Q^+$ ∪ {0}.

*Example 4.4*: Let

$$V = \left\{ \begin{pmatrix} \begin{array}{cc|ccc} [0,a_1] & [0,a_2] & [0,a_3] & [0,a_4] & [0,a_5] \\ [0,a_6] & [0,a_7] & [0,a_8] & [0,a_9] & [0,a_{10}] \\ \hline [0,a_{11}] & [0,a_{12}] & [0,a_{13}] & [0,a_{14}] & [0,a_{15}] \\ [0,a_{16}] & [0,a_{17}] & [0,a_{18}] & [0,a_{19}] & [0,a_{20}] \\ [0,a_{21}] & [0,a_{22}] & [0,a_{23}] & [0,a_{24}] & [0,a_{25}] \\ \hline [0,a_{26}] & [0,a_{27}] & [0,a_{28}] & [0,a_{29}] & [0,a_{30}] \\ [0,a_{31}] & [0,a_{32}] & [0,a_{33}] & [0,a_{34}] & [0,a_{35}] \end{array} \end{pmatrix} \middle| a_i \in Q^+ \cup \{0\}; 1 \le i \le 35 \right\}$$

be a semivector space super interval matrices over the semifield $S = Z^+ \cup \{0\}$ of type I.

Clearly it is impossible to define any product on V so V is only a semivector space and not a semilinear algebra of super interval matrices of type I over the semifield $S = Z^+ \cup \{0\}$.

*Example 4.5*: Let

$$V = \left\{ \begin{pmatrix} [0,a_1] \\ [0,a_2] \\ [0,a_3] \\ \hline [0,a_4] \\ [0,a_5] \\ [0,a_6] \\ [0,a_7] \\ \hline [0,a_8] \\ [0,a_9] \\ [0,a_{10}] \end{pmatrix} \middle| a_i \in Z^+ \cup \{0\}; 1 \le i \le 10 \right\}$$

be a semivector space of super interval column matrices over the semifield $S = Z^+ \cup \{0\}$ of type I. V is not a semi linear algebra of super interval column matrices of type I over the semifield S.

Now having seen examples we can define semi linear subalgebra of a semi linear algebra which is carried out as a matter of routine. The reader is expected to give examples of this special substructure. Now we define a new concept called pseudo super linear subalgebra of a semivector space of super interval matrices over the semifield S of type I.

**DEFINITION 4.2**: *Let V be a semivector space of super interval matrices over a semifield S of type I. Suppose W $\subseteq$ V, W be a semivector subspace of V over S of type I. If on W we can define a product so that W becomes a semilinear algebra but V is not a semilinear algebra then we define W to be a pseudo semilinear subalgebra of super interval matrices of type I of V over the semifield S.*

We will illustrate this by some examples.

*Example 4.6*: Let

$$V = \left\{ \begin{bmatrix} [0,a_1] & [0,a_2] & [0,a_3] \\ [0,a_4] & [0,a_5] & [0,a_6] \\ [0,a_7] & [0,a_8] & [0,a_9] \end{bmatrix} \right.$$

where $a_i \cup Z^+ \cup \{0\}$; $1 \le i \le 9\}$ be a semivector space of super interval matrices of type I over the semifield $Z^+ \cup \{0\} = S$.
Choose

$$W = \left\{ \begin{pmatrix} [0,a_1] & 0 & 0 \\ 0 & 0 & 0 \\ 0 & 0 & 0 \end{pmatrix} \middle| a_i \in Z^+ \cup \{0\} \right\} \subseteq V,$$

W is a semivector subspace of V over the semifield S of type I. Clearly we can define product on W for if

$$A = \begin{bmatrix} [0,5] & 0 & 0 \\ 0 & 0 & 0 \\ 0 & 0 & 0 \end{bmatrix} \text{ and } B = \begin{bmatrix} [0,8] & 0 & 0 \\ 0 & 0 & 0 \\ 0 & 0 & 0 \end{bmatrix}$$

are in W then

$$AB = \begin{bmatrix} [0,40] & 0 & 0 \\ 0 & 0 & 0 \\ 0 & 0 & 0 \end{bmatrix}$$

is in W. Thus W is a semilinear subalgebra of type I over S.

Thus W is a pseudo semilinear subalgebra of super interval matrices over S of V of type I.

*Example 4.7*: Let

$$V = \left\{ \begin{bmatrix} [0,a_1] & [0,a_2] & [0,a_3] & [0,a_4] \\ [0,a_5] & [0,a_6] & [0,a_7] & [0,a_8] \\ \hline [0,a_9] & [0,a_{10}] & [0,a_{11}] & [0,a_{12}] \\ [0,a_{13}] & [0,a_{14}] & [0,a_{15}] & [0,a_{16}] \end{bmatrix} \right\}$$

where $a_i \cup Z^+ \cup \{0\}$; $1 \le i \le 16$} be a semivector space of super interval matrices of type I over the semifield $S = Z^+ \cup \{0\}$.

Now choose a proper subset

$$W = \left\{ \begin{pmatrix} [0,a_1] & [0,a_2] & 0 & 0 \\ [0,a_3] & [0,a_4] & 0 & 0 \\ \hline 0 & 0 & 0 & 0 \\ 0 & 0 & 0 & 0 \end{pmatrix} \middle| a_i \in Z^+ \cup \{0\}, 1 \le i \le 4 \right\} \subseteq V.$$

Clearly W is a super semi vector subspace of V of type I over the semifield S.

Consider

$$A = \begin{pmatrix} [0,3] & [0,10] & 0 & 0 \\ [0,2] & [0,1] & 0 & 0 \\ \hline 0 & 0 & 0 & 0 \\ 0 & 0 & 0 & 0 \end{pmatrix} \text{ and } B = \begin{pmatrix} [0,1] & [0,2] & 0 & 0 \\ [0,3] & [0,4] & 0 & 0 \\ \hline 0 & 0 & 0 & 0 \\ 0 & 0 & 0 & 0 \end{pmatrix}$$

in V.

We see

$$A.B = \begin{pmatrix} [0,33] & [0,46] & 0 & 0 \\ [0,5] & [0,8] & 0 & 0 \\ \hline 0 & 0 & 0 & 0 \\ 0 & 0 & 0 & 0 \end{pmatrix} \in W.$$

Thus W is a semilinear algebra over the semifield S of super interval matrices. So W is a pseudo semi linear subalgebra of V of super interval matrices of type I of V over the semifield S.

If V has no pseudo semilinear subalgebra we call the semi vector space to be pseudo simple.

*Example 4.8*: Let

$$V = \left\{ \begin{pmatrix} [0,a_1] & [0,a_2] & [0,a_3] & [0,a_4] & [0,a_5] \\ [0,a_6] & [0,a_7] & [0,a_8] & [0,a_9] & [0,a_{10}] \\ [0,a_{11}] & [0,a_{12}] & [0,a_{13}] & [0,a_{14}] & [0,a_{15}] \\ [0,a_{16}] & [0,a_{17}] & [0,a_{18}] & [0,a_{19}] & [0,a_{20}] \\ [0,a_{21}] & [0,a_{22}] & [0,a_{23}] & [0,a_{24}] & [0,a_{25}] \end{pmatrix} \middle| a_i \in Q^+ \cup \{0\}, 1 \le i \le 25 \right\}$$

be a semivector space of super interval matrices of type I over the semifield $S = Z^+ \cup \{0\}$.

Consider

$$W = \left\{ \begin{pmatrix} 0 & 0 & 0 & 0 & 0 \\ 0 & 0 & [0,a_1] & [0,a_2] & 0 \\ 0 & 0 & [0,a_3] & [0,a_4] & 0 \\ 0 & 0 & 0 & 0 & 0 \\ 0 & 0 & 0 & 0 & 0 \end{pmatrix} \middle| a_i \in Q^+ \cup \{0\}, 1 \le i \le 4 \right\} \subseteq W.$$

It is easily verified W is a semivector subspace of V of super interval matrices of type I over the semifield $S = Z^+ \cup \{0\}$.

Consider

$$A = \begin{pmatrix} 0 & 0 & 0 & 0 & 0 \\ 0 & 0 & [0,3] & [0,1] & 0 \\ 0 & 0 & [0,8] & [0,2] & 0 \\ 0 & 0 & 0 & 0 & 0 \\ 0 & 0 & 0 & 0 & 0 \end{pmatrix} \text{ and } B = \begin{pmatrix} 0 & 0 & 0 & 0 & 0 \\ 0 & 0 & [0,4] & [0,5] & 0 \\ 0 & 0 & [0,7] & [0,7] & 0 \\ 0 & 0 & 0 & 0 & 0 \\ 0 & 0 & 0 & 0 & 0 \end{pmatrix}$$

in W.

Consider

$$A.B = \begin{pmatrix} 0 & 0 & 0 & 0 & 0 \\ 0 & 0 & [0,19] & [0,22] & 0 \\ 0 & 0 & [0,46] & [0,54] & 0 \\ 0 & 0 & 0 & 0 & 0 \\ 0 & 0 & 0 & 0 & 0 \end{pmatrix} \in W.$$

Thus we see a product '.' is well defined on W. Hence W is a pseudo semilinear subalgebra of V of type I over S.

*Example 4.9*: Let

$$V = \left\{ \begin{pmatrix} [0,a_1] & [0,a_2] \\ [0,a_3] & [0,a_4] \\ [0,a_5] & [0,a_6] \\ [0,a_7] & [0,a_8] \\ [0,a_9] & [0,a_{10}] \end{pmatrix} \middle| a_i \in Q^+ \cup \{0\}, 1 \le i \le 10 \right\}$$

be a semivector space of super interval matrices. V has no pseudo semilinear subalgebra so V is a pseudo simple semivector space of super interval matrices over the semifield $S = Z^+ \cup \{0\}$ of type I.

All these results can be easily worked out for type II semivector spaces with appropriate modifications. This task is left as an exercise to the reader.

Another way to study semivector spaces of super interval matrices is to introduce the notion of set semivector spaces of super interval matrices or set super semivector spaces. For these are very strong algebraic structures to derive properties and for it to find applications is we need to make them less strong algebraically. This is done by introducing set semivector spaces of super interval matrices.

**DEFINITION 4.3**: *Let S be any set and V be a set of super interval matrices. We say V is a set vector space (set super vector space) of super interval matrices over the set S if for all $v \in V$ and $s \in S$, vs and sv are in V.*

We will illustrate this by some simple examples.

***Example 4.10***: Let $V = \{([0, a_1] | [0, a_2] [0, a_3] | [0, a_4]), ([0, a_1] [0, a_2] | [0, a_3]), ([0, a_1] [0, a_2] [0, a_3] | [0, a_4] | [0, a_5] [0, a_6]) | a_i \in Z^+ \cup \{0\}; 1 \le i \le 6\}$ be a set of super interval row matrices. V is a set vector space over the set $S = \{0, 1, 2, 3, 4, 5\}$.

***Example 4.11***: Let

$$V = \left\{ \begin{bmatrix} [0,a_1] \\ \overline{[0,a_2]} \\ [0,a_3] \\ [0,a_4] \end{bmatrix}, ([0,a_1] \ [0,a_2] \ [0,a_3] | [0,a_4]), \begin{pmatrix} [0,a_1] | [0,a_2] \\ [0,a_3] | [0,a_4] \\ [0,a_5] | [0,a_6] \\ [0,a_7] | [0,a_8] \end{pmatrix} \middle| a_i \in Z^+ \cup \{0\}; 1 \le i \le 8 \right\}$$

be a set of super interval matrices. V is a set vector space of super interval matrices or set super vector space over the set $S = \{0, 1\}$.

***Example 4.12***: Let

$$V = \left\{ \begin{bmatrix} [0,a_1] & [0,a_2] \\ \overline{[0,a_3]} & \overline{[0,a_4]} \\ [0,a_5] & [0,a_6] \\ [0,a_7] & [0,a_8] \end{bmatrix}, \begin{bmatrix} [0,a_1] | [0,a_2] \\ [0,a_3] | [0,a_4] \\ [0,a_5] | [0,a_6] \\ [0,a_7] | [0,a_8] \end{bmatrix}, \begin{bmatrix} [0,a_1] | [0,a_2] \\ \overline{[0,a_3]} | \overline{[0,a_4]} \\ [0,a_5] | [0,a_6] \\ [0,a_7] | [0,a_8] \end{bmatrix}, \right.$$

$$\begin{bmatrix}[0,a_1]\\ [0,a_2]\\ [0,a_3]\\ [0,a_4]\\ [0,a_5]\end{bmatrix}, ([0,a_1]\ [0,a_2]\ [0,a_3]\ [0,a_4]\ [0,a_5]\,|\,[0,a_6]\ [0,a_7]),$$

$$\left\{\begin{bmatrix}[0,a_1]\,|\,[0,a_5] & [0,a_6]\,|\,[0,a_7]\\ [0,a_2]\,|\,[0,a_8] & [0,a_9]\,|\,[0,a_{10}]\\ [0,a_3]\,|\,[0,a_{11}] & [0,a_{12}]\,|\,[0,a_{13}]\\ [0,a_4]\,|\,[0,a_{14}] & [0,a_{15}]\,|\,[0,a_{16}]\end{bmatrix}\,\bigg|\,a_i \in Q^+ \cup \{0\}, 1 \le i \le 16\right\}$$

be a set of super interval matrices. Choose the set $S = \{5Z^+ \cup \{0\}, 2Z^+ \cup \{0\}\}$.

Clearly V is a set super vector space of super interval matrices over the set S. Now from these examples it is clear that we do not demand the perfect algebraic structure in V or in the set S just as we like we can induct any number of super interval matrices depending on demand and use a set S so that V is a set super vector space over S.

This sort of flexibility cannot be enjoyed by the semivector spaces of super interval matrices defined in earlier chapters of this book.

*Example 4.13*:

$$\left\{\begin{bmatrix}[0,a_1]\,|\,[0,a_2] & [0,a_3]\\ [0,a_4]\,|\,[0,a_5] & [0,a_6]\\ [0,a_7]\,|\,[0,a_8] & [0,a_9]\end{bmatrix},\right.$$

$$\begin{bmatrix}[0,a_1]\\ [0,a_2]\\ [0,a_3]\\ [0,a_4]\\ [0,a_5]\end{bmatrix}, ([0,a_1]\,|\,[0,a_2]\ [0,a_3]\,|\,[0,a_4]\ [0,a_5]\,|\,[0,a_6]),$$

$$\left\{ \left[ \begin{array}{c|c} [0,a_1] & [0,a_2] \\ \hline [0,a_3] & [0,a_4] \end{array} \right], \left[ \begin{array}{cc|cc} [0,a_1] & [0,a_1] & [0,a_1] & [0,a_1] \\ [0,a_1] & [0,a_1] & [0,a_1] & [0,a_1] \\ [0,a_1] & [0,a_1] & [0,a_1] & [0,a_1] \end{array} \right] \middle| a_i \in 3Z^+ \cup 5Z^+ \cup 17Z^+ \cup \{0\}, 1 \le i \le 9 \right\}$$

be a set of super interval matrices. We see V is a set super vector space over the set S = {0, 1, 3, 5, 17, 6, 10, 15}. We see the number of elements in V is infinite. So V is a super vector space of super interval matrices of infinite cardinality.

*Example 4.14*: Let

$$V = \{([0, 7] \mid 0 \mid [0, 2] \; [0, 5] \; [0, 6] \mid [0, 9]), (0 \mid 0 \mid 0 \; 0 \; 0 \mid 0),$$

$$\left[ \begin{array}{c|c} [0,5] & [0,2] \\ \hline [0,8] & [0,4] \\ \hline [0,12] & [0,7] \end{array} \right], \left[ \begin{array}{c|c} 0 & 0 \\ 0 & 0 \\ 0 & 0 \end{array} \right],$$

$$\left[ \begin{array}{c|ccccc|cc} [0,9] & [0,1] & [0,8] & [0,18] & [0,5] & 0 & [0,5] \\ [0,7] & 0 & [0,11] & [0,7] & [0,26] & [0,3] & [0,2] \\ [0,6] & [0,4] & [0,13] & [0,21] & [0,1] & [0,7] & [0,11] \end{array} \right],$$

$$\left[ \begin{array}{c|ccccc|cc} 0 & 0 & 0 & 0 & 0 & 0 & 0 \\ 0 & 0 & 0 & 0 & 0 & 0 & 0 \\ 0 & 0 & 0 & 0 & 0 & 0 & 0 \end{array} \right] \right\}$$

be a set of super interval matrices. V is a set super vector space of dimension or cardinality 4 over the set S = {0, 1}.

*Example 4.15*: Let V = {([0, $a_1$] [0, $a_2$] | [0, $a_3$] [0, $a_4$] | [0, $a_5$]), ([0, $a_1$] | [0, $a_2$] | [0, $a_3$] [0, $a_4$] | [0, $a_5$]), ([0, $a_1$] | [0, $a_2$] [0, $a_3$] [0, $a_4$] [0, $a_5$]), ([0, $a_1$] [0, $a_2$] [0, $a_3$] | [0, $a_4$] [0, $a_5$]), ([0, $a_1$] [0, $a_2$] | [0, $a_3$] | [0, $a_4$] [0, $a_5$]) | $a_i$ ∈ $Z^+ \cup \{0\}$} is a set of super interval row matrices.

Clearly V is set super vector space over the set S = {0, 1, 2, 4, 8}.

**Example 4.16**: Let

$$V = \left\{ \begin{bmatrix} [0,a_1] \\ \overline{[0,a_2]} \\ [0,a_3] \\ [0,a_4] \\ [0,a_5] \end{bmatrix}, ([0,a_1] \mid [0,a_2] \; [0,a_3] \; [0,a_4]), \right.$$

$$\left. \begin{bmatrix} [0,a_1] & [0,a_4] & [0,a_7] & [0,a_{10}] & [0,a_{13}] & [0,a_{16}] \\ [0,a_2] & [0,a_5] & [0,a_8] & [0,a_{11}] & [0,a_{14}] & [0,a_{17}] \\ [0,a_3] & [0,a_6] & [0,a_9] & [0,a_{12}] & [0,a_{15}] & [0,a_{18}] \end{bmatrix} \middle| a_i \in Q^+ \cup \{0\}, 1 \leq i \leq 18 \right\}$$

be the collection of super set interval matrices. Clearly V is a set super vector space over the set $S = 5Z^+ \cup \{0\}$.

Now having seen examples of set vector spaces of super interval matrices we now proceed onto describe their substructures.

**DEFINITION 4.4**: *Let V be a set super interval matrix vector space over the set S. Suppose $P \subseteq V$ (P a proper subset of V) and if P is also a set super vector space over the set S then we define P to be a super set vector interval matrix subspace of V over the set S.*

We will first illustrate this situation by some examples.

**Example 4.17**: Let

$$V = \left\{ \begin{bmatrix} [0,a_1] \\ \overline{[0,a_2]} \\ [0,a_3] \\ [0,a_4] \end{bmatrix}, \begin{bmatrix} [0,a_1] & [0,a_2] & [0,a_3] \\ [0,a_4] & [0,a_5] & [0,a_6] \\ [0,a_7] & [0,a_8] & [0,a_9] \\ [0,a_{10}] & [0,a_{11}] & [0,a_{12}] \end{bmatrix}, \right.$$

$([0, a_1] \mid [0, a_2] \; [0, a_3] \mid [0, a_4] \; [0, a_5] \; [0, a_6]) \mid a_i \in Z^+ \cup \{0\}; 1 \leq i \leq 12\}$ be a set super vector space of super interval matrices over the set $S = 3Z^+ \cup \{0\}$.

Consider

$$P = \left\{ \begin{bmatrix} [0,a_1] \\ \overline{[0,a_2]} \\ [0,a_3] \\ [0,a_4] \end{bmatrix}, ([0,a_1] \,|\, [0,a_2] \;\; [0,a_3] \,|\, [0,a_4] \;\; [0,a_5] \;\; [0,a_6]) \,\Big|\, a_i \in Z^+ \cup \{0\}; 1 \leq i \leq 6 \right\}$$

$\subseteq V$, P is a set super vector subspace of super interval matrices over the set $S = 3Z^+ \cup \{0\}$.

*Example 4.18*: Let

$$V = \left\{ ([0,a_1] \,|\, [0,a_2] \;\; [0,a_3] \;\; [0,a_4] \,|\, [0,a_5]), \right.$$

$$\begin{bmatrix} [0,a_1] & [0,a_2] & [0,a_3] & [0,a_4] \\ \hline [0,a_5] & [0,a_9] & [0,a_{10}] & [0,a_{11}] \\ [0,a_6] & [0,a_{12}] & [0,a_{13}] & [0,a_{14}] \\ [0,a_7] & [0,a_{15}] & [0,a_{16}] & [0,a_{17}] \\ [0,a_8] & [0,a_{18}] & [0,a_{19}] & [0,a_{20}] \end{bmatrix},$$

$$\left. \begin{pmatrix} [0,a_1] & [0,a_3] & [0,a_5] & [0,a_6] & [0,a_7] & [0,a_8] \\ \hline [0,a_2] & [0,a_4] & [0,a_9] & [0,a_{10}] & [0,a_{11}] & [0,a_{12}] \end{pmatrix} \,\Big|\, a_i \in Z^+ \cup \{0\}; 1 \leq i \leq 20 \right\}$$

be a set super vector space of super interval matrices over the set $S = 13Z^+ \cup \{0\}$.

Consider

$$W = \left\{ ([0,a_1] \,|\, [0,a_2] \;\; [0,a_3] \;\; [0,a_4] \,|\, [0,a_5]), \right.$$

$$\left. \begin{pmatrix} [0,a_1] & [0,a_3] & [0,a_5] & [0,a_6] & [0,a_7] & [0,a_8] \\ \hline [0,a_2] & [0,a_4] & [0,a_9] & [0,a_{10}] & [0,a_{11}] & [0,a_{12}] \end{pmatrix} \,\Big|\, a_i \in 13Z^+ \cup \{0\}; 1 \leq i \leq 12 \right\}$$

$\subseteq V$, W is a set super interval vector subspace of V over the set S.

***Example 4.19***: Let

$$V = \left\{ ([0, a_1] \mid [0, a_2] [0, a_3] \mid [0, a_4] [0, a_5] [0, a_6] [0, a_7]), \right.$$

$$\begin{bmatrix} [0,a_1] \\ \overline{[0,a_2]} \\ [0,a_3] \\ \overline{[0,a_4]} \\ [0,a_5] \\ [0,a_6] \\ [0,a_7] \end{bmatrix}, \begin{bmatrix} [0,a_1] & [0,a_2] & [0,a_3] & [0,a_4] & [0,a_5] \\ [0,a_2] & [0,a_7] & [0,a_8] & [0,a_9] & [0,a_{10}] \\ [0,a_3] & [0,a_{11}] & [0,a_{12}] & [0,a_{13}] & [0,a_{14}] \\ [0,a_4] & [0,a_{15}] & [0,a_{16}] & [0,a_{17}] & [0,a_{18}] \\ [0,a_5] & [0,a_{19}] & [0,a_{20}] & [0,a_{21}] & [0,a_{22}] \\ [0,a_6] & [0,a_{23}] & [0,a_{24}] & [0,a_{25}] & [0,a_{26}] \end{bmatrix},$$

$$\left. \begin{bmatrix} [0,a_1] \mid [0,a_2] & [0,a_3] \\ [0,a_4] \mid [0,a_5] & [0,a_6] \\ [0,a_7] \mid [0,a_8] & [0,a_9] \end{bmatrix} \right| a_i \in 3Z^+ \cup \{0\}; 1 \le i \le 26 \right\}$$

be a set super vector space over the set $S = 15Z^+ \cup \{0\}$ of super interval matrices.

Consider a proper subset

$$W = \left\{ \begin{bmatrix} [0,a_1] \\ \overline{[0,a_2]} \\ [0,a_3] \\ \overline{[0,a_4]} \\ [0,a_5] \\ [0,a_6] \\ [0,a_7] \end{bmatrix}, \begin{bmatrix} [0,a_1] \mid [0,a_2] & [0,a_3] \\ [0,a_4] \mid [0,a_5] & [0,a_6] \\ [0,a_7] \mid [0,a_8] & [0,a_9] \end{bmatrix} \middle| a_i \in 3Z^+ \cup \{0\}; 1 \le i \le 9 \right\} \subseteq V.$$

W is a set super vector subspace of super interval matrices over the set $S = 15Z^+ \cup \{0\}$.

Now having seen just subspace we now proceed onto define the notion of pseudo super set semivector subspace of a super set vector space over the set S.

**DEFINITION 4.5**: *Let V be a set super vector space over the set S. Suppose W ⊆ V be a proper subset such that W is a set super semivector space over the set S (S is also a semifield) then W is defined as the set super semivector subspace of V over the set S.*

We will illustrate this by some examples.

*Example 4.20*: Let

$$V = \left\{ \begin{bmatrix} [0,a_1] & [0,a_2] \\ \hline [0,a_3] & [0,a_4] \\ [0,a_5] & [0,a_6] \end{bmatrix}, ([0, a_1]\ [0, a_2]\ [0, a_3] \mid [0, a_4]\ [0, a_5] \mid [0, a_6]), \right.$$

$$\left. \begin{bmatrix} [0,a_1] & [0,a_2] & [0,a_3] \\ [0,a_4] & [0,a_5] & [0,a_6] \\ \hline [0,a_7] & [0,a_8] & [0,a_9] \end{bmatrix} \middle| a_i \in Z^+ \cup \{0\}; 1 \le i \le 9 \right\}$$

be a set super vector space of super interval matrices over the set $S = Z^+ \cup \{0\}$.

Consider $W = \{(([0, a_1]\ [0, a_2]\ [0, a_3] \mid [0, a_4]\ [0, a_5] \mid [0, a_6]) \mid a_i \in Z^+ \cup \{0\}; 1 \le i \le 6\} \subseteq V$, W is a set super vector subspace of super interval matrices over the set $S = Z^+ \cup \{0\}$.

*Example 4.21*: Let

$$V = \left\{ \begin{bmatrix} [0,a_1] \\ \hline [0,a_2] \end{bmatrix}, ([0, a_1]\ [0, a_2] \mid [0, a_3]\ [0, a_4]\ [0, a_5] \mid [0, a_6]), \right.$$

$$\left. \begin{bmatrix} [0,a_1] & [0,a_2] & [0,a_3] & [0,a_4] & [0,a_5] \\ [0,a_6] & [0,a_7] & [0,a_8] & [0,a_9] & [0,a_{10}] \\ \hline [0,a_{11}] & [0,a_{12}] & [0,a_{13}] & [0,a_{14}] & [0,a_{15}] \\ [0,a_{16}] & [0,a_{17}] & [0,a_{18}] & [0,a_{19}] & [0,a_{20}] \end{bmatrix} \middle| a_i \in 7Z^+ \cup \{0\}, 1 \le i \le 20 \right\}$$

be a set super vector space of super interval matrices over the set $S = \{0, 1\}$. Consider

$$P = \left\{ \begin{bmatrix} [0, a_1] \\ [0, a_2] \end{bmatrix}, ([0, a_1]\ [0, a_2] \mid [0, a_3]\ [0, a_4]\ [0, a_5] \mid [0, a_6]) \mid \right.$$

$a_i \in 7Z^+ \cup \{0\}; 1 \leq i \leq 6\} \subseteq V$, P is only a set super vector subspace of V over the set $S = \{0, 1\}$.

Now if we take

$$A = \left\{ \begin{bmatrix} [0, a_1] \\ [0, a_2] \end{bmatrix} \middle| a_i \in 7Z^+ \cup \{0\}; 1 \leq i \leq 2 \right\} \subseteq V.$$

A is still only a set super vector subspace of V over the set $S = \{0, 1\}$. V can never be a semivector subspace of V as $S = \{0, 1\}$ is only a set never a subfield.

**Remark**: If we replace $S = \{0, 1\}$ by $Z^+ \cup \{0\}$ then A will be a pseudo semivector subspace of V over the semifield $Z^+ \cup \{0\}$.

Now we proceed onto give examples of set super linear algebra of super interval matrices.

**Example 4.22**: Let

$$V = \left\{ \begin{bmatrix} [0, a_1] & [0, a_2] \\ [0, a_3] & [0, a_4] \\ [0, a_5] & [0, a_6] \\ [0, a_7] & [0, a_8] \end{bmatrix} \middle| a_i \in Z^+ \cup \{0\}; 1 \leq i \leq 8 \right\}$$

be a semigroup of super interval matrices. V is a set super linear algebra of super interval matrices over the set $S = \{0, 1, 3, 4, 8\}$.

We will illustrate this situation by some more examples.

*Example 4.23*:

$$\left\{ \begin{bmatrix} [0,a_1] & [0,a_2] \\ \hline [0,a_3] & [0,a_4] \\ [0,a_5] & [0,a_6] \\ \hline [0,a_7] & [0,a_8] \\ \hline [0,a_9] & [0,a_{10}] \\ [0,a_{11}] & [0,a_{12}] \end{bmatrix} \middle| a_i \in Z^+ \cup \{0\}; 1 \le i \le 12 \right\}$$

is a set super linear algebra (set linear algebra of super interval matrices or set super linear interval algebra) over the set $S = \{0, 2, 4, 6, 8, 10\}$.

*Example 4.24*: Let $V = \{([0, a_1] [0, a_2] [0, a_3] | [0, a_4] [0, a_5] | [0, a_6] [0, a_7] [0, a_8] | a_i \in Q^+ \cup \{0\}; 1 \le i \le 8\}$ is a set linear algebra of super row interval matrices over the set $S = 2Z^+ \cup 5Z^+ \cup 73Z^+ \cup \{0\}$.

*Example 4.25*: Let

$$V = \left\{ \begin{bmatrix} [0,a_1] & [0,a_2] & | & [0,a_9] & [0,a_{10}] & [0,a_{11}] & | & [0,a_{12}] \\ [0,a_3] & [0,a_4] & | & [0,a_{13}] & [0,a_{14}] & [0,a_{15}] & | & [0,a_{16}] \\ [0,a_5] & [0,a_6] & | & [0,a_{17}] & [0,a_{18}] & [0,a_{19}] & | & [0,a_{20}] \\ [0,a_7] & [0,a_8] & | & [0,a_{21}] & [0,a_{22}] & [0,a_{23}] & | & [0,a_{24}] \end{bmatrix} \middle| a_i \in R^+ \cup \{0\}; 1 \le i \le 24 \right\}$$

is a set super linear algebra of super row interval vectors over the set $S = \{\sqrt{3}, 0, \sqrt{7}, \sqrt{5}+1, 3, 2, 7, 19, 7/2, 5/9\}$.

Now we will proceed onto give examples of set linear subalgebra of super interval matrices, subset linear subalgebra of super interval matrices, pseudo set vector subspaces of set linear algebra and subset vector subspace of set vector spaces of super interval matrices. The task of giving definitions are left as exercises to the reader.

*Example 4.26*: Let $V = \{([0, a_1] [0, a_2] | [0, a_3] [0, a_4] [0, a_5] [0, a_6] | [0, a_7] [0, a_8] [0, a_9]) | a_i \in R^+ \cup \{0\}; 1 \le i \le 9\}$ be a set linear algebra over the set $S = 3Z^+ \cup 2Z^+ \cup 5Z^+ \cup 7Z^+ \cup \{0\}$.

Consider W = {([0, $a_1$] [0, $a_2$] | [0, $a_3$] [0, $a_4$] [0, $a_5$] [0, $a_6$] | [0, $a_7$] [0, $a_8$] [0, $a_9$]) | $a_i \in Z^+ \cup \{0\}$; $1 \le i \le 9$} $\subseteq$ V, W is a set linear subalgebra of super interval row matrices of V over the set S.

Consider P = {([0, $a_1$] [0, $a_2$] | 0 0 0 0 | [0, $a_7$] [0, $a_8$] [0, $a_9$]) | $a_i \in R^+ \cup \{0\}$; i=1, 2, 7, 8 and 9} $\subseteq$ V. P is a subset linear subalgebra of the set linear algebra V over the subset T = $3Z^+ \cup 2Z^+ \subseteq$ S.

Consider S = {([0, $a_1$] [0, $a_2$] | 0 0 0 0 | 0 0 0), (0 0 | [0, $b_1$] 0 0 [0, $b_2$] | 0 0 0) | $a_1$ $a_2$ $\in 3Z^+ \cup \{0\}$ and $b_1$ $b_2$ $\in 2Z^+ \cup \{0\}$} $\subseteq$ V. G is only a pseudo set vector subspace of V over the set S.

***Example 4.27***: Let

$$V = \left\{ \begin{bmatrix} [0,a_1] & [0,a_2] \\ [0,a_3] & [0,a_4] \\ [0,a_5] & [0,a_6] \\ \hline [0,a_7] & [0,a_8] \\ [0,a_9] & [0,a_{10}] \\ \hline [0,a_{11}] & [0,a_{12}] \\ [0,a_{13}] & [0,a_{14}] \\ [0,a_{15}] & [0,a_{16}] \end{bmatrix} \middle| a_i \in Z^+ \cup \{0\}; 1 \le i \le 16 \right\}$$

be a set linear algebra super interval column matrices over the set S = {0, 1}.

Consider

$$W = \left\{ \begin{bmatrix} 0 & 0 \\ 0 & 0 \\ [0,a_1] & [0,a_2] \\ \hline 0 & 0 \\ [0,a_3] & [0,a_4] \\ \hline 0 & 0 \\ 0 & 0 \\ [0,a_5] & [0,a_6] \end{bmatrix} \middle| a_i \in Z^+ \cup \{0\}; 1 \le i \le 6 \right\} \subseteq V;$$

W is a set linear subalgebra of super column interval vectors over the set S = {0, 1}. However we see V does not contain any subset linear subalgebra of super column interval vectors.

Consider

$$K = \left\{ \begin{bmatrix} [0,5] & [0,7] \\ 0 & 0 \\ 0 & 0 \\ \hline [0,2] & 0 \\ 0 & [0,1] \\ \hline [0,3] & [0,12] \\ [0,1] & [0,5] \\ 0 & 0 \end{bmatrix}, \begin{bmatrix} [0,3] & [0,1] \\ [0,2] & [0,1] \\ [0,5] & 0 \\ \hline 0 & [0,7] \\ [0,11] & 0 \\ \hline 0 & [0,101] \\ [0,1] & 0 \\ 0 & [0,3] \end{bmatrix}, \begin{bmatrix} 0 & 0 \\ 0 & 0 \\ 0 & 0 \\ \hline [0,7] & [0,9] \\ [0,2] & [0,4] \\ \hline 0 & 0 \\ 0 & 0 \\ [0,7] & [0,1] \end{bmatrix} \right\} \subseteq V,$$

K is pseudo set vector subspace of V over the set S = {0, 1}.

*Example 4.28*: Let

$$V = \left\{ \begin{bmatrix} [0,a_1] & [0,a_4] & [0,a_7] & [0,a_{10}] & [0,a_{13}] & [0,a_{16}] & [0,a_{19}] \\ [0,a_2] & [0,a_5] & [0,a_8] & [0,a_{11}] & [0,a_{14}] & [0,a_{17}] & [0,a_{20}] \\ [0,a_3] & [0,a_6] & [0,a_9] & [0,a_{12}] & [0,a_{15}] & [0,a_{18}] & [0,a_{21}] \end{bmatrix} \middle| a_i \in Q^+ \cup \{0\}; 1 \le i \le 21 \right\}$$

be a set linear algebra of super interval row vectors over the set S = {0, 1}. Consider P =

$$\left\{ \begin{bmatrix} [0,a_1] & [0,a_4] & [0,a_7] & [0,a_{10}] & [0,a_{13}] & [0,a_{16}] & [0,a_{19}] \\ [0,a_2] & [0,a_5] & [0,a_8] & [0,a_{11}] & [0,a_{14}] & [0,a_{17}] & [0,a_{20}] \\ [0,a_3] & [0,a_6] & [0,a_9] & [0,a_{12}] & [0,a_{15}] & [0,a_{18}] & [0,a_{21}] \end{bmatrix} \middle| a_i \in Z^+ \cup \{0\}; 1 \le i \le 21 \right\}$$

⊆ V is a set linear subalgebra of super interval row vectors over the set S of V.

Clearly V has no subset linear subalgebras over the subset of S.
However V has pseudo set vector subspaces.

Consider

$$H = \left\{ \begin{bmatrix} [0,8] & [0,7] & 0 & [0,4] & 0 & [0,8] & 0 \\ [0,12] & [0,4] & 0 & 0 & [0,8] & 0 & [0,11] \\ [0,3] & [0,9] & 0 & 0 & [0,3] & [0,4/7] & 0 \end{bmatrix}, \right.$$

$$\begin{bmatrix} 0 & 0 & [0,3] & [0,7/5] & 0 & 0 & [0,1/2] \\ 0 & [0,11/2] & 0 & 0 & [0,9/7] & [0,1/3] & 0 \\ 0 & 0 & [0,5] & 0 & [0,7/2] & 0 & [0,9/4] \end{bmatrix},$$

$$\left. \begin{bmatrix} 0 & 0 & 0 & 0 & 0 & [0,9] & [0,7] \\ 0 & 0 & 0 & 0 & 0 & [0,12] & 0 \\ 0 & 0 & 0 & 0 & 0 & 0 & [0,27/19] \end{bmatrix} \right\} \subseteq V$$

is a pseudo set vector subspace of V of order 3 over the set $\{0, 1\}$.

Now we proceed onto give examples of subset vector subspaces of a set vector space.

*Example 4.29*: Let

$$V = \left\{ \begin{bmatrix} [0,a_1] & [0,a_2] & [0,a_3] \\ [0,a_4] & [0,a_5] & [0,a_6] \\ [0,a_7] & [0,a_8] & [0,a_9] \end{bmatrix}, ([0, a_1] \right.$$

$$[0, a_2] \mid [0, a_3]\ [0, a_4]\ [0, a_5]\ [0, a_6] \mid [0, a_7]), \begin{bmatrix} [0,b_1] \\ [0,b_2] \\ [0,b_3] \\ [0,b_4] \\ [0,b_5] \\ [0,b_6] \\ [0,b_7] \\ [0,b_8] \end{bmatrix}$$

where $a_i \in \{0, 2, 4, 8, 20, 3, 45\}\ 1 \leq i \leq 9\}$ and $b_i \in \{0, 3, 7, 5, 15, 42, 2, 8, 4, 71, 85\}; 1 \leq i \leq 8\}$ be a set vector space over the set $\{0, 1\} = S$.

Clearly V has no subset vector subspace. Thus we can have set vector spaces which have no subset vector subspaces. Such set vector spaces we call as simple set vector spaces. The V given in example 4.29 is a set simple set vector spaces.

*Example 4.30*: Let

$$V = \left\{ \begin{pmatrix} [0,a_1] & [0,a_2] \\ [0,a_3] & [0,a_4] \\ \hline [0,a_5] & [0,a_6] \\ [0,a_7] & [0,a_8] \\ [0,a_9] & [0,a_{10}] \\ [0,a_{11}] & [0,a_{12}] \end{pmatrix} \right. ,$$

$([0, b_1] [0, b_2] | [0, b_3] [0, b_4] [0, b_5] | [0, b_6] [0, b_7] [0, b_8] [0, b_9]) | a_i \in 3Z^+ \cup \{0\}, 1 \le i \le 12\} \subseteq V$. H is a subset vector subspace of V over the subset $T = 6Z^+ \cup \{0\} \subseteq S = Z^+ \cup \{0\}$.

Interested reader can give more examples of this structure. Now we will just illustrate how set linear transformation of set vector spaces of super interval matrices are defined. Let V and W be any two set vector spaces of super interval matrices defined over the same set S.

Let T be a map from V to W, T is said to be a set linear transformation of V to W if T (v) = w and T (sv) = sT (v) for all v ∈ V, w ∈ W and s ∈ S.

*Example 4.31*: Let

$$V = \left\{ \begin{bmatrix} [0,a_1] \\ [0,a_2] \\ \hline [0,a_3] \\ [0,a_4] \\ \hline [0,a_5] \\ [0,a_6] \end{bmatrix} \right. |$$

$a_i \in 3Z^+ \cup \{0\}, 1 \le i \le 6$ and $([0, a_1] | [0, a_2] [0, a_3] [0, a_4] | [0, a_5] [0, a_6])$ where $a_i \in 5Z^+ \cup \{0\}, 1 \le i \le 6\}$ be a set vector space of super interval matrices defined over the set $S = Z^+ \cup \{0\}$.

$$W = \left\{ \left. \left( \begin{array}{cc|} [0,a_1] & [0,a_2] \\ [0,a_3] & [0,a_4] \\ \hline [0,a_5] & [0,a_6] \\ [0,a_7] & [0,a_8] \\ [0,a_9] & [0,a_{10}] \\ [0,a_{11}] & [0,a_{12}] \end{array} \right) \right| a_i \in 5Z^+ \cup \{0\}; 1 \le i \le 12 \right\}$$

$$\left( \begin{array}{cc|cc|cc|c} [0,a_1] & [0,a_3] & [0,a_5] & [0,a_6] & [0,a_7] & [0,a_8] \\ [0,a_2] & [0,a_4] & [0,a_9] & [0,a_{10}] & [0,a_{11}] & [0,a_{12}] \end{array} \right)$$

where $a_i \in 3Z^+ \cup \{0\}$, $1 \le i \le 12$} be a set vector space of super interval matrices defined over the set $S = Z^+ \cup \{0\}$.

Define a map $T : V \to W$ by

$$T \left( \begin{bmatrix} [0,a_1] \\ [0,a_2] \\ \hline [0,a_3] \\ [0,a_4] \\ [0,a_5] \\ \hline [0,a_6] \end{bmatrix} \right) = \left( \begin{array}{cc|cc|cc|c} [0,a_1] & [0,a_2] & [0,a_3] & [0,a_4] & [0,a_5] & [0,a_6] \\ [0,a_1] & [0,a_2] & [0,a_3] & [0,a_4] & [0,a_5] & [0,a_6] \end{array} \right)$$

and

$$T \left( ([0, a_1] \mid [0, a_2]\ [0, a_3]\ [0, a_4] \mid [0, a_5]\ [0, a_6]) \right) = \left( \begin{array}{cc} [0,a_4] & [0,a_4] \\ \hline [0,a_1] & [0,a_1] \\ [0,a_2] & [0,a_2] \\ [0,a_3] & [0,a_3] \\ \hline [0,a_5] & [0,a_6] \\ [0,a_6] & [0,a_6] \end{array} \right)$$

is a set linear transformation of V to W. Interested reader can give more examples of set linear transformation of super set vector spaces (set vector spaces of super interval matrices).

If in the set linear transformation of V to W; W is replaced by V itself then we call the set linear transformation as the set linear operator.

Let us illustrate the notion of set linear operator of set vector spaces of super interval matrices.

*Example 4.32*: Let

$$V = \left\{ \begin{bmatrix} [0,1] & [0,4] \\ [0,2] & [0,5] \\ [0,3] & [0,6] \end{bmatrix}, \begin{bmatrix} [0,2] & [0,8] \\ [0,4] & [0,10] \\ [0,6] & [0,12] \end{bmatrix}, \begin{bmatrix} 0 & 0 \\ 0 & 0 \\ 0 & 0 \end{bmatrix}, \right.$$

([0, 4] | [0, 6] [0, 8] | [0, 16] [0, 2] [0, 10] [0, 14] | [0, 8]), (0 | 0 0 | 0 0 0 0 | 0), ([0, 2] | [0, 3] [0, 4] | [0, 8] [0, 1] [0, 5] [0, 7] | [0, 4])} be set vector space of super interval matrices over the set S = {0, 1}.

Define a map T : V → V by

$$T \left( \begin{bmatrix} 0 & 0 \\ 0 & 0 \\ 0 & 0 \end{bmatrix} \right) = \begin{bmatrix} 0 & 0 \\ 0 & 0 \\ 0 & 0 \end{bmatrix},$$

T ([0 | 0 0 | 0 0 0 0 | 0]) = [0 | 0 0 | 0 0 0 0 | 0],

$$T \left( \begin{bmatrix} [0,1] & [0,4] \\ [0,2] & [0,5] \\ [0,3] & [0,6] \end{bmatrix} \right) = \begin{bmatrix} [0,2] & [0,8] \\ [0,4] & [0,10] \\ [0,6] & [0,12] \end{bmatrix}$$

and

T (([0, 4] | [0, 6] [0, 8] | [0, 16] [0, 2] [0, 10] [0, 14] | [0, 8])
= ([0, 2] | [0, 3] [0, 4] | [0, 8] [0, 1] [0, 5] [0, 7] | [0, 4]))

$$T \begin{bmatrix} [0,2] & [0,8] \\ [0,4] & [0,10] \\ [0,6] & [0,12] \end{bmatrix} = \left( \begin{bmatrix} [0,1] & [0,4] \\ [0,2] & [0,5] \\ [0,3] & [0,6] \end{bmatrix} \right)$$

and

$$T(([0, 2] \mid [0, 3] [0, 4] \mid [0, 8] [0, 1] [0, 5] [0, 7] \mid [0, 4])$$
$$= (([0, 4] \mid [0, 6] [0, 8] \mid [0, 16] [0, 2] [0, 10] [0, 14] \mid [0, 8])).$$

T is a set linear operator on V.

Interested reader is expected to construct more examples.

Now we leave it as an exercise for the reader to describe and define set linear transformation which is invertible.

We now proceed onto describe the notion of set linear functional of a set vector space of super interval matrices.

Let V be a set vector space over the set S. A set linear transformation of V to S in defined as linear functional on V if $f(c\alpha) = cf(\alpha)$

But since set S is reals and elements of super interval matrices are intervals we cannot define a linear functional.

Next we proceed onto describe semigroup vector space of super interval matrices.

**DEFINITION 4.6**: *Let V be a set of super interval matrices, S any semigroup under addition with zero. We call V to be a semigroup vector space of super interval matrices over S if the following conditions hold good.*

  *(1) $sv \in V$ for all $v \in V$ and $s \in S$.*
  *(2) $0v = 0$ for all $v \in V$ and $0 \in S$*
  *(3) $(s_1 + s_2) v = s_1 v + s_2 v$ for all $s_1, s_2 \in S$ and $v \in V$.*

We will illustrate this situation by some simple examples.

***Example 4.33***: Let

$$V = \left\{ \begin{bmatrix} [0, a_1] \\ [0, a_2] \\ \hline [0, a_3] \\ [0, a_4] \\ \hline [0, a_5] \\ [0, a_6] \end{bmatrix}, \right.$$

([0, $a_1$] [0, $a_2$] [0, $a_3$] [0, $a_4$] | [0, $a_5$] [0, $a_6$] [0, $a_7$] | [0, $a_8$] [0, $a_9$] | [0, $a_{10}$]),

$$\begin{bmatrix} [0,a_1] & [0,a_2] & | [0,a_3] & | [0,a_4] & [0,a_5] \\ [0,a_6] & [0,a_7] & | [0,a_8] & | [0,a_9] & [0,a_{10}] \\ \hline [0,a_{11}] & [0,a_{12}] & | [0,a_{13}] & | [0,a_{14}] & [0,a_{15}] \\ [0,a_{16}] & [0,a_{17}] & | [0,a_{18}] & | [0,a_{19}] & [0,a_{20}] \\ [0,a_{21}] & [0,a_{22}] & | [0,a_{23}] & | [0,a_{24}] & [0,a_{25}] \\ \hline [0,a_{26}] & [0,a_{27}] & | [0,a_{28}] & | [0,a_{29}] & [0,a_{30}] \\ [0,a_{31}] & [0,a_{32}] & | [0,a_{33}] & | [0,a_{34}] & [0,a_{35}] \end{bmatrix}$$  $a_i \in Z^+ \cup \{0\}$,

$1 \le i \le 35$} be a set vector space over $S = Z^+ \cup \{0\}$. Since S is a semigroup under addition we define V to be a semigroup vector space of super interval matrices over the semigroup S.

*Example 4.34*: Let

$$V = \left\{ \begin{bmatrix} [0,a_1] \\ \hline [0,a_2] \\ [0,a_3] \\ \hline [0,a_4] \\ \hline [0,a_5] \\ [0,a_6] \end{bmatrix}, ([0, a_1]\ [0, a_2]\ [0, a_3]\ |\ [0, a_4]\ [0, a_5]\ |\ [0, a_6]\ [0, a_7]), \right.$$

$$\begin{bmatrix} [0,a_1] & | [0,a_2] & [0,a_3] & | [0,a_4] \\ [0,a_5] & | [0,a_6] & [0,a_7] & | [0,a_8] \\ [0,a_9] & | [0,a_{10}] & [0,a_{11}] & | [0,a_{12}] \\ [0,a_{13}] & | [0,a_{14}] & [0,a_{15}] & | [0,a_{16}] \\ [0,a_{17}] & | [0,a_{18}] & [0,a_{19}] & | [0,a_{20}] \end{bmatrix}$$

$a_i \in 3Z^+ \cup \{0\}$, $1 \le i \le 20$} be a semigroup vector space of super interval matrices over the semigroup $S = 3Z^+ \cup \{0\}$.

**Example 4.35**: Let

$$M = \left\{ \begin{bmatrix} [0,a_1] & [0,a_2] \\ \hline [0,a_3] & [0,a_4] \\ \hline [0,a_5] & [0,a_6] \\ [0,a_7] & [0,a_8] \end{bmatrix}, \right.$$

$$\begin{bmatrix} [0,a_1] & [0,a_3] & [0,a_5] & [0,a_6] & [0,a_7] & [0,a_8] & [0,a_{13}] \\ [0,a_2] & [0,a_4] & [0,a_9] & [0,a_{10}] & [0,a_{11}] & [0,a_{12}] & [0,a_{14}] \end{bmatrix} \Big|$$

$a_i \in Q^+ \cup \{0\}\}$ be a semigroup vector space of super matrices over the semigroup $S = Z^+ \cup \{0\}$.

**Example 4.36**: Let

$$M = \left\{ \begin{bmatrix} [0,a_1] & [0,a_2] & [0,a_3] \\ \hline [0,a_4] & [0,a_5] & [0,a_6] \\ [0,a_7] & [0,a_8] & [0,a_9] \end{bmatrix}, \right.$$

$$\begin{bmatrix} [0,a_1] & [0,a_2] & [0,a_3] \\ \hline [0,a_4] & [0,a_5] & [0,a_6] \\ \hline [0,a_7] & [0,a_8] & [0,a_9] \end{bmatrix}, \begin{bmatrix} [0,a_1] & [0,a_2] & [0,a_3] \\ [0,a_4] & [0,a_5] & [0,a_6] \\ \hline [0,a_7] & [0,a_8] & [0,a_9] \end{bmatrix},$$

$$\begin{bmatrix} [0,a_1] & [0,a_2] & [0,a_3] \\ [0,a_4] & [0,a_5] & [0,a_6] \\ [0,a_7] & [0,a_8] & [0,a_9] \end{bmatrix}, \begin{bmatrix} [0,a_1] & [0,a_2] & [0,a_3] \\ \hline [0,a_4] & [0,a_5] & [0,a_6] \\ [0,a_7] & [0,a_8] & [0,a_9] \end{bmatrix},$$

$$\left. \begin{bmatrix} [0,a_1] & [0,a_2] & [0,a_3] \\ [0,a_4] & [0,a_5] & [0,a_6] \\ [0,a_7] & [0,a_8] & [0,a_9] \end{bmatrix} \Big| a_i \in Z_{44}; 1 \le i \le 9 \right\}$$

be a semigroup vector space of super interval matrices over the semigroup $Z_{44} = S$. Clearly M is of finite order.

**Example 4.37**: Let V = {(([0, $a_1$] | [0, $a_2$] [0, $a_3$] | [0, $a_4$] [0, $a_5$]), ([0, $a_1$] [0, $a_2$] | [0, $a_3$] [0, $a_4$] [0, $a_5$]), ([0, $a_1$] [0, $a_2$] [0, $a_3$] | [0, $a_4$] [0, $a_5$]), ([0, $a_1$] | [0, $a_2$] [0, $a_3$] [0, $a_4$] | [0, $a_5$]), ([0, $a_1$] [0, $a_2$] | [0, $a_3$] [0, $a_4$] [0, $a_5$]), ([0, $a_1$] [0, $a_2$] [0, $a_3$] | [0, $a_4$] [0, $a_5$]), ([0, $a_1$] | [0, $a_2$] | [0, $a_3$] [0, $a_4$] [0, $a_5$]), ([0, $a_1$] [0, $a_2$] | [0, $a_3$] [0, $a_4$] | [0, $a_5$]), ([0, $a_1$] | [0, $a_2$] [0, $a_3$] [0, $a_4$] | [0, $a_5$]) | $a_i \in Z_{13}$, $1 \le i \le 5$} be a semigroup vector space of super interval matrices over the semigroup S = ($Z_{13}$, +). Clearly V is of finite order.

Now having seen examples of semigroup vector spaces of super interval matrices we now proceed onto describe their substructures with examples.

**Example 4.38**: Let

$$V = \left\{ \left( \begin{array}{cc} [0,a_1] & [0,a_2] \\ \hline [0,a_3] & [0,a_4] \\ [0,a_5] & [0,a_6] \\ [0,a_7] & [0,a_8] \\ \hline [0,a_9] & [0,a_{10}] \end{array} \right) \middle| a_i \in Z_{25}, 1 \le i \le 10 \right.$$

and ([0, $a_1$] | [0, $a_2$] [0, $a_3$] [0, $a_4$] | [0, $a_5$]) | $a_i \in Z_{25}$, $1 \le i \le 5$} be a semigroup vector space of super interval matrices over the semigroup S = $Z_{25}$. Let H = {([0, $a_1$] | [0, $a_2$] [0, $a_3$] [0, $a_4$] | [0, $a_5$]) | $a_i \in Z_{25}$; $1 \le i \le 5$} $\subseteq$ V; H is a semigroup vector subspace of V over the semigroup S = $Z_{25}$.

**Example 4.39**: Let

$$V = \{([0, a_1] | [0, a_2] [0, a_3] [0, a_4] | [0, a_5]) \begin{bmatrix} [0,a_1] & [0,a_2] & [0,a_3] \\ [0,a_4] & [0,a_5] & [0,a_6] \\ [0,a_7] & [0,a_8] & [0,a_9] \end{bmatrix} |$$

$a_i \in Z^+ \cup \{0\}$, $1 \le i \le 9$} be a semigroup vector space of super interval matrices over the semigroup S=$5Z^+\cup\{0\}$.

Consider

$$P = \left\{ \left. \begin{bmatrix} [0,a_1] & [0,a_2] & | & [0,a_3] \\ [0,a_4] & [0,a_5] & | & [0,a_6] \\ [0,a_7] & [0,a_8] & | & [0,a_9] \end{bmatrix} \right| a_i \in 5Z^+ \cup \{0\}, 1 \leq i \leq 9 \right\} \subseteq V,$$

P is a semigroup vector subspace of super interval matrices over the semigroup $S = 5Z^+ \cup \{0\}$.

Now we proceed onto give examples of subsemigroup vector subspace of a semigroup vector subspace.

*Example 4.40*: Let

$$V = \left\{ \left( \begin{array}{cc} [0,a_1] & [0,a_2] \\ [0,a_3] & [0,a_4] \end{array} \right), \begin{bmatrix} [0,a_1] & | & [0,a_2] & | & [0,a_3] & [0,a_4] \\ [0,a_5] & | & [0,a_6] & | & [0,a_7] & [0,a_8] \\ [0,a_9] & | & [0,a_{10}] & | & [0,a_{11}] & [0,a_{12}] \end{bmatrix}, \right.$$

$$\left. \begin{bmatrix} [0,a_1] & | & [0,a_2] \\ [0,a_3] & | & [0,a_4] \end{bmatrix}, ([0, a_1] [0, a_2] [0, a_3] | [0, a_4]) \, \middle| \, a_i \in Z_{42}, 1 \leq i \leq 12 \right\}$$

be a semigroup vector space over the semigroup $S = Z_{42}$.
Consider

$$P = \left\{ \left. \begin{bmatrix} [0,a_1] & | & [0,a_2] \\ [0,a_3] & | & [0,a_4] \end{bmatrix}, \left( \begin{array}{cc} [0,a_1] & [0,a_2] \\ [0,a_3] & [0,a_4] \end{array} \right) \right| a_i \in Z_{42}, 1 \leq i \leq 4 \right\} \subseteq V$$

and $T = \{0, 2, 4, 6, 8, \ldots, 40\} \subseteq Z_{42}$ be a subsemigroup of $S = Z_{42}$. P is a subsemigroup vector subspace of V over the subsemigroup P of S.

*Example 4.41*: Let

$$V = \{([0, a_1] [0, a_2] | [0, a_3] [0, a_4] [0, a_5] | [0, a_6] \ldots [0, a_{12}]),$$

$$\begin{bmatrix} [0,a_1] & [0,a_2] & | & [0,a_5] \\ [0,a_3] & [0,a_4] & | & [0,a_6] \\ [0,a_7] & [0,a_8] & | & [0,a_{12}] \\ [0,a_9] & [0,a_{10}] & | & [0,a_{11}] \\ \hline [0,a_{13}] & [0,a_{14}] & | & [0,a_{15}] \\ [0,a_{16}] & [0,a_{17}] & | & [0,a_{18}] \end{bmatrix}, \begin{bmatrix} [0,a_1] & | & [0,a_4] \\ [0,a_2] & | & [0,a_5] \\ [0,a_3] & | & [0,a_6] \end{bmatrix}$$

$a_i \in Z^+ \cup \{0\}$, $1 \leq i \leq 18\}$ be a semigroup vector space over the semigroup $S = 3Z^+ \cup \{0\}$.

Consider

$H = \{([0, a_1] [0, a_2] | [0, a_3] [0, a_4] [0, a_5] | [0, a_6] \dots [0, a_{12}])$,

$$\begin{bmatrix} [0,a_1] & | & [0,a_4] \\ [0,a_2] & | & [0,a_5] \\ [0,a_3] & | & [0,a_6] \end{bmatrix} \mid a_i \in 3Z^+ \cup \{0\}, 1 \leq i \leq 12\} \subseteq V,$$

H is a subsemigroup vector subspace of V over the subsemigorup $T = 6Z^+ \cup \{0\} \subseteq S$.

Now having seen examples of substructures of a semigroup vector space we now proceed onto define the notion of linearly independent set in the semigroup vector space.

**DEFINITION 4.7**: *Let V be a semigroup vector space of super interval matrices over the semigroup S. A set of vectors $\{v_1, v_2, \dots, v_n\}$ in V is said to be a semigroup linearly independent set if $v_i \neq sv_j$ for any $s \in S$ for $i \neq j$; $1 \leq i, j \leq n$.*

We will first illustrate this situation by some examples.

*Example 4.42*: Let

$V = \{([0, a_1] [0, a_2] | [0, a_3] [0, a_4] [0, a_5] [0, a_6])$,

$$\begin{bmatrix} [0,a_1] & [0,a_2] & [0,a_3] \\ [0,a_4] & [0,a_5] & [0,a_6] \\ [0,a_7] & [0,a_8] & [0,a_9] \\ \hline [0,a_{10}] & [0,a_{11}] & [0,a_{12}] \end{bmatrix} \Bigg|$$

$a_i \in Z^+ \cup \{0\}$, $1 \le i \le 12\}$ be a semigroup vector space over the semigroup $S = Z^+ \cup \{0\}$. Consider $X = \{([0, 1] \ 0 \ | \ 0 \ 0 \ 0 \ 0)$, $(0, [0, 1] \ | \ 0 \ 0 \ 0 \ 0)$, $(0 \ 0 \ | \ 0 \ [0, 1] \ 0 \ 0\ )$, $(0 \ 0 \ | \ 0 \ 0 \ 0 \ [0, 1])$,

$$\begin{bmatrix} [0,1] & 0 & 0 \\ 0 & 0 & 0 \\ 0 & 0 & 0 \\ \hline 0 & 0 & 0 \end{bmatrix} \begin{bmatrix} 0 & 0 & 0 \\ [0,1] & 0 & 0 \\ 0 & 0 & 0 \\ \hline 0 & 0 & 0 \end{bmatrix} \begin{bmatrix} 0 & 0 & [0,1] \\ 0 & 0 & 0 \\ 0 & 0 & 0 \\ \hline 0 & 0 & 0 \end{bmatrix},$$

$$\begin{bmatrix} 0 & 0 & 0 \\ 0 & 0 & 0 \\ 0 & 0 & 0 \\ \hline 0 & 0 & [0,1] \end{bmatrix} \begin{bmatrix} 0 & 0 & 0 \\ 0 & 0 & [0,1] \\ 0 & 0 & 0 \\ \hline 0 & 0 & 0 \end{bmatrix} \Bigg\} \subseteq V$$

be a set of vectors in V. X is a linearly independent set of V.

*Example 4.43*: Let

$$M = \left\{ \begin{bmatrix} [0,a_1] \\ [0,a_2] \\ \hline [0,a_3] \\ [0,a_4] \\ [0,a_5] \\ [0,a_6] \end{bmatrix}, ([0, a_1] \ [0, a_2] \ [0, a_3] \ [0, a_4] \ | \ [0, a_5] \ [0, a_6]), \right.$$

$$V = \left\{ \begin{bmatrix} [0,a_1] & [0,a_4] \\ [0,a_2] & [0,a_5] \\ [0,a_3] & [0,a_6] \end{bmatrix} \mid a_i \in Z_{45},\ 1 \leq i \leq 6 \right\}$$

be a semigroup vector space of super interval matrices over the semigroup $S = Z_{45}$. Let

$$y = \left\{ \begin{bmatrix} 0 \\ [0,1] \\ 0 \\ 0 \\ 0 \\ 0 \end{bmatrix}, \begin{bmatrix} 0 \\ 0 \\ [0,7] \\ 0 \\ 0 \\ 0 \end{bmatrix}, \begin{bmatrix} 0 \\ 0 \\ 0 \\ [0,4] \\ 0 \\ 0 \end{bmatrix}, \begin{bmatrix} 0 \\ 0 \\ 0 \\ 0 \\ 0 \\ [0,13] \end{bmatrix}, \right.$$

$([0, 1]\ 0\ 0\ 0\ |\ 0\ 0)$, $(0\ 0\ 0\ [0, 1]\ |\ 0\ 0 )$, $(0\ 0\ [0, 1]\ 0\ |\ 0\ 0)$,
$(0\ 0\ 0\ 0\ |\ [0, 1]\ 0)$, $(0\ 0\ 0\ 0\ |\ 0\ [0, 1])$,

$$\left. \begin{bmatrix} [0,1] & 0 \\ 0 & 0 \\ 0 & 0 \end{bmatrix}, \begin{bmatrix} 0 & 0 \\ 0 & 0 \\ 0 & [0,1] \end{bmatrix}, \begin{bmatrix} 0 & 0 \\ 0 & 0 \\ 0 & [0,1] \end{bmatrix}, \begin{bmatrix} 0 & [0,1] \\ [0,1] & 0 \\ 0 & 0 \end{bmatrix} \right\} \subseteq V,$$

X is a linearly independent set of the semigroup super interval vector space V over the semigroup S.

*Example 4.44*: Let

$$V = \left\{ \begin{bmatrix} [0,a_1] \\ [0,a_2] \\ [0,a_3] \\ [0,a_4] \\ [0,a_5] \\ [0,a_6] \\ [0,a_7] \\ [0,a_8] \end{bmatrix}, \right.$$

$([0, a_1]\ [0, a_2]\ |\ [0, a_3]\ [0, a_4]\ |\ [0, a_5])\ |\ a_i \in Z_{12},\ 1 \leq i \leq 8\}$ be a semigroup vector space over the semigroup $S = Z_{12}$.

Consider

$$X = \left\{ \begin{bmatrix} 0 \\ 0 \\ 0 \\ \overline{[0,3]} \\ 0 \\ 0 \\ 0 \\ 0 \end{bmatrix}, \begin{bmatrix} 0 \\ 0 \\ 0 \\ \overline{[0,6]} \\ 0 \\ 0 \\ 0 \\ 0 \end{bmatrix}, \begin{bmatrix} 0 \\ 0 \\ 0 \\ 0 \\ \overline{[0,9]} \\ 0 \\ 0 \\ 0 \end{bmatrix}, \begin{bmatrix} 0 \\ 0 \\ 0 \\ 0 \\ \overline{[0,3]} \\ 0 \\ 0 \\ 0 \end{bmatrix}, \right.$$

$([0, 1]\ 0\ |\ 0\ 0\ 0)\ (0\ 0\ |\ 0\ 0\ [0, 4]),\ (0\ 0\ |\ 0\ 0\ [0, 2]),\ (0\ 0\ |\ 0\ 0\ [0, 1])\ ([0, 9]\ 0\ |\ 0\ 0\ 0)\} \subseteq V$ is not a linearly independent set of the semigroup super interval vector space V over the semigroup $Z_{12}$.

Now we will proceed onto define the notion of generating subset of a semigroup vector space over the semigroup S.

**DEFINITION 4.8**: *Let V be a semigroup vector space over the semigroup S under addition. Let $T = \{v_1, v_2, ..., v_n\} \subseteq V$ be a subset of V. We say T generates the semigroup vector space V over S if every element $v \in V$ can be got as $v = sv_i$, $v_i \in T$ and $s \in S$.*

We will illustrate this situation by some examples.

***Example 4.45***: Let V = {([0, a] | [0, a] [0, a] [0, a] [0, a] | [0, a]),

$$\begin{bmatrix} [0,a] \\ \overline{[0,a]} \\ \overline{[0,a]} \\ [0,a] \\ [0,a] \end{bmatrix}, \begin{bmatrix} [0,a] & [0,a] & 0 & 0 \\ [0,a] & [0,a] & 0 & 0 \\ \hline 0 & 0 & [0,a] & [0,a] \\ 0 & 0 & [0,a] & [0,a] \end{bmatrix}$$

$a \in Z^+ \cup \{0\}\}$ be a semigroup super interval matrix vector space over the semigroup $S = Z^+ \cup \{0\}$.

Consider

$$H = \{([0, 1] \mid [0, 1] \ [0, 1] \ [0, 1] \ [0, 1] \mid [0, 1]),$$

$$\begin{bmatrix} [0,1] \\ \overline{[0,1]} \\ \overline{[0,1]} \\ [0,1] \\ [0,1] \end{bmatrix}, \begin{bmatrix} [0,1] & [0,1] & 0 & 0 \\ [0,1] & [0,1] & 0 & 0 \\ 0 & 0 & [0,1] & [0,1] \\ 0 & 0 & [0,1] & [0,1] \end{bmatrix} \}$$

$\subseteq V$ generates V over the semigroup $S = Z^+ \cup \{0\}$.

*Example 4.46*: Let

$$V = \left\{ \left( \begin{array}{c|c} [0,a] & 0 \\ \hline 0 & [0,a] \end{array} \right), \begin{pmatrix} [0,a] & [0,a] & [0,a] & [0,a] \\ [0,a] & [0,a] & [0,a] & [0,a] \\ [0,a] & [0,a] & [0,a] & [0,a] \end{pmatrix}, \begin{bmatrix} [0,a] & [0,a] \\ [0,a] & [0,a] \\ 0 & [0,a] \\ \hline [0,a] & 0 \\ 0 & [0,a] \\ [0,a] & [0,a] \\ [0,a] & [0,a] \end{bmatrix} \right| $$

$a \in Z_{140}\}$ be a semigroup super interval matrix vector space over the semigroup $S = Z_{140}$.

Consider

$$X = \left\{ \left( \begin{array}{c|c} [0,1] & 0 \\ \hline 0 & [0,1] \end{array} \right), \begin{pmatrix} [0,1] & [0,1] & [0,1] & [0,1] \\ [0,1] & [0,1] & [0,1] & [0,1] \\ [0,1] & [0,1] & [0,1] & [0,1] \end{pmatrix}, \right.$$

$$\left.\begin{bmatrix} [0,1] & [0,1] \\ [0,1] & [0,1] \\ 0 & [0,1] \\ \hline [0,1] & 0 \\ 0 & [0,1] \\ [0,1] & [0,1] \\ [0,1] & [0,1] \end{bmatrix}\right\} \subseteq V,$$

X generates the semigroup vector space over the semigroup $S = Z_{140}$.

Now we proceed onto define the notion of semigroup linear algebra V over a semigroup S.

**DEFINITION 4.9**: *Let V be a semigroup super matrix interval vector space over the semigroup S. Suppose V is a semigroup under addition then we call V to be a semigroup super matrix interval linear algebra over the semigroup S.*

We will illustrate this situation by some examples.

*Example 4.47*: Let

$$V = \left\{ \begin{bmatrix} [0,a_1] & [0,a_2] & [0,a_3] \\ \hline [0,a_4] & [0,a_8] & [0,a_9] \\ [0,a_5] & [0,a_{10}] & [0,a_{11}] \\ [0,a_6] & [0,a_{12}] & [0,a_{13}] \\ [0,a_7] & [0,a_{14}] & [0,a_{15}] \end{bmatrix} \right.$$

where $a_i \in Z^+ \cup \{0\}$, $1 \le i \le 15\}$ is a semigroup super interval matrix linear algebra over the semigroup $S = Z^+ \cup \{0\}$.

*Example 4.48*: Let $V = \{([0, a_1] [0, a_2] \mid [0, a_3] [0, a_4] [0, a_5] [0, a_6] \mid [0, a_7] [0, a_8]) \mid a_i \in Z_{240}, 1 \le i \le 8\}$ is a semigroup interval super matrix linear algebra over the semigroup $S = Z_{240}$.

*Example 4.49*: Let

$$P = \left\{ \begin{bmatrix} [0,a_1] \\ [0,a_2] \\ \hline [0,a_3] \\ [0,a_4] \\ [0,a_5] \\ \hline [0,a_6] \\ [0,a_7] \\ \hline [0,a_8] \end{bmatrix} \mid a_i \in Z_{29}, 1 \le i \le 8 \right\}$$

be a semigroup interval super matrix linear algebra over the semigroup $S = Z_{29}$.

*Example 4.50*: Let

$$V = \left\{ \begin{bmatrix} [0,a_1] & [0,a_6] & [0,a_{11}] & [0,a_{16}] & [0,a_{21}] & [0,a_{26}] \\ [0,a_2] & [0,a_7] & [0,a_{12}] & [0,a_{17}] & [0,a_{22}] & [0,a_{27}] \\ [0,a_3] & [0,a_8] & [0,a_{13}] & [0,a_{18}] & [0,a_{23}] & [0,a_{28}] \\ [0,a_4] & [0,a_9] & [0,a_{14}] & [0,a_{19}] & [0,a_{24}] & [0,a_{29}] \\ [0,a_5] & [0,a_{10}] & [0,a_{15}] & [0,a_{20}] & [0,a_{25}] & [0,a_{30}] \end{bmatrix} \mid a_i \in Z_{25}, 1 \le i \le 30 \right\}$$

be a semigroup super interval matrix linear algebra over the semigroup $S = Z_{25}$.

It is easily verified that every semigroup linear algebra over a semigroup is a semigroup vector space over the semigroup. But a semigroup vector space in general is not a semigroup linear algebra.

Now we proceed onto give examples of semigroup linear subalgebra over a semigroup S and subsemigroup linear subalgebra over a semigroup of a semigroup linear algebra.

*Example 4.51*: Let

$$V = \left\{ \begin{bmatrix} \dfrac{[0,a_1]}{[0,a_2]} \\ [0,a_3] \\ \dfrac{[0,a_4]}{[0,a_5]} \\ [0,a_6] \\ [0,a_7] \\ [0,a_8] \end{bmatrix} \,\middle|\, a_i \in Z_{24},\ 1 \leq i \leq 8 \right\}$$

be a semigroup interval super matrix linear algebra over the semigroup $Z_{24} = S$.

Consider

$$H = \left\{ \begin{bmatrix} \dfrac{[0,a_1]}{0} \\ [0,a_2] \\ \dfrac{0}{[0,a_3]} \\ 0 \\ [0,a_4] \\ 0 \end{bmatrix} \,\middle|\, a_i \in Z_{24},\ 1 \leq i \leq 4 \right\} \subseteq V,$$

is a semigroup super interval matrix linear subalgebra over the semigroup $S = Z_{24}$. Take

$$P = \left\{ \begin{bmatrix} \dfrac{[0,a_1]}{[0,a_2]} \\ 0 \\ \dfrac{[0,a_3]}{0} \\ 0 \\ 0 \\ [0,a_4] \end{bmatrix} \,\middle|\, a_i \in Z_{24};\ 1 \leq i \leq 4 \right\} \subseteq V,$$

is a subsemigroup interval super matrix linear subalgebra of V over the subsemigroup $T = \{0, 2, 4, 6, 8, 10, 12, \ldots, 20, 22\} \subseteq Z_{24} = S$.

***Example 4.52***: Let $V = \{([0, a_1] [0, a_2] | [0, a_3] [0, a_4] [0, a_5] | [0, a_6] | [0, a_7] [0, a_8] | [0, a_9] [0, a_{10}]) | a_i \in Z^+ \cup \{0\}, 1 \leq i \leq 10\}$ be a semigroup linear algebra over the semigroup $S = Z^+ \cup \{0\}$.

Consider $W = \{(0\ 0 | [0, a_1] [0, a_2] [0, a_3] | 0\ 0\ 0 | [0, a_4] [0, a_5]) | a_i \in Z^+ \cup \{0\}, 1 \leq i \leq 5\} \subseteq V$; W is a semigroup linear subalgebra of V over the semigroup $S = Z^+ \cup \{0\}$.

***Example 4.53***: Let

$$V = \left\{ \begin{bmatrix} [0,a_1] & [0,a_2] \\ [0,a_3] & [0,a_4] \\ \hline [0,a_4] & [0,a_8] \\ [0,a_5] & [0,a_9] \\ [0,a_6] & [0,a_{10}] \\ \hline [0,a_7] & [0,a_{11}] \end{bmatrix} \middle| a_i \in Z_{25}; 1 \leq i \leq 11 \right\}$$

be a semigroup linear algebra over the semigroup $S = Z_{25}$.

$$P = \left\{ \begin{bmatrix} 0 & [0,b] \\ [0,a] & 0 \\ \hline 0 & [0,c] \\ [0,d] & 0 \\ \hline 0 & [0,e] \\ [0,f] & 0 \end{bmatrix} \middle| a,b,c,d,e,f \in Z_{25} \right\} \subseteq V;$$

P is a subsemigroup super matrix interval linear subalgebra of V over the subsemigroup $T = \{0, 5, 10, 15, 20\} \subseteq Z_{25} = S$.

*Example 4.54*: Let

$$A = \left\{ \begin{bmatrix} [0,a_1] & [0,a_2] & [0,a_3] \\ [0,a_4] & [0,a_5] & [0,a_6] \\ \hline [0,a_7] & [0,a_8] & [0,a_9] \\ [0,a_{10}] & [0,a_{11}] & [0,a_{12}] \\ [0,a_{13}] & [0,a_{14}] & [0,a_{15}] \\ \hline [0,a_{16}] & [0,a_{17}] & [0,a_{18}] \\ [0,a_{19}] & [0,a_{20}] & [0,a_{21}] \\ [0,a_{22}] & [0,a_{23}] & [0,a_{24}] \\ [0,a_{25}] & [0,a_{26}] & [0,a_{27}] \\ [0,a_{28}] & [0,a_{29}] & [0,a_{30}] \end{bmatrix} \right.$$

where $a_i \in Z_{13}$, $1 \le i \le 30$} be a semigroup super matrix interval linear algebra over the semigroup $S = Z_{13}$.
Let

$$P = \left\{ \begin{bmatrix} 0 & [0,a] & 0 \\ 0 & 0 & 0 \\ \hline 0 & 0 & 0 \\ 0 & 0 & 0 \\ [0,d] & [0,b] & [0,c] \\ \hline 0 & [0,g] & 0 \\ 0 & [0,h] & 0 \\ [0,e] & 0 & [0,i] \\ [0,f] & 0 & [0,j] \\ [0,k] & 0 & [0,l] \end{bmatrix} \middle| a,b,c,\ldots,l \in Z_{13} \right\} \subseteq V$$

is a semigroup super matrix interval linear subalgebra of V over the semigroup $S = Z_{13}$. However it is pertinent to mention here that S has no subsemigroups hence V has no subsemigroup linear subalgebras.

Inview of this we have the following theorem.

**THEOREM 4.2**: *Let V be a semigroup super matrix interval linear algebra over the semigroup S = $Z_p$, p a prime. V does not contain any subsemigroup super matrix interval linear subalgebra of V over the semigroup S.*

The proof is direct and hence is left as an exercise to the reader.

**THEOREM 4.3**: *Let V be a semigroup super matrix interval linear algebra of super interval matrices over the semigroup S = $Z_n$ (n < ∞, n a composite number) V has non trivial subsemigroup super matrix interval linear subalgebras over a subsemigroup T of S = $Z_n$.*

This proof is also direct and hence is left as an exercise to the reader.

We define the notion of pseudo semigroup super matrix interval subvector space of a semigroup super matrix interval linear algebra of super internal matrices.

**DEFINITION 4.10**: *Let V be a semigroup super matrix interval linear algebra over the semigroup S. Let P be a proper subset of V which is only a semigroup super matrix interval vector subspace of V over the semigroup S and not a semigroup linear subalgebra of V then we define P to be a pseudo semigroup vector subspace of V over the semigroup S.*

We will illustrate this situation by some examples.

*Example 4.55*: Let

$V = \{([0, a_1] | [0, a_2] [0, a_3] [0, a_4] | [0, a_5]) | a_i \in Z_6, 1 \leq i \leq 5\}$

be a semigroup linear algebra over the semigroup S = $Z_6$.
Consider

$H = \{([0, a] | 0\ 0\ 0 | [0, b]), ([0, 0] | [0, a] [0, b] [0, c] | [0, 0]) |$

a, b, c ∈ {0, 3} ⊆ $Z_6$} ⊆ V, H is only a semigroup vector space over S = $Z_6$. Thus H is a pseudo semigroup vector subspace of V over the semigroup S = $Z_6$.

*Example 4.56*: Let

$$V = \left\{ \begin{bmatrix} [0,a_1] \\ [0,a_2] \\ [0,a_3] \\ [0,a_4] \\ \hline [0,a_5] \\ [0,a_6] \\ [0,a_7] \\ \hline [0,a_8] \\ [0,a_9] \\ \hline [0,a_{10}] \\ [0,a_{11}] \end{bmatrix} \;\middle|\; a_i \in Z_{10}; 1 \leq i \leq 11 \right\}$$

be a semigroup linear algebra of super column interval matrices over the semigroup $S = Z_{10}$.

Consider

$$H = \left\{ \begin{bmatrix} [0,a_1] \\ 0 \\ 0 \\ 0 \\ \hline 0 \\ [0,a_2] \\ [0,a_3] \\ [0,a_4] \\ \hline 0 \\ 0 \\ 0 \end{bmatrix}, \begin{bmatrix} 0 \\ [0,a_1] \\ 0 \\ [0,a_2] \\ \hline 0 \\ 0 \\ 0 \\ 0 \\ \hline [0,a_3] \\ 0 \\ 0 \end{bmatrix}, \begin{bmatrix} 0 \\ 0 \\ 0 \\ 0 \\ \hline 0 \\ 0 \\ 0 \\ 0 \\ \hline 0 \\ [0,a_1] \\ [0,a_2] \end{bmatrix} \;\middle|\; a_1, a_2, a_3, a_4 \in \{0,5\} \subseteq Z_{10} \right\} \subseteq V,$$

H is a semigroup super interval matrix vector space over the semigroup $S = Z_{10}$. We see addition is not defined on H. So H is only a pseudo semigroup semi vector subspace of V over the semigroup $S = Z_{10}$. As in case of set linear algebras we can define semigroup linear transformation of a semigroup linear algebra V into a semigroup linear algebra W defined over the same semigroup S. If T is a semigroup linear transformation from a semigroup linear algebra V into V then we call T to be a semigroup linear operator on V.

We can define pseudo semigroup linear operator on V as in case of semigroup linear algebras [38, 46, 50].

Let V be a semigroup linear algebra over the semigroup S. Let W ⊆ V be a subsemigroup linear subalgebra of V over the subsemigroup P ⊆ S. (P a proper subsemigroup and W a proper linear subalgebra of V) Let T : V → W be a map such that
$$T(\alpha v + u) = T(\alpha) T(v) + T(u)$$
where u, v ∈ V, and for α ∈ S, T(α) ∈ P. We call T a pseudo semigroup linear operator on V.

The reader is expected to provide examples.
Now we define projection of V on a semigroup linear algebra.

Let V be the semigroup linear algebra of super interval matrices over the semigroup S. Let W ⊆ V be a semigroup linear subalgebra of V over the semigroup S. Let T be a semigroup linear operator on V. T is said to be a semigroup linear projection on V if T(v) = w for all v ∈ V and w ∈ W and
$$T(\alpha v + u) = \alpha T(u) + T(v)$$
T(v) and T(u) ∈ W for all α ∈ S and u, v ∈ V.

We will illustrate this situation by some examples.

*Example 4.57*: Let

$$V = \left\{ \begin{bmatrix} \overline{[0,a_1]} \\ \overline{[0,a_2]} \\ \overline{[0,a_3]} \\ \overline{[0,a_4]} \\ \overline{[0,a_5]} \\ \overline{[0,a_6]} \\ \overline{[0,a_7]} \end{bmatrix} \,\middle|\, a_i \in Z^+ \cup \{0\}, 1 \le i \le 7 \right\}$$

be a semigroup super matrix interval linear algebra of super interval column matrices over the semigroup $S = Z^+ \cup \{0\}$.

Let

$$W = \left\{ \begin{bmatrix} \overline{[0,a_1]} \\ 0 \\ 0 \\ \overline{[0,a_2]} \\ [0,a_3] \\ [0,a_4] \\ [0,a_5] \end{bmatrix} \middle| a_i \in 3Z^+ \cup \{0\}, 1 \le i \le 5 \right\} \subseteq V,$$

W is a semigroup linear algebra of super interval column matrices over the semigroup S.

Define a map $T : V \to W$ by

$$T \begin{bmatrix} \overline{[0,a_1]} \\ [0,a_2] \\ [0,a_3] \\ \overline{[0,a_4]} \\ [0,a_5] \\ [0,a_6] \\ [0,a_7] \end{bmatrix} = 3 \begin{bmatrix} \overline{[0,a_1]} \\ 0 \\ 0 \\ \overline{[0,a_4]} \\ [0,a_5] \\ [0,a_6] \\ [0,a_7] \end{bmatrix}.$$

It is easily verifield T is a pseudo semigroup linear operator on V.

Suppose $P : V \to W$ such that

$$P \left( \begin{bmatrix} \overline{[0,a_1]} \\ [0,a_2] \\ [0,a_3] \\ \overline{[0,a_4]} \\ [0,a_5] \\ [0,a_6] \\ [0,a_7] \end{bmatrix} \right) = \begin{bmatrix} \overline{[0,3a_1]} \\ 0 \\ 0 \\ \overline{[0,3a_4]} \\ [0,3a_5] \\ [0,3a_6] \\ [0,3a_7] \end{bmatrix}$$

then P is a semigroup linear projection on W.
Interested reader can construct examples of these.

We can write the semigroup vector space of super interval matrices as a direct union of semigroup vector subspaces.

We define this notion.

**DEFINITION 4.11**: *Let V be a semigroup vector space of super interval matrices over the semigroup S. Let $W_1, W_2, ..., W_n$ be semigroup vector subspaces of V over S if $V = \cup W_i$ and $W_i \cap W_j = (0)$ or $\phi$ if $i \neq j$ then we say V is the direct union of the semigroup vector subspaces of the semigroup vector space V over the semigroup S.*

We will illustrate this situation by some examples.

***Example 4.58***: Let $V = \{([0, a_1]\ [0, a_2]\ |\ [0, a_3]\ [0, a_4]\ [0, a_5]\ |\ [0, a_6])$,

$$\begin{bmatrix}[0,a_1]\\ [0,a_2]\\ [0,a_3]\\ [0,a_4]\\ [0,a_5]\\ [0,a_6]\end{bmatrix}, \begin{pmatrix}[0,a_1] & [0,a_2] & [0,a_3] & [0,a_4]\\ [0,a_5] & [0,a_6] & [0,a_7] & [0,a_8]\\ [0,a_9] & [0,a_{10}] & [0,a_{11}] & [0,a_{12}]\\ [0,a_{13}] & [0,a_{14}] & [0,a_{15}] & [0,a_{16}]\end{pmatrix}$$

$a_i \in Z^+ \cup \{0\}; 1 \le i \le 16\}$ be a semigroup semivector space of super interval matrices over the semigroup S.

Consider $W_1 = \{([0, a_1]\ [0, a_2]\ |\ [0, a_3]\ [0, a_4]\ [0, a_5])\ |\ [0, a_6]\ |\ a_i \in Z^+ \cup \{0\}\} \subseteq V$

$$W_2 = \left\{ \begin{bmatrix}[0,a_1]\\ [0,a_2]\\ [0,a_3]\\ [0,a_4]\\ [0,a_5]\\ [0,a_6]\end{bmatrix} \mid a_i \in Z^+ \cup \{0\}; 1 \le i \le 6 \right\} \subseteq V$$

and

$$W_3 = \left\{ \begin{pmatrix} [0,a_1] & [0,a_2] & [0,a_3] & [0,a_4] \\ [0,a_5] & [0,a_6] & [0,a_7] & [0,a_8] \\ [0,a_9] & [0,a_{10}] & [0,a_{11}] & [0,a_{12}] \\ [0,a_{13}] & [0,a_{14}] & [0,a_{15}] & [0,a_{16}] \end{pmatrix} \right.$$

$a_i \in Z^+ \cup \{0\}$; $1 \le i \le 16\} \subseteq V$ be three semigroup vector subspace of V over the semigroup $S = Z^+ \cup \{0\}$. Clearly $V = \cup W_i = W_1 \cup W_2 \cup W_3$ and $W_i \cap W_j = \phi$ if $i \neq j$. Thus V is the direct union of semigroup vector subspaces of V over S.

*Example 4.59*: Let

$$V = \left\{ \begin{bmatrix} [0,a_1] & [0,a_2] \\ [0,a_3] & [0,a_4] \\ \vdots & \vdots \\ [0,a_{11}] & [0,a_{12}] \end{bmatrix}, \right.$$

$([0, a_1] \ldots [0, a_5] \mid [0, a_6]\ [0, a_7]\ [0, a_8] \ldots [0, a_{16}])$,

$$\begin{bmatrix} [0,a_1] & [0,a_2] & \ldots & [0,a_{10}] & [0,a_{11}] & [0,a_{12}] \\ [0,a_{13}] & [0,a_{14}] & \ldots & [0,a_{22}] & [0,a_{23}] & [0,a_{24}] \\ \vdots & \vdots & \ldots & \vdots & \vdots & \vdots \\ [0,a_{49}] & [0,a_{50}] & \ldots & [0,a_{58}] & [0,a_{59}] & [0,a_{60}] \end{bmatrix}$$

$a_i \in Q^+ \cup \{0\}$; $1 \le i \le 60\}$ be a semigroup vector space of super interval matrices over the semigroup $S = 3Z^+ \cup \{0\}$.

Consider

$$W_1 = \left\{ \begin{bmatrix} [0,a_1] & [0,a_2] \\ [0,a_3] & [0,a_4] \\ \vdots & \vdots \\ [0,a_{11}] & [0,a_{12}] \end{bmatrix} \middle| a_i \in Q^+ \cup \{0\}, 1 \le i \le 12 \right\} \subseteq V,$$

$$W_2 = ([0, a_1] \ldots [0, a_5] \mid [0, a_6]\, [0, a_7]\, [0, a_8] \ldots [0, a_{16}]) \mid$$
$$a_i \in Q^+ \cup \{0\}; 1 \le i \le 60\}$$

and

$$W_3 = \left\{ \begin{bmatrix} [0,a_1] & [0,a_2] & \ldots & [0,a_{10}] & [0,a_{11}] & [0,a_{12}] \\ [0,a_{13}] & [0,a_{14}] & \ldots & [0,a_{22}] & [0,a_{23}] & [0,a_{24}] \\ \vdots & \vdots & \ldots & \vdots & \vdots & \vdots \\ [0,a_{49}] & [0,a_{50}] & \ldots & [0,a_{58}] & [0,a_{59}] & [0,a_{60}] \end{bmatrix} \middle| a_i \in Q^+ \cup \{0\}, 1 \le i \le 60 \right\}$$

$\subseteq V$ be semigroup vector subspaces of V over the semigroup S. Clearly $V = W_1 \cup W_2 \cup W_3$ and $W_i \cap W_j = \phi$ if $i \ne j$; $1 \le i, j \le 3$.

Thus V is the direct union of semigroup vector subspaces $W_1$, $W_2$ and $W_3$ of V over S.

Next we proceed onto define pseudo direct union and illustrate it with examples.

**DEFINITION 4.12**: *Let V be a semigroup vector space over the semigroup S. Let $W_1$, $W_2$, ..., $W_n$ be semigroup vector subspaces of V of super interval matrices over the semigroup S. If $V = \bigcup_{i=1}^{n} W_i$ with $W_i \cap W_j \ne \phi$ or $\{0\}$ if $i \ne j$ then we call V to be a pseudo direct union of semigroup vector subspaces of V over the semigroup S.*

We will illustrate this by some an example.

***Example 4.60***: Let $V = \{([0, a_1]\, [0, a_2] \mid [0, a_3]\, [0, a_4]\, [0, a_5]) \mid a_i \in 13Z^+ \cup 12Z^+ \cup 15Z^+ \cup 11Z^+ \cup \{0\}\}$ be the semigroup vector space of super interval row matrices over the semigroup $S = Z^+ \cup \{0\}$.

Consider $W_1 = \{([0, a_1]\, [0, a_2] \mid [0, a_3]\, [0, a_4]\, [0, a_5]) \mid a_i \in 13Z^+ \cup \{0\}\} \subseteq V$, $W_2 = \{([0, a_1]\, [0, a_2] \mid [0, a_3]\, [0, a_4]\, [0, a_5]) \mid a_i \in 12Z^+ \cup \{0\}\} \subseteq V$, $W_3 = \{([0, a_1]\, [0, a_2] \mid [0, a_3]\, [0, a_4]\, [0, a_5]) \mid a_i \in 15Z^+ \cup \{0\}\} \subseteq V$, and $W_4 = \{([0, a_1]\, [0, a_2] \mid [0, a_3]\, [0, a_4]\, [0, a_5]) \mid a_i \in 11Z^+ \cup \{0\}, 1 \le i \le 5\} \subseteq V$ be semigroup vector subspaces of V over the semigroup S.

Clearly $V = W_1 \cup W_2 \cup W_3 \cup W_4$ but $W_i \cap W_j \ne (0)$ or $\phi$. Thus V is the pseudo direct union of semigroup vector subspaces of V over the semigroup S.

We now proceed onto prove yet another new definition in case of semigroup linear algebra of super interval matrices over a semigroup S.

**DEFINITION 4.13**: *Let V be a semigroup linear algebra of super interval matrices over a semigroup S. We say V is a direct sum of semigroup linear subalgebras $W_1, W_2, \ldots, W_n$ of V if*

1. $V = W_1 + \ldots + W_n$
2. $W_i \cap W_j = \{0\}$ or $\phi$ if $i \neq j$
   $1 \leq i, j \leq n$.

We will illustrate this by some simple examples.

*Example 4.61*: Let

$$V = \left\{ \begin{bmatrix} [0,a_1] & [0,a_2] & [0,a_3] & [0,a_4] & [0,a_5] \\ [0,a_6] & [0,a_7] & [0,a_8] & [0,a_9] & [0,a_{10}] \\ [0,a_{11}] & [0,a_{12}] & [0,a_{13}] & [0,a_{14}] & [0,a_{15}] \end{bmatrix} \,\middle|\, a_i \in Z^+ \cup \{0\} \right\}$$

be the semigroup linear algebra of super interval row vectors over the semigroup $S = Z^+ \cup \{0\}$.

Consider

$$W_1 = \left\{ \begin{bmatrix} [0,a_1] & 0 & 0 & 0 & 0 \\ [0,a_2] & 0 & 0 & 0 & 0 \\ [0,a_3] & 0 & 0 & 0 & 0 \end{bmatrix} \,\middle|\, a_i \in Z^+ \cup \{0\}; 1 \leq i \leq 3 \right\} \subseteq V,$$

$$W_2 = \left\{ \begin{bmatrix} 0 & [0,a_1] & [0,a_2] & 0 & 0 \\ 0 & [0,a_3] & [0,a_4] & 0 & 0 \\ 0 & [0,a_5] & [0,a_6] & 0 & 0 \end{bmatrix} \,\middle|\, a_i \in Z^+ \cup \{0\}; 1 \leq i \leq 6 \right\} \subseteq V$$

and

$$W_3 = \left\{ \begin{bmatrix} 0 & 0 & 0 & [0,a_1] & [0,a_2] \\ 0 & 0 & 0 & [0,a_3] & [0,a_4] \\ 0 & 0 & 0 & [0,a_5] & [0,a_6] \end{bmatrix} \middle| a_i \in Z^+ \cup \{0\}; 1 \leq i \leq 6 \right\} \subseteq V$$

semigroup linear subalgebra of V over the semigroup S.

Clearly $V = W_1 \oplus W_2 \oplus W_3$ and

$$W_i \cap W_j = \begin{bmatrix} 0 & 0 & 0 & 0 & 0 \\ 0 & 0 & 0 & 0 & 0 \\ 0 & 0 & 0 & 0 & 0 \end{bmatrix} \text{ if } i \neq j, 1 \leq i, j \leq 3.$$

Thus V is the direct sum of semigroup linear subalgebra of V over the semigroup linear subalgebras of V over the semigroup $S = Z^+ \cup \{0\}$.

*Example 4.62*: Let

$$V = \left\{ \begin{bmatrix} [0,a_1] & [0,a_2] & [0,a_3] \\ [0,a_4] & [0,a_5] & [0,a_6] \\ \hline [0,a_7] & [0,a_8] & [0,a_9] \\ [0,a_{10}] & [0,a_{11}] & [0,a_{12}] \\ [0,a_{13}] & [0,a_{14}] & [0,a_{15}] \\ \hline [0,a_{16}] & [0,a_{17}] & [0,a_{18}] \end{bmatrix} \middle| a_i \in Z^+ \cup \{0\}; 1 \leq i \leq 18 \right\}$$

be the semigroup linear algebra of super interval matrices over the semigroup $Z^+ \cup \{0\} = S$.

Consider

$$W_1 = \left\{ \begin{bmatrix} 0 & 0 & [0,a_3] \\ 0 & 0 & [0,a_4] \\ \hline 0 & 0 & 0 \\ 0 & 0 & 0 \\ 0 & 0 & 0 \\ \hline [0,a_1] & [0,a_2] & 0 \end{bmatrix} \middle| a_i \in Z^+ \cup \{0\}; 1 \leq i \leq 4 \right\} \subseteq V,$$

$$W_2 = \left\{ \begin{bmatrix} 0 & 0 & 0 \\ 0 & 0 & 0 \\ \hline [0,a_3] & 0 & [0,a_1] \\ 0 & 0 & [0,a_2] \\ 0 & 0 & 0 \\ \hline 0 & 0 & 0 \end{bmatrix} \middle| a_i \in Z^+ \cup \{0\}; 1 \le i \le 3 \right\} \subseteq V,$$

$$W_3 = \left\{ \begin{bmatrix} [0,a_1] & [0,a_2] & 0 \\ [0,a_3] & [0,a_4] & 0 \\ \hline 0 & 0 & 0 \\ 0 & 0 & 0 \\ 0 & 0 & 0 \\ \hline 0 & 0 & [0,a_5] \end{bmatrix} \middle| a_i \in Z^+ \cup \{0\}; 1 \le i \le 5 \right\} \subseteq V$$

and

$$W_4 = \left\{ \begin{bmatrix} 0 & 0 & 0 \\ 0 & 0 & 0 \\ \hline 0 & [0,a_1] & 0 \\ [0,a_3] & [0,a_4] & 0 \\ [0,a_5] & [0,a_6] & [0,a_2] \\ \hline 0 & 0 & 0 \end{bmatrix} \middle| a_i \in Z^+ \cup \{0\}; 1 \le i \le 6 \right\} \subseteq V$$

be the set of semigroup linear subalgebras of V over the semigroup S.

Clearly $V = W_1 \oplus W_2 \oplus W_3 \oplus W_4$ with

$$W_i \cap W_j = \begin{bmatrix} 0 & 0 & 0 \\ 0 & 0 & 0 \\ 0 & 0 & 0 \\ 0 & 0 & 0 \\ 0 & 0 & 0 \\ 0 & 0 & 0 \end{bmatrix} \text{ if } i \ne j, \ 1 \le i, j \le 4.$$

Thus V is the direct sum of semigroup linear subalgebras of V over the semigroup S.

The study of semigroup linear operator, semigroup linear transformation and semigroup linear projection of semigroup linear algebras can be derived as a matter of routine and hence is left as an exercise for the reader.

Next we proceed onto define group vector spaces and group linear algebras of super interval matrices.

**DEFINITION 4.14**: *Let V be a collection of super interval matrices with zero super interval matrices and G be a group under addition. We call V to be a group vector space of super interval matrices over the group, if the following conditions are true.*

*(1) For every v ∈ V and g ∈ G g.v and v.g are in G.*
*(2) 0.v = 0 for every v ∈ V and 0 the additive identity of G.*

We will illustrate this situation by some simple examples.

*Example 4.63*: Let

$$V = \{([0, a_1] [0, a_2] | [0, a_3] [0, a_4] [0, a_5] | [0, a_6]), (0\ 0 | 0\ 0\ 0 | 0),$$

$$\begin{bmatrix} [0,a_1] & [0,a_2] & [0,a_3] \\ [0,a_4] & [0,a_5] & [0,a_6] \\ \hline [0,a_7] & [0,a_8] & [0,a_9] \\ [0,a_{10}] & [0,a_{11}] & [0,a_{12}] \\ [0,a_{13}] & [0,a_{14}] & [0,a_{15}] \\ [0,a_{16}] & [0,a_{17}] & [0,a_{18}] \end{bmatrix}, \begin{bmatrix} 0 & 0 & 0 \\ 0 & 0 & 0 \\ \hline 0 & 0 & 0 \\ 0 & 0 & 0 \\ 0 & 0 & 0 \\ 0 & 0 & 0 \end{bmatrix}, \begin{bmatrix} [0,a_1] \\ [0,a_2] \\ \hline [0,a_3] \\ [0,a_4] \\ [0,a_5] \\ \hline [0,a_6] \\ [0,a_7] \\ [0,a_8] \\ [0,a_9] \\ [0,a_{10}] \end{bmatrix}, \begin{bmatrix} 0 \\ 0 \\ \hline 0 \\ 0 \\ 0 \\ \hline 0 \\ 0 \\ 0 \\ 0 \\ 0 \end{bmatrix} |$$

$a_i \in Z_{45}$, $1 \leq i \leq 18$} be a group vector space over the group $G = Z_{45}$ under addition modulo 45.

*Example 4.64*: Let

$$V = \left\{ \begin{bmatrix} [0,a_1] & [0,a_2] \\ \overline{[0,a_3]} & [0,a_4] \\ \vdots & \vdots \\ [0,a_{13}] & [0,a_{14}] \\ \overline{[0,a_{15}]} & [0,a_{16}] \\ [0,a_{17}] & [0,a_{18}] \end{bmatrix}, \begin{bmatrix} 0 & 0 \\ 0 & 0 \\ \vdots & \vdots \\ 0 & 0 \\ 0 & 0 \\ 0 & 0 \end{bmatrix}, \right.$$

$$\begin{bmatrix} [0,a_1] & [0,a_2] & [0,a_3] & [0,a_4] \\ [0,a_5] & [0,a_6] & [0,a_7] & [0,a_8] \\ [0,a_9] & [0,a_{10}] & [0,a_{11}] & [0,a_{12}] \\ \overline{[0,a_{13}]} & [0,a_{14}] & [0,a_{15}] & [0,a_{16}] \\ [0,a_{17}] & [0,a_{18}] & [0,a_{19}] & [0,a_{20}] \\ [0,a_{21}] & [0,a_{22}] & [0,a_{23}] & [0,a_{24}] \end{bmatrix}, \begin{bmatrix} 0 & 0 & 0 & 0 \\ 0 & 0 & 0 & 0 \\ 0 & 0 & 0 & 0 \\ 0 & 0 & 0 & 0 \\ 0 & 0 & 0 & 0 \\ 0 & 0 & 0 & 0 \end{bmatrix},$$

([0, $a_1$] [0, $a_2$] [0, $a_3$] [0, $a_4$] | [0, $a_5$] | [0, $a_6$] [0, $a_7$] | [0, $a_8$] [0, $a_9$] [0, $a_{10}$]), (0 0 0 0 | 0 | 0 0) | $a_i \in Z_{27}$, $1 \leq i \leq 24$} be a group vector space of super interval matrices over the group $Z_{27}$ under addition modulo 27.

It is pertinent to mention here that we cannot have group vector spaces over a group using $Z^+ \cup \{0\}$ or $Q^+ \cup \{0\}$ or $R^+ \cup \{0\}$. The only group vector spaces we get are using modulo integers $Z_n$, $n < \infty$ under addition.

We now proceed onto define the notion of group vector subspace of a group vector space over a group G.

**DEFINITION 4.15**: *Let V be a group vector space of super interval matrices defined over the group G. Let $P \subseteq V$ (P a proper subset of V). If P itself is a group vector space over the group G then we call P to be a group vector subspace of V over the group G.*

We will illustrate this situation by some examples.

*Example 4.65*: Let

$$V = \left\{ \begin{bmatrix} [0,a_1] & [0,a_2] & [0,a_3] \\ [0,a_4] & [0,a_6] & [0,a_7] \\ [0,a_5] & [0,a_8] & [0,a_9] \end{bmatrix}, \begin{bmatrix} [0,a_1] \\ [0,a_2] \\ [0,a_3] \\ [0,a_4] \\ [0,a_5] \\ [0,a_6] \end{bmatrix}, \begin{bmatrix} 0 \\ 0 \\ 0 \\ 0 \\ 0 \\ 0 \end{bmatrix}, \begin{bmatrix} 0 & 0 & 0 \\ 0 & 0 & 0 \\ 0 & 0 & 0 \end{bmatrix}, \right.$$

$([0, a_1] \mid [0, a_2] [0, a_3] [0, a_4] [0, a_5] \mid [0, a_6] [0, a_7]), (0 \mid 0\ 0\ 0\ 0 \mid 0\ 0) \mid a_i \in Z_7, 1 \leq i \leq 9\}$ be a group vector space of super interval matrices over the group $Z_7 = G$ under addition modulo 7.

Consider $P = \{([0, a_1] \mid [0, a_2] [0, a_3] [0, a_4] [0, a_5] \mid [0, a_6] [0, a_7]), (0 \mid 0\ 0\ 0\ 0 \mid 0\ 0) \mid a_i \in Z_7, 1 \leq i \leq 7\} \subseteq V$, P is a group vector subspace of V over the group $G = Z_7$.

Consider

$$T = \left\{ \begin{bmatrix} [0,a_1] & [0,a_2] & [0,a_3] \\ [0,a_4] & [0,a_6] & [0,a_7] \\ [0,a_5] & [0,a_8] & [0,a_9] \end{bmatrix} \begin{bmatrix} 0 & 0 & 0 \\ 0 & 0 & 0 \\ 0 & 0 & 0 \end{bmatrix} \bigg| a_i \in Z_7; 1 \leq i \leq 9 \right\} \subseteq V,$$

T is a group vector subspace of V over the group $G = Z_7$.

*Example 4.66*: Let

$$V = \{([0, a_1] [0, a_2] \mid [0, a_3] \ldots [0, a_{12}] \mid [0, a_{13}] \ldots [0, a_{20}]),$$
$$(0\ 0 \mid 0 \ldots 0 \mid 0 \ldots 0),$$

$$\begin{bmatrix} [0,a_1] & [0,a_2] & [0,a_3] & [0,a_4] & [0,a_5] & [0,a_6] \\ [0,a_7] & [0,a_8] & [0,a_9] & [0,a_{10}] & [0,a_{11}] & [0,a_{12}] \end{bmatrix},$$

$$\begin{pmatrix} 0 & 0 & 0 & 0 & 0 & 0 \\ 0 & 0 & 0 & 0 & 0 & 0 \end{pmatrix},$$

$$\begin{bmatrix} [0,a_1] \\ \overline{[0,a_2]} \\ \vdots \\ [0,a_{11}] \\ \overline{[0,a_{12}]} \\ \vdots \\ [0,a_{16}] \end{bmatrix}, \begin{bmatrix} 0 \\ \overline{0} \\ \vdots \\ 0 \\ \overline{0} \\ \vdots \\ 0 \end{bmatrix} \Big|$$

$a_i \in Z_{240}$, $1 \le i \le 20$} be a group vector space of super interval matrices over the group $G = Z_{240}$.

Consider

$H = \{([0, a_1] [0, a_2] \mid [0, a_3] \ldots [0, a_{12}] \mid [0, a_{13}] \ldots [0, a_{20}])$,

$$(0\ 0 \mid 0 \ldots 0 \mid 0 \ldots 0), \begin{bmatrix} [0,a_1] \\ \overline{[0,a_2]} \\ \vdots \\ [0,a_{11}] \\ \overline{[0,a_{12}]} \\ \vdots \\ [0,a_{16}] \end{bmatrix}, \begin{bmatrix} 0 \\ \overline{0} \\ \vdots \\ 0 \\ \overline{0} \\ \vdots \\ 0 \end{bmatrix} \Big|$$

$a_i \in Z_{240}$, $1 \le i \le 16\} \subseteq V$, is a group vector subspace of V over the group $G = Z_{240}$.

Now we can define also the notion of subgroup vector subspace of a group vector space over a group G.

**DEFINITION 4.16**: *Let V be a group vector space of super interval matrices over the group G. Let $H \subseteq G$ be a proper subgroup of G. Let $W \subseteq V$, W be a group vector space over the group H; then we call W to be a subgroup vector subspace of V over the subgroup H of G. If V has no subgroup vector subspace we call V to be a pseudo simple group vector space over the group G.*

We will illustrate this situation by some examples.

**Example 4.67**: Let

$$V = \left\{ \begin{bmatrix} [0,a_1] & [0,a_2] & [0,a_3] \\ [0,a_4] & [0,a_5] & [0,a_6] \\ [0,a_7] & [0,a_8] & [0,a_9] \\ [0,a_{10}] & [0,a_{11}] & [0,a_{12}] \\ [0,a_{13}] & [0,a_{14}] & [0,a_{15}] \\ [0,a_{16}] & [0,a_{17}] & [0,a_{18}] \\ [0,a_{19}] & [0,a_{20}] & [0,a_{21}] \\ [0,a_{22}] & [0,a_{23}] & [0,a_{24}] \end{bmatrix}, \begin{bmatrix} 0 & 0 & 0 \\ 0 & & 0 \\ 0 & 0 & 0 \\ 0 & 0 & 0 \\ 0 & 0 & 0 \\ 0 & 0 & 0 \\ 0 & 0 & 0 \\ 0 & 0 & 0 \end{bmatrix}, \right.$$

$([0, a_1]\ [0, a_2]\ |\ [0, a_3]\ |\ [0, a_4]\ [0, a_5])\ (0\ 0\ |\ 0\ |\ 0\ 0)\ |\ a_i \in Z_{24}, 1 \le i \le 24\}$ be group vector space of super interval matrices over the group $G = Z_{24}$ under addition modulo 24. Consider $H = \{0, 2, 4, 6, 8, \ldots, 20, 22\} \subseteq Z_{24} = G$ be a subgroup of G under addition modulo 24.

Choose $W = \{([0, a_1]\ [0, a_2]\ |\ [0, a_3]\ |\ [0, a_4]\ [0, a_5])\ (0\ 0\ |\ 0\ |\ 0\ 0)\ |\ a_i \in \{0, 2, 4, 6, 8, \ldots, 18, 20, 22\} \subseteq Z_{24}, 1 \le i \le 5\} \subseteq V$, W is a subgroup super vector subspace of V over the subgroup H of G.

W is also a group vector subspace of V over G.

**Example 4.68**: Let

$$V = \left\{ \begin{pmatrix} [0,a_1] & [0,a_2] & | & [0,a_3] & [0,a_4] & [0,a_5] & [0,a_6] & | & [0,a_7] \\ [0,a_8] & [0,a_9] & | & [0,a_{10}] & [0,a_{11}] & [0,a_{12}] & [0,a_{13}] & | & [0,a_{14}] \end{pmatrix}, \right.$$

$$\begin{pmatrix} 0 & 0 & | & 0 & 0 & 0 & 0 & | & 0 \\ 0 & 0 & | & 0 & 0 & 0 & 0 & | & 0 \end{pmatrix}, \begin{pmatrix} [0,a_1] & [0,a_2] & [0,a_3] \\ [0,a_4] & [0,a_5] & [0,a_6] \\ [0,a_7] & [0,a_8] & [0,a_9] \\ [0,a_{10}] & [0,a_{11}] & [0,a_{12}] \\ [0,a_{13}] & [0,a_{14}] & [0,a_{15}] \\ [0,a_{16}] & [0,a_{17}] & [0,a_{18}] \\ [0,a_{19}] & [0,a_{20}] & [0,a_{21}] \end{pmatrix}, \begin{pmatrix} 0 & 0 & 0 \\ 0 & 0 & 0 \\ 0 & 0 & 0 \\ 0 & 0 & 0 \\ 0 & 0 & 0 \\ 0 & 0 & 0 \\ 0 & 0 & 0 \end{pmatrix} \mid$$

$a_i \in Z_{45}$, $1 \le i \le 21$} be a group vector space of super interval matrices over the group $G = Z_{45}$. Choose $H = \{0, 5, 10, 15, 20, 25, 30, 35, 40\} \subseteq Z_{45} = G$ is a subgroup of G.

Consider

$$W = \left\{ \begin{pmatrix} 0 & 0 & | & 0 & 0 & 0 & 0 & | & 0 \\ 0 & 0 & | & 0 & 0 & 0 & 0 & | & 0 \end{pmatrix}, \right.$$

$$\begin{pmatrix} [0,a_1] & [0,a_2] & | & [0,a_3] & [0,a_4] & [0,a_5] & [0,a_6] & | & [0,a_7] \\ [0,a_8] & [0,a_9] & | & [0,a_{10}] & [0,a_{11}] & [0,a_{12}] & [0,a_{13}] & | & [0,a_{14}] \end{pmatrix} |$$

$a_i \in Z_{45}$, $1 \le i \le 14\} \subseteq V$, W is a subgroup vector subspace of V over the subgroup H of G.

**Example 4.69**: Let

$$V = \left\{ \begin{bmatrix} 0 & | & 0 & 0 \\ 0 & | & 0 & 0 \\ 0 & | & 0 & 0 \end{bmatrix}, \begin{bmatrix} [0,a_1] & | & [0,a_2] & [0,a_3] \\ [0,a_4] & | & [0,a_6] & [0,a_7] \\ [0,a_5] & | & [0,a_8] & [0,a_9] \end{bmatrix}, \right.$$

$([0, a_1] [0, a_2] [0, a_3] | [0, a_4] [0, a_5] [0, a_6] | [0, a_7] [0, a_8] | [0, a_9])$, $(0\ 0\ 0\ |\ 0\ 0\ 0\ |\ 0\ 0\ |\ 0)\ |\ a_i \in Z_{13}$, $1 \le i \le 9\}$ be a group vector space of super interval matrices over the group $G = Z_{13}$ under addition modulo 13.

Now

$$W = \left\{ \begin{bmatrix} 0 & | & 0 & 0 \\ 0 & | & 0 & 0 \\ 0 & | & 0 & 0 \end{bmatrix} \begin{bmatrix} [0,a_1] & | & [0,a_2] & [0,a_3] \\ [0,a_4] & | & [0,a_6] & [0,a_7] \\ [0,a_5] & | & [0,a_8] & [0,a_9] \end{bmatrix} \middle| a_i \in Z_{13}, 1 \le i \le 9 \right\} \subseteq V$$

is a group vector subspace of V. However V has no subgroup vector subspaces as G has no proper subgroup under addition.

In view of this we have the following theorem the proof of which is straight forward.

**THEOREM 4.4**: *Let V be a group vector space of super interval matrices over the group $G = Z_p$, p a prime. G is a pseudo simple group vector space. We can define linearly dependent set and linearly independent* set in case of group vector spaces also. This task is a matter of routine and hence is left as an exercise for the reader.

We as in case of semigroup vector spaces define for any two group vector spaces V and W of super interval matrices over the group G, the group linear transformation.
T : V → W as T (αv) = α T (v) for all v ∈V and α ∈ G.

If V = W then the group linear transformation is defined as group linear operator. This study is interesting and the reader is expected to work with it.

Now we proceed onto define the notion of group linear algebra of interval super matrices over a group G.

**DEFINITION 4.17**: *Let V be a group vector space over the additive group G. If V is a group under addition then we define V to be a group interval super matrices of linear algebra over the group G.*

We will illustrate this situation by some examples.

*Example 4.70*: Let V = {([0, $a_1$] [0, $a_2$] | [0, $a_3$] [0, $a_4$] [0, $a_5$] | [0, $a_6$] | [0, $a_7$] [0, $a_8$]) | $a_i$ ∈ $Z_{42}$, 1 ≤ i ≤ 8}, V is a group under addition. V is a group linear algebra over the additive group G = $Z_{42}$.

*Example 4.71*: Let

$$V = \left\{ \begin{bmatrix} [0,a_1] & [0,a_2] & [0,a_3] \\ [0,a_4] & [0,a_5] & [0,a_6] \\ \hline [0,a_7] & [0,a_8] & [0,a_9] \\ [0,a_{10}] & [0,a_{11}] & [0,a_{12}] \\ [0,a_{13}] & [0,a_{14}] & [0,a_{15}] \\ [0,a_{16}] & [0,a_{17}] & [0,a_{18}] \end{bmatrix} \right.$$

where $a_i$ ∈ $Z_{29}$, 1 ≤ i ≤ 18} be a group linear algebra over the additive group G = $Z_{29}$.

*Example 4.72*: Let

$$V = \left\{ \begin{bmatrix} [0,a_1] & [0,a_2] \\ [0,a_3] & [0,a_4] \\ [0,a_5] & [0,a_6] \end{bmatrix} \right.$$

$$\left. \begin{vmatrix} [0,a_7] & [0,a_8] & [0,a_9] & [0,a_{10}] & [0,a_{11}] & [0,a_{12}] \\ [0,a_{13}] & [0,a_{14}] & [0,a_{15}] & [0,a_{16}] & [0,a_{17}] & [0,a_{18}] \\ [0,a_{19}] & [0,a_{20}] & [0,a_{21}] & [0,a_{22}] & [0,a_{23}] & [0,a_{24}] \end{vmatrix} \right.$$

$a_i \in Z_{240}$, $1 \leq i \leq 24$} be a group linear algebra of super row vectors over the group $G = Z_{240}$.

*Example 4.73*: Let

$$P = \left\{ \begin{bmatrix} [0,a_1] & [0,a_2] & [0,a_3] \\ \hline [0,a_4] & [0,a_5] & [0,a_6] \\ [0,a_7] & [0,a_8] & [0,a_9] \\ \hline [0,a_{10}] & [0,a_{11}] & [0,a_{12}] \\ [0,a_{13}] & [0,a_{14}] & [0,a_{15}] \\ [0,a_{16}] & [0,a_{17}] & [0,a_{18}] \\ [0,a_{19}] & [0,a_{20}] & [0,a_{21}] \\ \hline [0,a_{22}] & [0,a_{23}] & [0,a_{24}] \end{bmatrix} \right. \quad a_i \in Z_{14}, 1 \leq i \leq 24 \}$$

be a group linear algebra of super interval column vectors defined over the group $G = Z_{14}$.

It is important and interesting to note the following. Using super interval matrices with special type of intervals of the form [0, a] we see we cannot get group linear algebras built using $Z^+ \cup \{0\}$ or $Q^+ \cup \{0\}$ or $R^+ \cup \{0\}$ under addition.

*Example 4.74*: Let

$$V = \left\{ \begin{bmatrix} [0,a_1] & [0,a_2] & [0,a_3] \\ [0,a_4] & [0,a_5] & [0,a_6] \\ \hline [0,a_7] & [0,a_8] & [0,a_9] \\ [0,a_{10}] & [0,a_{14}] & [0,a_{12}] \end{bmatrix} \right.$$

$a_i \in Z_{43}$, $1 \le i \le 12$} be a group linear algebra of super interval matrices over the group G.

All group linear algebras constructed are of finite order hence we will have a finite basis or a finite generating set.

Another interesting feature is that using a single interval matrix one can derive several group linear algebras over the same group G. This is a marked difference between interval matrices and interval super matrices. Thus using interval super matrices helps us to get a variety of algebraic structures.

Now we proceed onto describe group linear subalgebra of super interval matrices and subgroup linear subalgebra of super interval matrices and give examples of them.

**DEFINITION 4.18**: *Let V be a group linear algebra of super interval matrices over the group G. Let $W \subseteq V$; if W is a proper subset of V and W itself is a group linear algebra over the group G then we define W to be a group linear subalgebra of super interval matrices of V over the group G.*

We will illustrate this situation by some examples.

*Example 4.75*: Let

$$V = \left\{ \begin{bmatrix} [0,a_1] & [0,a_2] \\ [0,a_3] & [0,a_4] \\ \hline [0,a_5] & [0,a_6] \\ [0,a_7] & [0,a_8] \\ [0,a_9] & [0,a_{10}] \\ \hline [0,a_{11}] & [0,a_{12}] \\ [0,a_{13}] & [0,a_{14}] \\ [0,a_{15}] & [0,a_{16}] \\ [0,a_{17}] & [0,a_{18}] \end{bmatrix} \middle| a_i \in Z_{25}; 1 \le i \le 18 \right\}$$

be a group linear algebra of the super interval row vectors over the group $G = Z_{25}$ under addition modulo 25.

Consider

$$W = \left\{ \begin{bmatrix} \begin{array}{cc} [0,a_1] & [0,a_2] \\ [0,a_3] & [0,a_4] \\ \hline 0 & 0 \\ 0 & 0 \\ 0 & 0 \\ \hline [0,a_5] & [0,a_5] \\ 0 & 0 \\ \hline [0,a_7] & [0,a_8] \\ 0 & 0 \end{array} \end{bmatrix} \middle| a_i \in Z_{25};\ 1 \leq i \leq 8 \right\} \subseteq V; $$

W is a group linear subalgebra of V over the group $G = Z_{25}$.

***Example 4.76***: Let V = {([0, $a_1$] [0,$a_2$] [0,$a_3$] | [0,$a_4$] [0,$a_5$] [0,$a_6$] [0,$a_7$] | [0,$a_8$] [0,$a_9$] | [0,$a_{10}$]) | $a_i \in Z_7$; $1 \leq i \leq 10$} be a group linear algebra over the group G = $Z_7$ under addition modulo 7. V is a set of super row interval matrices. Consider M = {([0,$a_1$] [0,$a_2$] [0,$a_3$] | 0 [0,$a_4$] 0 [0,$a_5$] | 0 [0,$a_6$] 0) | $a_i \in Z_7$; $1 \leq i \leq 6$} ⊆ V is a group linear subalgebra of super row interval matrices over the group G = $Z_7$.

***Example 4.77***: Let

$$P = \left\{ \begin{pmatrix} \begin{array}{c|ccc} [0,a_1] & [0,a_2] & [0,a_3] & [0,a_4] \\ [0,a_5] & [0,a_6] & [0,a_7] & [0,a_8] \\ \hline [0,a_9] & [0,a_{10}] & [0,a_{11}] & [0,a_{20}] \\ [0,a_{12}] & [0,a_{13}] & [0,a_{14}] & [0,a_{15}] \\ [0,a_{16}] & [0,a_{17}] & [0,a_{18}] & [0,a_{19}] \end{array} \end{pmatrix} \middle| a_i \in Z_{15};\ 1 \leq i \leq 20 \right\} $$

be the group linear algebra of super interval matrices over the group G = $Z_{15}$ under addition modulo 15.

Consider

$$T = \left\{ \begin{pmatrix} \begin{array}{c|ccc} [0,a_1] & 0 & 0 & 0 \\ [0,a_2] & 0 & 0 & 0 \\ \hline 0 & [0,a_3] & [0,a_4] & [0,a_5] \\ 0 & [0,a_6] & [0,a_7] & [0,a_8] \\ 0 & [0,a_9] & [0,a_{10}] & [0,a_{11}] \end{array} \end{pmatrix} \right.$$

$a_i \in Z_{15}$; $1 \le i \le 11\} \subseteq P$ is a group linear subalgebra of P over the group $G = Z_{15}$.

Now we proceed onto show every group linear algebra of subinterval matrices over G in general need not have group linear subalgebras over G.

We will illustrate this situation by some examples.

*Example 4.78*: Let

$$V = \left\{ \begin{bmatrix} [0,a] \\ [0,a] \\ \overline{[0,a]} \\ [0,a] \\ \overline{[0,a]} \\ [0,a] \\ [0,a] \\ \overline{[0,a]} \\ [0,a] \end{bmatrix} \quad a_i \in Z_{23} \right\}$$

be a group linear algebra of super column interval matrices over the group $G = Z_{23}$ under addition modulo 23. V has no nontrivial group linear subalgebra. In fact order of V is 23 and V is generated by

$$X = \left\{ \begin{bmatrix} [0,1] \\ \overline{[0,1]} \\ [0,1] \\ \overline{[0,1]} \\ [0,1] \\ [0,1] \\ \overline{[0,1]} \\ [0,1] \\ [0,1] \end{bmatrix} \right\} \subseteq V.$$

Thus dimension of V as a group linear algebra over G is one.

*Example 4.79*: Let

$$V = \left\{ \begin{bmatrix} [0,a] & [0,a] & [0,a] \\ [0,a] & [0,a] & [0,a] \\ \overline{[0,a]} & \overline{[0,a]} & \overline{[0,a]} \\ [0,a] & [0,a] & [0,a] \\ [0,a] & [0,a] & [0,a] \end{bmatrix} \middle| \; a_i \in Z_{17} \right\}$$

be a group linear algebra of super interval matrices over the group $G = Z_{17}$ under addition modulo 17. Clearly V has no group linear subalgebras. Order of V is 17 and dimension of V is one. V is generated by

$$X = \left\{ \begin{bmatrix} [0,1] & [0,1] & [0,1] \\ [0,1] & [0,1] & [0,1] \\ \overline{[0,1]} & \overline{[0,1]} & \overline{[0,1]} \\ [0,1] & [0,1] & [0,1] \\ [0,1] & [0,1] & [0,1] \end{bmatrix} \right\} \subseteq V$$

over G.

***Example 4.80***: Let V = {([0, a] [0,a] | [0,a] | [0,a] [0,a] [0,a] [0,a] [0,a] | [0,a] [0,a] [0,a] | [0,a]) | $a_i \in Z_{43}$} be a group linear algebra of super row interval matrices over the group $G = Z_{43}$ under addition modulo 43.

Clearly number of elements in V is 43. V has no proper group linear subalgebras over the group G. V is generated by X = {([0, 1] [0,1] | [0,1] | [0,1] [0,1] [0,1] [0,1] [0,1] | [0,1] [0,1] [0,1] | [0,1]|) $\subseteq$ V. Thus dimension of V as a group linear algebra over G is one.

Thus we have seen examples of group linear algebras of super interval matrices over the group G which does not contain any proper group linear subalgebras over the group G. We call these group linear algebras of super interval matrices as simple group linear algebras. We give a theorem the proof of which is direct and hence left as an exercise to the reader.

**THEOREM 4.5**: *Let V be a set of super interval matrices where every interval in that super interval matrix is [0, a] with a $\in Z_p$, p a prime. V be the group linear algebra over $Z_p = G$. V is a simple group linear algebra and has p elements and is of dimension one over $Z_p = G$.*

Next we proceed onto define subgroup linear subalgebra of a group linear algebra.

**DEFINITION 4.19**: *Let V be a group linear algebra of super interval matrices over the group G under addition. Suppose W $\subseteq$ V be a proper subset of V. Let H $\subseteq$ G be a proper subgroup of G. If W is a group linear algebra of super interval matrices over the group H, then we define W to be subgroup linear subalgebra of V over the subgroup H of the group G.*

We will illustrate this situation by some examples.

***Example 4.81***: Let

$$V = \left\{ \begin{pmatrix} [0,a_1] & [0,a_4] & [0,a_7] & [0,a_8] & [0,a_9] & [0,a_{10}] & [0,a_{19}] \\ [0,a_2] & [0,a_5] & [0,a_{11}] & [0,a_{12}] & [0,a_3] & [0,a_{14}] & [0,a_{20}] \\ [0,a_3] & [0,a_6] & [0,a_{15}] & [0,a_{16}] & [0,a_{17}] & [0,a_{18}] & [0,a_{21}] \end{pmatrix} \right.$$

$a_i \in Z_{100}$; $1 \le i \le 21$} be a group linear algebra of super row interval vectors over the group $G = Z_{100}$ under addition modulo 100. Consider H $= \{0, 10, 20, 30, \ldots, 90\} \subseteq G = Z_{100}$ is a subgroup of G under addition modulo 100.

Choose

$$W = \left\{ \begin{pmatrix} 0 & 0 & [0,a_1] & [0,a_2] & [0,a_3] & [0,a_4] & 0 \\ 0 & 0 & [0,a_5] & [0,a_6] & [0,a_7] & [0,a_{11}] & 0 \\ 0 & 0 & [0,a_8] & [0,a_9] & [0,a_{10}] & [0,a_{12}] & 0 \end{pmatrix} \right.$$

$a_i \in Z_{100}$; $1 \le i \le 12\} \subseteq V$; W is a subgroup linear subalgebra of V over the subgroup H of the group G.

*Example 4.82*: Let

$$V = \left\{ \begin{bmatrix} \begin{array}{cc} [0,a_1] & [0,a_2] \\ [0,a_3] & [0,a_4] \\ [0,a_5] & [0,a_6] \\ \hline [0,a_7] & [0,a_8] \\ [0,a_9] & [0,a_{10}] \\ \hline [0,a_{11}] & [0,a_{12}] \\ [0,a_{13}] & [0,a_{14}] \\ [0,a_{15}] & [0,a_{16}] \end{array} \end{bmatrix} \right| a_i \in Z_{45}; 1 \le i \le 16 \right\}$$

be a group linear algebra of super column vectors over the group $G = Z_{45}$. Consider $H = \{0, 5, 10, 15, \ldots, 35, 40\} \subseteq Z_{45}$ be a subgroup of G under addition modulo 45.

Choose

$$W = \left\{ \begin{bmatrix} \begin{array}{cc} [0,a_1] & [0,a_2] \\ 0 & 0 \\ 0 & 0 \\ \hline 0 & 0 \\ [0,a_3] & [0,a_4] \\ \hline [0,a_5] & [0,a_6] \\ 0 & 0 \\ 0 & 0 \end{array} \end{bmatrix} \right| a_i \in Z_{45}; 1 \le i \le 6 \right\} \subseteq V;$$

W is a subgroup linear subalgebra of super column interval vectors of V over H, the subgroup of G.

**Example 4.83**: Let

$$V = \left\{ \begin{bmatrix} [0,a_1] & [0,a_2] \\ \vdots & \vdots \\ [0,a_{15}] & [0,a_{16}] \end{bmatrix} \;\middle|\; a_i \in Z_{43};\; 1 \leq i \leq 16 \right\}$$

be a group linear algebra of super interval matrices over the group $G = Z_{43}$ under addition modulo 43.

Clearly G has no proper subgroups under addition. Hence V cannot contain any subgroup linear subalgebra. We call those group linear algebras which has no subgroup linear subalgebra as a pseudo simple group linear algebra of super interval matrices.

We will give a theorem which guarantees a class of pseudo simple group linear algebra of super interval matrices over a group G.

**THEOREM 4.6**: *Let V be a collection of super interval matrices with same partition with entries from $Z_p$, p a prime be the group linear algebra over the group $G = Z_p$ under addition modulo p. V is a pseudo simple group linear algebra over $G = Z_p$.*

We will give more examples before we proceed further.

**Example 4.84**: Let $G = Z_{43}$ be a group under addition modulo 43. Choose

$$V = \left\{ \begin{bmatrix} [0,a_1] & \cdots & [0,a_7] & [0,a_8] \\ [0,a_9] & \cdots & [0,a_{15}] & [0,a_{16}] \\ \hline [0,a_{17}] & \cdots & [0,a_{23}] & [0,a_{24}] \\ \vdots & \vdots & \vdots & \vdots \\ \hline [0,a_{65}] & \cdots & [0,a_{71}] & [0,a_{72}] \\ [0,a_{73}] & \cdots & [0,a_{79}] & [0,a_{80}] \\ [0,a_{81}] & \cdots & [0,a_{87}] & [0,a_{88}] \end{bmatrix} \right.$$

where $a_i \in Z_{43}$; $1 \leq i \leq 88$} be a group linear algebra of super interval matrices over the group $G = Z_{43}$. V has several group linear subalgebras but no subgroup linear subalgebra as G is a group of prime order 43.

As in case of semigroup linear algebras we can define the notion of direct sum in case of group linear algebras which is a matter of routine. Pseudo direct sum can also defined in case of group linear algebras.

We will give an example of each and the interested reader can provide more examples.

***Example 4.85***: Let V be a group linear algebra over the group $G = Z_{20}$. V = {($[0, a_1]$ $[0,a_2]$ $[0,a_3]$ | $[0,a_4]$ $[0,a_5]$ $[0,a_6]$ $[0,a_7]$ | $[0,a_8]$ $[0,a_9]$ | $[0,a_{10}]$) | $a_i \in Z_{20}$, $1 \leq i \leq 10$} be the group linear algebra of super interval row matrices. Consider $W_1$ = {($[0, a_1]$ $[0,a_2]$ $[0,a_3]$ | 0 0 0 0 | 0 0 | 0) | $a_i \in Z_{20}$, $1 \leq i \leq 3$}, $W_2$ = {( 0 0 0 | 0 0 0 0 | $[0, a_1]$ $[0,a_2]$ | 0)} $a_i \in Z_{20}$, $1 \leq i \leq 2$} and $W_3$ = {(0 0 0 | $[0, a_1]$ $[0,a_2]$ $[0,a_3]$ $[0,a_4]$ | 0 0 | $[0,a_6]$) | $a_i \in Z_{20}$, $1 \leq i \leq 5$} be group linear subalgebras of the group linear algebra V.`

Clearly $V = W_1 + W_2 + W_3$ and $W_i \cap W_j = (0)$ if $i \neq j$, $1 \leq i \leq 3$. Thus V is the direct sum of group linear subalgebras of V.

***Example 4.86***: Let

$$V = \left\{ \begin{bmatrix} [0,a_1] & [0,a_2] & [0,a_3] \\ [0,a_4] & [0,a_5] & [0,a_6] \\ \hline [0,a_7] & [0,a_8] & [0,a_9] \\ [0,a_{10}] & [0,a_{11}] & [0,a_{12}] \\ \hline [0,a_{13}] & [0,a_{14}] & [0,a_{15}] \end{bmatrix} \middle| a_i \in Z_{12}; \right.$$

$1 \leq i \leq 15$} be a group linear algebra of super interval matrices over the group $G = Z_{12}$ under addition modulo 12.
Consider

$$W_1 = \left\{ \begin{bmatrix} [0,a_1] & 0 & 0 \\ [0,a_2] & 0 & 0 \\ [0,a_3] & 0 & 0 \\ 0 & 0 & 0 \\ [0,a_4] & 0 & 0 \end{bmatrix} \;\middle|\; a_i \in Z_{12};\; 1 \le i \le 4 \right\},$$

$$W_2 = \left\{ \begin{bmatrix} 0 & 0 & 0 \\ 0 & 0 & 0 \\ [0,a_1] & 0 & 0 \\ [0,a_2] & 0 & 0 \\ [0,a_4] & [0,a_5] & [0,a_6] \end{bmatrix} \;\middle|\; a_i \in Z_{12};\; 1 \le i \le 6 \right\},$$

$$W_3 = \left\{ \begin{bmatrix} 0 & [0,a_1] & [0,a_2] \\ 0 & [0,a_3] & [0,a_4] \\ 0 & 0 & 0 \\ 0 & 0 & 0 \\ 0 & 0 & [0,a_5] \end{bmatrix} \;\middle|\; a_i \in Z_{12};\; 1 \le i \le 5 \right\},$$

$$W_4 = \left\{ \begin{bmatrix} 0 & 0 & 0 \\ 0 & 0 & 0 \\ 0 & [0,a_1] & [0,a_2] \\ 0 & [0,a_2] & 0 \\ 0 & 0 & [0,a_4] \end{bmatrix} \;\middle|\; a_i \in Z_{12};\; 1 \le i \le 4 \right\}$$

and

$$W_5 = \left\{ \begin{bmatrix} 0 & [0,a_2] & 0 \\ 0 & 0 & 0 \\ [0,a_1] & 0 & 0 \\ 0 & [0,a_3] & [0,a_4] \\ 0 & [0,a_5] & [0,a_6] \end{bmatrix} \;\middle|\; a_i \in Z_{12};\; 1 \le i \le 6 \right\},$$

be group linear subalgebras of V over the group $G = Z_{12}$. Clearly $V = W_1 + W_2 + W_3 + W_4 + W_5$ but $W_i \cap W_j \ne (0)$ or $\phi$, $1 \le i, j \le 5$. Thus V is only a pseudo direct sum of group linear subalgebras $W_1$, $W_2$, $W_3$, $W_4$ and $W_5$.

Now we proceed onto define the notion of pseudo group vector subspace of a group linear algebra.

**DEFINITION 4.20**: *Let V be a group linear algebra of super interval matrices over the group G. Let $W \subseteq V$ be a proper subset of V such that W is only a group vector space over the group G then we define W to be a pseudo group vector subspace of V over the group G.*

We will illustrate this situation by an example or two.

*Example 4.87*: Let

$$V = \left\{ \begin{bmatrix} [0,a_1] & [0,a_2] & [0,a_3] \\ [0,a_4] & [0,a_5] & [0,a_6] \\ \hline [0,a_7] & [0,a_8] & [0,a_9] \\ [0,a_{14}] & [0,a_{11}] & [0,a_{12}] \\ [0,a_{13}] & [0,a_{14}] & [0,a_{13}] \end{bmatrix} \right.$$

$a_i \in Z_{15}$; $1 \leq i \leq 15\}$, be a group linear algebra of super interval matrices over the group $G = Z_{15}$ under addition modulo 15.
Consider

$$W = \left\{ \begin{bmatrix} [0,a] & [0,a] & 0 \\ 0 & 0 & 0 \\ \hline 0 & 0 & [0,a] \\ [0,a] & 0 & 0 \\ 0 & [0,a] & 0 \end{bmatrix}, \begin{bmatrix} [0,a] & 0 & [0,a] \\ 0 & 0 & [0,a] \\ \hline 0 & [0,a] & 0 \\ 0 & 0 & [0,a] \\ 0 & 0 & [0,a] \end{bmatrix}, \right.$$

$$\begin{bmatrix} 0 & 0 & 0 \\ 0 & [0,a] & 0 \\ \hline 0 & 0 & 0 \\ 0 & [0,a] & 0 \\ [0,a] & [0,a] & [0,a] \end{bmatrix}$$

where $a \in Z_{15}\} \subseteq V$. Clearly W is only a group vector space over the group $G = Z_{15}$. Clearly W is not a group linear subalgebra over G hence W is a pseudo group vector subspace of V over the group $G = Z_{15}$.

**Example 4.88**: Let

$$V = \left\{ \begin{pmatrix} [0,a_1] & [0,a_4] & |[0,a_7]| & [0,a_{10}] & [0,a_{11}] & [0,a_{12}] & [0,a_{13}] \\ [0,a_2] & [0,a_5] & |[0,a_8]| & [0,a_{14}] & [0,a_{15}] & [0,a_{16}] & [0,a_{17}] \\ [0,a_3] & [0,a_6] & |[0,a_9]| & [0,a_{18}] & [0,a_{19}] & [0,a_{20}] & [0,a_{21}] \end{pmatrix} \right.$$

$a_i \in Z_8\}$; $1 \le i \le 21\}$ be a group linear algebra of super interval row vectors over the group $G = Z_8$.

Consider

$$W = \left\{ \begin{pmatrix} [0,a] & [0,a] & |0| & 0 & 0 & 0 & [0,a] \\ [0,a] & [0,a] & |0| & 0 & 0 & 0 & [0,a] \\ [0,a] & [0,a] & |0| & 0 & 0 & 0 & [0,a] \end{pmatrix}, \right.$$

$$\begin{pmatrix} 0 & [0,a] & |[0,a]| & 0 & 0 & [0,a] & 0 \\ 0 & [0,a] & |[0,a]| & 0 & 0 & 0 & [0,a] \\ 0 & [0,a] & |[0,a]| & 0 & 0 & [0,a] & 0 \end{pmatrix},$$

$$\begin{pmatrix} [0,a] & 0 & |[0,a]|[0,a] & 0 & 0 & 0 \\ 0 & [0,a] & | 0 & |[0,a] & 0 & 0 & 0 \\ [0,a] & [0,a] & |[0,a]|[0,a] & 0 & 0 & 0 \end{pmatrix}$$

$a_i \in Z_8 \subseteq V$; W is only a group vector space of super interval matrices over the group G. Hence W is the pseudo group vector subspace of V over the group $G = Z_8$.

**Chapter Five**

# SUPER FUZZY INTERVAL MATRICES

In this chapter we introduce the notion of super fuzzy linear algebras of super interval fuzzy matrices. We first introduce the notion of super interval fuzzy matrices before we proceed to define fuzzy set super interval vector spaces, fuzzy semigroup super interval vector spaces and fuzzy super interval group vector spaces. We will be defining two types of fuzzy set super interval sector spaces, fuzzy group super interval vector spaces and fuzzy semigroup super interval vector spaces.

The chapter has two sections in the first section super interval fuzzy matrices are introduced for the first time. Further in this section we define fuzzy super interval set linear algebra or fuzzy set super linear algebra. Special fuzzy set linear algebras are introduced in section two.

## 5.1 Super Fuzzy Interval Matrices

In this section we for the first time introduce the notion of super fuzzy interval matrices and describe some operations on them.

**DEFINITION 5.1.1**: *Let X = {([0, $a_1$] [0,$a_2$] | [0,$a_3$] [0,$a_4$] | [0,$a_5$] . . . [0,$a_n$] | [0,$a_1$] [0,$a_{n-1}$]) | $a_i$ ∈ [0,1]}. X is defined as the super fuzzy interval row matrix.*

So a super interval row matrix in which the entries are from [0,1] will be known as the super fuzzy interval row matrix.

***Example 5.1.1***: Let P = ([0, 0.7] [0,1] | 0 [0,0.2] [0,0.132] | [0,0.9312] [0,0.571] [0,0.3102] [0,0.014]) is a super fuzzy interval row matrix.

***Example 5.1.2***: Let X = {([0, 1] [0,0.5] [0,0.371] | [0,0.123] [0,0.3012] | [0,0.015])} is a super fuzzy interval row matrix.

**DEFINITION 5.1.2**: *Let*

$$T = \left\{ \begin{bmatrix} [0,a_1] \\ [0,a_2] \\ [0,a_3] \\ \hline [0,a_4] \\ \vdots \\ [0,a_r] \\ \hline [0,a_{r+1}] \\ \vdots \\ [0,a_{r-1}] \\ \hline [0,a_n] \end{bmatrix}, a_i \in [0,1]; 1 \le i \le n \right\}$$

*be a collection of super interval column matrices. T is a fuzzy super interval column matrix or super fuzzy interval column matrix.*

We will illustrate this situation by some examples.

***Example 5.1.3***: Let

$$P = \left\{ \begin{bmatrix} [0,0.5] \\ [0,0.2] \\ \hline [0,0.3] \\ [0,0.12] \\ [0,0.0123] \\ \hline [0,0.051] \\ [0,0.0125] \\ [0,0.051] \\ [0,0.0012] \\ [0,0.412] \\ [0,0.302] \end{bmatrix} \right\},$$

P is a super interval fuzzy column matrix.

***Example 5.1.4***: Let

$$X = \begin{bmatrix} [0,1] \\ [0,0.21] \\ [0,0.514] \\ \hline [0,0.35] \\ [0,0.021] \\ [0,0.101] \\ \hline [0,0.2] \\ [0,0.1] \\ [0,1] \end{bmatrix}$$

be a super fuzzy interval column matrix.

***Example 5.1.5***: Let

$$X = \begin{bmatrix} [0,1] \\ [0,0.2] \\ \hline [0,0.3] \\ \hline [0,0.4] \end{bmatrix}$$

be a super fuzzy interval column matrix.

**DEFINITION 5.1.3**: *Let*

$$V = \left\{ \begin{pmatrix} a_1 & a_2 & a_7 & a_{10} & a_{13} & \cdots & a_{n-2} \\ a_3 & a_4 & a_8 & a_{11} & a_{14} & \cdots & a_{n-1} \\ a_5 & a_6 & a_9 & a_{12} & a_{15} & \cdots & a_n \end{pmatrix} \,\middle|\, a_i \in [0,1];\ 1 \le i \le n \right\},$$

*V is defined as the super fuzzy interval row vector.*

We will illustrate this by some examples.

***Example 5.1.6***: Let
$$V = \left\{ \begin{pmatrix} [0,1] & [0,0.2] & [0.1] & [0,0.51] & [0.031] & [0,072] \\ 0 & [0,0.1] & [0.03] & [0,0.69] & [0.0.04] & [0,0.21] \end{pmatrix} \right\},$$
V is a super fuzzy interval row vector.

***Example 5.1.7***: Let
$$M = \begin{pmatrix} [0,0.03] & [0.1] & [0,0.72] & [0,0.38] & [0,0.11] & [0,0.2] \\ [0,0.16] & 0 & [0,0.103] & [0,0.105] & [0,0.167] & [0,0.7] \end{pmatrix}$$
be the super fuzzy interval row vector.

Now we proceed onto describe the notion of super fuzzy interval column vector.

**DEFINITION 5.1.4**: *Let*

$$V = \left\{ \begin{bmatrix} [0,a_1] & [0,a_2] & [0,a_3] \\ [0,a_4] & [0,a_5] & [0,a_6] \\ \vdots & \vdots & \vdots \\ \hline \vdots & \vdots & \vdots \\ [0,a_5] & [0,a_{n-4}] & [0,a_{n-3}] \\ [0,a_{n-2}] & [0,a_{n-1}] & [0,a_n] \end{bmatrix} \right.$$

where $a_i \in [0,1]$; $Z_{43}$ $1 \leq i \leq n$}, V is defined as the super fuzzy interval column vector.

We can have different partitions and different number of columns the only criteria being $0 \leq i \leq 1$.

We will illustrate this situation by some examples.

*Example 5.1.8*: Let

$$V = \begin{bmatrix} [0,1] & [0,0.3] & 0 \\ [0,0.4] & [0,1] & [0,1] \\ [0,1] & 0 & [0,0.4] \\ \hline [0,1] & [0,0.1] & [0,0.4] \\ 0 & [0,1] & [0,0.1] \\ 0 & [0,0.7] & 0 \\ 0 & [0,0.8] & [0,1] \\ \hline [0,1] & [0,0.3] & [0,0.1] \\ 0 & [0,0.7] & 0 \end{bmatrix}$$

be the super fuzzy interval column vector.

*Example 5.1.9*: Let

$$V = \begin{bmatrix} [0,a_1] & [0,a_2] & [0,a_3] & [0,a_{12}] \\ \hline [0,a_5] & [0,a_6] & [0,a_7] & [0,a_{14}] \\ [0,a_9] & [0,a_{10}] & [0,a_{11}] & [0,a_{18]} \\ \hline [0,a_{13}] & [0,a_{14}] & [0,a_{13}] & [0,a_{16}] \\ [0,a_{17}] & [0,a_{18}] & [0,a_{15}] & [0,a_{20}] \\ [0,a_{21}] & [0,a_{23}] & [0,a_{25}] & [0,a_{22}] \\ [0,a_{22}] & [0,a_{24}] & [0,a_{26}] & [0,a_{21}] \end{bmatrix}$$

where $0 \leq a_i \leq 1$ be a super fuzzy interval column vector.

Now we proceed onto define super fuzzy interval matrices or fuzzy super interval matrices.

**DEFINITION 5.1.5**: *Let*

$$V = \left\{ \begin{bmatrix} [0,a_1] & \cdots & [0,a_n] \\ [0,a_{n+1}] & \cdots & [0,a_{2n}] \\ \vdots & \vdots & \vdots \\ \vdots & & \\ [0,a_{(n-1)n+1}] & \cdots & [0,a_{mn}] \end{bmatrix} \right\}$$

*where* $0 \leq a_i \leq 1$. $i = 1, \ldots 2mn$. *V is defined as the super fuzzy interval matrix.*

It can be rectangular or square. Now we see we can on a given fuzzy interval matrix have several partitions which lead to different super fuzzy interval matrices. For instance consider the fuzzy interval row matrix.

X = ([0,0.5] [0, 0.3] [0, 0.1] [0, 1] 0 [0, 0.4]); how many super fuzzy interval row matrices can be obtained from X.

$X_1$ = ([0,0.5] | [0, 0.3] [0, 0.1] [0, 1] 0 [0, 0.4])
$X_2$ = ([0,0.5] [0, 0.3] | [0, 0.1] [0, 1] 0 [0, 0.4])
$X_3$ = ([0,0.5] [0, 0.3] [0, 0.1] | [0, 1] 0 [0, 0.4])
$X_4$ = ([0,0.5] [0, 0.3] [0, 0.1] [0, 1] | 0 [0, 0.4])
$X_5$ = ([0,0.5] [0, 0.3] [0, 0.1] [0, 1] 0 | [0, 0.4])
$X_6$ = ([0,0.5] | [0, 0.3] | [0, 0.1] [0, 1] | 0 [0, 0.4])
$X_7$ = ([0,0.5] [0, 0.3] | [0, 0.1] | …) and so on.

Thus by partitioning we get super interval fuzzy matrices which helps us to get several super interval fuzzy matrices. These structures can certainly find their applications in mathematical models.

*Example 5.1.10*: Let

$$V = \begin{bmatrix} [0,1] & [0,0.3] & [0,0.2] \\ [0,0.5] & [0,0.7] & [0,0.1] \\ [0,1] & [0,0.1] & [0,0.7] \\ [0,0.7] & [0,0.6] & [0,1] \\ [0,0.9] & 0 & [0,0.9] \end{bmatrix}$$

be a super interval fuzzy matrix we can have several partitions of the interval fuzzy matrix of order $5 \times 3$.

***Example 5.1.11***: Let

$$V = \left[\begin{array}{c|c} [0,1] & [0,0.3] \\ [0,0.4] & [0,0.1] \\ \hline [0,1] & [0,0.7] \\ [0,0.9] & [0,1] \end{array}\right]$$

be a super fuzzy interval matrix of natural order 4 × 2.

$$V = \begin{bmatrix} a_1 & a_2 \\ a_3 & a_4 \end{bmatrix}$$

where

$$a_1 = \begin{bmatrix} [0,1] \\ [0,0.4] \\ [0,1] \end{bmatrix}, \ a_2 = \begin{bmatrix} [0,0.3] \\ [0,1] \\ [0,0.7] \end{bmatrix},$$

$a_3 = ([0.9])$ and $a_4 = ([0, 1])$ are fuzzy interval submatrices of the super fuzzy interval matrix V.

***Example 5.1.12***: Let

$$V = \begin{bmatrix} a_1 & a_2 & a_9 \\ a_3 & a_4 & a_{10} \\ a_5 & a_6 & a_{11} \\ a_7 & a_8 & a_{12} \end{bmatrix}$$

be a super fuzzy interval matrix where

$$a_1 = \begin{pmatrix} [0,1] \\ [0,0.2] \end{pmatrix}, \ a_2 = \begin{pmatrix} 0 \\ [0,0.09] \end{pmatrix}$$

$$a_9 = \begin{pmatrix} [0,1] & 0 \\ 0 & [0,0.7] \end{pmatrix}, \ a_4 = \begin{pmatrix} 0 \\ [0,1] \\ [0,0.6] \end{pmatrix}, \ a_3 = \begin{pmatrix} [0,1] \\ 0 \\ [0,0.3] \end{pmatrix},$$

$$a_{10} = \begin{pmatrix} [0,1] & 0 \\ 0 & [0,0.3] \\ [0,0.7] & [0,0.1] \end{pmatrix}, \ a_5 = ([0, 1]), \ a_6 = ([0, 0.1101]),$$

$$a_{11} = ([0, 1] \ [0, 0.1 \ 01]), \ a_7 = \begin{pmatrix} [0,1] \\ [0,0.3] \\ 0 \end{pmatrix}, \ a_8 = \begin{pmatrix} 0 \\ [0,1] \\ 0 \end{pmatrix}$$

and

$$a_{12} = \begin{pmatrix} 0 & [0,1] \\ [0,1] & 0 \\ [0,0.3] & [0,0.1] \end{pmatrix}$$

are interval fuzzy submatrices of the super interval fuzzy matrix V. In fact V has 12 interval fuzzy submatrices. Further natural order of V is $9 \times 4$.

*Example 5.1.13*. Let

$$V = \begin{pmatrix} [0,1] & [0,0.2] \\ [0,0.5] & [0,0.3] \\ 0 & [0.1] \\ [0,0.4] & [0,0.8] \end{pmatrix} = \begin{pmatrix} a_1 & a_5 \\ a_2 & a_6 \\ a_3 & a_7 \\ a_4 & a_8 \end{pmatrix}$$

where $a_1, \ldots, a_8$ are the interval fuzzy submatrices of V. Further natural order of V is $4 \times 2$.

Now having seen examples of super fuzzy interval matrices we just give one more example of a super fuzzy interval square matrix whose natural order is a square matrix.

*Example 5.1.14*: Let

$$V = \begin{bmatrix} [0,1] & 0 & [0,0.3] & [0,1] \\ [0,0.3] & [0,1] & 0 & [0,0.2] \\ [0,0.2] & 0 & [0,1] & [0,0.3] \\ [0,1] & [0,1] & [0,0.2] & 0 \end{bmatrix}$$

be a super interval fuzzy matrix. $\begin{bmatrix} a_1 & a_2 & a_3 \\ a_4 & a_5 & a_6 \end{bmatrix}$ where $a_1, a_2, \ldots, a_6$ are fuzzy interval submatrices of the super fuzzy interval matrix V.

The natural order of V is 4 × 4 however its super interval matrix order is 2 × 3.

Now we proceed onto define min or max function on V a super fuzzy interval matrix.

Let V = ([0, $a_1$] [0,$a_2$] | ... | [0,$a_{12}$]) and W = ([0, $b_1$] [0, $b_2$] | . . . | [0, $b_{12}$]) be two super fuzzy interval row matrices where $0 \leq a_i, b_i \leq 1$ of same type and i = 1, 2, ..., 12.
Define

V + W = min (V, W)
   = (min ([0, $a_1$] [$a_1$,$b_1$]), min ([0, $a_2$] [0,$b_2$]) | ... | min ([0,$a_{12}$], [0,$b_{12}$])
   = ([0, min (a, $b_1$)] [0, min ($a_2$,$b_2$)] | . . . | [0, min ($a_{12,}$ $b_{12}$)]).

Thus we can define min operation on two super interval fuzzy row matrices of same order with identical partition.
Even if the partitions of the same natural order matrix is different we cannot define the min operation them.

We will illustrate this situation by some examples.

Let x = ([0, 1] | [0,0.3] [0, 0.4] | 0 [0,0.8], [0,0.2] | [0, 0.9]) and y = (0 | [0, 1] [0,0.5] | [0, 0.2] [0,0.7] [0,0.2] | [0, 0.1]) be two super fuzzy interval row matrices of same order with identical partition.
min (x, y) = (0 | [0, 0.3] [0,0.4] | 0 [0, 0.7] [0,0.2] | [0, 0.1]).
Thus when we have a collection of P super fuzzy interval row matrices with identical partition on them, then (P, min) is a semigroup known as the semigroup of super fuzzy interval row matrices.

We can show this by some examples.

**Example 5.1.15**: Let V = {([0, $a_1$] [0,$a_2$] | [0, $a_3$] [0,$a_4$] [0,$a_5$] | [0, $a_6$]) | $a_i$ ∈ [0,1]; 1 ≤ i ≤ 6} be the collection of all super fuzzy interval row matrices of order 1 × 6 with min operation defined on it. (V, min) is a semigroup of super fuzzy interval row matrices. Clearly cardinality of V is infinite.
If we in the definition of min replace by max still this will give a semigroup.

**Example 5.1.16**: Let V = {([0, $a_1$] | [0,$a_2$] ... [0, $a_{10}$] | [0,$a_{12}$] ... [0,$a_{11}$]) | $a_i \in$ [0,1]; $1 \le i \le 15$}; V under max operation is a semigroup of super interval fuzzy row matrices. Clearly order of V is infinite.

**Example 5.1.17**: Let V = {([0, $a_1$] | [0,$a_2$] | [0, $a_3$] [0,$a_4$] | [0,$a_5$]) | $a_i \in$ [0,1]; $1 \le i \le 5$} (V, max) is a fuzzy semigroup of super interval row matrices as well as (V, min) is a fuzzy super interval semigroup of row matrices.

Now we proceed onto define for a super fuzzy interval matrix multiplication by elements in [0,1]. We have three ways of finding product of any element [0,1] with a super fuzzy interval matrix.

Let A = ([0, $a_1$] [0,$a_2$] [0, $a_3$] | [0,$a_4$] | [0,$a_5$] [0, $a_6$]) where $a_i \in$ [0,1]; $1 \le i \le 6$}. Take $b_i \in$ [0,1]. min {$b_i$A} = ([0, min {$b_i$ $a_1$}] [0, min {$b_i$, $a_2$}] [0, min {$b_i$ $a_3$}] | [0, min {$b_i$ $a_4$}] | [0, min {$b_i$, $a_5$}] [0, min {$b_i$ $a_6$}]). We will illustrate this situation by an example.

Let A = ([0, 0.5] [0,0.2] [0, 0.1] | [0,0.7] [0,0.1] | [0, 0.31] | [0, 0.302] [0, 0.251] [0, 0.87]) be a super fuzzy interval row matrix. Let x = 0.31 $\in$ [0,1]. To find min (x, A); min (x, A) = ([min {0.31, 0.5}] [0, min {0.31, 0.2}] [0, min {0.31, 0.1}] | [0, min {0.31, 0.7}] [0, min {0.31, 0.1}] | [0, min {0.31, 0.31}] | [0, min {0.31, 0.302 }] [0, min {0.31, 0.251}] [0, min {0.31, 0.87}]) = ([0, 0.31] [0, 0.2] [0, 0.1] | [0, 0.31] [0, 0.1] | [0, 0.31] | [0, 0.302] [0, 0.251] [0, 0.31]).

This is min product method with a fuzzy number. Next we define the max product method with a fuzzy element from [0, 1] with a fuzzy super row interval matrix.

Let A = ([0, 0.2] [0,0.7] | [0, 0.3] 0 [0,1] [0,0.42] | [0, 1] 0 [0, 0.201] | [0, 0.001]) be a super fuzzy interval row matrix. Let x = 0.1 $\in$ [0,1] to find max (xA) = ([0, max{0.1, 0.2}] [0, max {0.1, 0.7}] | [0, max {0.1, 0.3}] ([0, max{0.1, 0}] [0, max {0.1, 1} [0, max {0.1, 0.42}] | [0, max {0.1, 1}] [0, max {0.1, 0}] [0, max {0.1, 0.201}] | [0, max {0.1, 0.001}]) = ([0, 0.2] [0, 0.7] | [0,0.3] [0, 0.1] [0, 1] [0, 0.42] | [0, 1] [0, 0.1] [0, 0.201] [0,0.1]).

Now similarly if the column super interval matrix A the entries are from the interval [0, 1] then we define A to be a super fuzzy column interval matrix.

We will first illustrate it by some examples.

*Example 5.1.18:* Let

$$V = \begin{bmatrix} [0,1] \\ \overline{[0,0.5]} \\ \overline{[0,0.2]} \\ [0,0.3] \\ [0,0.8] \end{bmatrix}$$

is the super fuzzy interval column matrix of natural order $5 \times 1$.

*Example 5.1.19*: Let

$$P = \begin{bmatrix} [0,0.7] \\ [0,0.1] \\ [0,1] \\ \overline{[0,0.3]} \\ \overline{[0,0.4]} \\ [0,0.5] \\ \overline{[0,0.6]} \\ \overline{[0,0.1]} \\ [0,0.5] \\ [0,0.7] \\ [0,0.8] \end{bmatrix}$$

be the super fuzzy interval column matrix of natural order $11 \times 1$.

*Example 5.1.20*: Let

$$T = \begin{bmatrix} \overline{[0,1]} \\ [0,0.2] \\ [0,0.3] \\ [0,0.4] \\ \overline{[0,0.5]} \end{bmatrix}$$

is the super fuzzy interval column matrix of natural order $5 \times 1$.

Now we can define the max or min only on the super fuzzy interval matrices identical partition and same natural order.

Thus if

$$x = \begin{bmatrix} [0,1] \\ 0 \\ \overline{[0,0.3]} \\ \overline{[0,0.7]} \\ [0,0.5] \\ \overline{[0,0.1]} \\ [0,0.8] \\ [0,0.5] \end{bmatrix} \text{ and } y = \begin{bmatrix} [0,0.3] \\ [0,0.8] \\ \overline{[0,1]} \\ \overline{[0,0.9]} \\ [0,0.2] \\ \overline{[0,0.4]} \\ [0,0.5] \\ [0,0.9] \end{bmatrix}$$

be two super fuzzy column interval matrices. We define min

$$\{x, y\} = \begin{bmatrix} [0,\min\{1,0.3\}] \\ [0,\min\{0,0.8\}] \\ \overline{[0,\min\{0.3,1\}]} \\ \overline{[0,\min\{0.7,0.9\}]} \\ [0,\min\{0.5,0.2\}] \\ \overline{[0,\min\{0.1,0.4\}]} \\ [0,\min\{0.8,0.5\}] \\ [0,\min\{0.5,0.9\}] \end{bmatrix} = \begin{bmatrix} [0,0.3] \\ [0,0] \\ \overline{[0,0.3]} \\ \overline{[0,0.7]} \\ [0,0.2] \\ \overline{[0,0.1]} \\ [0,0.5] \\ [0,0.5] \end{bmatrix}.$$

Now we can say the collection of all super fuzzy interval matrices with same natural order and same type of partition under min operation is a semigroup known as the semigroup of super interval fuzzy column matrices or fuzzy semigroup of super interval column matrices with min operation.

Now we can give one or two examples of fuzzy semigroups before we proceed to define max operation on them.

***Example 5.1.21***: Let

$$V = \left\{ \begin{bmatrix} \overline{[0,a_1]} \\ \overline{[0,a_2]} \\ [0,a_3] \\ \overline{[0,a_4]} \\ [0,a_5] \\ \overline{[0,a_6]} \\ [0,a_7] \\ [0,a_8] \\ [0,a_9] \\ [0,a_{10}] \end{bmatrix} \middle| \ a_i \in [0,1];\ 1 \le i \le 10 \right\}$$

be the collection of all super fuzzy interval column matrices of natural order $10 \times 1$ with the same type of partition on it $\{V, \min\}$; V with min operation is a fuzzy semigroup.

***Example 5.1.22***: Let

$$M = \left\{ \begin{bmatrix} [0,a_1] \\ [0,a_2] \\ \overline{[0,a_3]} \\ [0,a_4] \\ [0,a_5] \\ [0,a_6] \\ \overline{[0,a_7]} \\ \overline{[0,a_{10}]} \\ [0,a_{11}] \\ [0,a_{12}] \\ [0,a_{13}] \end{bmatrix} \middle| \ a_i \in [0,\tfrac{1}{2}] \subseteq [0,1];\ 1 \le i \le 13 \right\}$$

be a super fuzzy interval column matrix of same natural order $13 \times 1$ and similar type of partition. M is a fuzzy semigroup of super interval column matrices with min operation defined on it.

One of the advantages of using super fuzzy interval matrices is we get several fuzzy semigroups of same natural order depending upon the partition defined on it.

For instance take

$$X = \begin{bmatrix} [0, a_1] \\ [0, a_2] \\ [0, a_3] \\ [0, a_4] \\ [0, a_5] \end{bmatrix}$$

a $5 \times 1$ fuzzy interval column matrix. We can have only one fuzzy semigroup by defining the min operation on it. But on the other hand if we partition X we have the following partitions of X.

$$X_1 = \begin{bmatrix} [0, a_1] \\ \overline{[0, a_2]} \\ \vdots \\ [0, a_5] \end{bmatrix}, X_2 = \begin{bmatrix} [0, a_1] \\ [0, a_2] \\ \overline{[0, a_3]} \\ [0, a_4] \\ [0, a_5] \end{bmatrix}, X_3 = \begin{bmatrix} [0, a_1] \\ [0, a_2] \\ [0, a_3] \\ \overline{[0, a_4]} \\ [0, a_5] \end{bmatrix}, X_4 = \begin{bmatrix} [0, a_1] \\ \vdots \\ [0, a_4] \\ \overline{[0, a_5]} \end{bmatrix},$$

$$X_5 = \begin{bmatrix} [0, a_1] \\ \overline{[0, a_2]} \\ [0, a_3] \\ [0, a_4] \\ [0, a_5] \end{bmatrix}, X_6 = \begin{bmatrix} [0, a_1] \\ [0, a_2] \\ \overline{[0, a_3]} \\ [0, a_4] \\ [0, a_5] \end{bmatrix}, X_7 = \begin{bmatrix} [0, a_1] \\ [0, a_2] \\ [0, a_3] \\ \overline{[0, a_4]} \\ [0, a_5] \end{bmatrix}, X_8 = \begin{bmatrix} [0, a_1] \\ [0, a_2] \\ \overline{[0, a_3]} \\ [0, a_4] \\ [0, a_5] \end{bmatrix},$$

$$X_9 = \begin{bmatrix} [0, a_1] \\ [0, a_2] \\ \overline{[0, a_3]} \\ [0, a_4] \\ [0, a_5] \end{bmatrix}, X_{10} = \begin{bmatrix} [0, a_1] \\ [0, a_2] \\ \overline{[0, a_3]} \\ [0, a_4] \\ [0, a_5] \end{bmatrix}, X_{11} = \begin{bmatrix} [0, a_1] \\ [0, a_2] \\ \overline{[0, a_3]} \\ [0, a_4] \\ [0, a_5] \end{bmatrix}, X_{12} = \begin{bmatrix} [0, a_1] \\ [0, a_2] \\ \overline{[0, a_3]} \\ [0, a_4] \\ [0, a_5] \end{bmatrix},$$

$$X_{13} \begin{bmatrix} [0,a_1] \\ [0,a_2] \\ \hline [0,a_3] \\ [0,a_4] \\ [0,a_5] \end{bmatrix} \text{ and so on.}$$

We get more than 13 such fuzzy semigroups which is one of the advantages of using super interval fuzzy matrices instead of interval fuzzy matrices.

Now we can define max operation on fuzzy super interval column matrices which will be described by the following example.

Consider

$$X = \begin{bmatrix} [0,0.3] \\ [0,1] \\ \hline [0,0.7] \\ [0,0.1] \\ 0 \\ \hline [0,0.8] \\ [0,0.9] \\ [0,0.4] \\ [0,0.3] \end{bmatrix} \text{ and } Y = \begin{bmatrix} [0,0.2] \\ 0 \\ \hline [0,0.9] \\ 0 \\ [0,0.2] \\ \hline [0,0.7] \\ [0,0.1] \\ [0,1] \\ [0,0.5] \end{bmatrix}$$

two super interval fuzzy column matrices.
Now

$$\max\{X, Y\} = \begin{bmatrix} [0,\max\{0.3,0.2\}] \\ [0,\max\{1,0\}] \\ \hline [0,\max\{0.7,0.9\}] \\ [0,\max\{0.1,0\}] \\ [0,\max\{0,0.2\}] \\ \hline [0,\max\{0.8,0.7\}] \\ [0,\max\{0.9,0.1\}] \\ [0,\max\{0.4,1\}] \\ [0,\max\{0.3,0.5\}] \end{bmatrix} = \begin{bmatrix} [0,0.3] \\ [0,1] \\ \hline [0,0.9] \\ [0,0.1] \\ [0,0.2] \\ \hline [0,0.8] \\ [0,0.9] \\ [0,0.1] \\ [0,0.5] \end{bmatrix}.$$

Thus max for super interval fuzzy column matrices can be defined.

Let

$$V = \left\{ \begin{bmatrix} [0,a_1] \\ [0,a_2] \\ \hline [0,a_3] \\ [0,a_4] \\ \hline [0,a_5] \\ [0,a_6] \\ \vdots \\ \hline [0,a_r] \\ \hline [0,a_{r+1}] \\ \vdots \\ \hline [0,a_n] \end{bmatrix} \right.$$

$a_i \in [0,1]$; $1 \leq i \leq n$} be the collection of super interval fuzzy column matrices. V with max operation is a semigroup defined as the semigroup of super interval fuzzy column matrices.

Now we can define for any $a \in [0, 1]$ and for any super interval fuzzy column matrix in two ways using max or min principle. We will describe it with examples.

Consider

$$x = \begin{bmatrix} [0,0.5] \\ [0,0.7] \\ [0,1] \\ \hline 0 \\ \hline [0,0.8] \\ [0,0.7] \\ \hline [0,1] \\ 0 \\ [0,0.9] \\ \hline [0,0.3] \\ [0,0.2] \end{bmatrix}$$

and y = 0.6. To find y x we use now the min principle

$$\min(yx) = \begin{bmatrix} [0,\min\{0.6,0.5\}] \\ [0,\min\{0.6,0.7\}] \\ [0,\min\{0.6,1\}] \\ \overline{[0,\min\{0.6,0\}]} \\ \overline{[0,\min\{0.6,0.8\}]} \\ [0,\min\{0.6,0.7\}] \\ \overline{[0,\min\{0.6,1\}]} \\ [0,\min\{0.6,0\}] \\ [0,\min\{0.6,0.9\}] \\ \overline{[0,\min\{0.6,0.3\}]} \\ [0,\min\{0.6,0.2\}] \end{bmatrix} = \begin{bmatrix} [0,0.5] \\ [0,0.6] \\ [0,0.6] \\ [0,0] \\ \overline{[0,0.6]} \\ [0,0.6] \\ \overline{[0,0.6]} \\ [0,0] \\ [0,0.6] \\ \overline{[0,0.3]} \\ [0,0.2] \end{bmatrix}.$$

Thus we see with min product on a super fuzzy interval column matrix the structure is compatible. Now we can use also the product max on super fuzzy interval matrix.

Let

$$P = \begin{bmatrix} [0,1] \\ [0,0.3] \\ [0,0.2] \\ [0,0.9] \\ [0,0.7] \\ \overline{[0,0.1]} \\ [0,0.5] \\ [0,0.8] \\ \overline{0} \\ [0,0.5] \\ [0,0.6] \\ [0,0.7] \\ \overline{[0,0.9]} \end{bmatrix}.$$

Let x = 0.1 ∈ [0, 1].

$$\text{Max}(x, P) = \begin{bmatrix} [0,1] \\ [0,0.3] \\ [0,0.2] \\ [0,0.9] \\ [0,0.7] \\ \overline{[0,0.1]} \\ [0,0.5] \\ [0,0.8] \\ \overline{[0,0.1]} \\ [0,0.5] \\ [0,0.6] \\ [0,0.7] \\ \overline{[0,0.9]} \end{bmatrix}.$$

Now we have another operation called product which is compatible. If $X = ([0,a_1] \, [0,a_2] \mid [0,a_3] \ldots \mid [0,a_r] \ldots [0,a_n])$ with $a_i \in [0,1]$; $1 \leq i \leq n\}$ and b in [0, 1] then $b \times X = ([0,ba_1] \, [0,ba_2] \mid ([0,ba_3] \ldots \mid ([0,ba_r] \ldots ([0,ba_n])$; $ba_i \in [0,1]$; $1 \leq i \leq n$.

Likewise in case of super interval column fuzzy matrix we have for

$$Y = \begin{bmatrix} [0,a_1] \\ [0,a_2] \\ [0,a_3] \\ \overline{[0,a_4]} \\ \vdots \\ [0,a_t] \\ \vdots \\ [0,a_m] \end{bmatrix}$$

$a_i \in [0,1]$; $1 \leq i \leq m$ be a super fuzzy interval column matrix. Let $b \in [0,1]$,

$$bY = \begin{bmatrix} [0, ba_1] \\ [0, ba_2] \\ [0, ba_3] \\ \hline [0, ba_4] \\ \vdots \\ [0, ba_t] \\ \vdots \\ [0, ba_m] \end{bmatrix}$$

where $ba_i \in [0,1]$; $1 \leq i \leq m$. So we can use this product also and thus we have three types of products used on these matrices. Now if in a super interval row vector X, we take the elements only from the interval [0,1] we get X to be a super fuzzy interval row vector.

We will just illustrate it by an example. Let

$$X = \begin{pmatrix} [0,a_1] & [0,a_6] & | & [0,a_{11}] & [0,a_{16}] & [0,a_{21}] & | & [0,a_{26}] \\ [0,a_2] & [0,a_7] & | & [0,a_{12}] & [0,a_{17}] & [0,a_{22}] & | & [0,a_{27}] \\ [0,a_3] & [0,a_8] & | & [0,a_{13}] & [0,a_{18}] & [0,a_{23}] & | & [0,a_{28}] \\ [0,a_4] & [0,a_9] & | & [0,a_{14}] & [0,a_{19}] & [0,a_{24}] & | & [0,a_{29}] \\ [0,a_5] & [0,a_{10}] & | & [0,a_{15}] & [0,a_{20}] & [0,a_{25}] & | & [0,a_{30}] \end{pmatrix}$$

where $a_i \in [0,1]$; $1 \leq i \leq 30$ is a super fuzzy interval row vector of natural order $5 \times 6$ but $X = (b_1 | b_2 | b_3)$ where $b_1$, $b_2$ and $b_2$ are the interval fuzzy submatrices of the super interval fuzzy matrix X.

Now if in the super interval column vector

$$Y = \begin{bmatrix} [0,a_1] & [0,a_2] & [0,a_3] \\ [0,a_2] & [0,a_5] & [0,a_6] \\ [0,a_7] & [0,a_8] & [0,a_9] \\ \hline [0,a_{10}] & [0,a_{11}] & [0,a_{12}] \\ [0,a_{13}] & [0,a_{14}] & [0,a_{15}] \\ [0,a_{16}] & [0,a_{17}] & [0,a_{18}] \\ [0,a_{19}] & [0,a_{20}] & [0,a_{21}] \\ [0,a_{22}] & [0,a_{23}] & [0,a_{24}] \\ [0,a_{25}] & [0,a_{26}] & [0,a_{27}] \end{bmatrix}$$

we replace the values $a_i$ from the interval [0,1] then we get the super fuzzy interval column vector where $a_i \in [0,1]$; $1 \le i \le 27$.

Thus

$$y = \begin{bmatrix} \underline{c_1} \\ \underline{c_2} \\ \underline{c_3} \\ c_4 \end{bmatrix} = \begin{bmatrix} c_1 \\ c_2 \\ c_3 \\ c_4 \end{bmatrix}$$

where each $c_i$ is a fuzzy interval matrix, $1 \le i \le 4$. $c_i$'s are the fuzzy interval submatrices of the super interval column vector Y. Now having seen examples of super fuzzy column interval vector and super fuzzy interval row vector we now proceed onto define operations on them. Suppose V be a collection of super interval fuzzy row vectors of same natural order and same type of partition then on V we can define the min operation so that V with the min operation becomes a semigroup, known as the fuzzy super interval row vector semigroup with min operation. In fact (V, min) is a commutative semigroup.

Now if on the collection V of super fuzzy interval matrices (row vectors) we define the max operation then also V with max operation is a semigroup which is commutative we will illustrate this by some examples.

***Example 5.1.23***: Let

$$V = \left\{ \begin{pmatrix} [0,a_1] & [0,a_5] & [0,a_9] & [0,a_{13}] & [0,a_{12}] & [0,a_{21}] \\ [0,a_2] & [0,a_6] & [0,a_{10}] & [0,a_{14}] & [0,a_{18}] & [0,a_{22}] \\ [0,a_3] & [0,a_7] & [0,a_{11}] & [0,a_{15}] & [0,a_{19}] & [0,a_{23}] \\ [0,a_4] & [0,a_8] & [0,a_{12}] & [0,a_{16}] & [0,a_{20}] & [0,a_{24}] \end{pmatrix} \right.$$

$a_i \in [0,1]; \cup \{0\}$; $1 \le i \le 24\}$ is a semigroup of super fuzzy interval row vector under min operation.

This same V if max operation is defined then also V is a fuzzy semigroup of super interval fuzzy row vectors. Clearly the natural order of superinterval matrices in V is $4 \times 24$ but V = $\{(b_1\ b_2\ b_3)\}$ where $b_i$'s are interval fuzzy matrices that $b_1$ $b_2$ and $b_2$ are fuzzy interval submatrices of V.

Now if

$$A = \begin{pmatrix} [0,a_1] & [0,a_2] & | & [0,a_3] & | & [0,a_4] & [0,a_7] \\ [0,a_6] & [0,a_7] & | & [0,a_8] & | & [0,a_9] & [0,a_{10}] \\ [0,a_{11}] & [0,a_{12}] & | & [0,a_{13}] & | & [0,a_{14}] & [0,a_{15}] \\ \hline [0,a_{16}] & [0,a_{17}] & | & [0,a_{18}] & | & [0,a_{19}] & [0,a_{20}] \\ [0,a_{21}] & [0,a_{22}] & | & [0,a_{23}] & | & [0,a_{24}] & [0,a_{25}] \\ \hline [0,a_{26}] & [0,a_{27}] & | & [0,a_{28}] & | & [0,a_{29}] & [0,a_{29}] \end{pmatrix} = \begin{pmatrix} b_1 & b_4 & b_7 \\ b_2 & b_5 & b_8 \\ b_3 & b_6 & b_9 \end{pmatrix}$$

where $1 \leq i \leq 1$; $1 \leq i \leq 30$ is a super fuzzy interval matrix of natural order $6 \times 5$. We have $b_i$'s to be the interval fuzzy submatrices of A. Thus we can say a super fuzzy interval matrix is a super interval matrix in which entries are from [0,1] and submatrices are also fuzzy interval matrices.

We will give one or two examples of them before we proceed to give some structures on them.

Let

$$V = \begin{pmatrix} [0,0.3] & [0,0.8] & | & [0,0.5] & 0 & | & [0.0.1] \\ [0,0.7] & [0,1] & | & [0,0.8] & [0,0.9] & | & 0 \\ \hline [0,0.8] & 0 & | & [0,0.3] & 0 & | & [0,0.125] \\ 0 & [0,0.9] & | & [0,1] & [0,0.4] & | & 9 \\ \hline [0,1] & 0 & | & [0,0.1] & 0 & | & [0,0.15] \\ [0,0.5] & [0,0.7] & | & 0 & [0,0.1] & | & 0 \end{pmatrix}$$

be a super fuzzy interval matrix.

$$V = \begin{pmatrix} c_1 & c_2 & c_3 \\ c_4 & c_5 & c_6 \\ c_7 & c_8 & c_9 \end{pmatrix},$$

$c_i$'s are the submatrices of the super fuzzy interval matrix $1 \leq i \leq = 9$.

Thus

$$c_1 = \begin{pmatrix} [0,0.3] & [0,0.8] \\ [0,0.7] & [0,1] \end{pmatrix}, c_2 = \begin{pmatrix} [0,0.5] & 0 \\ [0,0.8] & [0,0.9] \end{pmatrix}, c_3 = \begin{pmatrix} [0,0.1] \\ [0,0] \end{pmatrix},$$

$$c_4 = \begin{pmatrix} [0,0.8] & [0,0] \\ [0,0] & [0,0.9] \\ [0,1] & [0,0] \end{pmatrix}, c_5 = \begin{pmatrix} [0,0.3] & 0 \\ [0,1] & [0,0.4] \\ [0,0.1] & 0 \end{pmatrix}, c_7 = ([0,0.5]\ [0,0.07]),$$

$$c_8 = (0, [0, 0.1]), c_9 = ([0, 0] = (0) \text{ and } c_6 = \begin{pmatrix} [0,0.125] \\ 0 \\ [0,0.15] \end{pmatrix}$$

are the interval submatrices of V.
Suppose

$$V = \left( \begin{array}{cc|cc|cc} [0,1] & 0 & [0,1] & [0,0.5] & [0,0.9] \\ [0,0.2] & [0,1] & [0] & [0,0.7] & 0 \\ [0,0.3] & 0 & [0,0.9] & [0,0.92] & [0,0.112] \\ [0,0.4] & [0,0.1] & 0 & [0,0.1] & 0 \\ 0 & [0,0.1] & [0,0.2] & 0 & [0,0.01] \\ \hline [0,0.5] & [0,0.1] & [0,0.7] & [0,0.1] & 0 \\ [0,0.6] & 0 & [0,0.9] & 0 & [0,0.02] \\ \hline [0,0.7] & [0,0.4] & [0,0.11] & [0,0.13] & 0 \\ 0 & [0,1] & [0,0.15] & 0 & [0,0.75] \\ [0,0.8] & 0 & [0,0.61] & [0,0.7] & [0,0.1] \end{array} \right)$$

$$= \begin{pmatrix} B_1 & B_4 & B_7 \\ B_2 & B_5 & B_8 \\ B_3 & B_6 & B_9 \end{pmatrix}$$

where $B_i$'s are subinterval fuzzy matrices or fuzzy interval submatrices of the super interval fuzzy matrix V.

Now having seen examples we can as in case of superinterval fuzzy row vector define min operation or max operation with scalars from [0,1], however we cannot define max or min operation among themselves. We can as in case of other super interval matrices define the notion of set super fuzzy interval matrices vector space of set fuzzy vector space or set fuzzy interval supervector space which is only a matter of routine. Hence we give only examples of these structures.

**Example 5.1.24**: Let

$$V = \left\{ \begin{pmatrix} [0,a_1] & [0,a_3] & | & [0,a_5] & [0,a_7] & | & [0,a_9] \\ [0,a_2] & [0,a_4] & | & [0,a_6] & [0,a_8] & | & [0,a_{10}] \end{pmatrix}, \right.$$

$([0,a_1] | [0, a_2]), [0, a_3] | [0, a_4])$; $a_i \in [0,1]$; $1 \leq i \leq 10$}, V is a set fuzzy vector space of super interval fuzzy matrices over the set $S = \{0, 1\}$.

**Example 5.1.25**: Let

$$V = \left\{ \begin{bmatrix} [0,a_1] \\ [0,a_2] \\ \overline{[0,a_3]} \\ [0,a_4] \\ \overline{[0,a_5]} \\ [0,a_6] \end{bmatrix}, \right.$$

$([0,a_1] | [0, a_2] [0, a_3] [0, a_4] | [0, a_5] | [0, a_6] [0, a_7]) | a_i \in [0,1]$; $1 \leq i \leq 7$}, be a set fuzzy vector space of super interval fuzzy matrices over the set $S = \{0, 0.02, 0.07, 0, 0.5, 1\}$.

**Example 5.1.26**: Let

$$V = \left\{ \begin{pmatrix} [0,a_1] & | & [0,a_5] \\ [0,a_2] & | & [0,a_6] \\ [0,a_3] & | & [0,a_7] \\ [0,a_4] & | & [0,a_8] \end{pmatrix}, \begin{pmatrix} \overline{[0,a_1]} & [0,a_5] \\ [0,a_2] & [0,a_6] \\ [0,a_3] & [0,a_7] \\ [0,a_4] & [0,a_8] \end{pmatrix}, \begin{pmatrix} [0,a_1] & | & [0,a_5] \\ [0,a_2] & | & [0,a_6] \\ [0,a_3] & | & [0,a_7] \\ [0,a_4] & | & [0,a_8] \end{pmatrix}, \right.$$

$$\begin{pmatrix} [0,a_1] & | & [0,a_5] \\ [0,a_2] & | & [0,a_6] \\ [0,a_3] & | & [0,a_7] \\ [0,a_4] & | & [0,a_8] \end{pmatrix}$$

$a_i \in [0,0.6]$; $1 \leq i \leq 8\}$, be the collection of super interval fuzzy matrices. V is a set fuzzy vector space of super fuzzy interval matrices over the set $S = \{0, \frac{1}{2}, 1\}$.

*Example 5.1.27*: Let

$$S = \left\{ \begin{pmatrix} [0,a_1] & [0,a_3] & [0,a_5] & [0,a_7] & [0,a_9] & [0,a_{11}] \\ [0,a_2] & [0,a_4] & [0,a_6] & [0,a_8] & [0,a_{10}] & [0,a_{12}] \end{pmatrix}, \right.$$

$$\left. \left( \begin{array}{c|c} [0,a_1] & [0,a_2] \\ \hline [0,a_3] & [0,a_4] \end{array} \right), \begin{bmatrix} [0,a_1] & [0,a_7] \\ [0,a_2] & [0,a_8] \\ [0,a_3] & [0,a_9] \\ [0,a_4] & [0,a_{10}] \\ [0,a_5] & [0,a_{11}] \\ [0,a_6] & [0,a_{12}] \end{bmatrix} \right\}$$

$a_i \in [0,1]$; $1 \leq i \leq 12\}$ be a collection of super fuzzy interval matrices. S is a set vector of fuzzy super interval matrices over the set $T = \{0, 1\}$.

*Example 5.1.28*: Let

$$W = \left\{ \left( \begin{array}{c|c} [0,1] & 0 \\ \hline [0,0.2] & [0,0.3] \\ \hline [0,0.4] & 0 \end{array} \right), \right.$$

$$\begin{pmatrix} [0,0.2] & [0,0.4] & [0,0.6] & [0,0.8] & [0,1] \\ 0 & [0,1] & 0 & [0,0.2] & 0 \\ [0,0.4] & 0 & [0,0.7] & 0 & [0,0.7] \\ [0,0.5] & [0,0.4] & [0,0.3] & [0,0.4] & [0,0.2] \end{pmatrix},$$

$$\left. \left( \begin{array}{c|c} 0 & 0 \\ 0 & 0 \\ \hline 0 & 0 \end{array} \right), \begin{pmatrix} 0 & 0 & 0 & 0 & 0 \\ 0 & 0 & 0 & 0 & 0 \\ 0 & 0 & 0 & 0 & 0 \\ 0 & 0 & 0 & 0 & 0 \end{pmatrix}, \right.$$

$$\left\{ \begin{bmatrix} 0 \\ 0 \\ \overline{0} \\ 0 \\ 0 \\ \overline{0} \\ 0 \\ 0 \end{bmatrix}, \begin{bmatrix} [0,0.5] \\ [0,0.7] \\ \overline{[0,0.3]} \\ 0 \\ [0,1] \\ [0,0.5] \\ \overline{0} \\ [0,0.4] \\ [0,0.7] \end{bmatrix} \right\}$$

be a set of super fuzzy interval matrices. Clearly W is a set vector space of super fuzzy interval matrices over the set $T = \{0, 1\}$ under min operator. It is pertinent to mention that max product of T over W is not defined.

Now having seen examples of set vector spaces of super fuzzy interval matrices we now proceed onto give examples of their substructures.

*Example 5.1.29*: Let

$$V = \left\{ \begin{pmatrix} [0,a_1] & [0,a_3] & [0,a_5] & [0,a_7] & [0,a_9] & [0,a_{11}] \\ [0,a_2] & [0,a_4] & [0,a_6] & [0,a_8] & [0,a_{10}] & [0,a_{12}] \end{pmatrix}, \begin{bmatrix} [0,a_1] \\ [0,a_2] \\ [0,a_3] \\ [0,a_4] \\ [0,a_5] \\ [0,a_6] \end{bmatrix}, \right.$$

$$\left. \begin{pmatrix} \overline{[0,a_1]} & [0,a_4] \\ [0,a_2] & [0,a_5] \\ [0,a_3] & [0,a_6] \end{pmatrix} \right|$$

$a_i \in [0,1]$; $1 \le i \le 12\}$ be the collection of super interval fuzzy matrices; V is a set vector of space of super interval fuzzy matrices over the set $S = \{0, 1\}$.

Consider

$$W = \left\{ \left( \begin{array}{cc} \overline{[0,a_1]} & [0,a_4] \\ \overline{[0,a_2]} & [0,a_5] \\ [0,a_3] & [0,a_6] \end{array} \right), \left| \begin{array}{c} [0,a_1] \\ \overline{[0,a_2]} \\ [0,a_3] \\ \overline{[0,a_4]} \\ [0,a_5] \\ [0,a_6] \end{array} \right| \right\} = a_i \in [0,1];$$

$1 \leq i \leq 6\} \subseteq V$; W is a set vector subspace of super fuzzy interval matrices of V over the set $S = \{0, 1\}$.

*Example 5.1.30*: Let

$$V = \left\{ \left| \left( \begin{array}{c} [0,a_1] \\ [0,a_2] \\ \overline{[0,a_3]} \\ [0,a_4] \\ [0,a_5] \\ [0,a_6] \\ \overline{[0,a_7]} \\ [0,a_8] \\ [0,a_9] \end{array} \right) \right|, \left( \begin{array}{cc|cc|cc} [0,a_1] & [0,a_3] & [0,a_5] & [0,a_7] & [0,a_9] & [0,a_{11}] \\ [0,a_2] & [0,a_4] & [0,a_6] & [0,a_8] & [0,a_{10}] & [0,a_{12}] \end{array} \right) \right\}$$

$a_i \in [0,1]$; $1 \leq i \leq 12\}$ be a set vector space of super interval fuzzy matrices over the set $S = \{0, 1\}$.

Consider

$$W = \left\{ \left( \begin{array}{cc|cc|cc} [0,a_1] & [0,a_3] & [0,a_5] & [0,a_7] & [0,a_9] & [0,a_{11}] \\ [0,a_2] & [0,a_4] & [0,a_6] & [0,a_8] & [0,a_{10}] & [0,a_{12}] \end{array} \right) \right\}$$

$a_i \in [0,1]$; $1 \leq i \leq 12\} \subseteq V$, W is a set vector subspace of V of super interval fuzzy matrices over the set $S = \{0, 1\}$.

Now having seen substructures we now proceed onto define set linear algebra of super fuzzy interval matrices or set fuzzy linear algebra of super fuzzy interval matrices.

A set vector space V of super fuzzy interval matrices over the set S is defined to be a set linear algebra over the sets if (i) V is closed under the binary operation '+' that is an additive abelian group that is for all a, b ∈ V; a + b ∈ V (+ is min or max operation only) (ii) s · (a + b) = s · a + s · b for all a, b ∈ V and s ∈ S, '·' is the min or max operation only. This will be illustrated by the following examples.

*Example 5.1.31*: Let

$$V = \left\{ \begin{pmatrix} \underline{[0,a_1] \quad [0,a_8]} \\ \underline{[0,a_2] \quad [0,a_9]} \\ \underline{[0,a_3] \quad [0,a_{10}]} \\ \underline{[0,a_4] \quad [0,a_{11}]} \\ \underline{[0,a_5] \quad [0,a_{12}]} \\ \underline{[0,a_6] \quad [0,a_{13}]} \\ [0,a_7] \quad [0,a_{14}] \end{pmatrix} \;\middle|\; a_i \in [0,1];\; 1 \le i \le 14 \right\}$$

be the collection of all super fuzzy interval matrices of same type. V is closed under min operation. Choose the set $S = \{0, \frac{1}{2n} ; n \in N\}$. V is a set linear algebra over the set S. Thus V is a set fuzzy linear algebra of super interval column vectors over the set S.

*Example 5.1.32*: Let

$$V = \left\{ \begin{pmatrix} [0,a_1] & [0,a_4] & [0,a_5] \\ [0,a_2] & [0,a_6] & [0,a_7] \\ [0,a_3] & [0,a_8] & [0,a_4] \end{pmatrix} \right.$$

where $a_i \in [0,1]$; $1 \leq i \leq 9$} be the fuzzy set linear algebra of super fuzzy interval matrices over the set $S = \{0, \frac{1}{3n} | n \in N\}$. Fox $x, y \in V$ we have the operation to be min $\{x, y\}$ and for $x \in V$ and $s \in S$ $sx = \max (s,x)$ carried out in the following way.

Let

$$x = \left\{ \begin{pmatrix} [0,0.3] & 0 & [0,1] \\ 0 & [0,0.1] & [0,0.8] \\ [0,0.9] & [0,0.2] & [0,0.4] \end{pmatrix} \right.$$

and

$$y = \begin{pmatrix} [0,1] & [0,0.3] & 0 \\ [0,0.4] & [0,0.2] & [0,0.1] \\ [0,0.7] & [0,1] & 0 \end{pmatrix}$$

and $s = 0.3$.

Now

$$x + y = \min (x, y) = \begin{pmatrix} [0,0.3] & 0 & 0 \\ 0 & [0,0.1] & [0,0.1] \\ [0,0.7] & [0,0.2] & 0 \end{pmatrix} \qquad (1)$$

$$sx = \max sx = \begin{pmatrix} [0,0.3] & [0,0.3] & [0,1] \\ [0,0.3] & [0,0.3] & [0,0.8] \\ [0,0.9] & [0,0.3] & [0,0.4] \end{pmatrix};$$

$$sy = \max (s, y) = \begin{pmatrix} [0,1] & [0,0.3] & [0,0.3] \\ [0,0.4] & [0,0.3] & [0,0.3] \\ [0,0.7] & [0,1] & [0,0.3] \end{pmatrix}$$

$$\min (sx, sy) = \begin{pmatrix} [0,3] & [0,0.3] & [0,0.3] \\ [0,0.3] & [0,0.3] & [0,0.3] \\ [0,0.7] & [0,0.3] & [0,0.3] \end{pmatrix} \qquad (2)$$

Now

$$\max(s, x+y) = \begin{pmatrix} [0,0.3] & | & [0,0.3] & [0,0.3] \\ [0,0.3] & | & [0,0.3] & [0,0.3] \\ [0,0.7] & | & [0,0.3] & [0,0.3] \end{pmatrix} \qquad (3)$$

2 and 3 are equal. Thus V is a set linear algebra of super fuzzy interval matrices over the set S.

**Example 5.1.33**: Let

$$V = \left\{ \begin{pmatrix} [0,a_1] & [0,a_2] \\ \overline{[0,a_3]} & \overline{[0,a_4]} \\ [0,a_4] & [0,a_5] \\ [0,a_6] & [0,a_7] \\ \overline{[0,a_8]} & \overline{[0,a_9]} \\ [0,a_{10}] & [0,a_{11}] \\ [0,a_{12}] & [0,a_{13}] \\ [0,a_{14}] & [0,a_{15}] \end{pmatrix} \right.$$

where $a_i \in [0,1]$; $1 \le i \le 15$} be a set linear algebra of super fuzzy interval matrices over the set S = {[0, 1] that is all elements in [0,1]} where on V, + is min operation and for $s \in S$ and $v \in V$ $sv = \min(s, v)$.

**Example 5.1.34**: Let V = {([0,$a_1$] [0, $a_2$] [0, $a_3$] [0, $a_4$] | [0, $a_5$] [0, $a_6$] [0, $a_7$] | [0, $a_8$] [0, $a_9$]) | $a_i \in [0,1]$; $1 \le i \le 9$} be a set linear algebra of super interval fuzzy row matrices over the set S = [0, 1] = {$a_i$ | $a_i \in [0,1]$}. The operation on V is max operation and for $s \in S$ and $v \in V$ $sv = \max(sv)$.

That is if x = ([0,0.1] [0,0.2] [0, 0.3] [0,0.7] | [0, 0.1] 0 [0, 0.5] | [0, 0.1] [0, 0.9]) and y = (0 [0, 1] [0, 0.1] [0, 0.2] | [0, 0.1] [0, 1] 0 | [0, 1] 0).
Then
xy =  max (x, y)
   = ([0, 0.1] [0, 1] [0, 0.3] [0, 0.7] | [0, 0.1] [0, 1] [0, 0.5] | [0, 1] [0, 0.9]).

Now for s = 0.06 and x ∈ V;

$$0.06 \cdot x = \max(\{0.06.x\})$$
$$= ([0, 0.6]\ [0, 0.6]\ [0, 0.6]\ [0, 0.7]\ |\ [0, 0.6]\ [0, 0.6]$$
$$[0, 0.6]\ |\ [0, 0.6]\ [0, 0.9]).$$

We shall say V is a set linear fuzzy algebra with (max, max) operation on it.

**Example 5.1.35**: Let

$$V = \left\{ \begin{pmatrix} [0,a_1] & | & [0,a_5] \\ [0,a_2] & | & [0,a_6] \\ [0,a_3] & | & [0,a_7] \\ [0,a_4] & | & [0,a_8] \end{pmatrix} \right.$$

where $a_i \in [0, 1]$, $1 \leq i \leq 8$} be a set linear algebra of super fuzzy interval matrices with max operation on it over the set $S = [0, 1] = \{x_i \mid x_i \in [0, 1]\}$ with min operation on it. Thus V is a set linear algebra with (max, min) operation on it.

Thus we have four types of operations of set linear algebra of super fuzzy interval matrices, viz, {min, min}, {max, min}, {min, max} and {max, max}. We will now proceed onto give substructures in set linear algebra of super fuzzy interval matrices.

The three types of substructures are set linear subalgebra, subset linear subalgebra and pseudosubset vector subspace. The definition is a matter of routine and hence is left as an exercise to the reader. However we will illustrate these three substructures by some examples.

**Example 5.1.36**: Let

$$V = \left\{ \begin{pmatrix} [0,a_1] & [0,a_2] & | & [0,a_3] & | & [0,a_4] \\ [0,a_5] & [0,a_6] & | & [0,a_7] & | & [0,a_8] \\ \hline [0,a_9] & [0,a_{10}] & | & [0,a_{11}] & | & [0,a_{15}] \\ [0,a_{12}] & [0,a_{13}] & | & [0,a_{14}] & | & [0,a_{16}] \end{pmatrix} \right.$$

$a_i \in [0,1]$; $1 \leq i \leq 16$} be a set linear algebra of super fuzzy interval matrices over the set $S = [0, ½] = \{a_i \mid a \in [0, 1/2]\}$ under {min, min} operation on it.

Choose

$$W \left\{ \left( \begin{array}{cc|c|c} [0,a_1] & [0,a_2] & [0,a_3] & [0,a_4] \\ [0,a_8] & [0,a_5] & [0,a_6] & [0,a_7] \\ \hline [0,a_9] & [0,a_{10}] & [0,a_{13}] & [0,a_{15}] \\ [0,a_{11}] & [0,a_{12}] & [0,a_{14}] & [0,a_{16}] \end{array} \right) \right|$$

$a_i \in [0, \frac{1}{2}]$; $1 \leq i \leq 16\} \subseteq V$; W is a set linear subalgebra of super interval fuzzy matrices over the set $S = [0, \frac{1}{2}]$ under the {min, min} operation on it.

*Example 5.1.37*: Let

$$V = \left\{ \left( \begin{array}{c|cc|ccc} [0,a_1] & [0,a_3] & [0,a_5] & [0,a_7] & [0,a_9] & [0,a_{11}] \\ [0,a_2] & [0,a_4] & [0,a_6] & [0,a_8] & [0,a_{10}] & [0,a_{12}] \end{array} \right) \right|$$

$a_i \in [0,1]$; $1 \leq i \leq 12\}$ be a set linear algebra of super fuzzy interval row vector over the set $S = \{a_i | a_i \in [0,1]\}$ with (max, max) operation on it. Consider

$$W = \left\{ \left( \begin{array}{c|cc|ccc} [0,a_1] & [0,a_3] & [0,1] & [0,1] & [0,a_7] & [0,1] \\ [0,a_2] & [0,a_4] & [0,1] & [0,1] & [0,a_8] & [0,a_6] \end{array} \right) \right|$$

where $a_i \in [0,1]$; $1 \leq i \leq 8\} \subseteq V$ is a set linear subalgebra of V of super fuzzy interval row vectors over the set S.

*Example 5.1.38*: Let

$$V = \left\{ \left( \begin{array}{cc} [0,a_1] & [0,a_2] \\ [0,a_3] & [0,a_4] \\ \hline [0,a_5] & [0,a_6] \\ [0,a_7] & [0,a_8] \\ [0,a_9] & [0,a_{10}] \\ \hline [0,a_{11}] & [0,a_{12}] \\ [0,a_{13}] & [0,a_{14}] \\ [0,a_{15}] & [0,a_{16}] \end{array} \right) \right.$$

$a_i \in [0,1]$; $1 \le i \le 16$} be a set linear algebra of super interval fuzzy column vectors over the set $S = \{0, 1\}$ with {min, min} operation on V.

Consider

$$W = \left\{ \begin{pmatrix} \begin{array}{|cc|} [0,a_1] & 0 \\ 0 & [0,a_2] \\ \hline [0,a_3] & 0 \\ 0 & [0,a_4] \\ 0 & 0 \\ \hline [0,a_5] & [0,a_6] \\ 0 & 0 \\ 0 & 0 \end{array} \end{pmatrix} \right.$$

$a_i \in [0,1]$; $1 \le i \le 6\} \subseteq V$; is a set linear subalgebra of super fuzzy interval column vectors of V over the set S.

*Example 5.1.39*: Let

$$V = \left\{ \begin{pmatrix} \begin{array}{|c|cc|} [0,a_1] & [0,a_2] & [0,a_3] \\ \hline [0,a_4] & [0,a_5] & [0,a_6] \\ [0,a_7] & [0,a_8] & [0,a_9] \\ \hline [0,a_{10}] & [0,a_{11}] & [0,a_{11}] \end{array} \end{pmatrix} \right.$$

where $a_i \in [0,1]$; $1 \le i \le 12$} be a set interval linear algebra with {min, max} operation on the set $S = \{0, 1\}$.

Consider

$$W = \begin{pmatrix} \begin{array}{|c|cc|} [0,1] & [0,a_1] & [0,a_2] \\ \hline [0,a_3] & [0,1] & [0,a_4] \\ 0 & [0,a_5] & [0,a_1] \\ [0,a_6] & 0 & [0,1] \end{array} \end{pmatrix}$$

$a_i \in [0,1]$; $1 \le i \le 12\} \subseteq V$; W is a set interval linear subalgebra of V with {min, max} on the set $S = \{0, 1\}$.

Now we will give examples of pseudo set vector space of super interval fuzzy matrices.

*Example 5.1.40*: Let

$$V = \left\{ \begin{pmatrix} [0,a_1] & [0,a_2] & [0,a_3] \\ \hline [0,a_4] & [0,a_5] & [0,a_6] \\ [0,a_7] & [0,a_8] & [0,a_9] \end{pmatrix} \right.$$

where $a_i \in [0,1]$; $1 \leq i \leq 9$} be a set linear algebra of super square interval fuzzy matrices over the set $S = \{0, 1\}$.

Consider

$$W = \left\{ \begin{pmatrix} 0 & [0,a_1] & [0,a_2] \\ \hline [0,a_3] & 0 & 0 \\ [0,a_4] & 0 & 0 \end{pmatrix}, \begin{pmatrix} [0,a_1] & 0 & 0 \\ \hline 0 & [0,a_3] & 0 \\ 0 & 0 & [0,a_2] \end{pmatrix} \right.$$

$0 \leq a_i \leq 1$, $i = 1, 2, 3, 4\} \subseteq V$; W is easily verified to be only a pseudo set vector subspace of super fuzzy interval matrices over the set $\{0, 1\}$. We can have several such pseudo set vector subspaces of V.

*Example 5.1.41*: Let

$$V = \begin{pmatrix} [0,a_1] & [0,a_3] \\ [0,a_2] & [0,a_4] \end{pmatrix}, ([0,a_1]\ [0,a_2]\ [0,a_3] \mid [0,a_4]\ [0,a_5]),$$

$$\begin{bmatrix} [0,a_1] \\ \hline [0,a_1] \\ [0,a_3] \\ \hline [0,a_4] \\ [0,a_5] \end{bmatrix}, \begin{pmatrix} [0,a_1] & [0,a_2] & [0,a_3] \\ [0,a_4] & [0,a_5] & [0,a_6] \\ \hline [0,a_7] & [0,a_8] & [0,a_9] \end{pmatrix}$$

$a_i \in [0,1]$; $1 \leq i \leq 9$} be a set vector space of super interval fuzzy matrices over the set $S = \{0, 1\}$. It is clear that a set vector space in general is not a set linear algebra but every set linear algebra is a set vector space from the example 5.1.29.

Now we proceed onto define semigroup vector space of super interval fuzzy matrices over a semigroup.

**DEFINITION 5.1.6**: *Let V be the collection of all super interval fuzzy matrices and S a semigroup such that V is a semigroup vector space over the semigroup S. We call V a semigroup fuzzy set vector space of super interval fuzzy matrices over the semigroup S.*

  i)     $sv \in V$ for all $v \in V$ and $s \in S$
  ii)    $0 \cdot V = 0 \in V$ for all $v \in V$ and $0 \in S$.
  iii)    $(s_1 + s_2) v = s_1 v + s_2 v$ for all $v \in V$ and $s_1, s_2 \in S$.

We will illustrate this situation by some examples.

***Example 5.1.42***: Let

$$V = \left\{ \begin{bmatrix} [0,a_1] \\ \overline{[0,a_2]} \\ [0,a_3] \\ \overline{[0,a_4]} \\ [0,a_5] \\ [0,a_6] \\ \overline{[0,a_7]} \end{bmatrix}, \right.$$

$([0,a_1]\ [0,a_2]\ [0,a_3]\ |\ [0,a_4]\ |\ [0,a_5]\ [0,a_6])$, $\left. \begin{pmatrix} [0,a_1] & [0,a_3] & [0,a_4] \\ \hline [0,a_4] & [0,a_5] & [0,a_6] \end{pmatrix} \right|$

$a_i \in [0,1]$; $1 \leq i \leq 7\}$ be the semigroup of vector space of super interval fuzzy matrices over the semigroup $S = \{a_i \mid a_i \in [0,1]$, with min operation on it$\}$.

***Example 5.1.43***: Let

$$V = \left\{ \begin{pmatrix} [0,a_1] & [0,a_2] & [0,a_3] \\ \hline [0,a_3] & [0,a_5] & [0,a_6] \\ [0,a_4] & [0,a_7] & [0,a_9] \end{pmatrix}, \right.$$

$$\begin{pmatrix} [0,a_1] & [0,a_5] & [0,a_9] & [0,a_{13}] & [0,a_{17}] & [0,a_{21}] \\ [0,a_2] & [0,a_6] & [0,a_{10}] & [0,a_{14}] & [0,a_{18}] & [0,a_{22}] \\ [0,a_3] & [0,a_7] & [0,a_{11}] & [0,a_{15}] & [0,a_{19}] & [0,a_{23}] \\ [0,a_4] & [0,a_8] & [0,a_{12}] & [0,a_{16}] & [0,a_{20}] & [0,a_{24}] \end{pmatrix},$$

$$\begin{bmatrix} [0,a_1] & [0,a_2] \\ [0,a_3] & [0,a_4] \\ \hline [0,a_5] & [0,a_6] \\ [0,a_7] & [0,a_8] \\ [0,a_9] & [0,a_{10}] \\ \hline [0,a_{11}] & [0,a_{12}] \end{bmatrix}$$

.$a_i \in [0,1]$; $1 \leq i \leq 24\}$ be a semigroup fuzzy vector space over the semigroup $S = \{a_i \in [0,1]$ with min operation on it$\}$ be the collection of all super interval fuzzy matrices.

***Example 5.1.44***: Let $V = \{([0,a_1]\ [0,a_2]\ [0,a_3]\ |\ [0,a_4]\ [0,a_5])$, $([0,a_1]\ [0,a_2]\ [0,a_2]\ [0,a_3]\ |\ [0,a_4]\ [0,a_5])$, $([0,a_1]\ [0,a_2]\ [0,a_3]\ [0,a_4]\ |\ (0,a_5])$, $([0,a_1]\ |\ [0,a_2]\ |\ [0,a_3]\ |\ [0,a_4]\ |\ [0,a_5])$, $([0,a_1]\ |\ [0,a_2]\ [0,a_3]\ [0,a_4]\ |\ [0,a_5])$, $([0,a_1])\ |\ [0,a_2]\ [0,a_3]\ [0,a_4]\ [0,a_5])\ |\ a_i \in [0,1]$; $1 \leq i \leq 5\}$ be the semigroup vector space of super interval fuzzy matrices with entries from [0,1] over the semigroup $S = \{a_i\ |\ a_i \in [0,1]\}$ with max operation on it.

Now having seen examples of semigroup vector spaces of super fuzzy interval matrices are now proceed onto define semigroup vector subspaces of super fuzzy interval matrices.

**DEFINITION 5.1.7**: *Let V be a semigroup vector space of super interval matrices over the fuzzy semigroup S. Let $W \subseteq V$ (W a proper subset of V). If W itself is a semigroup vector space of super interval matrices over the semigroup S then we define W to be the semigroup vector subspace of V over the semigroup S.*

We will illustrate this situation by an example.

**Example 5.1.45**: Let V = = {([0,a$_1$] [0,a$_2$] | [0,a$_3$] ... [0,a$_8$] | [0,a$_9$]),

$$\left( \begin{matrix} [0,a_1] & [0,a_3] \\ [0,a_2] & [0,a_4] \end{matrix} \right), \begin{bmatrix} [0,a_1] \\ [0,a_2] \\ \overline{[0,a_3]} \\ [0,a_4] \\ \overline{[0,a_5]} \\ [0,a_6] \end{bmatrix}$$

a$_i$ ∈ [0,1]; 1 ≤ i ≤ 9} be a semigroup vector space of super interval matrices over the semigroup S = {a$_i$| a$_i$ ∈ [0,1]}. Consider W = {([0,a$_1$] [0,a$_2$]| [0,a$_3$] ... [0,a$_8$] | [0,a$_9$]),

$$\begin{bmatrix} [0,a_1] \\ [0,a_2] \\ \overline{[0,a_3]} \\ [0,a_4] \\ \overline{[0,a_5]} \\ [0,a_6] \end{bmatrix}$$

a$_i$ ∈ [0,1]; 1 ≤ i ≤ 9} ⊆ V, W is a semigroup vector subspace of V over the semigroup S.

**Example 5.1.46**: Let

$$V = \left\{ \begin{bmatrix} [0,a_1] & [0,a_2] \\ [0,a_3] & [0,a_4] \\ \overline{[0,a_5]} & \overline{[0,a_6]} \\ [0,a_7] & [0,a_8] \\ \overline{[0,a_9]} & \overline{[0,a_{10}]} \\ [0,a_{11}] & [0,a_{12}] \\ [0,a_{13}] & [0,a_{14}] \end{bmatrix} \right. ,$$

$$\left( \begin{bmatrix} [0,a_1] & [0,a_3] & [0,a_5] \\ [0,a_2] & [0,a_4] & [0,a_6] \end{bmatrix} \middle| \begin{matrix} [0,a_4] & [0,a_9] \\ [0,a_8] & [0,a_{10}] \end{matrix} \middle| \begin{matrix} [0,a_{11}] & [0,a_{13}] \\ [0,a_{12}] & [0,a_{14}] \end{matrix} \right)$$

$a_i \in [0,1]$; $1 \leq i \leq 14$} be a semigroup vector space of super interval fuzzy matrices over the semigroup $S = \{a_i \mid a_i \in [0,1]\}$ under minimum operation. Consider

$$W = \left\{ \begin{bmatrix} [0,a_1] & [0,a_2] \\ [0,a_3] & [0,a_4] \\ \hline [0,a_5] & [0,a_6] \\ [0,a_7] & [0,a_8] \\ [0,a_9] & [0,a_{10}] \\ \hline [0,a_{11}] & [0,a_{12}] \\ [0,a_{13}] & [0,a_{14}] \end{bmatrix} \middle| a_i \in [0,1]; 1 \leq i \leq 14 \right\} \subseteq V,$$

W is a semigroup vector subspace of V over the semigroup S.

*Example 5.1.47*: Let V = {(([0,a_1] [0,a_2] | [0,a_3] … [0,a_{10}] | [0,a_{11}] [0,a_{12}] [0,a_{13}]), ([0,a_1] [0,a_2] [0,a_3] | [0,a_4] … [0,a_{11}] | [0,a_{12}] [0,a_{13}]), [0,a_1] | [0,a_2] … [0,a_{12}] | [0,a_{13}]) | $a_i \in [0,1]$; $1 \leq i \leq 13$} be a semigroup vector space of super interval fuzzy row matrices over the semigroup $S = \{a_i \in [0,1]$; under max operation}. Consider W = {(([0,a_1] [0,a_2] | [0,a_3] … [0,a_{10}] | [0,a_{11}] [0,a_{12}] [0,a_{13}]), ([0,a_1] [0,a_2] … [0,a_{12}] | [0,a_{13}]) | $a_i \in [0,1]$; $1 \leq i \leq 13$} $\subseteq$ V, W is a semigroup vector subspace of V over the semigroup S.

Now having seen semigroup vector subspaces we now proceed onto define subsemigroup vector subspace of a semigroup vector over S.

**DEFINITION 5.1.8**: *Let V be a semigroup vector space of super fuzzy interval matrices over the fuzzy semigroup S. Let W $\subseteq$ V be a proper subset of V and H a proper fuzzy subsemigroup of S. Suppose W is a semigroup vector space over the semigroup H then we define W to be a subsemigroup super fuzzy subvector space of V over the subsemigroup H of S. If V has no subsemigroup subvector space then we call V to be a pseudo simple super fuzzy semigroup vector subspace of V.*

We will illustrate this situation by some simple examples.

*Example 5.1.48*: Let

$$V = \left\{ \begin{bmatrix} [0,a_1] & [0,a_2] \\ \hline [0,a_3] & [0,a_4] \\ [0,a_5] & [0,a_6] \\ \hline [0,a_7] & [0,a_8] \\ \hline [0,a_9] & [0,a_{10}] \\ [0,a_{11}] & [0,a_{12}] \\ [0,a_{13}] & [0,a_{14}] \\ [0,a_{15}] & [0,a_{16}] \end{bmatrix}, ([0,a_1]\ [0,a_2]\ [0,a_3]\ |\ [0,a_4]\ [0,a_5]\ |\ [0,a_6])\ |\right.$$

$a_i \in [0,1]$; $1 \leq i \leq 16$} be a semigroup vector space of super fuzzy interval matrices over the semigroup $S = \{a_i\ |\ a_i \in [0,1]$ under min operation$\}$. Let $W = \{([0,a_1]\ [0,a_2]\ [0,a_3]\ |\ [0,a_4]\ [0,a_5]\ |\ [0,a_6])\ |\ a_i \in [0,1]\} \subseteq V$ be a proper subset of V and $H = \{a_i\ |\ a_i \in [0, \frac{1}{2}]\} \subseteq S$ be a fuzzy subsemigroup of S under min operation. W is a subsemigroup vector subspace of V over the subsemigroup H of S.

*Example 5.1.49*: Let

$$V = \left\{ \begin{bmatrix} [0,a_1] \\ \hline [0,a_2] \\ \hline [0,a_3] \\ [0,a_4] \\ \hline [0,a_5] \\ [0,a_6] \\ [0,a_7] \end{bmatrix}, ([0,a_1]\ [0,a_2]\ |\ [0,a_3]\ [0,a_4]\ [0,a_5]\ |\ [0,a_6])\ |\ a_i \in [0,1]; \right.$$

$1 \leq i \leq 7$} be a semigroup vector space of super fuzzy interval matrices over the fuzzy semigroup $S = \{a_i \in [0,1]\ |$ min operation$\}$.
Consider

$$W = \left\{ \begin{bmatrix} [0,a_1] \\ \overline{[0,a_2]} \\ [0,a_3] \\ [0,a_4] \\ \overline{[0,a_5]} \\ [0,a_6] \\ [0,a_7] \end{bmatrix} \middle| a_i \in [0,1]; 1 \leq i \leq 7 \right\} \subseteq V;$$

and H = {$a_i$ | $a_i \in$ [0, 1/3] | min operation} $\subseteq$ S. W is a subsemigroup vector subspace of V over the subsemigroup H of S.

*Example 5.1.50*: Let V = {([0,$a_1$] | [0,$a_2$] [0,$a_3$] ... [0,$a_8$] | [0,$a_9$] [0,$a_{10}$]) | $a_i \in$ [0,1]; 1 $\leq$ i $\leq$ 10} be a semigroup vector space over the semigroup S = {0,1} under product. Clearly V is a pseudo simple semigroup vector space over the semigroup S.

Now we can define semigroup linear algebra of super fuzzy interval matrices over the semigroup S.

**DEFINITION 5.1.9**: *Let V be a semigroup vector space of super interval fuzzy matrices over the semigroup S. If V itself is a semigroup closed under binary operation under and has the additive identity then we define V to be a semigroup linear algebra of super fuzzy interval matrices over the semigroup S.*

We will illustrate this by some examples.

*Example 5.1.51*: Let

$$V = \left\{ \begin{pmatrix} [0,a_1] & [0,a_2] & [0,a_3] \\ \overline{[0,a_4]} & [0,a_7] & [0,a_8] \\ [0,a_5] & [0,a_9] & [0,a_{10}] \\ [0,a_6] & [0,a_{11}] & [0,a_{12}] \end{pmatrix} \right\},$$

where $a_i \in$ [0,1]; 1 $\leq$ i $\leq$ 12} be a semigroup linear algebra of super interval fuzzy matrices over the semigroup [0,1] under min operation.

*Example 5.1.52*: Let

$$P = \left\{ \left( \begin{array}{c|cc} [0,a_1] & [0,a_2] & [0,a_3] \\ \hline [0,a_4] & [0,a_6] & [0,a_7] \\ [0,a_5] & [0,a_8] & [0,a_9] \end{array} \right) \right.$$

$a_i \in [0,1]$ where $1 \leq i \leq 9$} be a semigroup linear algebra of super interval matrices over the semigroup $S = \{0,1\}$ under min operation.

Now the definition of substructures can be taken as a matter of routine by the reader. We however supply some examples.

*Example 5.1.53*: Let

$$V = \left\{ \left( \begin{array}{cc} [0,a_1] & [0,a_2] \\ [0,a_2] & [0,a_3] \\ \hline [0,a_4] & [0,a_5] \\ [0,a_6] & [0,a_7] \\ [0,a_8] & [0,a_9] \\ \hline [0,a_{10}] & [0,a_{11}] \\ [0,a_{12}] & [0,a_{13}] \end{array} \right) \right.$$

where $a_i \in [0,1]$; $1 \leq i \leq 14$} be the semigroup linear algebra over the semigroup $S = \{0,1\}$ under product of super fuzzy interval column vectors. Choose

$$W = \left\{ \left( \begin{array}{cc} [0,a_1] & [0,a_2] \\ [0,a_3] & [0,a_4] \\ \hline [0,a_5] & [0,a_6] \\ [0,a_7] & [0,a_8] \\ [0,a_9] & [0,a_{10}] \\ [0,a_{11}] & [0,a_{12}] \\ [0,a_{13}] & [0,a_{14}] \end{array} \right) \right.$$

$a_i \in [0,1/2]$; $1 \leq i \leq 14$} $\subseteq V$ is the semigroup linear subalgebra of V over the semigroup S.

*Example 5.1.54*: Let

$$V = \left\{ \begin{pmatrix} [0,a_1] & [0,a_4] & [0,a_2] & [0,a_{10}] & [0,a_{13}] & [0,a_{17}] \\ [0,a_2] & [0,a_5] & [0,a_8] & [0,a_{11}] & [0,a_{14}] & [0,a_{18}] \\ [0,a_3] & [0,a_6] & [0,a_9] & [0,a_{12}] & [0,a_{15}] & [0,a_{16}] \end{pmatrix} \right.$$

$a_i \in [0,1]$; $1 \leq i \leq 16$} be a semigroup linear algebra of super fuzzy row interval vectors over the semigroup $S = \{0, 1\}$ under multiplication. Consider

$$W = \left\{ \begin{pmatrix} [0,a_1] & 0 & [0,a_4] & 0 & [0,a_7] & [0,a_{10}] \\ [0,a_2] & 0 & [0,a_5] & 0 & [0,a_8] & [0,a_{11}] \\ [0,a_3] & 0 & [0,a_6] & 0 & [0,a_9] & [0,a_{12}] \end{pmatrix} \right.$$

$a_i \in [0,1]$; $1 \leq i \leq 12\} \subseteq V$; W is a semigroup linear subalgebra of the super row fuzzy interval vector over the semigroup $S = \{0, 1\}$.

Now we will proceed onto describe subsemigroup linear subalgebra of a semigroup linear subalgebra of a semigroup linear algebra of super interval fuzzy matrices over the semigroup S.

*Example 5.1.55*: Let

$$V = \left\{ \begin{pmatrix} [0,a_1] \\ [0,a_2] \\ [0,a_3] \\ [0,a_4] \\ [0,a_5] \\ [0,a_6] \\ [0,a_7] \\ [0,a_8] \\ [0,a_9] \\ [0,a_{10}] \\ [0,a_{11}] \\ [0,a_{12}] \\ [0,a_{13}] \end{pmatrix} \right.$$

where $a_i \in [0,1]; 1 \le i \le 13$} be a semigroup super linear algebra of super fuzzy interval column matrix under min operation over the semigroup $S = \{a_i \mid a_i \in [0, 1]\}$ under min operation.

Consider

$$W = \left\{ \begin{pmatrix} \overline{[0,a_1]} \\ 0 \\ \overline{[0,a_2]} \\ 0 \\ \overline{[0,a_3]} \\ 0 \\ \overline{[0,a_4]} \\ 0 \\ \overline{[0,a_5]} \\ 0 \\ \overline{[0,a_6]} \\ 0 \\ [0,a_7] \end{pmatrix} \right.$$

$a_i \in [0,1]; 1 \le i \le 7\} \subseteq V$ is a subsemigroup super linear subalgebra of V over the subsemigroup $P = \{a_i \mid a_i \in [0, ½] \subseteq S$ under min operation.

If a semigroup linear algebra of super interval fuzzy matrices does not contain any proper subsemigroup linear subalgebra then we define V to be a pseudo simple semigroup super linear algebra. We will illustrate this situation by a example.

***Example 5.1.56***: Let

$$V = \left\{ \begin{bmatrix} [0,a_1] & [0,a_2] \\ [0,a_3] & [0,a_4] \\ [0,a_5] & [0,a_6] \\ [0,a_7] & [0,a_8] \end{bmatrix} \right.$$

where $a \in [0,1]; 1 \le i \le 8\}$ be a semigroup linear algebra of super interval fuzzy matrices over the semigroup $S = \{0, 1\}$ under multiplication. V is a pseudo simple semigroup linear algebra.

Now as in case of semigroup linear algebra of super fuzzy interval matrices we can define semigroup linear transformation and semigroup linear operation. This task is also left as an exercise to the reader. Other related results can also be derived as in case of other semigroup linear algebras. The notion of direct sum and pseudo direct sum is also direct and the reader is expected to study these structures.

*Example 5.1.57*: Let

$$V = \left\{ \begin{bmatrix} [0,a_1] \\ [0,a_2] \\ [0,a_3] \\ [0,a_4] \\ [0,a_5] \\ [0,a_6] \\ [0,a_7] \end{bmatrix} \middle| \, a \in [0,1]; \, 1 \leq i \leq 7 \right\}$$

be a group vector space over the fuzzy group G. That is for every (1) g ∈ G and v ∈ V g, v and vg ∈ V. (2) for $v_1, v_2 \in V$ and g ∈ G $g(v_1 + v_2) = gv_1 + gv_2$, G can be a fuzzy group under min or max operation.

*Example 5.1.58*: Let

$$V = \left\{ \begin{pmatrix} [0,a_1] & [0,a_3] & [0,a_5] \\ [0,a_2] & [0,a_4] & [0,a_6] \end{pmatrix}, \begin{bmatrix} [0,a_1] \\ [0,a_2] \\ [0,a_3] \\ [0,a_4] \\ [0,a_5] \\ [0,a_6] \\ [0,a_7] \end{bmatrix}, \right.$$

$$\left. (([0,a_1] \, [0,a_2] \, [0,a_3]), \begin{pmatrix} [0,a_1] & [0,a_4] & [0,a_7] \\ [0,a_2] & [0,a_5] & [0,a_8] \\ [0,a_3] & [0,a_6] & [0,a_9] \end{pmatrix} \right\}$$

a ∈ [0,1]; 1 ≤ i ≤ 9} be a group vector space over the fuzzy group G.

We can define the concepts of group vector subspaces and subgroup vector subspaces which is a matter of routine we can also derive all interesting properties in case of group vector space of super fuzzy interval matrices over a fuzzy group. Also the concept of group linear algebra can be studied. Notions of group linear transformation and group linear operator can be defined and described. However problems are given for the interested reader in the last chapter of this book.

## 5.2 Special Fuzzy Linear Algebras Using Super Fuzzy Interval matrices

In this section we for the first time we introduce the notion of special fuzzy linear algebra of super interval matrices.

**DEFINITION 5.2.1**: *Let V be a set vector space of super interval matrices A with entries from $Z^+ \cup \{0\}$ or $R^+ \cup \{0\}$ or $Q^+ \cup \{0\}$ or $Z_n$ over a set $S \subseteq Z^+ \cup \{0\}$ or $R^+ \cup \{0\}$ or $Q^+ \cup \{0\}$ or $Z_n$.*

*Let $\eta : V \to [0,1]$ such that $\eta(A) \mapsto$ (to be a super fuzzy interval matrix) and $\eta(s) = t \in [0,1]$ for all $s \in S$. Then we define $V_\eta$ or $(V, \eta)$ to be the special set fuzzy vector space of super fuzzy interval matrices or fuzzy special set vector space of super fuzzy interval matrices.*

We will illustrate this situation by some examples.

***Example 5.2.1***: Let V be a set super vector space over the set S, where

$$V = \left\{ \begin{bmatrix} [0,a_1] & [0,a_4] \\ [0,a_2] & [0,a_5] \\ [0,a_3] & [0,a_6] \end{bmatrix} \right\},$$

$([0,a_1]\ [0,a_2]\ [0,a_3]\ |\ [0,a_4]\ [0,a_5]) \mid a_i \in Z^+ \cup \{0\};\ 1 \le i \le 6\}$ over the set $S = \{0,1\}$.

Define $\eta: V \to [0,1]$ and $\eta : S \to [0,1]$ as follows.

$$\eta \begin{pmatrix} \begin{array}{c|c} [0,a_1] & [0,a_4] \\ \hline [0,a_2] & [0,a_5] \\ [0,a_3] & [0,a_6] \end{array} \end{pmatrix} = \begin{bmatrix} \begin{array}{c|c} \left[0,\dfrac{1}{a_1}\right] & \left[0,\dfrac{1}{a_4}\right] \\ \hline \left[0,\dfrac{1}{a_3}\right] & \left[0,\dfrac{1}{a_5}\right] \\ \left[0,\dfrac{1}{a_3}\right] & \left[0,\dfrac{1}{a_6}\right] \end{array} \end{bmatrix}$$

if $a_i \neq 0$ then $\eta(a_i) = 1/a_i$ if $a_i = 0$ then $\eta(a_i)=1$, $\eta(0) = 0$ and $\eta(1) = 1$. It is easily verified $(V, \eta)$ is a special fuzzy set super linear algebra provided $\eta$ defined so is such that $\eta(\alpha v) = \eta(\alpha) \eta(v)$ for all $\alpha \in S$ and $v \in V$.

*Example 5.2.2*: Let

$$V = \left\{ \begin{pmatrix} \begin{array}{cc|c} [0,a_1] & [0,a_2] & [0,a_3] \\ [0,a_4] & [0,a_5] & [0,a_6] \\ \hline [0,a_7] & [0,a_8] & [0,a_9] \\ [0,a_{10}] & [0,a_{11}] & [0,a_{12}] \end{array} \end{pmatrix} \right.$$

where $a_i \in Z^+ \cup \{0\}$; $1 \leq i \leq 12\}$ be a set super linear algebra defined over the set $S = \{Z^+ \cup \{0\}\}$. Define $\eta: V \to [0,1]$ and $\eta : S \to [0,1]$

$$\eta(A) = \begin{bmatrix} \begin{array}{cc|c} \left[0,\dfrac{1}{a_1}\right] & \left[0,\dfrac{1}{a_2}\right] & \left[0,\dfrac{1}{a_3}\right] \\ \left[0,\dfrac{1}{a_4}\right] & \left[0,\dfrac{1}{a_5}\right] & \left[0,\dfrac{1}{a_6}\right] \\ \hline \left[0,\dfrac{1}{a_8}\right] & \left[0,\dfrac{1}{a_8}\right] & \left[0,\dfrac{1}{a_9}\right] \\ \left[0,\dfrac{1}{a_9}\right] & \left[0,\dfrac{1}{a_n}\right] & \left[0,\dfrac{1}{a_{12}}\right] \end{array} \end{bmatrix}$$

where

$$\eta(a_i) = \begin{cases} \dfrac{1}{a_i} & \text{if } a_i \neq 0 \\ 1 & \text{if } a_i = 0 \end{cases}$$

It is easily verified $\eta(V)$ is a special set fuzzy super linear algebra under the operation min or max.

*Example 5.2.3*: Let

$$V = \left\{ \begin{bmatrix} [0,a_1] \\ \overline{[0,a_2]} \\ [0,a_3] \\ \overline{[0,a_4]} \\ [0,a_5] \end{bmatrix} \right\},$$

$([0,a_1]\ [0,a_2]\ |\ [0,a_3]\ [0,a_4]\ [0,a_5])\ |\ a_i \in Z^+ \cup \{0\}\}$ be a set vector space over the set $S = Z^+ \cup \{0\}$. Define

$$\eta(x) = \begin{cases} \dfrac{1}{x} & \text{if } x \neq 0 \\ 1 & \text{if } x = 0 \end{cases} \quad \text{for } x \in Z^+$$

$\eta(V)$ becomes a special set fuzzy super vector space over the fuzzy set $\eta(S)$.

Likewise if V is semigroup vector space of super interval matrices with entries from the semigroup S, then we define $\eta(V)$ or $(V, \eta)$ as the special map on V and S such that $\eta(V)$ is the set of all super interval fuzzy matrices and $\eta(S)$ is the special fuzzy semigroup, that is it contains elements from [0,1] with min or max or product operation on it. We define $\eta(V)$ or $(V, \eta)$ as the special fuzzy semigroup of super vector space over the special fuzzy semigroup S.

We will describe this situation by some examples.

*Example 5.2.4*: Let

$$V = \left\{ \begin{bmatrix} [0,a_1] & [0,a_2] \\ [0,a_3] & [0,a_4] \\ [0,a_5] & [0,a_6] \\ [0,a_7] & [0,a_8] \\ [0,a_9] & [0,a_{10}] \end{bmatrix}, \right.$$

$([0,a_1]\ [0,a_2]\ |\ [0,a_3]\ [0,a_4]\ [0,a_5])\ |\ [0,a_6]\ |\ [0,a_7]\ |\ a_i \in Z^+ \cup \{0\}\}$ be a semigroup vector space of super interval matrices over the semigroup S. Let

$$\eta(a) = \begin{cases} \dfrac{1}{a} & \text{for all } a \in Z^+ \\ 1 & \text{for } 0 \end{cases}.$$

With this definition of $\eta$; $\eta(V)$ is the collection of super interval fuzzy matrices with entries from the fuzzy semigroup $\eta(S)$. Thus $\eta(V)$ is the special fuzzy interval matrix semigroup vector space over the fuzzy semigroup $\eta(S)$.

*Example 5.2.5*: Let

$$V = \left\{ \begin{bmatrix} [0,a_1] & [0,a_2] & [0,a_3] \\ [0,a_4] & [0,a_5] & [0,a_6] \\ [0,a_7] & [0,a_{10}] & [0,a_{11}] \\ [0,a_8] & [0,a_{13}] & [0,a_{12}] \\ [0,a_9] & [0,a_{14}] & [0,a_{15}] \end{bmatrix}, \right.$$

$$\begin{bmatrix} [0,a_1] & [0,a_5] & [0,a_6] & [0,a_{13}] & [0,a_{17}] & [0,a_{21}] \\ [0,a_2] & [0,a_7] & [0,a_8] & [0,a_{14}] & [0,a_{18}] & [0,a_{22}] \\ [0,a_3] & [0,a_9] & [0,a_{10}] & [0,a_{15}] & [0,a_{19}] & [0,a_{23}] \\ [0,a_4] & [0,a_{11}] & [0,a_{12}] & [0,a_{16}] & [0,a_{20}] & [0,a_{24}] \end{bmatrix},$$

$$\left\{ \begin{bmatrix} [0,a_1] \\ [0,a_2] \\ \hline [0,a_3] \\ [0,a_4] \\ \hline [0,a_5] \\ [0,a_6] \\ [0,a_7] \\ \hline [0,a_8] \\ [0,a_9] \\ [0,a_{10}] \\ [0,a_{11}] \end{bmatrix} \middle| \; a_i \in Z^+ \cup \{0\}; \; 1 \leq i \leq 24 \right\}$$

a semigroup vector space of super interval matrices over the semigroup $S = Q^+ \cup \{0\}$. Define

$$\eta(x) = \begin{cases} \dfrac{1}{3} & \text{if } x \in Z^+ \\ \dfrac{1}{4} & \text{if } x \in Q^+ \setminus Z^+ \\ 1 & \text{if } x = 0 \end{cases}.$$

Clearly $\eta(V)$ is the set of all special fuzzy interval super matrices or super interval matrices with entries from [0,1]. Further $\eta(S)$ is a special fuzzy semigorup. Thus $\{\eta(V), \eta(S)\}$ is the special fuzzy semigroup super matrix interval vector space over the fuzzy semigroup $\eta(S)$. Now in a similar way we can define the notion of special fuzzy group vector space over a special fuzzy group using the notion of group vector space of super interval matrices over the group G. This is left as an exercise to the reader. Further we can define using set linear algebra of super interval matrices over the set S special fuzzy set super matrix interval linear algebra over the fuzzy set $\eta(S)$.

Definition and study of fuzzy substructure is also a matter of routine interested reader can do this job and illustrate these concepts with examples.

Now if we have η(V) and η(W) to be two special fuzzy set super matrix interval vector spaces over the same fuzzy set η(S) we can define fuzzy set linear transformation η(V) to η(W).

Similar study can be carried out in case of special fuzzy semigroup vector spaces and special fuzzy group vector spaces when the structures are defined over the same fuzzy special semigroup and fuzzy special group respectively.

When the range space coincides with the domain space in all the three cases we can get the corresponding fuzzy set linear operator super matrix interval of vector spaces / special fuzzy semigroup linear operator of super matrix interval vector spaces or special fuzzy group linear operator of super matrix interval vector spaces defined over a special fuzzy set / fuzzy special semigroup / special fuzzy group respectively.

For more about such type of usual spaces and their fuzzy analogue refer [38, 46-7, 49, 52].

Thus one can keep this is mind and define these fuzzy structures.

If $I(Z^+ \cup \{0\}) = \{[0,a] \mid a \in Z^+ \cup \{0\}\}$ = all integer intervals of this special type we can get the fuzzy analogue

$\eta : I(Z^+ \cup \{0\}) \to I([0,1])$ where $I([0,1]) = \{[0,a] \mid a \in [0,1]\}$ is the collection of all fuzzy intervals of the special form $[0,a]$; $a \in [0,1]$. $\eta : I(Q^+ \cup \{0\}) \to I([0,1])$ or $\eta : I(R^+ \cup \{0\}) \to I([0,1]$ gives fuzzy intervals and we can define min or max or product and make them special fuzzy semigroups.

We will illustrate this situation by one or two examples.

**Example 5.2.6**: Let

$$V = \left\{ \begin{bmatrix} [0,a_1] & [0,a_2] \\ \hline [0,a_3] & [0,a_4] \\ [0,a_5] & [0,a_5] \end{bmatrix} \right.$$

where $a_i \in Z_{12}$; $1 \leq i \leq 6$} be a group linear algebra over the group $G = Z_{12}$ under addition modulo 12.

Define for every $A \in V$, $\eta(A)$ a fuzzy interval super matrix where $\eta(A)$ has its entries from $I([0,1])$ such that

$$\eta(A+B) \geq \min\{\eta(A), \eta(B)\}$$

(that is for every $\eta([0,a_i], [0,b_i]) \geq \min\{\eta([0,a_i]), \eta([0,b_i]))$ where $[0,a_i]$ is an entry in A and $[0,b_i]$ is an entry in B) to be more exact.

$$\eta\left(\begin{bmatrix} [0,a_1] & [0,a_2] \\ \hline [0,a_3] & [0,a_4] \\ \hline [0,a_5] & [0,a_6] \end{bmatrix}, \begin{bmatrix} [0,b_1] & [0,b_2] \\ \hline [0,b_3] & [0,b_4] \\ \hline [0,b_5] & [0,b_6] \end{bmatrix}\right) \geq$$

$$\begin{bmatrix} [0,\min(\eta([0,a_1],\eta[0,b_1]))] & [0,\min(\eta([0,a_2],\eta[0,b_2]))] \\ \hline [0,\min(\eta([0,a_3],\eta[0,b_3]))] & [0,\min(\eta([0,a_4],\eta[0,b_4]))] \\ \hline [0,\min(\eta([0,a_5],\eta[0,b_5]))] & [0,\min(\eta([0,a_6],\eta[0,b_6]))] \end{bmatrix}$$

$$\eta\left(\begin{bmatrix} [0,a_1] & [0,a_2] \\ \hline [0,a_3] & [0,a_4] \\ \hline [0,a_5] & [0,a_5] \end{bmatrix}\right) > \eta\left(\begin{bmatrix} [0,-a_1] & [0,-a_2] \\ \hline [0,-a_3] & [0,-a_4] \\ \hline [0,-a_5] & [0,-a_5] \end{bmatrix}\right)$$

and

$$\eta\begin{bmatrix} 0 & 0 \\ \hline 0 & 0 \\ \hline 0 & 0 \end{bmatrix} = \begin{bmatrix} 1 & 1 \\ \hline 1 & 1 \\ \hline 1 & 1 \end{bmatrix},$$

$\eta(rA) = \eta(r)\eta(A)$ where $\eta(r) = 1/r$, $r \in Z_{12} \setminus \{0\}$ $\eta(0) = 1$; for all $A, B \in V$ and $r \in G$. $V\eta$ is the special fuzzy group linear algebra.

***Example 5.2.7***: Let

$$V = \left\{ \begin{bmatrix} [0,a_1] \\ [0,a_2] \\ \overline{[0,a_3]} \\ [0,a_4] \\ \overline{[0,a_5]} \\ [0,a_6] \end{bmatrix}, \right.$$

$([0,a_1]\ [0,a_2]\ [0,a_3]\ |\ [0,a_4]\ |\ [0,a_5]\ [0,a_6])\ |\ a_i \in Z_{20}\}$ be a group vector space of super interval matrices over the group $G = Z_{20}$.

Define

$$\eta(A) = \begin{bmatrix} \left[0, \dfrac{1}{a_1}\right] \\ \overline{\left[0, \dfrac{1}{a_2}\right]} \\ \left[0, \dfrac{1}{a_3}\right] \\ \overline{\left[0, \dfrac{1}{a_4}\right]} \\ \overline{\left[0, \dfrac{1}{a_5}\right]} \\ \left[0, \dfrac{1}{a_6}\right] \end{bmatrix}$$

if $a_i \neq 0$, $1 \leq i \leq 6$ if $a_i = 0$ replace by 1.

$$\eta(B) = \left(\left[0, \dfrac{1}{a_1}\right]\left[0, \dfrac{1}{a_2}\right]\left[0, \dfrac{1}{a_3}\right]\left[0, \dfrac{1}{a_4}\right]\left[0, \dfrac{1}{a_5}\right]\left[0, \dfrac{1}{a_6}\right]\right)$$

if $a_i \neq 0$; $1 \leq i \leq 6$ if $a_i = 0$ replace by 1. $\eta(rA) > \eta(A)$ for all $r \in Z_{20}$.

$$\eta : Z_{20} \to [0,1]$$

is such that

$$\eta(a) = \begin{cases} \dfrac{1}{a} & \text{if } a \neq 0 \\ \dfrac{1}{4} & \text{if } a = 0 \end{cases}.$$

Thus $V\eta$ is a special group fuzzy vector space of fuzzy super interval matrices over the fuzzy group $\eta(Z_{20})$.

**Chapter Six**

# APPLICATION OF SUPER INTERVAL MATRICES AND SET LINEAR ALGEBRAS BUILT USING SUPER INTERVAL MATRICES

Super interval matrices are obtained by replacing the entries in super matrices by intervals of the form [0, a] where a $\in$ $Z_n$ or $Z^+ \cup \{0\}$ or $Q^+ \cup \{0\}$ or $R^+ \cup \{0\}$. Also it can equally be realized as usual interval matrices which are partitioned in different ways to get different types of super interval matrices. Thus given any m × n super interval matrix M by partitioning M we can get several m × n super interval matrices. Also for each type of partition we have a collection and that collection can have nice algebraic structure. So in the place of one algebraic structure we get several of them with same natural order which is not possible in usual interval matrices. Hence in applying them to fields like mathematical / fuzzy modeling and in finite element methods using intervals one can get better results saving economy and time. Further when in super interval matrices if we replace the real intervals by fuzzy intervals we get super fuzzy interval matrices which can find its application in all problems where fuzzy interval solution is expected.

Certainly these new structures will be a boon to finite interval analysis were usual interval matrices can be replaced by super interval matrices and the results can be accurate time saving and the expert can experiment with different partitions of the interval matrix to obtain the expected or the best solution.

**Chapter Seven**

# SUGGESTED PROBLEMS

In this chapter we suggest around 110 problems. Some problems are simple, some difficult and some can be treated as a research problem. These problems will certainly help the reader in understanding the concepts mentioned in the book.

1. Obtain some interesting properties enjoyed by super interval column vectors.

2. How many types of super interval row vectors can be got from $X = \{([0,0]\ [0,a_2]\ [0,a_3]\ \ldots\ [0,a_{10}]) \mid a_i \in Z_{12}\}$?

3. Prove the types of super interval column vector got from
$$Y = \left\{ \begin{bmatrix} [0,a_1] \\ [0,a_2] \\ [0,a_3] \\ \vdots \\ [0,a_{10}] \end{bmatrix} \middle| a_i \in Z_{12};\ 1 \leq i \leq 10 \right\} \text{ and } X \text{ in (2) are isomorphic.}$$

4. Find the number of super interval matrix got from the interval matrix

$$A = \begin{bmatrix} [0,a_1] & [0,a_2] & [0,a_3] & [0,a_4] & [0,a_5] \\ [0,a_6] & [0,a_7] & [0,a_8] & [0,a_9] & [0,a_{10}] \\ [0,a_{11}] & [0,a_{12}] & [0,a_{13}] & [0,a_{14}] & [0,a_{15}] \\ [0,a_{16}] & [0,a_{17}] & [0,a_{14}] & [0,a_{20}] & [0,a_{21}] \\ [0,a_{22}] & [0,a_{23}] & [0,a_{24}] & [0,a_{25}] & [0,a_{26}] \end{bmatrix}.$$

Hence or otherwise obtain the number of super matrices obtained from any n × n interval matrix.

5. Find the order of the X; where X = {([0, $a_1$] [0,$a_2$] [0,$a_3$] | [0,$a_4$] [0,$a_5$]) | $a_i \in Z_{10}$; 1 ≤i≤5}.

   (a) Is X a group under interval row matrix addition?
   (b) Is X a group under multiplication?
   (c) Can X have S-subsemigroups if X is taken under multiplication?

6. Let $P = \left\{ \begin{bmatrix} [0,a_1] \\ [0,a_2] \\ \hline [0,a_3] \\ [0,a_4] \\ \hline [0,a_5] \\ [0,a_6] \end{bmatrix} \middle| a_i \in Z_{20}; 1 \le i \le 6 \right\}$ be the collection of all super interval column vectors of the same type.

   (a) Is P a group?
   (b) What is the order of P?
   (c) How many types of super interval column vector groups can be got using the interval column matrix
   $$\begin{bmatrix} [0,a_1] \\ [0,a_2] \\ [0,a_3] \\ [0,a_4] \\ [0,a_5] \\ [0,a_6] \end{bmatrix} ?$$

7. Let $B = \left\{ \begin{pmatrix} [0,a_1] & [0,a_5] & | & [0,a_8] & [0,a_9] & [0,a_{10}] \\ [0,a_3] & [0,a_6] & | & [0,a_{11}] & [0,a_{12}] & [0,a_{13}] \\ [0,a_2] & [0,a_{14}] & | & [0,a_{14}] & [0,a_{15}] & [0,a_{17}] \end{pmatrix} \middle| a_i \in Z_{15}; 1 \leq i \leq 15 \right\}$ be
the super interval row vector of the same type.
   (a) Prove B is a group under addition.
   (b) What is the order of B?
   (c) Can B have subgroups?
   (d) How many types of super interval row vectors can be got using the interval matrix $X = \begin{pmatrix} [0,a_1] & [0,a_2] & [0,a_3] & [0,a_4] & [0,a_5] \\ [0,a_6] & [0,a_7] & [0,a_8] & [0,a_9] & [0,a_{10}] \\ [0,a_{11}] & [0,a_{12}] & [0,a_{13}] & [0,a_{14}] & [0,a_{15}] \end{pmatrix}$
   where $a_i \in Z_{15}$; $1 \leq i \leq 15$.

8. Give an example of a Smarandache super interval row semigroup.

9. Does there exists a S-semigroup of super interval row semigroup of order 27?

10. Let $S = \{(a_1\ a_2\ a_3\ |\ a_4\ a_5\ a_6\ a_7\ |\ a_8) \mid a_i \in Z^+ \cup \{0\}\ ;\ 1 \leq i \leq 8\}$ be the collection of all super row matrices of the same type
    (a) Prove S is a semigroup under super row matrix addition
    (b) Can S have subsemigroups?
    (c) Is S a Smarandache semigroup?

11. Let $S = \{([0, a_1] \mid [0, a_2]\ [0, a_3] \mid [0, a_4]\ [0, a_5]\ [0, a_6]\ [0, a_7]) \mid a_i \in Z^+ \cup \{0\}; 1 \leq i \leq 7\}$ be a semivector space of super row interval matrices defined over the semifield $S = Z^+ \cup \{0\}$.
    (a) Find a basis of V.
    (b) Find a linearly independent subset of V with cardinality 12.
    (c) Find a semivector subspace of V.
    (d) Find a linearly dependent subset of V of order 4.
    (e) Can V be a semilinear algebra of super interval row matrices over S?

12. Let $V = \left\{ \begin{bmatrix} [0,a_1] & [0,a_2] & [0,a_3] \\ [0,a_4] & [0,a_5] & [0,a_6] \\ [0,a_7] & [0,a_8] & [0,a_7] \\ \hline [0,a_{10}] & [0,a_{11}] & [0,a_{12}] \\ [0,a_{13}] & [0,a_{14}] & [0,a_{15}] \\ [0,a_{16}] & [0,a_{17}] & [0,a_{18}] \end{bmatrix} \right.$ where $a_i \in Z^+ \cup \{0\}$; $1 \leq i \leq 7\}$ be a semivector space of super interval matrices over the semifield $T = Z^+ \cup \{0\}$ of type I.

(a) Find a basis B of V.

(b) Prove there can be linearly independent subsets of V whose cardinality exceeds the cardinality of the basis B of V.

(c) Prove V can have infinite number of semivector subspaces of type I over T.

(d) Define a linear operator on V.

(e) Can V be made into a semilinear algebra of type I over T?

13. Obtain some interesting properties enjoyed by the semivector space of super interval square matrices of type I over a semifield S.

14. Find the difference between type I and type II semivector spaces.

15. Prove in general semiring of interval super matrices are not semifields.

16. Obtain a necessary and sufficient condition for a super interval semiring to be a semifield.

17. Determine some important and interesting features enjoyed by semifield of super interval matrices.

18. Let $S = \{([0, a_1]\ [0,a_2] \mid [0,a_3] \mid [0,a_4]\ [0,a_5]\ [0,a_6]\ [0,a_7]) \mid a_i \in Z^+, 1 \leq i \leq 7\} \cup \{(0\ 0 \mid 0 \mid 0\ 0\ 0\ 0)\}$. Is S a semifield?
(Product in S is defined component wise).

19. Let $S = \left\{ \left( \begin{array}{c|c} [0,a] & [0,b] \\ \hline [0,b] & [0,c] \end{array} \right) \middle| a, b, c \in R^+ \right\} \cup \begin{pmatrix} 0 & 0 \\ 0 & 0 \end{pmatrix}$. Can every S be a semifield? Can S be a strict semiring? What is the possible algebraic structure enjoyed by S?

20. Give an example of a semivector space of super interval matrices of type I over a semifield which has exactly 4 elements in its basis B.

21. What is the advantage of using semivector spaces of super interval matrices in the place of vector space of super matrices?

22. Give some applications of semivector space of super interval matrices over the semifield S of type I.

23. Give an example of a semivector space of super interval matrices of type I which has infinite basis.

24. Does there exist a semi vector space of super interval matrices of type I with finite number of elements in it?

25. Let $V = \{([0, a_1] \mid [0,a_2] [0,a_3] \mid [0,a_4] [0,a_5] [0,a_6] [0,a_7] [0,a_8] [0,a_9]) \mid a_i \in Z^+ \cup \{0\}\ 1 \leq i \leq 9\}$ be a semivector space of super interval row matrices of type I over $S = Z^+ \cup \{0\}$ and

$$W = \left\{ \left[ \begin{array}{cccc} [0,a_1] & [0,a_2] & [0,a_3] & [0,a_4] \\ [0,a_5] & [0,a_6] & [0,a_7] & [0,a_8] \\ [0,a_9] & [0,a_{10}] & [0,a_{11}] & [0,a_{12}] \\ [0,a_{13}] & [0,a_{14}] & [0,a_{15}] & [0,a_{16}] \\ \hline [0,a_{17}] & [0,a_{18}] & [0,a_{19}] & [0,a_{20}] \\ [0,a_{21}] & [0,a_{22}] & [0,a_{23}] & [0,a_{24}] \\ [0,a_{25}] & [0,a_{26}] & [0,a_{27}] & [0,a_{28}] \\ [0,a_{29}] & [0,a_{30}] & [0,a_{31}] & [0,a_{32}] \\ [0,a_{33}] & [0,a_{34}] & [0,a_{35}] & [0,a_{36}] \end{array} \right] \middle| a_i \in Z^+ \cup \{0\}; 1 \leq i \leq 36 \right\}$$

be a semivector space of super interval column vectors of type I over the semifield $S = Z^+ \cup \{0\}$.

(a) Determine two distinct linear transformations from V to W.
(b) Find a linear transformation T : W → V so that T is invertible.
(c) Suppose F = {set of all linear transformations from V to W}; what can be the algebraic structure enjoyed by F?

26. Let V = $\left\{ \begin{bmatrix} [0,a_1] & [0,a_2] & [0,a_3] & [0,a_4] \\ [0,a_5] & [0,a_6] & [0,a_7] & [0,a_8] \\ [0,a_9] & [0,a_{10}] & [0,a_{11}] & [0,a_{12}] \\ [0,a_{13}] & [0,a_{14}] & [0,a_{15}] & [0,a_{16}] \\ [0,a_{17}] & [0,a_{18}] & [0,a_{19}] & [0,a_{20}] \\ [0,a_{21}] & [0,a_{22}] & [0,a_{23}] & [0,a_{24}] \\ [0,a_{25}] & [0,a_{26}] & [0,a_{27}] & [0,a_{28}] \\ [0,a_{29}] & [0,a_{30}] & [0,a_{31}] & [0,a_{32}] \\ [0,a_{33}] & [0,a_{34}] & [0,a_{35}] & [0,a_{36}] \end{bmatrix} \right.$ where $a_i \in Z^+ \cup \{0\}$;

$1 \le i \le 36\}$ be a semivector space of super interval matrices of type I over the semifield $S = Z^+ \cup \{0\}$.

(a) Find three different linear operators on S.
(b) What is the algebraic structure enjoyed by the set of all linear operators on V?

27. Let V = $\left\{ \begin{bmatrix} [0,a_1] \\ [0,a_2] \\ [0,a_3] \\ \overline{[0,a_4]} \\ [0,a_5] \\ [0,a_6] \\ [0,a_7] \end{bmatrix} \right.$ $a_i \in Q^+ \cup \{0\}$ $1 \le i \le 7\}$ be a semivector space

over the semifield $S = Z^+ \cup \{0\}$ of super interval column matrices of type I. Let W = $\left\{ \begin{pmatrix} [0,a_1] & [0,a_2] & [0,a_3] \\ 0 & [0,a_4] & [0,a_5] \\ [0,a_6] & 0 & [0,a_7] \end{pmatrix} \right.$ $a_i \in Z^+ \cup \{0\}$ $1 \le i$

≤ 7} be a semivector space of super interval matrices of type I over the semifield $S = Z^+ \cup \{0\}$.
(a) Find a linear transformation from V to W?
(b) What is the structure enjoyed by the collection of all linear transformations from V to W?
(c) What is the algebraic structure satisfied by the set of all linear transformations from W to V?
(d) Compare the algebraic structure in (a) and (b).
(e) Can any linear transformation from W to V be invertible? Justify your claim.

28. Let $V = \{([0, a_1] \mid [0,a_2] [0,a_3] \mid [0,a_4] [0,a_5] [0,a_6]) \mid a_i \in Q^+ \cup 1 \leq i \leq 6\}$ be a semivector space of super interval row matrices of type I over the interval semifield $S = I(Z^+ \cup \{0\})$. Find a basis of V.

29. How many semivector spaces of super interval matrices of type I can be obtained by varying the partitions of a 5 × 5 interval matrix?

30. Let $T = \left\{ \begin{bmatrix} [0,a_1] & [0,a_2] & [0,a_3] & [0,a_4] \\ [0,a_5] & [0,a_8] & [0,a_9] & [0,a_{10}] \\ [0,a_6] & [0,a_{11}] & [0,a_{12}] & [0,a_{13}] \\ [0,a_7] & [0,a_{14}] & [0,a_{15}] & [0,a_{10}] \end{bmatrix} \right.$ where $a_i \in Z^+ \cup \{0\}$; 1 ≤ i ≤ 16} be a semivector space of super interval matrices of type I over $S = Z^+ \cup \{0\}$
(a) How many partitions can be made on a 4 × 4 interval matrix?
(b) Give the total number of distinct semivector spaces of super interval 4 × 4 matrices by varying the partitions.
(c) Find a basis of T?
(d) Will the basis have same cardinality even if partitions are different on the 4 × 4 interval matrices.
(e) Find at least 3 linear operators on V.

31. Can a semivector space of super interval 5 × 6 matrices over the semifield $S = R^+ \cup \{0\}$ have subspaces of finite dimensions?

32. Is P = $\left\{\begin{bmatrix} [0,a_1] & [0,a_2] & \cdots & [0,a_{12}] \\ [0,a_{13}] & [0,a_{14}] & \cdots & [0,a_{24}] \\ \hline [0,a_{25}] & [0,a_{26}] & \cdots & [0,a_{36}] \\ \hline [0,a_{37}] & [0,a_{38}] & \cdots & [0,a_{48}] \end{bmatrix}\right\}$ where $a_i \in 3Z^+ \cup \{0\}$; $1 \le i \le 48$} a semivector space of super interval matrices of type I over the semifield $S = Z^+ \cup \{0\}$ is of finite dimension over S? Justify your answer.

33. Can there exist a semivector space of super interval matrices of dimension 23?

34. Let V = $\left\{ \begin{pmatrix} [0,a_1] & [0,a_2] & [0,a_3] \\ [0,a_4] & [0,a_5] & [0,a_6] \end{pmatrix}, \begin{bmatrix} [0,a_1] \\ [0,a_2] \\ \hline [0,a_3] \\ [0,a_4] \\ \hline [0,a_5] \\ [0,a_6] \end{bmatrix}, ([0,a_1] \mid [0,a_2]\; [0,a_3]\; [0,a_4] \mid [0,a_5]) \right.$ where $a_i \in Q^+ \cup \{0\}$; $1 \le i \le 6$} be a set vector space of super interval matrices over the set $S = \{0, 4, 5, 9, 28, 45\}$
(a) Find at least 3 set vector subspaces of V.
(b) Prove or disprove V can have infinite number of set vector subspaces of super interval matrices over S.
(c) Does there exists a finite generating subset for V?

35. Let V = $\left\{ (([0, a_1] \mid [0,a_2]\; [0,a_3] \mid [0,a_4]\; [0,a_5]\; [0,a_6]\; [0,a_7]), \begin{bmatrix} [0,a_1] \\ [0,a_2] \\ \hline [0,a_3] \\ [0,a_4] \\ \hline [0,a_5] \\ [0,a_6] \end{bmatrix}, \begin{pmatrix} [0,a_1] & [0,a_2] & \cdots & [0,a_{12}] \\ [0,a_{13}] & [0,a_{14}] & \cdots & [0,a_{24}] \end{pmatrix} \right.$ $\mid a_i \in Z^+ \cup \{0\}; 1 \le i \le 24\}$ be a set vector space of super interval matrices over the set $S = \{0, 1\}$.

(a) Find set vector subspaces of V.
(b) Can V have semivector subspace over S?
(c) Find a generating set of V over S.
(d) Prove V is an infinite set.

36. What are the special properties enjoyed by set supervector spaces of super interval matrices?

37. Give an example of a set super vector space which has only 5 elements in it.

38. Give an example of a set super vector space of super interval matrices which has no proper set super vector subspace.

39. Give an example of an infinite super set vector space of super interval matrices defined over the set $S = Z^+ \cup \{0\}$.

40. Prove

$$V = \left\{ \begin{bmatrix} [0,a_1] & [0,a_2] & [0,a_3] \\ [0,a_4] & [0,a_5] & [0,a_6] \\ [0,a_7] & [0,a_8] & [0,a_9] \\ [0,a_{10}] & [0,a_{11}] & [0,a_{12}] \end{bmatrix}, \begin{bmatrix} [0,a_1] & [0,a_2] & [0,a_3] \\ [0,a_4] & [0,a_5] & [0,a_6] \\ [0,a_7] & [0,a_8] & [0,a_9] \\ [0,a_{10}] & [0,a_{11}] & [0,a_{12}] \end{bmatrix}, \right.$$

$$\left. \begin{bmatrix} [0,a_1] & [0,a_2] & [0,a_3] \\ [0,a_4] & [0,a_5] & [0,a_6] \\ [0,a_7] & [0,a_8] & [0,a_9] \\ [0,a_{10}] & [0,a_{11}] & [0,a_{12}] \end{bmatrix}, \begin{bmatrix} [0,a_1] & [0,a_2] & [0,a_3] \\ [0,a_4] & [0,a_5] & [0,a_6] \\ [0,a_7] & [0,a_8] & [0,a_9] \\ [0,a_{10}] & [0,a_{10}] & [0,a_{12}] \end{bmatrix} \right\}$$

$a_i \in Q^+ \cup \{0\}; 1 \le i \le 12\}$ is a set vector space of super interval matrices over the set $S = Z^+ \cup \{0\}$.
(a) Can V have semivector super subspaces?
(b) Obtain atleast 3 set vector subspaces of V over S.
(c) Prove V can have infinitely many set vector subspaces.

41. Give an example of a semivector space of super interval matrices which is not a semi linear algebra of type I.

42. Give an example of a semilinear algebra of super interval matrices of type I which has a finite basis.

43. Let $V = \left\{ \begin{pmatrix} [0,a_1] & [0,a_2] & [0,a_5] \\ [0,a_3] & [0,a_4] & [0,a_6] \end{pmatrix}, ([0,a_1] \mid [0,a_2] \; [0,a_3] \; [0,a_4] \mid [0,a_5] \; [0,a_6])\right.$, $\left. \begin{bmatrix} [0,a_1] \\ \overline{[0,a_2]} \\ [0,a_3] \\ [0,a_4] \\ [0,a_5] \\ [0,a_6] \end{bmatrix} \right\}$ $a_i \in R^+ \cup \{0\}; 1 \leq i \leq 6$ be a set super vector space of super interval matrices over the set $S = Z^+ \cup \{0\}$.

   (a) Find set super vector subspaces of V.
   (b) Can V have set super semilinear subalgebra?
   (c) Can V have a finite generating subset?

44. Obtain some interesting properties enjoyed by set super vector spaces defined over a set S.

45. What are the advantages of using set super vector space of super interval matrices over usual set vector spaces?

46. Prove the set $V = \left\{ \begin{bmatrix} [0,a_1] \\ \overline{[0,a_2]} \\ [0,a_3] \\ [0,a_4] \\ [0,a_5] \\ [0,a_6] \end{bmatrix} \right.$ $\left. a_i \in Z^+ \cup \{0\}; 1 \leq i \leq 6 \right\}$ is a semigroup under addition. Can V be a S-semigroup?

47. Give an example of a semigroup of all super interval matrices (with respect to a fixed partition) which is S-semigroup.

48. Does there exist a set of super interval matrices (with respect a fixed partition) which is a group?

49. Let V = {([0,a_1] | [0,a_2] [0,a_3] [0,a_4] | [0,a_5] [0,a_6] [0,a_7] [0,a_8] [0,a_9]) | $a_i \in Q \cup \{0\}$; $1 \le i \le 9$}. Prove V is a group under component wise multiplication. Can V have subgroup?

50. Let V = $\left\{ \begin{bmatrix} [0,a_1] & [0,a_2] & [0,a_3] & [0,a_4] \\ [0,a_5] & [0,a_6] & [0,a_7] & [0,a_8] \\ \hline [0,a_9] & [0,a_{10}] & [0,a_{11}] & [0,a_{12}] \\ [0,a_{13}] & [0,a_{14}] & [0,a_{15}] & [0,a_{18}] \\ [0,a_{17}] & [0,a_{18}] & [0,a_{19}] & [0,a_{20}] \\ \hline [0,a_{21}] & [0,a_{22}] & [0,a_{23}] & [0,a_{24}] \\ [0,a_{25}] & [0,a_{26}] & [0,a_{27}] & [0,a_{28}] \\ \hline [0,a_{29}] & [0,a_{30}] & [0,a_{31}] & [0,a_{321}] \\ [0,a_{33}] & [0,a_{34}] & [0,a_{35}] & [0,a_{36}] \end{bmatrix} \right.$ where $a_i \in Z^+ \cup$

{0}; $1 \le i \le 36$} be semivector space of super interval matrices of type I over $Z^+ \cup \{0\} = S$.
(a) Find subspaces of V.
(b) Find a basis of V.
(c) Find a linearly dependent set in V.
(d) Prove even a two element subset of V can be linearly dependent set.
(e) Define a linear operator on V which is invertible.
(f) Find a linear operator on V which is not invertible.

51. Let V = {([0,5] | [0,7] [0,12] [0, $\sqrt{3}$ ] | [0, 5] [0, $\sqrt{5}$ ]), (0 0 0 | 0 | 0 0),

$\begin{bmatrix} [0,7] \\ [0,4] \\ [0,9] \\ [0,11] \end{bmatrix}$, $\begin{bmatrix} 0 \\ \overline{0} \\ 0 \\ \overline{0} \end{bmatrix}$, $\left( \begin{array}{c|c} [0,7] & [0,2] \\ \hline [0,5] & \left[ 0, \dfrac{3}{2} \right] \\ \hline [0,1] & \left[ 0, \dfrac{9}{5} \right] \end{array} \right)$, $\left( \begin{array}{c|c} 0 & 0 \\ \hline 0 & 0 \\ \hline 0 & 0 \end{array} \right)$,

$$\left\{ \begin{bmatrix} [0,8] & [0,3] & [0,9] & [0,11] \\ \hline [0,4] & [0,\sqrt{5}] & [0,2] & [0,1] \\ [0,2] & [0,\sqrt{7}] & [0,5] & [0,7] \end{bmatrix}, \begin{bmatrix} 0 & 0 & 0 & 0 \\ \hline 0 & 0 & 0 & 0 \\ 0 & 0 & 0 & 0 \end{bmatrix} \right\}$$

be a set super vector space over the set S = {0, 1} (a) Find set super vector subspaces of V. (b) Define a linear operator from V to V? Can it be non trivial?

52. Prove the set of all partitions of a 7 × 5 interval matrix with entries from $Z^+ \cup \{0\}$ is a set vector space over the set {0, 1}.

53. Find the number of partitions of the interval matrix

$$P = \left\{ \begin{bmatrix} [0,a_1] & [0,a_2] & [0,a_3] & [0,a_4] & [0,a_5] \\ [0,a_6] & \ldots & \ldots & \ldots & [0,a_{10}] \\ [0,a_{11}] & \ldots & \ldots & \ldots & [0,a_{15}] \\ [0,a_{16}] & \ldots & \ldots & \ldots & [0,a_{20}] \\ [0,a_{21}] & [0,a_{22}] & [[0,a_{23}] & [0,a_{24}] & [0,a_{25}] \end{bmatrix} \right\}$$

where $a_i \in Z^+ \cup \{0\}$}. Hence or otherwise what is number of super interval matrices one can get using a n × n interval matrix.

54. Give an example of a set linear algebra of super interval matrices of finite order.

55. Obtain some interesting properties about set super linear algebra of super interval matrices.

56. Prove a set vector space of super interval matrices in general is not a set linear algebra.

57. Let V = {([0,a₁] | [0,a₂] [0,a₃] [0,a₄] | [0,a₅] . . . [0,a₁₁]) | $a_i \in Q^+ \cup \{0\}$; 1 ≤ i ≤ 11}. Prove V is a set linear algebra of super interval matrices over the set S = {$3Z^+ \cup 2Z^+ \cup 5Z^+ \cup \{0\}$}.
   (a) Find super set sublinear algebra of V.
   (b) Find pseudo super set vector subspaces of V.

58. Define for super set vector spaces V and W of super interval matrices the concept of set linear transformation and illustrate it by examples.

59. Prove if V is a set super vector space over the set S. Suppose $W_1, \ldots, W_t$ be n proper set vector subspaces of V. Then $\bigcap_i W_i = \phi$ can also occur.

60. Define the concept of set linear operator of a set vector space of super interval matrices and illustrate it by an example.

61. Can we define in case of set super vector space V defined over the set S the concept of set linear functionals from V into S.

62. Let $V = \left\{ \begin{bmatrix} [0, a_1] \\ [0, a_2] \\ \underline{[0, a_3]} \\ [0, a_4] \\ \underline{[0, a_5]} \\ [0, a_6] \end{bmatrix} \middle| a_i \in Z^+ \cup \{0\}; 1 \le i \le 6 \right\}$ be a set super linear algebra (set linear algebra of super column matrices) over the set S = {0, 1, 2, 5, 10, 27}.

(a) Find set sublinear algebras (set linear subalgebras or set linear super subalgebras) of V.
(b) Define a linear operator on V.
(c) Find pseudo set vector subspaces of V.
(d) Find subspaces $W_1$ $W_2$ $W_3$ in V so that $V = \cup W_i$ with $W_i \cap W_j = \begin{bmatrix} 0 \\ 0 \\ \underline{0} \\ 0 \\ \underline{0} \\ 0 \end{bmatrix}$, if $i \ne j$.

63. Let V = {([0,a_1] | [0,a_2] ... [0,a_7] [0,a_8] ... [0,a_{16}]) | $a_i \in Q^+ \cup \{0\}$; $1 \le i \le 16$} be a set super linear algebra over the set S = {0, 2, 8, 25, 48, 64} and W = $\left\{ \begin{bmatrix} \begin{array}{cc} [0,a_1] & [0,a_2] \\ [0,a_3] & [0,a_4] \\ [0,a_5] & [0,a_6] \\ \hline [0,a_7] & [0,a_8] \\ [0,a_9] & [0,a_{10}] \\ \hline [0,a_{11}] & [0,a_{12}] \\ [0,a_{13}] & [0,a_{14}] \end{array} \end{bmatrix} \middle| a_i \in Z^+ \cup \{0\}; 1 \le i \le 14 \right\}$ be a set linear algebra of super column interval vectors over the set S = {0, 2, 8, 25, 48, 64}.

   (a) Find a set linear transformation from V to W.
   (b) Find a set linear transformation from W to V.
   (c) Does there exist an invertible set linear transformation from V to W?
   (d) Suppose G = {set of all set linear transformations from V to W}, does G have any nice algebraic structure?

64. Let V = $\left\{ \begin{bmatrix} \begin{array}{cc|c} [0,a_1] & [0,a_2] & [0,a_3] \\ [0,a_4] & [0,a_5] & [0,a_6] \\ \hline [0,a_7] & [0,a_8] & [0,a_9] \\ [0,a_{10}] & [0,a_{11}] & [0,a_{12}] \\ [0,a_{13}] & [0,a_{14}] & [0,a_{15}] \end{array} \end{bmatrix} \middle| a_i \in Z^+ \cup \{0\}; 1 \le i \le 15 \right\}$

   be a set linear algebra of super interval matrices over the set S = {0, 1}.
   (a) Find at least two set linear subalgebras of V.
   (b) Find set linear subalgebras $W_i$ of V so that V = $\bigcup_i W_i$
   (c) Find at least two set linear operation on V.

65. Is the collection of a set linear operators on a set vector space of super interval matrices over a set S a set vector space over S?

66. Let $V = \left\{ \begin{pmatrix} [0,a_1] & [0,a_3] & [0,a_4] & [0,a_5] \\ [0,a_2] & [0,a_6] & [0,a_7] & [0,a_8] \end{pmatrix} \right.$ where $a_i \in R^+ \cup \{0\}\}$ be a set linear algebra of super interval row vectors over the set $S = Z^+ \cup \{0\}$. Find the set of set linear operators on V. Is it a set linear algebra over S?

67. Let $V = \{([0,a_1] \mid [0,a_2] \mid [0,a_3]\ [0,a_4]),\ ([0,a_1] \mid [0,a_2] \mid [0,a_3]\ [0,a_4]),\ ([0,a_1] \mid [0,a_2]\ [0,a_3]\ [0,a_4]),\ ([0,a_1]\ [0,a_2] \mid [0,a_3] \mid [0,a_4]),\ ([0,a_1]\ [0,a_2] \mid [0,a_3]\ [0,a_4]),\ ([0,a_1] \mid [0,a_2]\ [0,a_3]\ [0,a_4]),\ ([0,a_1]\ [0,a_2]\ [0,a_3] \mid [0,a_4]) \mid a_i \in 5Z^+ \cup \{0\}\}$ be a set vector space of super row interval matrices over the set $S = Z^+ \cup \{0\}$.
   (a) Find a generating set of V.
   (b) Find set vector subspaces of V.
   (c) Can V have pseudo set linear subalgebras?
   (d) Find linear operators on V which are invertible.
   (e) Find set subspaces $W_i$ of V so that $V = \cup W_i = \phi$.
   (f) Is it possible to find set linear operators which are not invertible?

68. Let $V = \{([0,a_1] \mid [0,a_2]\ [0,a_3]) \mid a_i \in Q^+ \cup \{0\};\ 1 \le i \le 3\}$ be a semigroup of super interval row matrices. Find transpose of V? Is $V^t$ a semigroup?

69. Let $V = \left\{ \begin{pmatrix} [0,a_1] & [0,a_2] & [0,a_3] \\ [0,a_4] & [0,a_6] & [0,a_7] \\ [0,a_5] & [0,a_8] & [0,a_9] \end{pmatrix} \right.$ $a_i \in Z^+ \cup \{0\};\ 1 \le i \le 9\}$ be a set linear algebra of super interval matrices over the set $S = 3Z^+ \cup \{0\} \cup 2Z^+$.
   (a) Find $W_1$ and $W_2$ two subspaces in V so that $V = W_1 \oplus W_2$.
   (b) Find at least one set linear operator on V.
   (c) Find at least two pseudo set vector subspaces of interval super matrix of V.

70. Define any other interesting properties about super set linear algebra of super interval matrices.

71. Does there exist a set linear algebra of super interval matrices which has no set linear subalgebra?

72. Define the concept of generating set in case of set vector space of super interval matrices and illustrate it by examples.

73. Define the concept of sectional subset vector sectional subspace in case of super set vector space (set vector space of super interval matrices).

74. Obtain some interesting relations described in problem 73.

75. Let $V = \{([0, a_1] [0,a_2] | [0,a_3]) | ([0,a_1] [0,a_2] [0,a_3] | [0,a_4] [0,a_5])$,
$\begin{bmatrix} [0,a_1] \\ \overline{[0,a_2]} \\ [0,a_3] \\ [0,a_4] \\ [0,a_5] \end{bmatrix}$ | $a_i \in 3Z^+ \cup \{0\} \cup 5Z^+$ $1 \leq i \leq 5\}$ be a set vector space

of super interval matrices defined over the set $S = \{0, 1\}$;
(a) Find a set linear operator on V.
(b) Find set vector subspaces of V.
(c) Write $V = \bigcup_i W_i$, $W_i$ set vector subspaces of V.

76. Obtain some interesting properties about semigroup vector space of super interval matrices over a semigroup S.

77. Let $V = \left\{ \begin{bmatrix} [0,a_1] & [0,a_2] \\ [0,a_3] & [0,a_4] \\ [0,a_5] & [0,a_6] \\ [0,a_7] & [0,a_8] \end{bmatrix} \right.$, $([0, a_1] [0,a_2] | [0,a_3])$ $([0,a_4] [0,a_5] |$

$[0,a_6]) | a_i \in Z^+ \cup \{0\}; 1 \leq i \leq 8\}$ be the semigroup vector space of super interval matrices over the semigroup $S = \in Z^+ \cup \{0\}$.

78. 
(a) Find at least two semigroup vector subspaces of V.
(b) Find at least two subsemigroup vector subspaces of V.
(c) Find atleast two semigroup linear operators on V.

Let $V = \left\{ \begin{bmatrix} [0,a_1] & [0,a_4] & [0,a_7] & [0,a_{10}] \\ [0,a_2] & [0,a_5] & [0,a_8] & [0,a_{11}] \\ [0,a_3] & [0,a_6] & [0,a_9] & [0,a_{12}] \end{bmatrix} \middle| a_i \in Z^+ \cup \{0\}; 1 \le i \le 12 \right\}$ be a semigroup super interval row vector linear algebra over the semigroup $S = Z^+ \oplus \cup \{0\}$. Let T = {all semigroup linear operators on V}. Find the algebraic structure enjoyed by T.

79. Give an example of a semigroup linear algebra of super interval matrices of order 28.

80. Let $V = \left\{ \begin{bmatrix} [0,a_1] & [0,a_6] & [0,a_{11}] & [0,a_{16}] \\ [0,a_2] & [0,a_7] & [0,a_{12}] & [0,a_{17}] \\ [0,a_3] & [0,a_8] & [0,a_{13}] & [0,a_{18}] \\ [0,a_4] & [0,a_9] & [0,a_{14}] & [0,a_{19}] \\ [0,a_5] & [0,a_{10}] & [0,a_{15}] & [0,a_{20}] \end{bmatrix} \middle| a_i \in Z_7 ; 1 \le i \le 20 \right\}$ be a semigroup super interval matrix linear algebra over the semigroup $S = Z_7$.
(a) Find the order of V.
(b) How many subsemigroup super interval linear subalgebra can V have?

81. Determine some important properties enjoyed by group linear algebra of superinterval matrices defined over the group G under addition.

82. Give an example of a pseudo simple group linear algebra.

83. Does there exist a group linear algebra of super interval matrices of order 23? Can this group linear algebra have group linear subalgebras?

84. Let $V = \left\{ \begin{bmatrix} [0,a_1] & [0,a_2] & [0,a_3] & [0,a_4] \\ [0,a_5] & [0,a_6] & [0,a_7] & [0,a_8] \\ [0,a_4] & [0,a_{10}] & [0,a_{11}] & [0,a_{12}] \\ [0,a_{13}] & [0,a_{14}] & [0,a_{15}] & [0,a_{16}] \end{bmatrix} \middle| a_i \in Z_{47} ; 1 \le i \le 16 \right\}$ be a group super interval matrix linear algebra over the group $G = Z_{47}$ under addition.
   (a) Find order of V.
   (b) Find a basis of V.
   (c) What is the dimension of V?
   (d) Can V have subgroup linear subalgebras?
   (e) Find atleast two invertible linear operators on V.

85. Obtain some interesting properties about superinterval fuzzy matrices.

86. How many super interval matrices can be constructed using a $7 \times 5$ interval matrix?

87. How many super interval row vectors can be constructed using a $4 \times 12$ interval matrix?

88. Give an example of a super $7 \times 6$ interval fuzzy row vector.

89. How many super interval matrices can be constructed using an interval matrix natural order $7 \times 42$?

90. Obtain some applications of super interval fuzzy matrices.

91. Give an example of a semigroup of super fuzzy interval matrices.

92. Given $A = \left\{ \begin{bmatrix} [0,a_1] & \cdots & [0,a_4] \\ \vdots & \cdots & \vdots \\ [0,a_{12}] & \cdots & [0,a_{16}] \end{bmatrix} \middle| a_i \in Z_{13} ; 1 \le i \le 16 \right\}$ be a group linear algebra over the group $G = Z_{13}$ under addition modulo 13.
   (a) Find group linear subalgebras of A.

(b) Prove A is pseudo simple.
(c) What is order A?
(d) What is dimension of A?
(e) Find a generating subset of A.
(f) Find a noninvertible linear operator on A.

93. Let $V = \left\{ \left( \begin{array}{cccc} [0,a_1] & [0,a_2] & [0,a_3] & [0,a_4] \\ \hline [0,a_5] & [0,a_6] & [0,a_7] & [0,a_8] \\ \hline [0,a_9] & [0,a_{10}] & [0,a_{11}] & [0,a_{12}] \\ \hline [0,a_{13}] & [0,a_{14}] & [0,a_{15}] & [0,a_{16}] \end{array} \right) \,\middle|\, a_i \in Z_{42};\ 1 \leq i \leq 16 \right\}$

be a group linear algebra of super interval column vectors over the group $G = Z_{42}$ under addition modulo 42.
(a) Find a basis of V.
(b) Prove V is finite dimensional
(c) Find 2 subgroup linear subalgebras.
(d) Find 3 pseudo group vector subspace of V.

94. Let $V = \left\{ \left[ \begin{array}{c|cc} [0,a_1] & [0,a_3] & [0,a_5] \\ \hline [0,a_2] & [0,a_4] & [0,a_6] \end{array} \right], \begin{bmatrix} [0,a_1] \\ \hline [0,a_2] \\ \hline [0,a_3] \\ [0,a_4] \\ [0,a_5] \\ [0,a_6] \\ [0,a_7] \end{bmatrix}, \left( \begin{array}{cc|c|c} [0,a_1] & [0,a_4] & [0,a_7] & [0,a_{10}] \\ [0,a_2] & [0,a_5] & [0,a_8] & [0,a_{11}] \\ [0,a_3] & [0,a_6] & [0,a_9] & [0,a_{12}] \end{array} \right) \,\middle|\, a_i \in Z_{15};\ 1 \right.$

$\leq i \leq 12$} be a group vector space of super interval matrices over the group $G = Z_{15}$ under addition modulo 15.
(a) Find order of V.
(b) Find at least 3 group vector subspaces of V.
(c) Find three group linear operators on V.
(d) Find the two subgroup vector subspaces of V.
(e) Obtain the special fuzzy group vector subspace of V.

95. Obtain any interesting property associated with super fuzzy interval matrices.

96. Define set vector space of super fuzzy interval matrices and illustrate it by an example.

97. How is special fuzzy set vector space different from the set vector space of super interval fuzzy matrices over a set S?

98. Show by an example the difference between these two structures mentioned in problem 97.

99. Prove every semigroup vector space of super interval matrices need not be a semigroup linear algebra.

100. Let $V = \left\{ \begin{bmatrix} [0,a_1] & [0,a_7] & [0,a_{13}] \\ [0,a_2] & [0,a_8] & [0,a_{14}] \\ \hline [0,a_3] & [0,a_9] & [0,a_{15}] \\ [0,a_4] & [0,a_{10}] & [0,a_{16}] \\ \hline [0,a_5] & [0,a_{11}] & [0,a_{17}] \\ [0,a_6] & [0,a_{121}] & [0,a_{18}] \end{bmatrix} \right.$ where $a_i \in Z_7$;

$1 \leq i \leq 18\}$ be the group linear algebra of super interval matrices over the group $G = Z_7$.
(a) Find the order of V.
(b) Find a generating subset of V.
(c) Find atleast two distinct linear operators on V.
(d) Find group linear subalgebras of V?
(e) Is V pseudo simple?
(f) Using V find the special fuzzy group linear algebra of super interval matrices over $Z_7 = G$.
(g) Find a fuzzy group linear subalgebra of V.

101. Let $V = \{([0, a] [0,a] | [0,a]) ([0,a] | [0,a] [0,a] [0,a]) | a \in Z_5 \}$ be a group linear algebra of super interval row matrices over the group $G = Z_5$.

(a) Is the order of V five?
    (b) Can V be generated by more than one element?
    (c) Find a basis of V.
    (d) Can V have group linear subalgebras?
    (e) Is V pseudo simple?
    (f) How many group linear operators on V can be defined?
    (g) How many special fuzzy group linear algebras can be obtained?

102. Obtain some interesting properties enjoyed by semigroup linear operators of a semigroup vector space of super interval matrices over a semigroup S.

103. Study the same problem in case of group linear algebras.

104. Let $V = \left\{ \begin{bmatrix} [0,a_1] & [0,a_6] & [0,a_{11}] & [0,a_{16}] \\ [0,a_2] & [0,a_7] & [0,a_{12}] & [0,a_{17}] \\ [0,a_3] & [0,a_8] & [0,a_{13}] & [0,a_{18}] \\ [0,a_4] & [0,a_9] & [0,a_{14}] & [0,a_{19}] \\ [0,a_5] & [0,a_{10}] & [0,a_{15}] & [0,a_{20}] \end{bmatrix} \right\}$ where $a_i \in Z_{24}$; $1 \leq i \leq 20$} be a group linear algebra of super interval matrices over the group $G = Z_{24}$ under addition.
    (a) Find order of V.
    (b) Find a basis of V.
    (c) How many subgroup linear subalgebras of V exist?
    (d) Find atleast five group linear subalgebras of V.
    (e) Define atleast two invertible linear operators on V.
    (f) Obtain atleast 6 special fuzzy group linear algebras using V.

105. Let $V = \{([0, a_1] \mid [0,a_2] [0,a_3] [0,a_4] \mid [0,a_5] [0,a_6]) \mid a_i \in Z_{20}$; $1 \leq i \leq 6\}$ be a group linear algebra over the group $G = Z_{20}$ under addition.

$$W = \left\{ \begin{bmatrix} [0, a_1] \\ [0, a_2] \\ [0, a_3] \\ [0, a_4] \\ [0, a_5] \\ [0, a_6] \\ [0, a_7] \end{bmatrix} \middle| \; a_i \in Z_{20}; \; 1 \le i \le 7 \right\}$$ be a group linear algebra over

the group $G = Z_{20}$

(a) Find atleast 2 invertible group linear transformations of V to W.

(b) If H = {all group linear transformations of V to W}, what is the algebraic structure enjoyed by H?

(c) How many distinct elements does H contain?

(d) What is the order of V?

(e) What is the order of W?

(f) If J = {all group linear transformations of W to V} what is the algebraic structure enjoyed by J?

(g) Is $H \cong J$?

(h) Find the special fuzzy group linear algebras of V and W.

(i) Find a group fuzzy linear transformation of $(V, \eta)$ with $(W, \mu)$.

106. Let $M = \left\{ \begin{bmatrix} [0, a_1] \\ [0, a_2] \\ [0, a_3] \\ [0, a_4] \\ [0, a_5] \end{bmatrix}, \; ([0, a_1] \, [0, a_2] \, | \, [0, a_3]) \; ([0, a_4] \, | \, [0, a_5]) \; \middle| \; a_i \in Z_{14} \, ;  \right.$

$1 \le i \le 5\}$ be a group vector space over the group $Z_{14}$ under addition.

(a) Find subgroup vector subspace of M.

(b) What is the order of M?

(c) Find a basis of M.

(d) Can basis of M exceed 10?

(e) Find atleast 6 (six) group vector subspaces of V.

(f) Find atleast five distinct group linear operator on V.

107. Let $V = \left\{ \begin{bmatrix} [0,a_1] \\ \overline{[0,a_2]} \\ [0,a_3] \\ \overline{[0,a_4]} \\ [0,a_5] \end{bmatrix}, ([0, a_1] \; [0,a_2] \mid [0,a_3] \; [0,a_4] \; [0,a_5]), \right.$

$\left. \begin{pmatrix} [0,a_1] & [0,a_3] & [0,a_5] \\ [0,a_2] & [0,a_4] & [0,a_6] \end{pmatrix} \right| a_i \in [0,1]; \; 1 \le i \le 6 \}$ be a set fuzzy vector space over the set $S = \{0,1\}$.

(a) Find three set fuzzy vector subspaces of V.
(b) Prove V is pseudo simple.
(c) Find five set linear operators on V.
(d) Find a basis of V.
(e) Is cardinality of V infinite?

108. Let $P = \left\{ \begin{bmatrix} [0,a_1] & [0,a_5] & [0,a_9] \\ [0,a_2] & [0,a_6] & [0,a_{10}] \\ [0,a_3] & [0,a_7] & [0,a_{11}] \\ [0,a_4] & [0,a_8] & [0,a_{12}] \end{bmatrix} \right| a_i \in [0,1]; 1 \le i \le 12 \}$ be a set fuzzy interval matrices. How many set of fuzzy super interval matrices can be obtained by partition the $4 \times 3$ interval matrix?

109. Let $W = \{([0, a_1] \mid [0,a_2] \; [0,a_3] \mid [0,a_4]) \mid a_i \in [0,1]; \; 1 \le i \le 4\}$ be a group fuzzy interval linear algebra over a fuzzy group G.
(a) Find group fuzzy interval linear subalgebras of W.
(b) Find atleast three group linear operators on W.
(c) Find a basis of V.

110. Let $P = \left\{ \begin{bmatrix} [0,a_1] \\ [0,a_2] \\ [0,a_3] \\ [0,a_4] \\ [0,a_5] \\ [0,a_6] \end{bmatrix} \middle| a_i \in Z^+ \cup \{0\}; 1 \leq i \leq 6 \right\}$ be a semigroup linear algebra over the semigroups $S = Z^+ \cup \{0\}$;

   (a) Find a basis of P.

   (b) Find the special fuzzy semigroup linear algebra of V.

   (c) Prove one can construct infinite number of special fuzzy semigroup linear algebras using V.

   (d) Find a linear operator on V.

   (e) Can a relation exist between semigroup linear operator on V and the special semigroup fuzzy linear operator on V.

111. Obtain some interesting properties enjoyed by special fuzzy interval group linear algebras.

112. Give some applications of group vector spaces of super interval matrices over a group in finite interval analysis.

113. Give some applications of super fuzzy interval matrices.

# FURTHER READING

# INDEX

# ABOUT THE AUTHORS

**Dr.W.B.Vasantha Kandasamy** is an Associate Professor in the Department of Mathematics, Indian Institute of Technology Madras, Chennai. In the past decade she has guided 13 Ph.D. scholars in the different fields of non-associative algebras, algebraic coding theory, transportation theory, fuzzy groups, and applications of fuzzy theory of the problems faced in chemical industries and cement industries. She has to her credit 646 research papers. She has guided over 68 M.Sc. and M.Tech. projects. She has worked in collaboration projects with the Indian Space Research Organization and with the Tamil Nadu State AIDS Control Society. She is presently working on a research project funded by the Board of Research in Nuclear Sciences, Government of India. This is her 56$^{th}$ book.

On India's 60th Independence Day, Dr.Vasantha was conferred the Kalpana Chawla Award for Courage and Daring Enterprise by the State Government of Tamil Nadu in recognition of her sustained fight for social justice in the Indian Institute of Technology (IIT) Madras and for her contribution to mathematics. The award, instituted in the memory of Indian-American astronaut Kalpana Chawla who died aboard Space Shuttle Columbia, carried a cash prize of five lakh rupees (the highest prize-money for any Indian award) and a gold medal.
She can be contacted at vasanthakandasamy@gmail.com
Web Site: http://mat.iitm.ac.in/home/wbv/public_html/
or http://www.vasantha.in

---

**Dr. Florentin Smarandache** is a Professor of Mathematics at the University of New Mexico in USA. He published over 75 books and 200 articles and notes in mathematics, physics, philosophy, psychology, rebus, literature.

In mathematics his research is in number theory, non-Euclidean geometry, synthetic geometry, algebraic structures, statistics, neutrosophic logic and set (generalizations of fuzzy logic and set respectively), neutrosophic probability (generalization of classical and imprecise probability). Also, small contributions to nuclear and particle physics, information fusion, neutrosophy (a generalization of dialectics), law of sensations and stimuli, etc. He can be contacted at smarand@unm.edu